\documentclass[12pt,a4paper]{article}
\usepackage[utf8]{inputenc}

\usepackage{mathrsfs,amsthm,xcolor,verbatim,bbm,amsmath,amsfonts,amssymb,nicefrac,enumitem,hyperref,bm,mathtools,xparse,etoolbox}
\usepackage{graphicx}
\usepackage{subcaption}  
\usepackage[capitalize,nameinlink,noabbrev]{cleveref}
\usepackage{cite}

\usepackage{thmtools}
\usepackage[margin=1in]{geometry}
\usepackage{tikz}
\usetikzlibrary{matrix,chains,positioning,decorations.pathreplacing,arrows,shapes,math,arrows.meta,calc, shapes.geometric}
\usepackage[most]{tcolorbox}

\makeatletter
\AtBeginEnvironment{tcolorbox}{%
   \def\@mpfn{footnote}%
}
\makeatother



\tikzset{
    neuron/.style={circle, draw, minimum size=0.5cm, text centered},
    input neuron/.style={neuron, fill=yellow!30},
    hidden neuron/.style={neuron, fill=blue!10, draw=blue, thick},
    output neuron/.style={neuron, fill=green!20},
    conn/.style={-stealth, thick},
    ellipsis/.style={font=\footnotesize},
}
\usepackage{capt-of}
\usepackage{stmaryrd} % for double brackets

\crefname{enumi}{item}{items}
\crefname{equation}{}{}
\crefname{subsection}{Subsection}{Subsections}
\crefname{figure}{Figure}{Figures}

\usepackage{xurl}

\usepackage{booktabs} % PHILIPPE
\usepackage{adjustbox}
\usepackage{array}          % column-type definitions for the optimizer table
\usepackage{colortbl}       % \cellcolor in the optimizer table
\usepackage{csvsimple-l3}   % CSV-driven optimizer table (see optimizer_table.tex)

\hypersetup{
    colorlinks,
    linkcolor={red!80!black},
    citecolor={green},
    urlcolor={blue!80!black}
}

\theoremstyle{plain}
\newtheorem{theorem}{Theorem}[section]

\newtheorem{prop}[theorem]{Proposition}

\newtheorem{definition}[theorem]{Definition}
\newtheorem{setting}[theorem]{Setting}



\crefname{lemma}{Lemma}{Lemmas}
\Crefname{lemma}{Lemma}{Lemmas}
\crefname{prop}{Proposition}{Propositions}
\Crefname{prop}{Proposition}{Propositions}
\crefname{cor}{Corollary}{Corollaries}
\Crefname{cor}{Corollary}{Corollaries}
\crefname{remark}{Remark}{Remarks}
\Crefname{remark}{Remark}{Remarks}
\crefname{definition}{Definition}{Definitions}
\Crefname{definition}{Definition}{Definitions}
\crefname{setting}{Setting}{Settings}
\Crefname{setting}{Setting}{Settings}
\crefname{conjecture}{Conjecture}{Conjectures}
\Crefname{conjecture}{Conjecture}{Conjectures}
\crefname{figure}{Figure}{Figures}
\Crefname{figure}{Figure}{Figures}

\usepackage{graphicx}

\usepackage{expl3,xparse} % expl3 utilities and xparse

\ExplSyntaxOn
\clist_new:N \l__my_fig_clist
\int_new:N  \l__my_fig_int

\NewDocumentCommand{\figref}{m}
  { \my_figref_process:n { #1 } }

\cs_new_protected:Npn \my_figref_process:n #1
  {
    \clist_set:Nn \l__my_fig_clist { #1 }
    \int_set:Nn \l__my_fig_int { \clist_count:N \l__my_fig_clist }

    \int_step_inline:nn { \l__my_fig_int }
      {
        \int_compare:nNnTF { ##1 } = { 1 }
          { }
          {
            \int_compare:nNnTF { ##1 } = { \l__my_fig_int }
              {, \space and\space}
              {,\space }
          }
        Figure~\hyperref[\clist_item:Nn \l__my_fig_clist { ##1 }]
                   {\ref*{\clist_item:Nn \l__my_fig_clist { ##1 }}}
      }
  }
\ExplSyntaxOff

\theoremstyle{definition}

\DeclareMathAlphabet{\mathpzc}{OT1}{pzc}{m}{it}

\DeclareFontEncoding{LS1}{}{}
\DeclareFontSubstitution{LS1}{stix}{m}{n}
\DeclareMathAlphabet{\mathscr}{LS1}{stixscr}{m}{n}

\newcommand{\tpd}{two-point distributed\cfadd{definition: two point distributed}}

\newcommand{\E}{\mathbb{E}}
\renewcommand{\P}{\mathbb{P}}

\newcommand{\R}{\mathbb{R}}
\newcommand{\N}{\mathbb{N}}
\newcommand{\Z}{\mathbb{Z}}

\newcommand{\smalll}{\mathbb{L}}
\newcommand{\bbL}{\mathbb{L}}

\newcommand{\setX}{K}

\renewcommand{\d}{ \mathrm{d}}

\renewcommand{\c}[1]{\mathfrak{c}^{#1}}

\newcommand{\cons}{c}
\newcommand{\adamvector}{f}

\newcommand{\cB}{\mathcal{B}}

\newcommand{\cF}{\mathcal{F}}

\newcommand{\cK}{\mathcal{K}}
\newcommand{\cL}{\mathcal{L}}
\newcommand{\cM}{\mathcal{M}}

\newcommand{\cZ}{\mathcal{Z}}

\newcommand{\bfc}{\mathbf{c}}

\newcommand{\bfm}{\mathbf{m}}
\newcommand{\bfn}{\mathbf{n}}

\newcommand{\bfx}{\mathbf{x}}

\newcommand{\globalmin}{\vartheta}
\newcommand{\roott}{\theta}
\newcommand{\enpointp}{p}
\newcommand{\enpointq}{q}

\newcommand{\bfB}{\mathbf{B}}
\newcommand{\bfC}{\mathbf{C}}

\newcommand{\bfM}{\mathbf{M}}

\newcommand{\bfV}{\mathbf{V}}

\newcommand{\bfX}{\mathbf{X}}

\newcommand{\scrc}{\mathscr{c}}

\newcommand{\scrM}{\mathscr{M}}

\newcommand{\vertiii}[1]{{\left\vert\kern-0.25ex\left\vert\kern-0.25ex\left\vert #1 
    \right\vert\kern-0.25ex\right\vert\kern-0.25ex\right\vert}}

\newcommand{\g}{\|\mathfrak{g}\|_{Lip}}

\newcommand{\fC}{\mathfrak{C}}

\newcommand{\fM}{\mathfrak{M}}
\newcommand{\fN}{\mathfrak{N}}

\newcommand{\fZ}{\mathfrak{Z}}
\newcommand{\bbX}{\mathbb{X}}

\newcommand{\fc}{\mathfrak{c}}
\newcommand{\fd}{\mathfrak{d}}

\newcommand{\fn}{\mathfrak{n}}

\newcommand{\fp}{\mathfrak{p}}
\newcommand{\fq}{\mathfrak{q}}

\newcommand{\scrd}{\mathscr{d}}
\newcommand{\f}{f}

\newcommand{\bbT}{\mathbb{T}}

\newcommand{\dimX}{\mathscr{d}}

\newcommand{\consc}{c}

\newcommand{\mandV}[4]{\mathbb{V}^{#1,#2,#3}_{#4}}
\newcommand{\mandVcom}[5]{\mathbb{V}^{#1,#2,#3}_{#4,#5}}
\newcommand{\Adamprocess}[3]{\Theta^{#1,#2}_{#3}}
\newcommand{\Adamprocesscom}[4]{\Theta^{#1,#2}_{#3,#4}}
\newcommand{\mandVV}[4]{\mathbf{V}^{#1,#2,#3}_{#4}}
\newcommand{\mandVcomV}[5]{\mathbf{V}^{#1,#2,#3}_{#4,#5}}
\newcommand{\AdamprocessV}[3]{\varTheta^{#1,#2}_{#3}}
\newcommand{\AdamprocesscomV}[4]{\varTheta^{#1,#2}_{#3,#4}}

\newcommand{\loss}{\mathcal{L}}
\newcommand{\para}{\beta}

\newcommand{\bbY}{\mathbb{Y}}

\newcommand{\stronglyconvex}{generalized convex\cfadd{definition: strongly convex}}
\newcommand{\stronglyconvexx}{strongly convex\cfadd{definition: strongly convex}}

\newcommand\restr[2]{{% we make the whole thing an ordinary symbol
  \left.\kern-\nulldelimiterspace % automatically resize the bar with \right
  #1 % the function
  \vphantom{|} % pretend it's a little taller at normal size
  \right|_{#2} % this is the delimiter
  }}

\DeclarePairedDelimiter{\spro}{\langle}{\rangle}

\newcommand{\qandq}{\quad \text{and} \quad }
\newcommand{\qqandqq}{\qquad\text{and}\qquad}
\newcommand{\andq}{\text{and}\qquad}

\newcommand{\bbM}{\mathbb{M}}
\newcommand{\democrator}{\mathbb{V}}

\ExplSyntaxOn

\seq_new:N \g_cflist_loaded
\seq_new:N \g_cflist_pending

\NewDocumentCommand{\cfadd} { m } {
  \seq_if_in:NnF \g_cflist_loaded { #1 } {
    \seq_if_in:NnF \g_cflist_pending { #1 } {
      \seq_gput_right:Nn \g_cflist_pending { #1 }
    }
  }
}

\NewDocumentCommand{\cfconsiderloaded} { m } {
  \seq_gput_right:Nn \g_cflist_loaded {#1}
}

\NewDocumentCommand{\cfremove} { m } {
  \seq_gremove_all:Nn \g_cflist_pending { #1 }
}

\NewDocumentCommand{\cfload} { o } {
  \seq_if_empty:NTF \g_cflist_pending {
    \IfValueTF{#1}{\ignorespaces}{\unskip}
  } {
    (cf.\ \cref{\seq_use:Nn \g_cflist_pending {,}})\IfValueTF{#1}{#1~}{\unskip}
    \seq_gconcat:NNN \g_cflist_loaded \g_cflist_loaded \g_cflist_pending
    \seq_gclear:N \g_cflist_pending
    \IfValueT{#1}{\ignorespaces}
  }
}

\NewDocumentCommand{\cfclear} {} {
  \seq_gclear:N \g_cflist_loaded
  \seq_gclear:N \g_cflist_pending
}

\NewDocumentCommand{\cfout} { o } {
  \seq_if_empty:NTF \g_cflist_pending {\unskip\IfValueT{#1}{\ignorespaces}} {
    (cf.\ \cref{\seq_use:Nn \g_cflist_pending {,}})\IfValueTF{#1}{#1~}{\unskip}
    \seq_gclear:N \g_cflist_pending
    \IfValueT{#1}{\ignorespaces}
  }
}

\NewDocumentCommand{\ifnocf} { m } {
  \seq_if_empty:NT \g_cflist_pending { #1 }
}

\ExplSyntaxOff

\ExplSyntaxOn

\bool_new:N \g_noteobserve

\NewDocumentCommand{\setnote}{}{
  \bool_gset_true:N \g_noteobserve
}

\NewDocumentCommand{\setobserve}{}{
  \bool_gset_false:N \g_noteobserve
}

\NewDocumentCommand{\nobs}{ o }{
  \IfValueT{#1}{
    \str_if_eq:noTF {note} {#1} {
      \bool_gset_true:N \g_noteobserve
    } {
      \str_if_eq:noTF {Note} {#1} {
        \bool_gset_true:N \g_noteobserve
      } {
        \bool_gset_false:N \g_noteobserve
      }
    }
  }
  \bool_if:nTF { \g_noteobserve } {
    \bool_gset_false:N \g_noteobserve
    note
  } {
    \bool_gset_true:N \g_noteobserve
    observe
  }
  \IfValueF{#1}{~}
}

\NewDocumentCommand{\Nobs}{ o }{
  \IfValueT{#1}{
    \str_if_eq:noTF {note} {#1} {
      \bool_gset_true:N \g_noteobserve
    } {
      \str_if_eq:noTF {Note} {#1} {
        \bool_gset_true:N \g_noteobserve
      } {
        \bool_gset_false:N \g_noteobserve
      }
    }
  }
  \bool_if:nTF { \g_noteobserve } {
    \bool_gset_false:N \g_noteobserve
    Note
  } {
    \bool_gset_true:N \g_noteobserve
    Observe
  }
  \IfValueF{#1}{~}
}

\ExplSyntaxOff

\ExplSyntaxOn

\bool_new:N \g_hencetherefore

\NewDocumentCommand{\hence}{ o }{
  \IfValueT{#1}{
    \str_if_eq:noTF {hence} {#1} {
      \bool_gset_true:N \g_hencetherefore
    } {
      \str_if_eq:noTF {Hence} {#1} {
        \bool_gset_true:N \g_hencetherefore
      } {
        \bool_gset_false:N \g_hencetherefore
      }
    }
  }
  \bool_if:nTF { \g_hencetherefore } {
    \bool_gset_false:N \g_hencetherefore
    hence
  } {
    \bool_gset_true:N \g_hencetherefore
    therefore
  }
  \IfValueF{#1}{~}
}

\NewDocumentCommand{\Hence}{ o }{
  \IfValueT{#1}{
    \str_if_eq:noTF {hence} {#1} {
      \bool_gset_true:N \g_hencetherefore
    } {
      \str_if_eq:noTF {Hence} {#1} {
        \bool_gset_true:N \g_hencetherefore
      } {
        \bool_gset_false:N \g_hencetherefore
      }
    }
  }
  \bool_if:nTF { \g_hencetherefore } {
    \bool_gset_false:N \g_hencetherefore
    Hence,~we~obtain
  } {
    \bool_gset_true:N \g_hencetherefore
    Therefore,~we~obtain
  }
  \IfValueF{#1}{~}
}

\ExplSyntaxOff

\ExplSyntaxOn

\seq_const_from_clist:Nn \g_prove_mru {
  establish,
  demonstrate,
  prove,
  show,
  imply,
  ensure
}

\prop_new:N \l__verbs
\prop_put:Nnn \l__verbs {show} {shows}
\prop_put:Nnn \l__verbs {imply} {implies}
\prop_put:Nnn \l__verbs {demonstrate} {demonstrates}
\prop_put:Nnn \l__verbs {prove} {proves}
\prop_put:Nnn \l__verbs {establish} {establishes}
\prop_put:Nnn \l__verbs {ensure} {ensures}
\prop_put:Nnn \l__verbs {assure} {assures}

\tl_new:N \g_wordtmp
\seq_new:N \l_mytmps

\cs_generate_variant:Nn \str_if_in:nnTF { nVTF }
\cs_generate_variant:Nn \str_if_in:nnTF { xVTF }

\NewDocumentCommand{\prove}{ o }{
  \IfValueTF{#1}{
    \seq_clear:N \l_mytmps
    \seq_map_inline:Nn \g_prove_mru {
      \str_if_eq:nnTF {##1} {ensure} {
        \str_set:Nn \l_temps {n}
      } {
        \str_set:Nx \l_temps {\str_head_ignore_spaces:n {##1}}
      }
      \str_if_in:xVTF {#1} \l_temps {
        \seq_put_right:Nn \l_mytmps {##1}
      } { }
    }
    \seq_get_right:NN \l_mytmps \g_wordtmp
  } {
    \seq_get_right:NN \g_prove_mru \g_wordtmp
  }
  \tl_use:N \g_wordtmp
  \IfValueTF{#1}{}{~}
  \seq_gput_left:NV \g_prove_mru \g_wordtmp
  \seq_gremove_duplicates:N \g_prove_mru
}

\NewDocumentCommand{\proves}{ o }{
  \IfValueTF{#1}{
    \seq_clear:N \l_mytmps
    \seq_map_inline:Nn \g_prove_mru {
      \str_if_eq:nnTF {##1} {ensure} {
        \str_set:Nn \l_temps {n}
      } {
        \str_set:Nx \l_temps {\str_head_ignore_spaces:n {##1}}
      }
      \str_if_in:xVTF {#1} \l_temps {
        \seq_put_right:Nn \l_mytmps {##1}
      } { }
    }
    \seq_get_right:NN \l_mytmps \g_wordtmp
  } {
    \seq_get_right:NN \g_prove_mru \g_wordtmp
  }
  \str_set:NV \l_tmpa_str \g_wordtmp
  \prop_get:NVN \l__verbs \l_tmpa_str \l_tmpa_tl
  \tl_use:N \l_tmpa_tl
  \IfValueTF{#1}{}{~}
  \seq_gput_left:NV \g_prove_mru \g_wordtmp
  \seq_gremove_duplicates:N \g_prove_mru
}

\newcommand{\llabel}[1]{\savelabel{#1}\label{\loc.#1}\ignorespaces}

\clist_new:N \l_localreflist
\clist_new:N \l_reflist

\NewDocumentCommand{\lref} { m } {
  \clist_set:No \l_localreflist {#1}
  \clist_clear:N \l_reflist
  \clist_map_inline:Nn \l_localreflist { \clist_put_right:Nn \l_reflist {\loc.##1} }
  \cref{\l_reflist}
}

\NewDocumentCommand{\Lref} { m } {
  \clist_set:No \l_localreflist {#1}
  \clist_clear:N \l_reflist
  \clist_map_inline:Nn \l_localreflist { \clist_put_right:Nn \l_reflist {\loc.##1} }
  \Cref{\l_reflist}
}

\NewDocumentCommand{\itref}{ m m }{
  \clist_set:No \l_localreflist {#2}
  \clist_clear:N \l_reflist
  \clist_map_inline:Nn \l_localreflist { \clist_put_right:Nn \l_reflist {#1.##1} }
  \cref{\l_reflist}~in~\cref{#1}
}

\seq_new:N \l_enum_seq
\int_new:N \l_num_items

\bool_new:N \g_commaused_bool

\providecommand{\comma}{}

\cs_new:Nn \enum_it:nn {
  \int_case:nnF {\l_num_items - #1} {
    {0} {
      \renewcommand{\comma}{}
      #2\space
    }
    {1} {
      \bool_gset_false:N \g_commaused_bool
      \renewcommand{\comma}{,~\bool_gset_true:N \g_commaused_bool}
      #2
      \bool_if:NTF \g_commaused_bool {} {,~}
      and~
    }
  } {
    \bool_gset_false:N \g_commaused_bool
    \renewcommand{\comma}{,~\bool_gset_true:N \g_commaused_bool}
    #2
    \bool_if:NTF \g_commaused_bool {} {,~}
  }
}

\cs_new:Nn \enum_it_U:nn {
  \int_case:nnF {\l_num_items - #1} {
    {0} {
      \renewcommand{\comma}{}
      #2
      \space
    }
    {1} {
      \bool_gset_false:N \g_commaused_bool
      \renewcommand{\comma}{,~\bool_gset_true:N \g_commaused_bool}
      #2
      \bool_if:NTF \g_commaused_bool {} {,~}
      and~
    }
  } {
    \bool_gset_false:N \g_commaused_bool
    \renewcommand{\comma}{,~\bool_gset_true:N \g_commaused_bool}
    \int_compare:nTF {#1=1} {
      \text_titlecase_first:n {#2}
    } {
      #2
    }
    \bool_if:NTF \g_commaused_bool {} {,~}
  }
}

\cs_generate_variant:Nn \tl_if_eq:nnTF {onTF}

\cs_new:Nn \enum:nnnn {
  \seq_set_split:Nnn \l_enum_seq ; {#1}
  \seq_remove_all:Nn \l_enum_seq { }
  \seq_remove_all:Nn \l_enum_seq {#2}
  \seq_log:N \l_enum_seq
  \int_set:Nn \l_num_items {\seq_count:N \l_enum_seq}
  \int_log:N \l_num_items
  \int_case:nnF {\l_num_items} {
    { 0 } { 0 }
    { 1 } {
      \IfBooleanTF{#4} {
        \tl_set:Nn \l_text_case_exclude_arg_tl {\cref}
        \bool_set_false:N \l_text_titlecase_check_letter_bool
        \text_titlecase_first:n {\seq_use:Nn \l_enum_seq {}}
      } {
        \seq_use:Nn \l_enum_seq {}
      }
      \space
      \tl_if_eq:onTF{#3}{-}{}{
        \bool_if:NTF \l_plural_bool {
          \prove[#3]~
        } {
          \proves[#3]~
        }
      }
    }
    { 2 } {
      \IfBooleanTF{#4} {
        \tl_set:Nn \l_text_case_exclude_arg_tl {\cref}
        \bool_set_false:N \l_text_titlecase_check_letter_bool
        \text_titlecase_first:n {\seq_use:Nn \l_enum_seq {\ and\ }}
      } {
        \seq_use:Nn \l_enum_seq {\ and\ }
      }
      \space
      \tl_if_eq:onTF{#3}{-}{}{
        \prove[#3]~
      }
    }
  } {
    \IfBooleanTF{#4} {
      \tl_set:Nn \l_text_case_exclude_arg_tl {\cref}
      \seq_indexed_map_function:NN \l_enum_seq \enum_it_U:nn
    } {
      \seq_indexed_map_function:NN \l_enum_seq \enum_it:nn
    }
    \tl_if_eq:onTF{#3}{-}{}{
      \prove[#3]~
    }
  }
}

\cs_generate_variant:Nn \enum:nnnn {nxnn}
\cs_generate_variant:Nn \enum:nnnn {nxxn}

\NewDocumentCommand{\enum}{O{} m O{-} s}{
  \IfBooleanTF{#4}{
    \enum:nxnn {#2} {#1} {sindep} \BooleanFalse
  } {
    \enum:nxxn {#2} {#1} {#3} \BooleanFalse
  }
}

\NewDocumentCommand{\dott}{}{\ifnocf{.}\space}

\bool_new:N \g_arg_start_bool
\bool_gset_true:N \g_arg_start_bool

\NewDocumentCommand{\startnewargseq}{}{\bool_gset_true:N \g_arg_start_bool \tl_set:Nn \g_label_tl {}}

\cs_generate_variant:Nn \seq_if_in:NnTF {NxTF}
\cs_generate_variant:Nn \seq_remove_all:Nn {Nx}

\int_new:N \l_random_int

\cs_generate_variant:Nn \tl_if_head_eq_catcode:nNTF {oNTF}
\cs_generate_variant:Nn \tl_if_head_eq_catcode:nNTF {VNTF}

\cs_generate_variant:Nn \tl_if_head_eq_catcode:nNTF {eNTF}
\cs_generate_variant:Nn \tl_log:n {o}
\cs_generate_variant:Nn \tl_log:n {f}
\cs_generate_variant:Nn \tl_log:n {x}
\cs_generate_variant:Nn \tl_log:n {e}

\cs_generate_variant:Nn \tl_if_in:nnTF {onTF}
\cs_generate_variant:Nn \tl_if_in:NnTF {NeTF}

\cs_generate_variant:Nn \tl_if_head_eq_meaning:nNTF {VNTF}

\seq_const_from_clist:Nn \g_arg_mru_this {
  Ahpr,
  Tapr,
  Ctapr,
  H
}

\seq_const_from_clist:Nn \g_arg_mru_nothis {
  Ia,
  Nwc,
  N,
  Itns,
  Fm,
  Itnswc,
  Mo
}

\seq_new:N \l_arg_seq
\tl_new:N \l_cons_tl
\tl_new:N \l_dummy_tl

\bool_new:N \g_debug_bool
\bool_gset_false:N \g_debug_bool

\bool_new:N \l_insidearg_bool

\bool_new:N \g_firstargletter_bool

\sys_gset_rand_seed:n {0903}

\bool_new:N \l_plural_bool
\tl_new:N \l_arg_verbs_tl

\NewDocumentCommand{\argument}{mom}{
\color{black}
  \bool_set_false:N \l_plural_bool
  \tl_set:Nn \l_arg_verbs_tl {sindep}
  \keys_define:nn { benno/argument } {
    plural .value_forbidden:n = true,
    plural .code:n = {\bool_set_true:N \l_plural_bool},
    verbs .value_required:n = false,
    verbs .tl_set:N = \l_arg_verbs_tl,
  }
  \IfValueT{#2}{
    \keys_set:nn { benno/argument } {#2}
  }
  \bool_log:N \l_plural_bool
  \bool_gset_true:N \l_insidearg_bool
  \seq_set_split:Nnn \l_arg_seq ; {#1}
  \seq_remove_all:Nn \l_arg_seq { }
  \seq_log:N \l_arg_seq
  \tl_set:Nn \l_cons_tl {#3}
  \tl_trim_spaces:N \l_cons_tl
  \seq_if_in:NxTF \l_arg_seq {\lref{\g_label_tl}} {
    \seq_remove_all:Nx \l_arg_seq {\lref{\g_label_tl}}
    \seq_get_left:NNTF \l_arg_seq \l_dummy_tl {
      \tl_trim_spaces:N \l_dummy_tl
      \bool_gset_false:N \g_firstargletter_bool
      \tl_if_head_eq_catcode:VNTF \l_dummy_tl a {
        \bool_gset_true:N \g_firstargletter_bool
      } {
        \tl_if_head_eq_meaning:VNTF \l_dummy_tl {\cref} {
          \tl_set:Nx \l_tmpa_tl {\tl_tail:N \l_dummy_tl}
          \tl_set:Nx \l_tmpb_tl {\tl_head:N \l_tmpa_tl}
          \bool_gset_true:N \g_firstargletter_bool
          \tl_if_in:NeTF \l_tmpb_tl {lem\c_colon_str} {} {
            \tl_if_in:NeTF \l_tmpb_tl {thm\c_colon_str} {} {
              \tl_if_in:NeTF \l_tmpb_tl {prop\c_colon_str} {} {
                \tl_if_in:NeTF \l_tmpb_tl {cor\c_colon_str} {} {
                  \bool_gset_false:N \g_firstargletter_bool
                }
              }
            }
          }
        } {
        }
      }
      \bool_if:NTF \g_firstargletter_bool {
        \seq_set_eq:NN \l_tmpa_seq \g_arg_mru_this
        \seq_remove_all:Nn \l_tmpa_seq {H}
        \seq_get_right:NN \l_tmpa_seq \l_tmpa_tl
        \int_case:nnF {\seq_count:N \l_arg_seq} {
          {1} {
            \str_case:VnF {\l_tmpa_tl} {
              {Ahpr} {
                \bool_if:NT \g_debug_bool {C1.1}
                \seq_gput_left:Nn \g_arg_mru_this {Ahpr}
                \seq_gremove_duplicates:N \g_arg_mru_this
                \enum:nxnn {#1} {\lref{\g_label_tl}} {-} {\BooleanTrue}
                \hence~
                \bool_if:NTF \l_plural_bool {
                  \prove[\l_arg_verbs_tl]~\ignorespaces #3
                } {
                  \proves[\l_arg_verbs_tl]~\ignorespaces #3
                }
              }
              {Tapr} {
                \bool_if:NT \g_debug_bool {C1.2}
                \seq_gput_left:Nn \g_arg_mru_this {Tapr}
                \seq_gremove_duplicates:N \g_arg_mru_this
                \enum[\lref{\g_label_tl}]{
                  This;
                  #1
                }[\l_arg_verbs_tl]\ignorespaces #3
              }
              {Ctapr} {
                \bool_if:NT \g_debug_bool {C1.3}
                \seq_gput_left:Nn \g_arg_mru_this {Ctapr}
                \seq_gremove_duplicates:N \g_arg_mru_this
                Combining~
                \enum[\lref{\g_label_tl}]{
                  this;
                  #1
                } \proves[\l_arg_verbs_tl]~\ignorespaces #3
              }
            } {}
          }
        } {
          \str_case:VnF {\l_tmpa_tl} {
             {Ahpr} {
              \bool_if:NT \g_debug_bool {C2.1}
              \seq_gput_left:Nn \g_arg_mru_this {Ahpr}
              \seq_gremove_duplicates:N \g_arg_mru_this
              \enum:nxnn {#1} {\lref{\g_label_tl}} {-} {\BooleanTrue}
              \hence~
              \prove[\l_arg_verbs_tl]~\ignorespaces #3
            }
            {Tapr} {
              \bool_if:NT \g_debug_bool {C2.2}
              \seq_gput_left:Nn \g_arg_mru_this {Tapr}
              \seq_gremove_duplicates:N \g_arg_mru_this
              \enum[\lref{\g_label_tl}]{
                This;
                #1
              }[\l_arg_verbs_tl]\ignorespaces #3
            }
            {Ctapr} {
              \int_case:nn {\int_rand:nn {0} {1}} {
                {0} {
                  \bool_if:NT \g_debug_bool {C2.3}
                  \seq_gput_left:Nn \g_arg_mru_this {Ctapr}
                  \seq_gremove_duplicates:N \g_arg_mru_this
                  Combining~
                  \enum[\lref{\g_label_tl}]{
                    this;
                    #1
                  } \proves[\l_arg_verbs_tl]~\ignorespaces #3
                }
                {1} {
                  \bool_if:NT \g_debug_bool {C2.4}
                  \seq_gput_left:Nn \g_arg_mru_this {Ctapr}
                  \seq_gremove_duplicates:N \g_arg_mru_this
                  Combining~
                  \enum:nxnn {#1} {\lref{\g_label_tl}} {-} {\BooleanFalse}
                  \hence~
                  \proves[\l_arg_verbs_tl]~\ignorespaces #3
                }
              }
            }
          } {}
        }
      } {
        \seq_set_eq:NN \l_tmpa_seq \g_arg_mru_this
        \seq_remove_all:Nn \l_tmpa_seq {H}
        \seq_remove_all:Nn \l_tmpa_seq {Ahpr}
        \seq_get_right:NN \l_tmpa_seq \l_tmpa_tl
        \int_case:nnF {\seq_count:N \l_arg_seq} {
          {1} {
            \str_case:VnF {\l_tmpa_tl} {
              {Tapr} {
                \bool_if:NT \g_debug_bool {C3.1}
                \seq_gput_left:Nn \g_arg_mru_this {Tapr}
                \seq_gremove_duplicates:N \g_arg_mru_this
                \enum[\lref{\g_label_tl}]{
                  This;
                  #1
                }[\l_arg_verbs_tl]\ignorespaces #3
              }
              {Ctapr} {
                \bool_if:NT \g_debug_bool {C3.2}
                \seq_gput_left:Nn \g_arg_mru_this {Ctapr}
                \seq_gremove_duplicates:N \g_arg_mru_this
                Combining~
                \enum[\lref{\g_label_tl}]{
                  this;
                  #1
                } \proves[\l_arg_verbs_tl]~\ignorespaces #3
              }
            } {}
          }
        } {
          \str_case:VnF {\l_tmpa_tl} {
            {Tapr} {
              \bool_if:NT \g_debug_bool {C4.1}
              \seq_gput_left:Nn \g_arg_mru_this {Tapr}
              \seq_gremove_duplicates:N \g_arg_mru_this
              \enum[\lref{\g_label_tl}]{
                This;
                #1
              }[\l_arg_verbs_tl]\ignorespaces #3		
            }
            {Ctapr} {
              \int_case:nn {\int_rand:nn {0} {1}} {
                {0} {
                  \bool_if:NT \g_debug_bool {C4.2}
                  \seq_gput_left:Nn \g_arg_mru_this {Ctapr}
                  \seq_gremove_duplicates:N \g_arg_mru_this
                  Combining~
                  \enum[\lref{\g_label_tl}]{
                    this;
                    #1
                  } \proves[\l_arg_verbs_tl]~\ignorespaces #3		
                }
                {1} {
                  \bool_if:NT \g_debug_bool {C4.3}
                  \seq_gput_left:Nn \g_arg_mru_this {Ctapr}
                  \seq_gremove_duplicates:N \g_arg_mru_this
                  Combining~
                  \enum:nxnn {#1} {\lref{\g_label_tl}} {-} {\BooleanFalse}
                  \hence~
                  \proves[\l_arg_verbs_tl]~\ignorespaces #3    
                }
              }
            }
          } {}
        }
      }
    } {
      \tl_if_head_eq_catcode:oNTF \l_cons_tl a {
        \seq_set_eq:NN \l_tmpa_seq \g_arg_mru_this
        \seq_remove_all:Nn \l_tmpa_seq {Ctapr}
        \seq_remove_all:Nn \l_tmpa_seq {Ahpr}
        \seq_get_right:NN \l_tmpa_seq \l_tmpa_tl
        \str_case:VnF {\l_tmpa_tl} {
          {H} {
            \bool_if:NT \g_debug_bool {C5.1}
            \seq_gput_left:Nn \g_arg_mru_this {H}
            \seq_gremove_duplicates:N \g_arg_mru_this
            Hence,~we~obtain~\ignorespaces #3
          }
          {Tapr} {
            \bool_if:NT \g_debug_bool {C5.2}
            \seq_gput_left:Nn \g_arg_mru_this {Tapr}
            \seq_gremove_duplicates:N \g_arg_mru_this
            This~\proves[\l_arg_verbs_tl]~\ignorespaces #3
          }
        } {}
      } {
        \bool_if:NT \g_debug_bool {C6.1}
        \seq_gput_left:Nn \g_arg_mru_this {Tapr}
        \seq_gremove_duplicates:N \g_arg_mru_this
        This~\proves[\l_arg_verbs_tl]~\ignorespaces #3
      }
    } 
  } {
    \int_compare:nNnTF {\seq_count:N \l_arg_seq} = {0} {
      \bool_if:NTF \g_arg_start_bool {
        \bool_if:NT \g_debug_bool {C7.1}
        \Nobs\unskip
        #3
      } {
        \bool_if:NT \g_debug_bool {C7.2}
        \Moreover~
        #3
      }
    } {
      \bool_if:NTF \g_arg_start_bool {
        \bool_if:NT \g_debug_bool {C8.1}
        \tl_log:N \l_arg_verbs_tl
        \Nobs~that~
        \enum{
          #1
        }[\l_arg_verbs_tl]\ignorespaces #3
      } {
        \int_compare:nNnTF {\seq_count:N \l_arg_seq} = {1} {
          \seq_set_eq:NN \l_tmpa_seq \g_arg_mru_nothis
          \seq_remove_all:Nn \l_tmpa_seq {Nwc}
          \seq_remove_all:Nn \l_tmpa_seq {Itnswc}
          \seq_get_right:NN \l_tmpa_seq \l_tmpa_tl
        } {
          \seq_get_right:NN \g_arg_mru_nothis \l_tmpa_tl
        }
        \str_case:VnF {\l_tmpa_tl} {
          {Mo} {
            \bool_if:NT \g_debug_bool {C9.1}
            \seq_gput_left:Nn \g_arg_mru_nothis {Mo}
            \seq_gremove_duplicates:N \g_arg_mru_nothis
            Moreover,~\nobs~that~
            \enum{
              #1
            }[\l_arg_verbs_tl]\ignorespaces #3		
          }
          {Fm} {
            \bool_if:NT \g_debug_bool {C9.2}
            \seq_gput_left:Nn \g_arg_mru_nothis {Fm}
            \seq_gremove_duplicates:N \g_arg_mru_nothis
            Furthermore,~\nobs~that~
            \enum{
              #1
            }[\l_arg_verbs_tl]\ignorespaces #3		
          }
          {Ia} {
            \bool_if:NT \g_debug_bool {C9.3}
            \seq_gput_left:Nn \g_arg_mru_nothis {Ia}
            \seq_gremove_duplicates:N \g_arg_mru_nothis
            In~addition,~\nobs~that~
            \enum{
              #1
            }[\l_arg_verbs_tl]\ignorespaces #3		
          }
          {N} {
            \bool_if:NT \g_debug_bool {C9.4}
            \seq_gput_left:Nn \g_arg_mru_nothis {N}
            \seq_gremove_duplicates:N \g_arg_mru_nothis
            Next,~\nobs~that~
            \enum{
              #1
            }[\l_arg_verbs_tl]\ignorespaces #3		
          }
          {Itns} {
            \bool_if:NT \g_debug_bool {C9.5}
            \seq_gput_left:Nn \g_arg_mru_nothis {Itnswc}
            \seq_gput_left:Nn \g_arg_mru_nothis {Itns}
            \seq_gremove_duplicates:N \g_arg_mru_nothis
            In~the~next~step~we~\nobs~that~
            \enum{
              #1
            }[\l_arg_verbs_tl]\ignorespaces #3		
          }
          {Nwc} {
            \bool_if:NT \g_debug_bool {C9.6}
            \seq_gput_left:Nn \g_arg_mru_nothis {Nwc}
            \seq_gremove_duplicates:N \g_arg_mru_nothis
            Next~we~combine~
            \enum{
              #1
            }to~obtain~\ignorespaces #3
          }
          {Itnswc} {
            \bool_if:NT \g_debug_bool {C9.7}
            \seq_gput_left:Nn \g_arg_mru_nothis {Itns}
            \seq_gput_left:Nn \g_arg_mru_nothis {Itnswc}
            \seq_gremove_duplicates:N \g_arg_mru_nothis
            In~the~next~step~we~combine~
            \enum{
              #1
            }to~obtain~\ignorespaces #3
          }
        } {}
      }
    }
  }
  \bool_gset_false:N \g_arg_start_bool
  \bool_gset_false:N \l_insidearg_bool
  \cfload[.]%~\unskip
  \color{black}
}

\tl_new:N \g_label_tl
\tl_gset:Nn \g_label_tl { }

\NewDocumentCommand{\savelabel}{m}{
  \bool_if:NTF \l_insidearg_bool {
    \tl_gset:Nn \g_label_tl {#1}
  } {
    \tl_gset:Nn \g_label_tl { }
  }
}

\ExplSyntaxOff

\ExplSyntaxOn

\NewDocumentEnvironment {athm} {m m o} {
\str_if_eq:noTF {example} {#1} {
  \bool_gset_true:N \g_example_bool
} {
  \bool_gset_false:N \g_example_bool
}
\cfclear
\IfNoValueTF{#3}{
\begin{#1}\label{#2}\global\def\loc{#2}
}{
\begin{#1}[#3]\label{#2}\global\def\loc{#2}
}
}{
\end{#1}
}

\NewDocumentEnvironment {adef} {m} {
\begin{definition}\label{#1}\global\def\loc{#1}
}{
\end{definition}
}

\NewDocumentEnvironment{aproof} {} {
\bool_if:NTF \g_example_bool {
  \bool_gset_true:N \g_arg_start_bool
  \begin{proof}[Proof~for~\cref{\loc}]
} {
  \bool_gset_true:N \g_arg_start_bool
  \begin{proof}[Proof~of~\cref{\loc}]
}
\bool_gset_false:N \g_finishproof_bool
}{
\bool_if:NTF \g_finishproof_bool {}
{\finishproofthus}
\end{proof}
}

\NewDocumentCommand{\finishproofthus} {} {
  \bool_gset_true:N \g_finishproof_bool 
  \bool_if:NTF \g_example_bool {
    The~proof~for~\cref{\loc}~is~thus~complete.
  } {
    The~proof~of~\cref{\loc}~is~thus~complete.
  }
}
\NewDocumentCommand{\finishproofthis} {} {
  \bool_gset_true:N \g_finishproof_bool 
  \bool_if:NTF \g_example_bool {
    This~completes~the~proof~for~\cref{\loc}.
  } {
    This~completes~the~proof~of~\cref{\loc}.
  }
}

\ExplSyntaxOff

\NewDocumentEnvironment{cproof}{m}
{\begin{proof}[Proof of \cref{#1}]}%
{\noindent The proof of \cref{#1} is thus complete.
\end{proof}}

\NewDocumentEnvironment{cproof2}{m}
{\begin{proof}[Proof of \cref{#1}]}%
{\noindent This completes the proof of \cref{#1}.
\end{proof}}

\hypersetup{
    colorlinks,
    linkcolor={blue!80!black},
    citecolor={green},
    urlcolor={blue!80!black}
}

\makeatletter
\ExplSyntaxOn
\seq_new:N \g_abbrs
\prop_new:N \g_abbr_counts
\tl_new:N \l_abbr_count_tl

\ExplSyntaxOn

\bool_new:N \g_forexample

\NewDocumentCommand{\eg}{ o }{
	\IfValueT{#1}{
		\str_if_eq:noTF {fe} {#1} {
			\bool_gset_true:N \g_forexample
		} {\bool_gset_false:N \g_forexample}
	}
	\bool_if:nTF { \g_forexample } {
		\bool_gset_false:N \g_forexample
		for~example
	}{
		\bool_gset_true:N \g_forexample
		for~instance
	}
}

\NewDocumentCommand{\abbr}{m m O{#1} m m O{#4} m}{
	\expandafter\newcommand\csname#3\endcsname[1][]{
		\seq_if_in:NnTF \g_abbrs {#1} {
			\prop_get:NnN \g_abbr_counts {#1} \l_abbr_count_tl
			\prop_gput:Nnx \g_abbr_counts {#1} {\int_eval:n {\l_abbr_count_tl + 1}}
			\hyperref[#1]{#7}
		} {
			\seq_gput_left:Nn \g_abbrs {#1}
			\prop_gput:Nnn \g_abbr_counts {#1} {1}
			\expandafter\gdef\csname#1@def\endcsname{#2}
			\phantomsection\label{#1}
			\str_if_eq:nnTF{##1}{}{\emph{#2}}{##1}~(\hyperref[#1]{#7})
		}
	}
	\expandafter\newcommand\csname#6\endcsname[1][]{
		\seq_if_in:NnTF \g_abbrs {#1} {
			\prop_get:NnN \g_abbr_counts {#1} \l_abbr_count_tl
			\prop_gput:Nnx \g_abbr_counts {#1} {\int_eval:n {\l_abbr_count_tl + 1}}
			\hyperref[#1]{#4}
		} {
			\expandafter\gdef\csname#1@def\endcsname{#5}
			\seq_gput_left:Nn \g_abbrs {#1}
			\prop_gput:Nnn \g_abbr_counts {#1} {1}
			\phantomsection\label{#1}
			\str_if_eq:nnTF{##1}{}{\emph{#5}}{##1}~(\hyperref[#1]{#4})
		}
	}
}

\ExplSyntaxOff
\makeatother
\abbr{ODEs}{ordinary differential equations}{ODE}{ordinary differential equation}{ODEs}
\abbr{SDEs}{stochastic differential equations}{SDE}{stochastic differential equation}{SDEs}
\abbr{GD}{gradient descent}{GDs}{gradient descents}{GD}
\abbr{ANN}{artificial neural network}{ANNs}{artificial neural networks}{ANN}
\abbr{RMSprop}{root mean square propagation}{RMSprops}{the root mean square propagation}{RMSprop}
\abbr{Adam}{adaptive moment estimation}{Adams}{adaptive moment estimation}{Adam}
\abbr{AdamW}{Adam with decoupled weight decay}{AdamWs}{Adams with decoupled weight decay}{AdamW}
\abbr{Adagrad}{adaptive gradient \SGD}{Adagrads}{adaptive gradient \SGD}{Adagrad}
\abbr{Muon}{momentum orthogonalized by Newton-Schulz}{Muons}{momentum orthogonalized by Newton-Schulz}{Muon}
\abbr{SOP}{stochastic optimization problem}{SOPs}{stochastic optimization problems}{SOP}
\abbr{OP}{optimization problem}{OPs}{optimization problems}{OP}
\abbr{DNN}{deep artificial neural network}{DNNs}{deep artificial neural networks}{DNN}
\abbr{SGD}{Stochastic gradient descent}{SGDs}{stochastic gradient descent}{SGD}
\abbr{iid}{independent and identically distributed}{i.i.d}{independent and identically distributed}{i.i.d.}
\abbr{DL}{Deep learning}{DLs}{Deep learning}{DL}
\abbr{AI}{artificial intelligence}{AIs}{artificial intelligence}{AI}
\abbr{LLM}{large language model}{LLMs}{large language models}{LLM}
\abbr{PDE}{partial differential equation}{PDEs}{partial differential equations}{PDE}
\abbr{as}{almost surely}{ass}{almost surely}{a.s.}
\abbr{pas}{$\P$-almost surely}{pass}{almost surely}{$\P$-a.s.}
\abbr{ReLUs}{rectified linear unit}{ReLU}{rectified linear unit}
\makeatother

\title{
Adam symmetry theorem: characterization of\\ the convergence of the stochastic Adam optimizer
}
\author{Steffen Dereich$^{1}$, Thang Do$^{2}$, Arnulf Jentzen$^{3,4}$, \\and Philippe von Wurstemberger$^{5}$
	\bigskip
	\\
	\small{$^1$ Institute for Mathematical Stochastics, Faculty of Mathematics and Computer Science,}\vspace{-0.1cm}\\
\small{University of M\"unster, Germany, e-mail: \texttt{steffen.dereich@uni-muenster.de}}
\smallskip
\\
    	\small{$^2$ School of Data Science, The Chinese University of Hong Kong, Shenzhen}
	\vspace{-0.1cm}\\
	\small{ (CUHK-Shenzhen), China, e-mail: \texttt{minhthangdo@link.cuhk.edu.cn}}
 \smallskip
	\\
	\small{$^3$ School of Data Science and School of Artificial Intelligence, The Chinese University}
	\vspace{-0.1cm}\\
	\small{of Hong Kong, Shenzhen (CUHK-Shenzhen), China, e-mail: \texttt{ajentzen@cuhk.edu.cn}}
	\smallskip
	\\
 \small{$^4$ Applied Mathematics: Institute for Analysis and Numerics, Faculty of Mathematics and}
	\vspace{-0.1cm}\\
	\small{Computer Science, University of M{\"u}nster, Germany, e-mail: \texttt{ajentzen@uni-muenster.de}}
	\smallskip
	\\
    \small{$^5$ School of Data Science, The Chinese University of Hong Kong, Shenzhen}
	\vspace{-0.1cm}\\
	\small{ (CUHK-Shenzhen), China, e-mail: \texttt{philippevw@cuhk.edu.cn}}
 \smallskip
	\\
}

\date{\today}
\allowdisplaybreaks 
\begin{document}
\maketitle
\newpage
\begin{abstract}
    Beside the standard stochastic gradient descent (SGD) optimization method, the Adam optimizer due to Kingma and Ba (2014) is currently probably the best-known optimization method for the training of deep neural networks in artificial intelligence (AI) systems. For the standard SGD method a rigorous convergence analysis is already known for a long time, at least for strongly convex stochastic optimization problems (SOPs). For the Adam optimizer, in turn, it remains -- despite the popularity and the success of Adam -- an \emph{open research problem} to provide a rigorous convergence analysis even for the class of strongly convex SOPs.
    
 In one of the main results of this work we establish convergence rates for Adam in terms of the number of gradient steps (convergence rate \nicefrac{1}{2} with respect to the size of the learning rate), the size of the mini-batches (convergence rate 1 with respect to the size of the mini-batches), and the size of the second moment parameter of Adam (convergence rate 1 with respect to the distance of the second moment parameter to one) for the class of strongly convex SOPs. 
 
 In a further main result of this work, which we refer to as \emph{Adam symmetry theorem}, we illustrate the optimality of the established convergence rates by proving for a special class of simple quadratic strongly convex SOPs that Adam converges as the number of gradient steps increases to infinity to the solution of the SOP (the unique minimizer of the strongly convex objective function) if and \emph{only} if the random variables in the SOP (the data in the SOP) are \emph{symmetrically distributed}. In particular, in the standard case where the random variables in the SOP are not symmetrically distributed we \emph{disprove} that Adam converges to the minimizer of the SOP as the number of Adam steps increases to infinity. 
 
 We also complement the conclusions of our convergence analysis and the Adam symmetry theorem by several numerical simulations that indicate the sharpness of the established convergence rates and that illustrate the practical appearance of the phenomena revealed in the \emph{Adam symmetry theorem}.
 \end{abstract}
\pagebreak
\tableofcontents
\newpage
\section{Introduction}
\renewcommand{\fd}{d}
\SGD\ optimization schemes are currently the standard methods to train deep \ANNs\ for \AI\ based systems. However, in many situations one does not employ the standard \SGD\ method as the training algorithm to minimize the risk function of the learning problem but instead suitable accelerated and adaptive variants of the standard \SGD\ method are considered \cite{ArBePhi2024,Ruder2016AnOO}. The most popular variant of such accelerated and adaptive \SGD\ optimization methods is presumably the \Adam\ optimizer proposed in 2014 by Kingma \& Ba \cite{KingmaBa2014_Adam}. 
In particular, \Adam\ and an extended variant of \Adam\ are among the most widely used optimizers in the training of state-of-the-art \LLMs\ such as DeepSeek-V3 \cite{LLMTable241219437} and OpenAI's GPT-3 \cite{LLMTable200514165}. We refer to \cref{tab:llm-optimizers} in \cref{sec:adam_training_llms} below for a detailed list of the optimizers employed to train \LLMs.

Despite the popularity and success of the \Adam\ optimizer in numerical computations, it remains, except of in a few special cases (cf., \eg, \cite[Theorem 1.1]{DereichAdamconvergence2024}), a fundamental open research problem to develop a convergence analysis for \Adam\ even in the setting of strongly convex \SOPs\ in which case essentially sharp convergence rates for the standard \SGD\ method are for a long time well-known (cf., \eg, \cite{Guro2023handbook,Jentzen_2020} and the references therein). The \Adam\ optimizer is, however, a much more advanced and complicated optimization algorithm than the standard \SGD\ method that involves averages of gradients evaluated at the previous steps of the optimization process (momentum \cite{Polyak1964SomeMO}) as well as adaptive modifications of the learning rate schedule (\RMSprop\ \cite{Hinton24_RMSprop}) which seems to make the error analysis for \Adam\ much more subtle than that of the standard \SGD\ method. The issue is also illustrated by the fact that it is has recently, in a very wide generality, been established (see \cite[Corollary 1.10]{DereichJentzenKassing2025} and, \eg, \cref{theo: convergence point of Adam copy 2} in \cref{sec: converge to wrong point} below)
that if the \Adam\ optimizer converges almost surely or in probability to a random variable, then this random limit point must be almost surely a zero of the \Adam\ vector field function that has been proposed in \cite[Definition 2.4]{DereichAdamconvergence2024} (cf., \eg, also \cref{theo: convergence point of Adam copy 3 pre} in \cref{sec: converge to wrong point} below). 

In this work we make use of the works \cite{DereichAdamconvergence2024,DereichJentzenKassing2025} to deliver a partial answer to the above described fundamental research problem for the \Adam\ optimizer. More specifically, for the general class of strongly convex \SOPs\ we develop in this work -- under the assumption that the \Adam\ optimizer does not blow up but stays bounded -- a comprehensive convergence analysis with convergence rates (see \cref{main theorem}, \cref{cor: convegence of Adam}, and \cref{cor: convegence of Adam2} below). In particular, we establish an upper bound for the overall strong $L^p$-distance with $p \in [1,\infty)$ between the \Adam\ optimization process and the global minimizer of the strongly convex \SOP\ that decays to zero with convergence rates as the size of the mini-batches is sufficiently large (convergence rate 1), as the second moment parameter in \Adam\ is sufficiently close to one (convergence rate 1), and as the number of \Adam\ steps is sufficiently large (convergence rate $\nicefrac{1}{2}$ with respect to the size of the learning rate); see \cref{main theorem} in \cref{subsec: main theorem 1} below. 

The well-known error analysis for the standard \SGD\ optimization method for strongly convex \SOPs\ is much less sophisticated and provides an upper bound for the strong $L^p$-distance with $p \in [1,\infty)$ between the standard \SGD\ optimization process and the unique global minimizer of the \SOP\ that decays to zero with convergence rates as the number of \SGD\ steps is sufficiently large (convergence rate \nicefrac{1}{2} with respect to the size of the learning rate); cf., \eg, \cite{Guro2023handbook,Jentzen_2020,MR4055054} and the references therein. In our first main result, which we refer to as \emph{\Adam\ symmetry theorem}, we reveal for a class of simple quadratic strongly convex \SOPs\ that in nearly all cases it is impossible to establish such a simple convergence analysis for the \Adam\ optimizer. Specifically, in \cref{main theorem 2} in \cref{subsec: main theorem 2} below we reveal for the considered class of \SOPs\ that \Adam\ converges to the solution of the \SOP\ as the number of \Adam\ steps increase to infinity if and \emph{only if} the stochastic data in the \SOP\ is \emph{\textbf{symmetrically} distributed}\footnote{Note that for every probability space $(\Omega,\cF,\P)$, every $d \in \N$, every $A \subseteq \R^d$, and every random variable $X \colon \Omega \to A$ it holds that $X$ is symmetrically distributed if and only if there exists $x \in \R^d$ such that for every measurable $B \subseteq \R^d$ it holds that $\P( X - x \in B ) = \P( x - X \in B )$.}.

\subsection{First main result: Adam symmetry theorem -- wrong limiting values of Adam}\label{subsec: main theorem 2}
The natural number $d \in  \N = \{1, 2, 3, \dots \}$ in \cref{main theorem 2} specifies the dimension of the considered \SOP\ (the number of parameters/degrees of freedom that need to be optimized) and the \iid\ random variables $X_{n,m}$, $(n, m) \in  \N^2$, represent the data in the considered \SOP. 

In \cref{main theorem 2} we prove that in the situation of a class of simple quadratic \SOPs\ we have that \Adam\ converges to the solution of the \SOP\ (the minimizer of the \SOP), which is well-known for standard \SGD, if and \emph{only} if the external data $X_{ n, m }$, $(n,m) \in \N^2$, in the \SOP\ is \emph{symmetrically distributed} (see \cref{item 4: main theorem 2} in \cref{main theorem 2} below). In particular, \cref{main theorem 2} discovers that in nearly all cases of external data (namely, all cases of not symmetrically distributed data) items \ref{item 2: main theorem 2} and \ref{item 1: main theorem 2} of \cref{main theorem 2} fail to hold and the \Adam\ optimizer does \emph{not} converge to the unique global minimizer $\globalmin$ of the \SOP. 
\renewcommand{\mandV}[3]{\Theta^{#1,#2}_{#3}}
\renewcommand{\mandVcom}[4]{\Theta^{#1,#2}_{#3,#4}}
\newpage
\begin{samepage}
\begin{tcolorbox}[colback=white!95!gray,
                  colframe=black,
                  boxrule=0.5pt,
                  sharp corners,
                  enhanced,
                  breakable,
                 ]
\begin{athm}{theorem}{main theorem 2}[\textcolor{red}{Adam symmetry theorem}]
     Let $(\Omega,\cF,\P)$ be a probability space, let $\fd\in \N$, $\xi,v,w\in \R^\fd$, let $X_{n,m}\colon \Omega\allowbreak\to \{v,w\}$, $(n,m) \in \N^2$, be \iid\ random variables, let $\smalll=(\smalll(\theta,x))_{(\theta,x)\in \R^\fd\times\R^{\fd}}\allowbreak\colon\R^\fd\times\R^{\fd}\to\R$ satisfy for all $\theta,x\in \R^\fd$ that $\smalll(\theta,x)=\|\theta-x\|^2$, let $\para_1\in (0,1)$, $\para_2\in ((\para_1)^2,1)$, $\varepsilon,p\in (0,\infty)$, let $(\gamma_n)_{n\in \N}\subseteq (0,\infty)$ be non-increasing, for every $k,M,n\in \N_0$ let 
    $
\mandV{k}{M}{n}=(\mandVcom{k}{M}{n}{1},\dots,\mandVcom{k}{M}{n}{\fd})\colon\allowbreak \Omega\to\R^\fd$ be a random variable, assume for all $k\in \{1,2\}$, $M,n\in \N$, $i\in \{1,2,\dots,\fd\}$ that
    \begin{gather}
    \textcolor{magenta}{\mandV{k}{M}{0}=0},\qquad \textcolor{magenta}{\mandVcom{k}{M}{n}{i}= \beta_k \mandVcom{k}{M}{n-1}{i}+(1-\beta_k)\bigl[\textstyle \frac 1M \sum_{m=1}^M(\nabla_{\theta_i} \smalll)(\mandV{0}{M}{n-1},X_{n,m})\bigr]^k}, \label{eq2: main theorem 2}\\
\textcolor{magenta}{\mandV{0}{M}{0}=\xi},\ \text{and}\ \textcolor{magenta}{\mandVcom{0}{M}{n}{i}=\mandVcom{0}{M}{n-1}{i} -\gamma_n [1-(\para_1)^n]^{-1}\bigl[\varepsilon+\bigl[[1-(\para_2)^n]^{-1}\mandVcom{2}{M}{n}{i}\bigr]^{\nicefrac{1}{2}}\bigr]^{-1}\mandVcom{1}{M}{n}{i}}\notag,
     \end{gather}
     assume $\limsup_{n\to\infty}(( \gamma_n )^{ - 2 } ( \gamma_n - \gamma_{ n + 1 } ) +\sum_{m=n}^\infty(\gamma_m)^p) =0$, and let $\vartheta \in \R^{\fd}$ satisfy $\E[\smalll(\vartheta,X_{1,1})]= \inf_{ \theta \in \R^\fd } \E[\smalll(\theta,X_{1,1})]$ \cfload.
      Then\footnotemark\ the following statements are \textbf{equivalent}:
     \begin{enumerate}[label=(\roman*)]
      \item \label{item 2: main theorem 2} There exists 
         \textcolor{magenta}{$\scrc \in \R$} such that for all \textcolor{magenta}{$M,n \in 
         \N$} it holds that \textcolor{magenta}{$\bigl( \E\bigl[ \| \mandV{0}{M}{n}- \vartheta \|^p \bigr] \bigr)^{ 1 / p } \leq \scrc \sqrt{ \gamma_n }$}.
     \item \label{item 1: main theorem 2} It holds for all \textcolor{magenta}{$M\in \N$} that \textcolor{magenta}{$\limsup_{ n \to \infty }\E\bigl[ \min\{ 1, \| \mandV{0}{M}{n} - \vartheta\|\} \bigr] = 0$}.
         \item \label{item 4: main theorem 2} It holds that \textcolor{magenta}{$X_{ 1, 1 }$} is \textcolor{red}{\textbf{symmetrically}} distributed.
     \end{enumerate}
\end{athm} 
\end{tcolorbox}
\end{samepage}
\footnotetext{Note that for all $n \in \N$, $v = ( v_1, \dots, v_n ) \in \R^n$ it holds that $\| v \| = [\sum_{ i = 1 }^n | v_i |^2]^{ 1 / 2 }$.}
\renewcommand{\mandV}[4]{\mathbb{V}^{#1,#2,#3}_{#4}}
\renewcommand{\mandVcom}[5]{\mathbb{V}^{#1,#2,#3}_{#4,#5}}
\cref{main theorem 2} is a direct consequence of \cref{main theorem 2'} in \cref{sec: convergence symmetric} below. The triple $( \Omega, \cF, \P )$ in \cref{main theorem 2} is the underlying probability space on which the data $X_{ n,m }$, $(n,m) \in \N^2$, of the considered \SOP\ is defined. 
The function
\begin{equation}\label{eq: introduction ***1}
  \R^d \ni \theta \mapsto \E[\smalll(\theta,X_{1,1})] \in \R
\end{equation}
in \cref{main theorem 2} is the objective function (the function that we intend to minimize in \cref{main theorem 2}) and the vector $\vartheta \in \R^d$ in \cref{main theorem 2} is the unique global minimizer of the function in \cref{eq: introduction ***1}. Specifically, in \cref{main theorem 2} we intend to employ the \Adam\ optimizer to approximately compute the unique global minimizer $\globalmin$ of the \SOP
\begin{equation}\label{eq: loss1}
  \textstyle\min_{ \theta \in \R^d } \E\bigl[ \smalll( \theta, X_{1,1} ) \bigr].
\end{equation}
In  \cref{eq2: main theorem 2} for every size of the mini-batches $M \in \N$ the \Adam\ optimization process $(\Theta^{0,M}_n )_{n\in\N_0}$ is specified (cf.\ \cite{KingmaBa2014_Adam} and, \eg, \cite[Definition 7.9.1]{ArBePhi2024}). Next we observe that
\begin{itemize}
\item the real number $\beta_1 \in (0,1)$ in \cref{main theorem 2} represents the momentum parameter for the \Adam\ optimizer,
\item the real number $\beta_2 \in (0,1)$ in \cref{main theorem 2} represents the second moment parameter for the \Adam\ optimizer,
\item the real number $\varepsilon\in (0,\infty)$ represents the regularization parameter in the \Adam\ optimizer that ensures that one does not divide by zero in \cref{eq2: main theorem 2}, and
\item the sequence of real numbers $\gamma_n \in (0,\infty)$, $n \in \N$, in \cref{main theorem 2} represents the sequence of learning rates/step sizes of the \Adam\ optimizer. 
\end{itemize}
In \cref{main theorem 2} we employ the popular assumption that $( \beta_1 )^2 < \beta_2$. This common assumption is not really restrictive as the common choices for the parameters of \Adam\ $\beta_1 = 0.9$, $\beta_2 = 0.95$, $\varepsilon = 10^{ - 1 }$ (cf.\ \cref{tab:llm-optimizers} below) and $\beta_1 = 0.9$, $\beta_2 = 0.999$, $\varepsilon= 10^{-8}$ (cf., \eg, \cite[Algorithm 1]{KingmaBa2014_Adam}, \cite{pytorch2025adam}, and \cite{tensorflow2025adam}) clearly satisfy this condition $( \beta_1 )^2 = ( 0.9 )^2 = 0.81 <\beta_2$ (these common choices for the parameters of \Adam\ even satisfy the condition $\beta_1 < \beta_2$). Furthermore, in \cref{main theorem} we also impose the assumption that
  $\limsup_{n\to\infty}( ( \gamma_n )^{ - 2 } ( \gamma_n - \gamma_{ n + 1 } )+ \sum_{m=n}^\infty(\gamma_m)^p) =0$. We note that this assumption is equivalent to the common assumptions that
\begin{equation}\label{assumption: gamma1}
\textstyle \sum_{n=1}^{ \infty } ( \gamma_n)^p < \infty\qqandqq\textstyle \limsup_{n\to\infty} (( \gamma_n )^{ - 2 } ( \gamma_n - \gamma_{ n + 1 } ) )  = 0.
\end{equation}
Loosely speaking, these well-known assumptions ensure that the learning rates converge to zero but do not converge sufficiently quickly to zero. These assumptions have been employed in the literature to analyze \SGD\ optimization methods (cf., \eg, \cite[Theorem 1.1]{DereichAdamconvergence2024}). We note that the assumption that $\textstyle \limsup_{n\to\infty} (( \gamma_n )^{ - 2 } ( \gamma_n - \gamma_{ n + 1 } ) )  = 0$ in \cref{assumption: gamma1} ensures that the learning rates fail to be summable in the sense that $\sum_{ n = 1 }^{ \infty } \gamma_n = \infty$ (see, \eg, \cref{lem: unsummable} in \cref{subsec: convergence with sufficient batch} below). Moreover, we observe that for every $c, r \in (0,\infty)$ it holds that the sequence of real numbers $c n^{ - r } \in \R$, 
$n\in \N$, satisfies the condition that $\textstyle \limsup_{n\to\infty} (( \gamma_n )^{ - 2 } ( \gamma_n - \gamma_{ n + 1 } ) )  = 0$ in \cref{assumption: gamma1} if and only if $r <1$ (see, \eg, \cref{lem: verify 3} in \cref{subsec: convergence with sufficient batch} below).

\cref{main theorem 2} proves that in many cases it is not possible to establish convergence of \Adam\ to the minimizer of the \SOP\ in the number of gradient steps as for the standard \SGD\ method (see items \ref{item 2: main theorem 2} and \ref{item 4: main theorem 2} in \cref{main theorem 2}) but instead we have in many cases that \Adam\ fails to converge to the minimizer of the \SOP\ while the standard \SGD\ method does, of course, converge to the minimizer of the \SOP\ as the number of gradient steps increases to infinity. Nonetheless, in \cref{main theorem} in \cref{subsec: main theorem 1} below we reveal another error analysis for \Adam\ which implies convergence as the number of \Adam\ steps increases to infinity, as the size of the mini-batches increase to infinity, and as the second moment parameter of \Adam\ increases to one, respectively (and due to \cref{main theorem 2} we know that in many cases the mini-batch and the second moment parameter terms, respectively, can not be omitted).

\subsection{Second main result: Convergence rates for the stochastic Adam optimizer}\label{subsec: main theorem 1}
The natural number $d \in \N $ in \cref{main theorem} specifies the dimension of the considered \SOP\ (the number of parameters/degrees of freedom that need to be optimized), the natural number $\dimX \in \N$ specifies the dimension of data of the considered \SOP, and the \iid\ random variables $X_{ n, m } \colon \Omega \to\R^\dimX$, $(n,m) \in \N^2$, represent the data in the considered \SOP. We now present \cref{main theorem} in a completely self-contained way with all mathematical details and, thereafter, we explain some of the mathematical objects in and the conclusion of \cref{main theorem} in words.
\begin{samepage}
\renewcommand{\mandV}[4]{\Theta^{#1,#2,#3}_{#4}}
\renewcommand{\mandVcom}[5]{\Theta^{#1,#2,#3}_{#4,#5}}
\begin{tcolorbox}[colback=white!95!gray, colframe=black, boxrule=0.5pt,
                  sharp corners, enhanced]
\begin{athm}{theorem}{main theorem}[\textcolor{red}{Adam error bounds}]
    Let $(\Omega,\cF,\P)$ be a probability space, let $\fd,\dimX,p\in \N$, let $X_{n,m}\colon \Omega\allowbreak\to \R^{\dimX}$, $(n,m) \in \N^2$, be bounded \iid\ random variables, let $\smalll=(\smalll(\theta,x))_{(\theta,x)\in \R^\fd\times\R^{\dimX}}\allowbreak\in C^2(\R^\fd\times\R^{\dimX},\R)$, $\varepsilon\in (0,\infty)$, let $(\gamma_n)_{n\in \N}\subseteq (0,\infty)$ be non-increasing, for every $k,M,n\in \N_0$, $\para\in \R^2$ let 
    $
\mandV{k}{M}{\para}{n}=(\mandVcom{k}{M}{\para}{n}{1},\dots,\mandVcom{k}{M}{\para}{n}{\fd})\colon\allowbreak \Omega\to\R^\fd$ be a random variable, assume for all $k\in \{1,2\}$, $M,n\in \N$, $i\in \{1,2,\dots,\fd\}$, $\para=(\para_1,\para_2)\in (0,1)^2$ that
    \begin{equation}\label{eq1: main theorem}
   \textcolor{magenta}{ \mandV{k}{M}{\para}{0}=0},\quad \textcolor{magenta}{\mandVcom{k}{M}{\para}{n}{i}= \beta_k \mandVcom{k}{M}{\para}{n-1}{i}+(1-\beta_k)\bigl[\textstyle \frac 1M \sum_{m=1}^M(\nabla_{\theta_i} \smalll)(\mandV{0}{M}{\beta}{n-1},X_{n,m})\bigr]^k},
        \end{equation}
        \begin{equation}\label{eq2: main theorem}
      \text{and} \ \textcolor{magenta}{\mandVcom{0}{M}{\para}{n}{i}=\mandVcom{0}{M}{\para}{n-1}{i} -\gamma_n [1-(\para_1)^n]^{-1}\bigl[\varepsilon+\bigl[[1-(\para_2)^n]^{-1}\mandVcom{2}{M}{\para}{n}{i}\bigr]^{\nicefrac{1}{2}}\bigr]^{-1}\mandVcom{1}{M}{\para}{n}{i}}, 
     \end{equation}
     assume 
       $\limsup_{n\to\infty}( ( \gamma_n )^{ - 2 } ( \gamma_n - \gamma_{ n + 1 } )+ \sum_{m=n}^\infty(\gamma_m)^p) =0$, assume for all $M\in \N$, $\para\in \R^2$ that $\mandV{0}{M}{\beta}{0}$ and $(X_{n,m})_{(n,m)\in \N^2}$ are independent, assume that $\R^d\ni \theta\mapsto \E[\smalll(\theta,X_{1,1})]\in \R$ is \stronglyconvexx, and let $\globalmin \in \R^\fd$ satisfy $\E[\smalll(\globalmin,X_{1,1})]=\inf_{\theta\in \R^\fd}\E[\smalll(\theta,X_{1,1})]$. Then there exist $(\scrc_{\para})_{\para\in \R^2}\subseteq \R$ such that for all $M,n \in \N$, $\para = ( \para_1, \para_2 ) \in (0,1) \times [ p^{ - 1 }, 1 )$ with $( \para_1 )^2 < \para_2$ and $\inf_{N,m\in \N}\P\bigl(\|\mandV{0}{N}{\beta}{m}\|\leq p\bigr)=1$ it holds that
\begin{equation} \label{conclude: main theorem}
     \textcolor{magenta}{\bigl( \E\bigl[ \| \mandV{0}{M}{\beta}{n} - \globalmin\|^p \bigr] \bigr)^{ 1/p } \leq  \scrc_0 M^{ - 1 } ( 1 - \beta_2 \mathbbm 1_{ [\scrc_0,\infty)} ( M ) ) + \scrc_{ \beta } \sqrt{\gamma_n}}.
\end{equation}
\end{athm}
\end{tcolorbox}
\end{samepage}
\cref{main theorem} is a direct consequence of \cref{cor: convegence of Adam2} in \cref{sec: convergence general} below. \cref{cor: convegence of Adam2}, in turn, follows from an application of \cref{cor: convegence of Adam}, which is one of the main results of this work. The twice continuously differentiable function $\smalll \colon \R^d\times\R^\dimX \to\R$ in \cref{main theorem} represents the loss function of the \SOP\ that we intend to solve in \cref{main theorem}, the strongly convex function
\begin{equation}\label{eq: introduction ***}
  \R^d \ni \theta \mapsto \E\bigl[\smalll(\theta,X_{1,1})\bigr] \in \R
\end{equation}
in \cref{main theorem} represents the objective function, and the vector $\vartheta \in \R^d$ in \cref{main theorem} is the unique global minimizer of the function in \cref{eq: introduction ***}. 
Strong convexity of the objective function in \cref{eq: introduction ***} means here that there exists a strictly positive real number $c \in (0,\infty)$ such that the function $\R^d \ni \theta \mapsto \E[ \smalll( \theta, X_{ 1, 1 } ) ]-c\|\theta\|^2\in \R$ is convex (see, \eg, \cref{item 2: strongly convex} in \cref{definition: strongly convex} in \cref{sec: convergence general} below).

In \cref{eq1: main theorem} and \cref{eq2: main theorem} for every size of the mini-batches $M \in \N$ and every first and second moment parameter vector $\beta = ( \beta_1,\beta_2 ) \in (0,1)^2$ the \Adam\ optimization process $( \Theta^{ 0, M, \beta }_n )_{ n \in \N_0 }$ is specified. In \cref{main theorem} we also employ the assumption that the norm of the \Adam\ optimization process does not diverge to infinity but stays bounded. In the situation of certain simple quadratic \SOPs\ it is known that the \Adam\ optimization process does indeed satisfy this assumption; see \cite{DereichGraArAd2025}.
 
In \cref{conclude: main theorem} in \cref{main theorem} it is revealed that the overall strong $L^p$-error between the unique global minimizer $\vartheta$ of the function in \cref{eq: introduction ***} and the \Adam\ optimization process $\Theta_n^{ 0, M, \beta }$ at step $n$ is as small as desired as the number of \Adam\ steps $n \in \N$ is sufficiently large (convergence order $\nicefrac{1}{2}$ with respect to the size of the learning rate), as the size of the mini-batches $M \in  \N$ is sufficiently large (convergence order 1 with respect to the size of the mini-batches), and as the second moment parameter $\beta_2 \in (0,1)$ is sufficiently close to one (convergence order 1 with respect to the distance of the second moment parameter to one), respectively. Not for \Adam\ but in the situation of the standard \SGD\ method applied to strongly convex \SOPs\ a strengthened version of \cref{conclude: main theorem} without the additional term $\scrc_0 M^{ - 1 } ( 1 - \beta_2 \mathbbm 1_{ [\scrc_0,\infty)} ( M ) )$ in \cref{conclude: main theorem} is well-known. More formally, in the setup of \cref{main theorem} but for the standard \SGD\ method $\varTheta^M_n \colon \Omega \to \R^d$, $(n,M) \in \N_0 \times \N$, it is well-known that there exists a real number $\scrc \in \R$ such that for all $M, n \in \N$ it holds that  
\begin{equation}\label{convegence of sgd}
\bigl( \E\bigl[ \| \varTheta^M_n - \globalmin\|^p \bigr] \bigr)^{ \nicefrac{1}{p} } 
\leq
\scrc\sqrt{ \gamma_{n} }.
\end{equation}
A natural question is whether an estimate of the form \cref{convegence of sgd}, which is strictly stronger than \cref{conclude: main theorem}, can also be established for the \Adam\ optimizer. It is precisely one of the consequences of \cref{main theorem 2} that this can \emph{not} be achieved.

In implementations the size of the mini-batches $M \in \N$ often satisfies $M = 32$ or $M = 64$ and the second moment parameter of \Adam\ $\beta_2 \in (0, 1)$ often satisfies $\beta_2 = 0.95$ (cf.\ \cref{tab:llm-optimizers} below) while $\gamma_n$ is often in between $6\times 10^{-7}$ and $10^{ -3 }$ (cf., \eg, \cref{tab:llm-optimizers} below, \cite[Algorithm 1]{KingmaBa2014_Adam}, \cite{pytorch2025adam}, and \cite{tensorflow2025adam}). Taking the fact that in the case $M = 32$ and $\beta_2 = 0.95$ we have that 
\begin{equation}
  \textstyle M^{ - 1 } ( 1-\beta_2) = (32)^{ - 1 } ( 1 - 0.95 ) =  (32)^{ - 1 } (20)^{ - 1 } = (640)^{ - 1 }
  = 
  \sqrt{ ( 409600 )^{ - 1 } }
\end{equation}
into account indicates that the first error term $M^{ - 1 } ( 1 -\beta_2 \mathbbm 1_{[\scrc_0,\infty)}(M) )$ on the right hand side of \cref{conclude: main theorem} is -- in the situation where the learning rate $\gamma_n$ is great than $\frac{ 1 }{ 409 600 }$ -- smaller than the second error term $\sqrt{ \gamma_n }$ on the right hand side of \cref{conclude: main theorem}. In particular, we note that \cref{conclude: main theorem} in \cref{main theorem} suggests that $n \in \N$ needs to be so large so that the learning rate $\gamma_n$ is smaller or of size $M^{ - 2 } ( 1 - \beta_2 )^2$ (which is equal to $32^{ - 2 } ( 1 - 0.95 )^2 = 2^{ - 10 } 20^{ - 2  } = 2^{ - 12 } 10^{ - 2 } = \frac{ 1 }{ 409600 }$ in the case $M = 32$ and $\beta_2 = 0.95$) so that
the first error term $ M^{ - 1 } ( 1 -\beta_2 \mathbbm 1_{[\scrc_0,\infty)}(M) )$ becomes relevant/dominant and can not be neglected anymore.

To sum it up, it is a key contribution of this work to reveal that in the convergence analysis for the \Adam\ optimizer there appears an additional additive term $M^{ - 1 } ( 1 -\beta_2 \mathbbm 1_{[\scrc_0,\infty)}(M) )$ in the error bound for \Adam\ (cf.\ \cref{main theorem}) that can in many cases not be omitted (cf.\ \cref{item 4: main theorem 2} in \cref{main theorem 2}) but this additional term $M^{ - 1 } ( 1 -\beta_2 \mathbbm 1_{[\scrc_0,\infty)}(M) )$ is often (but not always) smaller than the standard error term $\sqrt{ \gamma_n }$ (that also appears in the error bound for the standard \SGD\ method; cf.\ \cref{convegence of sgd}) so that the error term $\sqrt{ \gamma_n }$ dominates the additional error term $M^{ - 1 } ( 1 -\beta_2 \mathbbm 1_{[\scrc_0,\infty)}(M) )$ unless the learning rate $\gamma_n$ is quite small (for instance, of size $( 32 )^{ - 2 } ( 1 - 0.95 )^{ 2 } = 2^{ - 10 } 20^{ - 2  } = 2^{ - 12 } 10^{ - 2 } = \frac{ 1 }{ 409 600 }$).

We also refer to \cref{sec: numerial simulation} below for several numerical simulations that indicate the sharpness of the convergence rates established in \cref{main theorem} and that illustrate the practical appearance of the phenomena revealed in \cref{main theorem 2}.
\subsection{Literature review}
In this subsection we briefly review selected further results 
from the literature that theoretically analyze the \Adam\ optimizer 
and related \SGD\ optimization methods. 

 \subsubsection*{Upper bounds for Adam and related optimizers for stochastic optimization problems}

We refer, \eg, to \cite{Barakattopological2021,DereichAdamconvergence2024,ar2112.03459,ar2104.14840,Alexconvergence2025,ar2304.13972,ReddiKale2019,ar1808.03408,ar2403.15146,ar2208.09900,ar2404.01436,ZhangChen2022,ZouShen2019} for works providing upper bounds for the
error of \Adam\ (cf., \eg, \cite[Theorem 2.5]{DereichAdamconvergence2024}), for the regret of \Adam\ (cf., \eg, \cite[Section 2]{ReddiKale2019}), or for minima or averages of the norm of the gradient of the objective function
evaluated along the \Adam\ optimization process (cf., for instance, \cite[Theorem 1]{ar2112.03459}, \cite[Theorem 6]{ar2104.14840}, \cite[Theorem 1]{ar1808.03408}, \cite[Theorem 3]{ar2403.15146}, \cite[Theorem 4.1]{ar2208.09900}, \cite[Theorem 2]{ar2404.01436}, and \cite[Theorem 4 and Theorem 5]{ZouShen2019}) when applied to \SOPs.

For example, Theorem 4 in \cite{ZouShen2019} provides an upper bound for the $\nicefrac{4}{3}$-moment of the norm of the gradient of the objective function (the function one intends to minimize) evaluated at a random time of a generalized \Adam\ optimization process (cf., \eg, also \cite[Theorem 1]{ar1808.03408}). The random time is uniformly distributed on the set $\{ 1, 2, \dots, n \}$ where $n \in \N$ is the number of considered gradient/optimization steps. In the situation of \Adam, the upper bound in \cite[Theorem 4]{ZouShen2019} essentially contains the term $( 1 - \beta_2 )^{ 1 / 2 }$ multiplied by a constant where $\beta_2 \in (0,1)$ is the second moment parameter of \Adam\ in \cref{eq2: main theorem 2} and \cref{eq2: main theorem}, respectively. 

It should also be pointed out that, in the situation of strongly convex \OPs, we have that upper bounds for the gradient of the objective function evaluated at the optimization process imply upper bounds for the error of the optimization. Indeed, if the objective function $ \cL \in C^1( \R^{ \fd }, \R )$ is strongly convex with the strong convexity constant $c \in (0,\infty)$ (cf., \eg, \cref{item 1: strongly convex,item 2: strongly convex} in \cref{definition: strongly convex} below), then there exists a unique global minimizer $\vartheta \in \R^{ \fd }$ of $\cL$ and at $\vartheta$ the coercivity-type condition that for all $\theta \in \R^{ \fd }$ it holds\footnote{Note that for all $n \in \N$, $v = ( v_1, \dots,v_n )$, $w = ( w_1,\dots,w_n ) \in \R^n$ it holds that $\spro{v, w} = \sum_{i=1}^nv_iw_i$.} that
\begin{equation}\label{eq: strongly convex introduction}
  \spro{ \theta-\vartheta, \cL(\theta)-\cL(\vartheta)} \geq c\| \theta - \vartheta\|^2
\end{equation}
is fulfilled (see, \eg, \cite[Corollary 5.7.27]{ArBePhi2024}). Applying the Cauchy-Schwarz inequality to the left hand side of \cref{eq: strongly convex introduction}, in turn, implies that for all $\theta \in \R^{ \fd }$ it holds that 
\begin{equation}\label{eq: strongly convex introduction 2}
  \| \theta - \vartheta \| \leq c^{ - 1 } \| ( \nabla \cL )( \theta ) \|.
\end{equation} 
In the strong convex setting we thus have that upper bounds for the norm of the gradient of the objective function also imply upper bounds for the standard norm distance to the unique global minimizer of the objective function.

Some of the above cited works study \Adam\ when the hyperparameters $\beta_1$ and $\beta_2$, respectively, are not fixed (as in the default setting $\beta_1 = 0.9$ and $\beta_2 = 0.999$;  see \cite[Algorithm 1]{KingmaBa2014_Adam} and, \eg, \cite{pytorch2025adam,tensorflow2025adam}) but when they increase to 1 with an appropriate scaling/speed of convergence (cf., \eg, \cite[Theorem 5.2]{Barakattopological2021}, \cite[Theorem 2]{ar2112.03459}, and \cite[Theorem 7]{ar2104.14840}).

Many works also study \Adam\ when the objective function is given as a normalized sum of deterministic functions. In particular, the work \cite[Theorem 4.1]{ar2208.09900} basically provides upper bounds for the minimum over all \Adam\ steps until time $n$ of the norm of the gradient of the objective function evaluated at the \Adam\ optimization process under the so-called $(L_0 , L_1 )$-smoothness condition \cite[Definition 1]{ar1905.11881} (which is essentially a linear growth bound on the Hessian of the summands in the objective function with constants $L_0$ and $L_1$ \cite[Assumption 3.1]{ar2208.09900}) and under a suitable linear growth bound on the noise in the sense that, roughly speaking, the sum of the squared norms of the gradients of the summands in the objective function can be bounded by the sum of a constant and a constant multiplied by the squared norm of the gradient of the objective function \cite[Assumption 3.1]{ZhangChen2022}\,(cf., \eg, also \cite[Assumption 2.2]{ar2208.09900}).

We also point, \eg, to \cite{Barakattopological2021,ar2509.21614,DereichAdamconvergence2024,DereichJentzenKassing2025,MR4795573,MR4253710,ar1511.06251,MR3948080,ar2205.10287}  for works that study continuous-time limits of \SGD\ optimization methods in the form of \ODEs\ \cite{Barakattopological2021,DereichAdamconvergence2024,DereichJentzenKassing2025,MR4253710} (cf., \eg, also \cite{ar2009.06125}) and \SDEs\ \cite{ar2509.21614,MR4795573,ar1511.06251,MR3948080,ar2205.10287} and that establish upper bounds for the distance between the considered \SGD\ optimization method
and the solution process of the associated continuous-time differential equation.

We refer, \eg, to \cite{KuNaKim2023adamfamily,NaXiXinKim2023adamfamily} for works that analyze \Adam\ type optimization methods when applied to non-smooth objective functions. Non-smooth objective functions arise, for example, in the training of \ANNs\ with non-smooth activation functions such as the \ReLU\ $\R \ni x \mapsto \max\{ x, 0 \} \in \R$ and the leaky \ReLU\ $ \R \ni x \mapsto \max\{ x, \gamma x \} \in \R $, $ \gamma \in (0,1) $, activation functions (cf., \eg, \cite[Sections 1.2.3 and 1.2.11]{ArBePhi2024}) which are not continuously differentiable functions.

We also point the reader, \eg, to \cite{dereich2024learningrateadaptivestochastic}, \cite{GodichonBaggioni2023}, \cite{ar1706.06569}, \cite[Theorem 1, Theorem 2, Theorem 3, and Theorem 4]{ar1805.08114}, \cite{ar1706.05507}, \cite[Theorem 2, Theorem 6, and Theorem 8]{ar2305.18471}, and the references therein for works that do not study \Adam\ but provide upper bounds for other (related) adaptive gradient based optimization methods – such as the \RMSprop\ \cite{Hinton24_RMSprop} and the \Adagrad\ \cite{duchi2011adaptive} optimization methods – applied to \SOPs\ and we refer, \eg, to \cite{DereichKassingconvergence2021} and \cite[Theorem 2 and Corollary 1]{ar2007.07989} and the references therein for works that do not investigate \Adam\ but study other (related) accelerated gradient based optimization methods – such as the momentum \SGD\ method \cite{Polyak1964SomeMO} – applied to \SOPs. The \Adam\ optimizer can essentially be viewed as a combination of the \RMSprop\ optimizer and the momentum \SGD\ optimizer. However, the interplay of both of these optimization methods within the \Adam\ optimizer result in new phenomena such as the \Adam\ symmetry theorem in \cref{main theorem 2} above.

We also refer, for instance, to \cite{ar1808.02941,ar2106.08208,ar1805.07557,LuoAdaBoundAMSBound2019,ReddiKale2019} for works that provide upper bounds for the regret or for averages of the norm of the gradient of the objective function for modified/extended variants of \Adam\ such as the AMSGrad \cite{ReddiKale2019}, the NosAdam \cite[Theorem 3.2 and Corollary 3.2.1]{ar1805.07557}, the AdaFom \cite[Corollary 3.2]{ar1808.02941}, the AdaBound \cite[Theorem 4 and Corollary 4.1]{LuoAdaBoundAMSBound2019}, and the AMSBound \cite[Theorem 5]{LuoAdaBoundAMSBound2019} optimizers.

Some of the above cited works contain error analyses that apply to much more general \SOPs\ than the class of strongly convex \SOPs\ considered in \cref{main theorem} and \cref{cor: convegence of Adam}, respectively.
However, to the best of our knowledge, \cref{main theorem} and \cref{cor: convegence of Adam}, respectively, are the first results in the scientific literature that establish for the \emph{error of \Adam} the upper bound \cref{conclude: main theorem} with the asymptotic $M ^{-1} (1 - \beta_2 )$ for the non-learning rate error term. In this context we also note that the numerical simulations in \cref{sec: numerial simulation} strongly suggest that the convergence rate
$1$ with respect to the reciprocal of the size of the mini-batches $M \in \N$ ($ \nicefrac{ 1 }{ M } $ to the power of $ 1 $) and the convergence rate $1$ with respect to the distance of $\beta_2 $ to 1 ($ 1 - \beta_2 $ to the power of $ 1 $) in the asymptotic $M^{-1} (1 - \beta_2 )$ for the non-learning rate error term (see \cref{conclude: main theorem}) seems to be \emph{optimal}.

We also point out that in \cref{main theorem} (and \cref{main theorem 2}, respectively) we do neither employ $(L_0 , L_1 )$-smoothness assumptions \cite[Definition 1]{ar1905.11881} (cf., \eg, also \cite[Definition 3.2]{ar2304.13972} and \cite[Section 3.1]{ar2306.01264}) nor a suitable linear growth bound on the noise but instead the key assumption that we impose is that the objective function of the considered \SOP\ (the function that we intend to minimize) is \emph{strongly convex}. Of course, the objective functions (the functions one intends to minimize) that appear in \ANN\ training problems are in nearly all cases neither strongly convex nor convex (cf., \eg, \cite[Corollary 2.19]{ScarpaMultiLayers}). Nonetheless, the class of strongly convex \SOPs\ is maybe the most basic class of optimization problems for which one intends to understand an optimization algorithm and we believe that it is an issue of fundamental relevance to aim for developing a complete mathematical understanding of the \Adam\ optimizer in the situation of strongly convex \SOPs. Building up on this, it is of course also highly relevant to study \Adam\ in the situation of non-convex \SOPs\ and, in particular, in the training of \ANNs.

Moreover, under this strong convexity assumption, a further novelty of this work is that in \cref{main theorem 2} and \cref{main theorem} we do neither bound the regret of \Adam, nor the gradient of the objective function, and nor any minima or temporal averages of such quantities but instead we estimate the actual optimization error, that is, in \cref{main theorem 2} and \cref{main theorem} we bound the $L^p$-norm of the standard norm distance between the \Adam\ optimization process after $n \in \N$ gradient steps and the solution of the \SOP\ (the global minimizer of the \SOP) from above. There are several works (cf., \eg, \cite{ar1808.03408,ar2208.09900,ZouShen2019}) that provide upper bounds for minima or averages of the norm of the gradient of the objective function and, in the strongly convex setting, upper bounds for the norm of the gradient of the objective function imply upper bounds for the standard norm distance to the unique global minimizer of the objective function (see \cref{eq: strongly convex introduction 2} above) but such works only imply upper bounds for minima or averages of the standard norm distance to the unique global minimzer of the objective function while \cref{main theorem} (and \cref{main theorem 2}, respectively) provide upper bounds for the standard norm distance of the \Adam\ optimization process to the unique global minimizer $\vartheta$ after \emph{exactly} $n \in \N$ gradient/\Adam\ steps (instead of minima or averages over the first $n$ steps).

 We also refer, \eg, to the survey articles \cite{WeChaoSteLei2020,Ruder2016AnOO,Sun2019OptimizationFD} and the monographs \cite{bach2024learning,Guro2023handbook,ArBePhi2024,PhilipJa2024} for further illustrative reviews and references on \SGD\ optimization methods.

\subsubsection*{Upper bounds for Adam and related optimizers for deterministic optimization problems}

We also refer, \eg, to \cite{ar1911.07596}, \cite{Sebastianimprovement2022}, \cite{Bock2019APO}, \cite[Theorem 3.4]{ar1807.06766}, \cite{DeArAdsharp2025}, \cite{Polyak1964SomeMO}, \cite[Theorem 4.2]{Shi2021RMSpropCW}, and the references therein for works that analyze \Adam\ and related gradient based optimization methods when applied to deterministic \OPs\ (which is often also considered as the so-called full/complete batch mode \cite{Bock2019APO}).

In the situation of \Adam\ applied to deterministic \OPs\ the learning rates do not need to converge to zero (see \cref{assumption: gamma1}) but may stay constant but sufficiently small and a convergence speed for \Adam\ in the number of gradient steps can then be achieved that is higher than the error term $\sqrt{ \gamma_n }$ in \cref{item 2: main theorem 2} in the \Adam\ symmetry theorem in \cref{main theorem 2} and in the second summand on the right hand side of \cref{conclude: main theorem} in \cref{main theorem}. In particular, in the deterministic case with sufficiently small constant learning rates the error of \Adam\ after $n \in \N$ gradient steps can be bounded from above by the term $\scrc_0 \exp( - \scrc_1 n )$ where $\scrc_0, \scrc_1 \in (0,\infty)$ are strictly positive real numbers; see item (iv) in Theorem 1.1 in \cite{DeArAdsharp2025} for the optimal value of $\scrc_1$.

In this context we also want to point out that in the situation of deterministic \OPs\ the non-convergence phenomena with the respect to the number of gradient steps revealed in the \Adam\ symmetry theorem in \cref{main theorem 2} above does not show up. Indeed, in the case where the \SOP\ in \cref{main theorem 2} is deterministic (corresponding to the case $v=w$ in \cref{main theorem 2}) we have that the variance of the data vanishes in the sense that $\operatorname{Var}( X_{ 1, 1 } ) = 0$ in which case we have that the data $X_{ 1, 1 }$ is always symmetrically distributed (as it is even one-point distributed) and, thus, items \ref{item 2: main theorem 2} and \ref{item 1: main theorem 2} in \cref{main theorem 2} hold true and the first error term on the right hand side of \cref{conclude: main theorem} can be omitted.

\subsubsection*{Lower bounds for Adam and related optimizers for stochastic optimization problems}

There are also several works in the scientific literature that establish \emph{lower bounds and non-convergence/divergence results} for \Adam\ and related gradient based optimization methods (cf., \eg, \cite{Sebastiannonconvergence2025,CheriditoJentzenRossmanek2021,DeRoAr2024nonconvergence,HannibalJentzenThang2024,DoArAd2025,gallon2022blowphenomenagradientdescent,ArAd2024,LSSK2020,ReddiKale2019,ZhangChen2022}) and some of these works are also related to parts of the conclusions in the \Adam\ symmetry theorem in \cref{main theorem 2} above.

In particular, the work \cite[Theorem 3]{ReddiKale2019} constructs for every first moment \Adam\ parameter $\beta_1 \in (0,1)$ and every second moment \Adam\ parameter $\beta_2 \in ( ( \beta_1 )^2, 1) )$ an appropriate convex \SOP\ (depending on $\beta_1$ and $\beta_2$) such that \Adam\ does not converge to the optimal solution (cf. also \cite[Theorem 1 and Theorem 2]{ReddiKale2019}). 

We also refer, \eg, to \cite[Section 3.2 and Appendix E]{ZhangChen2022} for blow-up results (divergence to infinity) for \Adam\ applied to a specific example objective function (cf., \eg, also \cite{gallon2022blowphenomenagradientdescent}, \cite{Ziwei2020}, \cite{JuDaDeAr2025}, \cite{ar1906.05890}, \cite[Section 1.2]{Petersenfound2021}, and \cite{ar2110.02732}). In this context we also note that for the class of \SOPs\ considered in the \Adam\ symmetry theorem in \cref{main theorem 2} we know already that \Adam\ does not diverge to infinity but stays bounded according to \cite{DereichGraArAd2025} as long as the standard condition $( \beta_1 )^2 < \beta_2$ is satisfied and for the class of \SOPs\ considered in \cref{main theorem} we assume boundedness of the \Adam\ optimizer (see above \cref{conclude: main theorem} for details).

We also refer, \eg, to \cite{HannibalJentzenThang2024,DoArAd2025,ArAd2024} for works that show that in the training of \ANNs\ with the \ReLU\ and related activations we have for \Adam\ and many other \SGD\ optimization methods that the objective function evaluated at the optimization process fails to converge to the infimal value of the objective function (cf., \eg, also \cite{CheriditoJentzenRossmanek2021,LSSK2020}). In these works \cite{HannibalJentzenThang2024,DoArAd2025,ArAd2024} the non-convergence property is not a special issue of \Adam\ but also applies to the standard \SGD\ method and arises due to the highly non-convex optimization landscape in the training of \ANNs\ with the \ReLU\ and related activations.

We also refer, \eg, to  \cite{DeRoAr2024nonconvergence} for results that show that \Adam\ and other \SGD\ optimization methods do not converge to any random point if the sequence of learning rates does not converge to zero (cf., \eg, also \cite [Lemma 7.2.13]{ArBePhi2024}).
The non-convergence property in \cite{DeRoAr2024nonconvergence} is also not a special issue of \Adam\ but applies to the standard \SGD\ method too and arises due
to the fact that the learning rates do not asymptotically vanish. In \cref{main theorem 2} and \cref{main theorem} above the learning rates are assumed to converge to zero in the sense that $\lim_{ n \to \infty } \gamma_n = 0$. This follows from the assumption in \cref{assumption: gamma1} that $\sum_{ n = 1 }^{ \infty } ( \gamma_n )^p <\infty$.

To the best of our knowledge, the \Adam\ symmetry theorem in \cref{main theorem 2} is the first
 result in the scientific literature that provides a \emph{complete characterization} for the convergence of \Adam\ to the solution of the \SOP. Interestingly, according to \cref{main theorem 2} convergence of \Adam\ to the solution of the \SOP\ is \emph{equivalent} to the data in the \SOP\ being \emph{symmetrically distributed} and in this sense we propose to refer to \cref{main theorem 2} as \emph{\Adam\ symmetry theorem}.

\subsection{Adam in the training of large language models (LLMs)}
\label{sec:adam_training_llms}

\Adam\ and its extension \AdamW\ are among the most widely used optimizers in the training of state-of-the-art \LLMs.
In \cref{fig:llm-optimizers,tab:llm-optimizers} we list the optimizers employed in the pre-training of several prominent \LLMs\ for all cases in which this information has been publicly disclosed.
In every such case, the optimizer is a variant of \Adam, \AdamW, or of the \Muon\ optimizer.
We note that in the \Muon\ optimizer and its variants \Muon Clip and \Muon\ Split, the bias terms are trained according to the \AdamW\ optimizer.
Consequently, in all \LLMs\ listed in \cref{fig:llm-optimizers,tab:llm-optimizers} for which the optimizer has been publicly disclosed, at least the bias terms are trained according to a variant of the \Adam\ optimizer.
This demonstrates the popularity of the \Adam\ optimizer.

\begin{table}[h!]
\centering
\renewcommand{\arraystretch}{1.3}
\begin{adjustbox}{max width=0.8\textwidth}
\begin{tabular}{ | l l l | }
  \hline
  \textbf{Model} & \textbf{Optimizer} & \textbf{References} \\
  \hline
  \multicolumn{3}{|l|}{\textbf{OpenAI}} \\
  GPT-3
    & \Adam\ \((\beta_1{=}0.9,\;\beta_2{=}0.95,\;\gamma_{\text{peak}}{=}6\text{e-}5,\;\gamma_{\text{end}}{=}6\text{e-}6)\)
    & \cite{LLMTable200514165} % (checked)
    \\
  GPT-4, GPT-5
    & Not disclosed
    & \cite{LLMTable230308774,LLMTable260103267} % (checked) 
    \\
  \hline
  \multicolumn{3}{|l|}{\textbf{Meta}} \\
  LLaMA\,1 
    & \AdamW\ \((\beta_1{=}0.9,\;\beta_2{=}0.95,\;\gamma_{\text{peak}}{=}1.5\text{e-}4,\;\gamma_{\text{end}}{=}1.5\text{e-}5)\)
    & \cite{LLMTable230213971} % (checked)
    \\
  LLaMA\,2
    & \AdamW\ \((\beta_1{=}0.9,\;\beta_2{=}0.95,\;\gamma_{\text{peak}}{=}1.5\text{e-}4,\;\gamma_{\text{end}}{=}1.5\text{e-}5)\)
    & \cite{LLMTable230709288} % (checked)
    \\
  LLaMA\,3
    & \AdamW\ \((\gamma_{\text{peak}}{=}8\text{e-}5,\;\gamma_{\text{end}}{=}8\text{e-}7)\)
    & \cite{LLMTable240721783} % (checked)
    \\
  \hline
  \multicolumn{3}{|l|}{\textbf{DeepSeek}} \\
  DeepSeek-V3
    & \AdamW\ \((\beta_1{=}0.9,\;\beta_2{=}0.95,\;\gamma_{\text{peak}}{=}2.2\text{e-}4,\;\gamma_{\text{end}}{=}7.3\text{e-}6)\)
    & \cite{LLMTable241219437} % (checked)
    \\
  \hline
  \multicolumn{3}{|l|}{\textbf{Google DeepMind}} \\
  Gemini\,1.0, 1.5, 2.5, 3, 3.1
    & Not disclosed
    & \cite{LLMTable231211805,GeminiTeam2024,GeminiTeam2025,GoogleDeepMind2025} % (checked. make references nicer)
    \\
  \hline
  \multicolumn{3}{|l|}{\textbf{Alibaba}} \\
  Qwen\,2.5
    & Not disclosed
    & \cite{LLMTable241215115} % (checked)
    \\
  \hline
  \multicolumn{3}{|l|}{\textbf{Anthropic}} \\
  Claude\,3, 4
    & Not disclosed
    & \cite{LLMTableClaude3ModelCard,LLMTableClaude4ModelCard} % (checked)
    \\
  \hline
  \multicolumn{3}{|l|}{\textbf{BigScience}} \\
  BLOOM
    & \AdamW\ \((\beta_1{=}0.9,\;\beta_2{=}0.95,\;\gamma_{\text{peak}}{=}6\text{e-}5,\;\gamma_{\text{end}}{=}6\text{e-}6)\)
    & \cite{LLMTable221105100} % (checked)
    \\
  \hline
  \multicolumn{3}{|l|}{\textbf{EleutherAI}} \\
  Pythia
    & \Adam\ \((\beta_1{=}0.9,\;\beta_2{=}0.95,\;\gamma_{\text{peak}}{=}1.2\text{e-}4,\;\gamma_{\text{end}}{=}1.2\text{e-}5)\)
    & \cite{LLMTable230401373} % (checked)
    \\
  \hline
  \multicolumn{3}{|l|}{\textbf{TII}} \\
  Falcon
    & \AdamW\ \((\gamma_{\text{peak}}{=}1.25\text{e-}4,\;\gamma_{\text{end}}{=}1.25\text{e-}5)\)
    & \cite{LLMTable231116867} % (checked)
    \\
  Falcon2-11B
    & \AdamW\ \((\beta_1{=}0.9,\;\beta_2{=}0.999,\;\gamma_{\text{peak}}{=}3.7\text{e-}4,\;\gamma_{\text{end}}{=}1.89\text{e-}5)\)
    & \cite{LLMTable240714885} % (checked)
    \\
  \hline
  \multicolumn{3}{|l|}{\textbf{Allen Institute}} \\
  OLMo
    & \AdamW\ \((\beta_1{=}0.9,\;\beta_2{=}0.95,\;\gamma_{\text{peak}}{=}3\text{e-}4,\;\gamma_{\text{end}}{=}3\text{e-}5)\)
    & \cite{LLMTable240200838} % (checked)
    \\
  OLMo\,2
    & \AdamW\ \((\gamma_{\text{peak}}{=}6\text{e-}4,\;\gamma_{\text{end}}{=}6\text{e-}5)\)
    & \cite{OLMo2025} % (checked)
    \\
  \hline
  \multicolumn{3}{|l|}{\textbf{Microsoft}} \\
  phi-1
    & \AdamW\ \((\gamma_{\text{peak}}{=}1\text{e-}3)\)
    & \cite{Gunasekar2023} % (checked)
    \\
  phi-1.5
    & \AdamW\ \((\beta_1{=}0.9,\;\beta_2{=}0.98,\;\gamma_{\text{peak}}{=}\gamma_{\text{end}}{=}2\text{e-}4)\)
    & \cite{LLMTable230905463} % (checked)
    \\
  phi-2, phi-3, phi-4
    & Not disclosed
    & \cite{Abdin2024,Abdin2024a} % (checked)
    \\
  \hline
  \multicolumn{3}{|l|}{\textbf{Mistral AI}} \\
  Mistral\,7B
    & Not disclosed
    & \cite{LLMTable231006825} % (checked)
    \\
  \hline
  \multicolumn{3}{|l|}{\textbf{Moonshot AI}} \\
  Moonlight\,16B
    & \Muon\
    & \cite{LLMTable250216982} % (checked)
    \\
  Kimi K2
    & MuonClip
    & \cite{LLMTable250720534} % (checked)
    \\
  \hline
  \multicolumn{3}{|l|}{\textbf{Zhipu AI}} \\
  GLM-5
    & \Muon\ Split
    & \cite{GLM5Team2026} % (checked)
    \\
  \hline
  \multicolumn{3}{|l|}{\textbf{NVIDIA}} \\
  Nemotron\,3 Super
    & \AdamW\ \((\beta_1{=}0.9,\;\beta_2{=}0.95,\;\gamma_{\text{peak}}{=}4.5\text{e-}4,\;\gamma_{\text{end}}{=}4.5\text{e-}6)\)
    & \cite{NVIDIA2026} % (checked)
    \\
  \hline
  \multicolumn{3}{|l|}{\textbf{Arcee AI}} \\
  Arcee Trinity
    & \Muon\
    & \cite{LLMTable260217004} % (checked)
    \\
  \hline
  \multicolumn{3}{|l|}{\textbf{xAI}} \\
  Grok-1, 2, 3, 4
    & Not disclosed
    & \cite{GROK1,GROK2,GROK3,GROK4} % (checked)
    \\
  \hline
\end{tabular}
\end{adjustbox}
\caption{%
\label{tab:llm-optimizers}
  Stochastic optimization algorithms used during pre-training of prominent \LLMs.
  Here $\gamma_{\text{peak}}$ denotes the peak learning rate and $\gamma_{\text{end}}$ the final learning rate.
  Most models employ a schedule that linearly warms up to $\gamma_{\text{peak}}$ over a few hundred steps and then decays to $\gamma_{\text{end}}$.
  Learning rates are given in scientific notation, \eg, $1.5\mathrm{e}{-4} = 1.5 \times 10^{-4}$.
  % \newline
  % ${}^{\dagger}$The technical report explicitly withholds training details,
  % including the optimizer.
  % \newline
  % ${}^{\ddagger}$The Gemini\,1.0 technical report does not explicitly name the
  % optimizer; \AdamW\ is attributed to the Gemini family in related documentation.
  % \newline
  % ${}^{\S}$Anthropic's model cards and system cards describe Constitutional AI
  % and RLHF post-training but do not disclose the pre-training optimizer
  % for any Claude generation.
}
\end{table}

% =====================================================================
% optimizer_table.tex
% ---------------------------------------------------------------------
% Pre-training optimizers of the top-28 LLMs.  This block reproduces the
% table from the "2026_05_18 Optimizer slide" presentation so that it can
% be dropped into the main article via a single \include command.
%
% Requires (loaded in the preamble of the main document):
%   array, booktabs, xcolor, colortbl, graphicx, tikz, etoolbox,
%   csvsimple-l3
%
% Data file : llm_optimizers.csv  (copied from the slide folder)
% Logos     : logos/*.png         (copied from the slide folder)
% =====================================================================

% ── Access-type background colours ───────────────────────────────────
\definecolor{propBG}{RGB}{210,225,250}      % blue  – Proprietary
\definecolor{openBG}{RGB}{205,240,210}      % green – Open Weights (free)
\definecolor{restrBG}{RGB}{255,238,190}     % amber – Open Weights (restricted)

% ── Score-bar settings ───────────────────────────────────────────────
\newcommand{\barbaseline}{20}   % bar length proportional to (score - barbaseline)
\newcommand{\barmaxcm}{3}        % maximum bar length in cm (score = 60)
\definecolor{barlow}{RGB}{20,60,120}        % colour at score min (46)
\definecolor{barhigh}{RGB}{130,0,0}         % colour at score max (60)

% ── Muon-family optimizer highlight ──────────────────────────────────
\definecolor{muonBG}{RGB}{255,255,255}      % white – any Muon variant

% ── Column-type shortcuts ────────────────────────────────────────────
\newcolumntype{R}[1]{>{\raggedleft\arraybackslash}p{#1}}
\newcolumntype{L}[1]{>{\raggedright\arraybackslash}p{#1}}
\newcolumntype{C}[1]{>{\centering\arraybackslash}p{#1}}

% ── Static convenience macros ────────────────────────────────────────
\newcommand{\nd}{}
\newcommand{\nf}{\textit{\textcolor{gray}{?}}}

% ── Score bar + optimizer combined cell ──────────────────────────────
%   #1 = score (integer); #2 = optimizer name; #3 = ostatus
%   Confirmed -> bold label inside bar; Inferred -> bold label + dagger;
%   otherwise -> no label.  \pdfstrcmp is used because csvsimple-l3 reads
%   the CSV content with catcode-12 (other) characters.
\newcommand{\scorebar}[3]{%
  \begin{tikzpicture}[baseline=-0.5ex, inner sep=0pt, outer sep=0pt]
    \pgfmathsetmacro{\bw}{(#1 - \barbaseline) / (60 - \barbaseline) * \barmaxcm}%  bar width in cm
    \pgfmathtruncatemacro{\cpct}{(#1 - 46) * 100 / 14}%   colour 0–100
    \colorlet{bcol}{barhigh!\cpct!barlow}%
    \fill[bcol] (0,-0.09cm) rectangle (\bw cm, 0.11cm);%
    \node[anchor=east, inner xsep=1pt, inner ysep=0pt,
        font=\fontsize{4.8}{4.8}\sffamily\fontseries{l}\selectfont, text=white]
      at (\bw cm, 0.01cm) {#1};%
    \ifnum\pdfstrcmp{#3}{Confirmed}=0%
      \node[anchor=west, inner sep=0pt,
            font=\fontsize{4.8}{4.8}\selectfont, text=white]
        at (0.04cm, 0.01cm) {\textbf{#2}};%
    \fi%
    \ifnum\pdfstrcmp{#3}{Inferred}=0%
      \node[anchor=west, inner sep=0pt,
            font=\fontsize{4.5}{4.5}\selectfont, text=white]
        at (0.04cm, 0.01cm) {\textbf{#2}$^{\dagger}$};%
    \fi%
  \end{tikzpicture}%
}

% ── Table header (opens the tabular; closed by caller) ───────────────
\newcommand{\tableheadFirst}{%
  \renewcommand{\arraystretch}{0.75}%
  \setlength{\tabcolsep}{1.5pt}%
  \begin{tabular}{C{0.30cm} L{4.8cm} L{3cm}}
  \toprule
   & \textbf{Model} & \textbf{Score \& Optimizer} \\
  \midrule
}

% ── The slide itself ─────────────────────────────────────────────────
\begin{center}
\resizebox{0.8\linewidth}{!}{%
\begin{tcolorbox}[
  width=9.2cm, halign=center,
  colback=white, colframe=black!55, boxrule=0.6pt, arc=2pt,
  left=4mm, right=4mm, top=3mm, bottom=3.5mm,
]
  {\small\textbf{Pre-training optimizers of top-28 \LLMs}}\par
  \vspace{1.5mm}
  \tiny
  \tableheadFirst
  \begin{filecontents*}{llm_optimizers.csv}
rnk,model,mfr,score,atype,optimizer,ostatus,src,notes,logofile
1,GPT-5.5 (xhigh),OpenAI,60,Proprietary,—,Undisclosed,https://openai.com/index/introducing-gpt-5-5/,No technical report with training details released,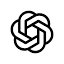
2,GPT-5.5 (high),OpenAI,59,Proprietary,—,Undisclosed,https://openai.com/index/introducing-gpt-5-5/,No technical report with training details released,logos/openai.png
3,Claude Opus 4.7 (Adaptive Reasoning),Anthropic,57,Proprietary,—,Undisclosed,https://www.anthropic.com/claude,No technical report with training details released,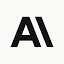
4,Gemini 3.1 Pro Preview,Google DeepMind,57,Proprietary,—,Undisclosed,https://deepmind.google/models/model-cards/gemini-3-1-pro/,No technical report with optimizer details,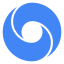
5,GPT-5.5 (medium),OpenAI,57,Proprietary,—,Undisclosed,https://openai.com/index/introducing-gpt-5-5/,No technical report with training details released,logos/openai.png
6,Kimi K2.6,Moonshot AI,54,Open Weights,MuonClip,Confirmed,https://arxiv.org/abs/2507.20534,MuonClip = Muon + QK-Clip stability mechanism; Modified MIT license (commercial use allowed),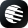
7,MiMo-V2.5-Pro,Xiaomi,54,Open Weights,AdamW,Confirmed,https://venturebeat.com/technology/open-source-xiaomi-mimo-v2-5-and-v2-5-pro-are-among-the-most-efficient-and-affordable-at-agentic-claw-tasks,Muon for matrix params; AdamW for embeddings/norms/LM-head; MIT license,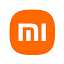
8,GPT-5.3 Codex (xhigh),OpenAI,54,Proprietary,—,Undisclosed,https://developers.openai.com/api/docs/models/gpt-5.5,No technical report with training details released,logos/openai.png
9,Grok 4.3 (high),xAI,53,Proprietary,—,Undisclosed,https://x.ai/news/grok-4,No model card or technical report published by xAI,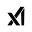
10,Muse Spark,Meta,52,Proprietary,—,Undisclosed,https://ai.meta.com/blog/introducing-muse-spark-msl/,Meta's first non-open-weights frontier model; no optimizer details disclosed,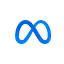
11,Claude Opus 4.7 (non-reasoning),Anthropic,52,Proprietary,—,Undisclosed,https://www.anthropic.com/claude,No technical report with training details released,logos/anthropic.png
12,Qwen3.6 Max Preview,Alibaba,52,Proprietary,—,Undisclosed,https://artificialanalysis.ai/models/qwen3-6-max,API-only closed-weights model; no optimizer info in available documentation,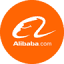
13,Claude Sonnet 4.6 (max),Anthropic,52,Proprietary,—,Undisclosed,https://www.anthropic.com/claude,No technical report with training details released,logos/anthropic.png
14,DeepSeek V4 Pro (Max),DeepSeek,52,Open Weights,Muon,Confirmed,https://fe-static.deepseek.com/chat/transparency/deepseek-V4-model-card-EN.pdf,Muon replaces AdamW for weight matrices; AdamW for embeddings/heads; MIT license,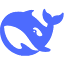
15,GLM-5.1,Z.AI (Zhipu AI),51,Open Weights,Muon Split,Confirmed,https://arxiv.org/abs/2602.15763,Muon Split adapts matrix orthogonalization per attention head for MLA; MIT license,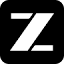
16,GPT-5.5 (low),OpenAI,51,Proprietary,—,Undisclosed,https://openai.com/index/introducing-gpt-5-5/,No technical report with training details released,logos/openai.png
17,Qwen3.6 Plus,Alibaba,50,Proprietary,—,Undisclosed,https://openrouter.ai/qwen/qwen3.6-plus-preview,API-only closed-weights model; no optimizer info in available documentation,logos/alibaba.png
18,DeepSeek V4 Pro (High),DeepSeek,50,Open Weights,Muon,Confirmed,https://fe-static.deepseek.com/chat/transparency/deepseek-V4-model-card-EN.pdf,Same base model as DeepSeek V4 Pro (Max); MIT license,logos/deepseek.png
19,GLM-5,Z.AI (Zhipu AI),50,Open Weights,Muon Split,Confirmed,https://arxiv.org/abs/2602.15763,Same paper as GLM-5.1; MIT license,logos/zhipuai.png
20,MiniMax-M2.7,MiniMax,50,Open Weights (Commercial Restricted),—,Undisclosed,https://www.minimax.io/news/minimax-m27-en,Pre-training optimizer undocumented (the reported optimizer is Inferred); self-evolution refers to RL post-training only; MiniMax proprietary license (requires written auth. for commercial use). ,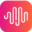
21,MiMo-V2.5,Xiaomi,49,Open Weights,AdamW,Confirmed,https://arxiv.org/abs/2601.02780,Claude infers the optimizer “Muon + AdamW” from MiMo-V2.5-Pro tech report (same team/training stack) - I can’t find it. What I can see is that AdamW is used to train V2 and it seems like this is the backbone for V2.5 where they just changed the post-training ; MIT license,logos/xiaomi.png
22,GPT-5.4 mini (xhigh),OpenAI,49,Proprietary,—,Undisclosed,https://developers.openai.com/api/docs/models/gpt-5.5,No technical report with training details released,logos/openai.png
23,GPT-5.4 (low),OpenAI,48,Proprietary,—,Undisclosed,https://developers.openai.com/api/docs/models/gpt-5.5,No technical report with training details released,logos/openai.png
24,GLM-5-Turbo,Z.AI (Zhipu AI),47,Proprietary,Muon Split,Confirmed,https://arxiv.org/abs/2602.15763,Proprietary fine-tune of GLM-5 base; optimizer inferred from base model paper,logos/zhipuai.png
25,DeepSeek V4 Flash (Max),DeepSeek,47,Open Weights,Muon,Confirmed,https://fe-static.deepseek.com/chat/transparency/deepseek-V4-model-card-EN.pdf,Smaller/faster variant in V4 family; same optimizer; MIT license,logos/deepseek.png
26,Gemini 3 Flash,Google DeepMind,46,Proprietary,—,Undisclosed,https://deepmind.google/technologies/gemini/,No technical report with optimizer details,logos/deepmind.png
27,DeepSeek V4 Flash (High),DeepSeek,46,Open Weights,Muon,Confirmed,https://fe-static.deepseek.com/chat/transparency/deepseek-V4-model-card-EN.pdf,Same base model as DeepSeek V4 Flash (Max); MIT license,logos/deepseek.png
28,Qwen3.6 27B,Alibaba,46,Open Weights,—,Undisclosed,https://huggingface.co/Qwen/Qwen3.6-27B,No optimizer details found in available documentation; Apache 2.0 license,logos/alibaba.png
\end{filecontents*}
\csvreader[
    head to column names,
    range={1-28},
    late after line=\\,
  ]{llm_optimizers.csv}{}{%
    \raisebox{-1pt}{\includegraphics[height=0.22cm,keepaspectratio]{\logofile}} & \model
         & \scorebar{\score}{\optimizer}{\ostatus}%
  }
  \bottomrule
  \end{tabular}
  \par\vspace{2mm}
  {\tiny\textcolor{gray}{Intelligence scores: artificialanalysis.ai, 18 May 2026}}
\end{tcolorbox}%
}
\par\vspace{1.5mm}
\begin{minipage}{\textwidth}
  \captionof{figure}{%
    The 28 highest-scoring \LLMs\ according to the \emph{Artificial Analysis Intelligence Index} in \cite{ArtificialAnalysis}, together with their pre-training
    optimizers in the cases in which it was disclosed.
  }
  \label{fig:llm-optimizers}
\end{minipage}
\end{center}

\subsection{Structure of this article}
The remainder of this work is organized as follows. In \cref{sec: convergence general} we employ the findings in \cite{DereichAdamconvergence2024} to establish \emph{strong and weak convergence rates} for the \Adam\ optimizer for the class of strongly convex \SOPs\ and, thereby, prove \cref{main theorem} in \cref{subsec: main theorem 1} above. In \cref{sec: converge to wrong point} we make use of the finding in \cite{DereichJentzenKassing2025} to \emph{disprove} for the class of \SOPs\ considered in \cref{main theorem 2} above that \Adam\ converges to a critical point of the \SOP\ if the data in the \SOP\ is not symmetrically distributed. In \cref{prop: convergence in symmetric case} in \cref{sec: convergence symmetric}, in turn, we establish for the class of \SOPs\ considered in \cref{main theorem 2} that \Adam\ converges to a critical point of the \SOP\ if the data of the \SOP\ is \emph{symmetrically distributed}. Furthermore, in \cref{sec: convergence symmetric} we also combine \cref{prop: convergence in symmetric case} with the findings of \cref{sec: converge to wrong point} to establish \cref{main theorem 2} above. Finally, in \cref{sec: numerial simulation} we complement the analytical findings of this work by means of several \emph{numerical simulations} that indicate the sharpness of the established convergence rates and that illustrate the practical appearance of the revealed theoretical phenomena.

\section{Convergence rates for Adam for strongly convex stochastic optimization}\label{sec: convergence general}
In this section we make use of the findings in our preliminary work \cite{DereichAdamconvergence2024} to establish in \cref{cor: convegence of Adam} in \cref{subsec: convergence with sufficient batch} below strong and weak convergence rates for the \Adam\ optimizer applied to the class of strongly convex \SOPs. In \cref{cor: convegence of Adam2} in \cref{subsec: convergence with sufficient batch} we present a special case of \cref{cor: convegence of Adam} in which we restrict ourselves to strong convergence rates for \Adam. \cref{main theorem} in the introduction is a direct application of \cref{cor: convegence of Adam2}. Some of the arguments in this section are inspired by the arguments in Sections 8 and 9 in our preliminary work \cite{DereichAdamconvergence2024}.
 
\subsection{Strong convexity}
Throughout this work we frequently employ the well-known concepts of strong and generalized convexity (cf., \eg, \cite[Definition 2.11]{Guro2023handbook} and  \cite[Definitions 5.7.18 and 5.7.19]{ArBePhi2024}). For convenience of the reader we briefly recall these concepts within the next notion.
\begin{definition}[Strong convexity]\label{definition: strongly convex}
    Let $d\in \N$ and let $f\colon \R^d\to\R$ be a function. Then 
    \begin{enumerate}[label=(\roman*)]
        \item \label{item 1: strongly convex} for every $c\in \R$ we say that $f$ is $c$-\stronglyconvex\ if and only if it holds that 
    $
\textstyle\R^d\ni\theta\mapsto f(\theta)-\frac{c}{2}\|\theta\|^2\in \R$
    is convex and
    \item \label{item 2: strongly convex} we say that $f$ \stronglyconvexx\ if and only if there exists $c\in (0,\infty)$ such that $f$ is $c$-\stronglyconvex.
    \end{enumerate}
\end{definition}
\subsection{Mathematical formulation of the Adam vector field}
\newcommand{\setTheta}{V}
\cfclear
\begin{setting}\label{setting: Adam vector field}
    Let $(\Omega,\mathcal F,\P)$ be a probability space, let $\fd,\dimX\in \N$,  let $\setX\subseteq \R^\dimX$ be compact, let $X_{n,m}\colon \Omega\allowbreak\to \setX$, $(n,m) \in (\N_0)^2$, be \iid\ random variables, let $\smalll=(\smalll(\theta,x))_{(\theta,x)\in \R^\fd\times\R^{\dimX}}\allowbreak\in C^2(\R^\fd\times\R^{\dimX},\R)$, $\alpha\in (0,1)$, $\beta\in (\alpha^2,1)$, $\varepsilon,\chi\in (0,\infty)$, let $(\gamma_n)_{n\in \N}\subseteq (0,\infty)$ be non-increasing,
     for every $M\in \N$, $n\in \N_0$ let $\bbX_M^n=(\bbX_M^{n,1},\dots,\bbX_M^{n,\fd})\colon \R^\fd\times\Omega\to \R^\fd$ satisfy for all $\theta=(\theta_1,\dots,\theta_\fd)\in \R^\fd$ that
    \begin{equation}\label{def: bbX}
        \bbX_M^{n,i}(\theta)=-\frac 1M \biggl[\textstyle \sum\limits_{m=1}^M(\frac{\partial}{\partial\theta_i}\smalll)(\theta,X_{n,m})\biggr],
    \end{equation}
    assume that $\R^\fd\ni \theta\mapsto \E[\smalll(\theta,X_{1,1})]\in \R$ is $\chi$-\stronglyconvex, and for every $M\in \N$ let $f_M=(f_M^1,\dots,f_M^\fd)\colon \R^\fd\to\R^\fd$ satisfy for all $i\in \{1,2,\dots,\fd\}$, $\theta\in \R^\fd$ that
     \begin{equation}\label{def: f}
         \textstyle f_M^i(\theta)=(1-\alpha)\,\E\Bigl[\frac{\sum_{n=0}^\infty\alpha^n\bbX_M^{n,i}(\theta)}{\varepsilon+\sqrt{(1-\beta)\sum_{n=0}^\infty\beta^n|\bbX_M^{n,i}(\theta)|^2}}\Bigr] \ifnocf.
     \end{equation}
     \cfout[.]
\end{setting}
\begin{athm}{lemma}{lem: explicit batch size}
   Let $(\Omega,\cF,\P)$ be a probability space, let $\alpha\in (0,1)$, $\beta\in (\alpha^2,1)$, $\varepsilon,\consc\in (0,\infty)$, let $\bbX_{k}\colon \Omega\to [-\consc,\consc]$, $k\in \N_0$, be \iid\ random variables, let $\bbY_k\colon\Omega\to [-\consc,\consc]$, $k\in \N_0$, be \iid\ random variables, assume that $(\bbX_k,\bbY_k)$, $k\in \N_0$, are identically distributed, and let $F,G\in \R$ satisfy
   \begin{equation}\llabel{def: F}
       \begin{split}
        F&=\textstyle\E\biggl[\frac{(1-\alpha)\sum_{k=0}^\infty\alpha^k\bbY_k}{\varepsilon+\sqrt{(1-\beta)\sum_{k=0}^\infty\beta^k|\bbX_k|^2}}\biggr]\\
        &-\E\biggl[\textstyle\frac{\bigl[(1-\alpha)\sum_{k=0}^\infty\alpha^k\bbX_{k}\bigr]\bigl[(1-\beta)\sum_{k=0}^\infty\beta^k\bbX_k\bbY_k\bigr]}{\bigl(\varepsilon+\sqrt{(1-\beta)\sum_{k=0}^\infty\beta^k|\bbX_k|^2}\bigr)^2\sqrt{\mathbbm 1_{\{0\}}(\sum_{k=0}^\infty|\bbX_k|)+(1-\beta)\sum_{k=0}^\infty\beta^k|\bbX_k|^2}}\biggr]
        \end{split}
   \end{equation}
   and
   \begin{equation} \llabel{def: G}
       \begin{split}
           G=\frac{\varepsilon\E[\bbY_0]}{\bigl(\varepsilon+|\E[\bbX_0]|\bigr)^2}.
       \end{split}
   \end{equation}
      Then
   \begin{equation}\llabel{conclude}
       |F-G|\leq\textstyle \bigl(2\consc+\frac{2\consc^2}{\varepsilon}+\frac{8\consc^3+4\consc}{\varepsilon^2}\bigr)\bigl(\operatorname{Var}(\bbX_0)\bigr)^{1/2}+\bigl(\frac{\consc}{\varepsilon^2}+\frac{1}{\varepsilon}\bigr)\bigl(\operatorname{Var}(\bbY_0)\bigr)^{1/2}.
   \end{equation}
\end{athm}
\begin{aproof}
    Throughout this proof assume without loss of generality that $\P(\bbX_0=0)<1$ (otherwise note that the fact that
    \begin{equation}
        F=\E\biggl[\frac{(1-\alpha)\sum_{k=0}^\infty\alpha^k\bbY_k}{\varepsilon}\biggr]=\frac{\E[\bbY_0]}{\varepsilon}=\frac{\E[\bbY_0]}{\varepsilon+|\E[\bbX_0]|}=\frac{\varepsilon\E[\bbY_0]}{\bigl(\varepsilon+|\E[\bbX_0]|\bigr)^2}=G
    \end{equation}
    establishes \lref{conclude}) and let $X\colon \Omega\to\R$, $Y\colon \Omega\to\R$, $Z\colon \Omega\to\R$, $T\colon \Omega\to\R$, and $\bbT\in \R$ satisfy
    \begin{equation}\llabel{def: XYZ}
        \textstyle X=(1-\alpha)\sum\limits_{k=0}^\infty\alpha^k\bbX_{k},\qquad Y=(1-\alpha)\sum\limits_{k=0}^\infty\alpha^k\bbY_k, \qquad Z=(1-\beta)\sum\limits_{k=0}^\infty\beta^k|\bbX_k|^2,
    \end{equation}
    \begin{equation}\llabel{def: T}
      \textstyle  T=(1-\beta)\sum\limits_{k=0}^\infty\beta^k\bbX_k\bbY_k,\qqandqq \bbT=\E[\bbX_0]\E[\bbY_0]=(1-\beta)\sum\limits_{k=0}^\infty\beta^k\E[\bbX_k]\E[\bbY_k].
    \end{equation}
    \argument{the assumption that $\P(\bbX_0=0)<1$; the fact that $\bbX_k$, $k\in \N_0$, are \iid}{that for all $v \in \{0,1\}$ it holds that
    \begin{equation}\llabel{eqcpre}
        \begin{split}
            &\P\bigl( \textstyle\sum_{ k = 0 }^{ \infty } \beta^{ v k } | \bbX_k |^{ v + 1 } = 0 \bigr)
=
\P\bigl( \forall\, k \in \N_0 \colon| \bbX_k |^{ v + 1 } = 0 \bigr)
\\
&=
\P\bigl( \forall\, k \in \N_0 \colon \bbX_k  = 0 \bigr)
=
\P\bigl( \cap_{ k = 0 }^{ \infty }  \{\bbX_k  = 0 \} \bigr)
\\
&=
\P\bigl( \cap_{ K \in \N_0 } ( \cap_{ k = 0 }^{ K } \{  \bbX_k  = 0 \} ) \bigr)
=
\lim_{ K \to \infty }
\bigl[\textstyle
  \prod_{ k = 0 }^{ K } \P( \bbX_k = 0 )
\bigr]
\\
&=\textstyle
\lim_{ K \to \infty }
\bigl( | \P( \bbX_0 = 0 ) |^K \bigr)
=
0.
        \end{split}
    \end{equation}}
    \argument{\lref{eqcpre};the assumption that $\P( \bbX_0 = 0 ) < 1$;}{that 
    \begin{equation}\llabel{eqc}
\P\bigl(\textstyle\sum_{k=0}^{\infty}|\bbX_k|>0\bigr)=\P(Z>0)=1\qqandqq \E[|\bbX_0|^2]>0.
\end{equation}}
    \argument{\lref{def: G}; \lref{def: XYZ};\lref{def: T};the assumption that $\bbX_k$, $k\in \N_0$, are \iid;the assumption that $\bbY_k$, $k\in \N_0$, are \iid}{that
    \begin{equation}\llabel{def: f2}
         \begin{split}
        G &=\E[ \bbY_0 ] 
\left[ 
  \frac{ 
    \varepsilon + | \E[ \bbX_0 ] | 
  }{ ( \varepsilon + | \E[ \bbX_0 ] | )^2 } 
  - 
  \frac{ 
    | \E[ X_0 ] | 
  }{ ( \varepsilon + | \E[ \bbX_0 ] | )^2 } 
\right]\\
&=\frac{\E[\bbY_0]}{\varepsilon+|\E[\bbX_0]|}-\frac{|\E[\bbX_0]|\E[\bbY_0]}{(\varepsilon+|\E[\bbX_0]|)^2}=\frac{\E[\bbY_0]}{\varepsilon+|\E[\bbX_0]|}-\frac{\E[\bbX_0]\E[\bbX_0]\E[\bbY_0]}{(\varepsilon+|\E[\bbX_0]|)^2\bigl(|\E[\bbX_0]|+\mathbbm 1_{\{0\}}(\E[\bbX_0])\bigr)}\\
        &
        =\frac{(1-\alpha)\sum_{k=0}^\infty\alpha^k\E[\bbY_k]}{\varepsilon+|\E[\bbX_0]|}-\frac{\bigl[(1-\alpha)\sum_{k=0}^\infty\alpha^k\E[\bbX_{k}]\bigr]\bigl[(1-\beta)\sum_{k=0}^\infty\beta^k\E[\bbX_k]\E[\bbY_k]\bigr]}{(\varepsilon+|\E[\bbX_0]|)^2|\bigl(|\E[\bbX_0]|+\mathbbm 1_{\{0\}}(\E[\bbX_0])\bigr)}\\
        &=\frac{\E[Y]}{\varepsilon+|\E[\bbX_0]|}-\frac{\E[X]\bbT}{(\varepsilon+|\E[\bbX_0]|)^2(|\E[\bbX_0]|+\mathbbm 1_{\{0\}}(\E[\bbX_0]))}.
        \end{split}
    \end{equation}}
    \argument{\lref{def: F};\lref{def: XYZ};\lref{def: T};\lref{eqc}}{that
    \begin{equation}\llabel{eqtg1}
       \textstyle F=\E\bigl[\frac{Y}{\varepsilon+\sqrt{Z}}\bigr]-\E\bigl[\frac{XT}{(\varepsilon+\sqrt{Z})^2\sqrt{Z+\mathbbm 1_{\{0]}(Z)}}\bigr]
    \end{equation}}
    \argument{\lref{eqtg1};\lref{def: f2}}{that
    \begin{equation}\llabel{def: F-G}
        \begin{split}
           & \textstyle|F-G|\\
           &\leq\textstyle \E\bigl[\bigl|\frac{Y}{\varepsilon+\sqrt{Z}}-\frac{\E[Y]}{\varepsilon+|\E[\bbX_0]|}\bigr|\bigr]+\E\bigl[\big|\frac{XT}{(\varepsilon+\sqrt{Z})^2\sqrt{Z+\mathbbm 1_{\{0]}(Z)}}-\frac{\E[X]\bbT}{(\varepsilon+|\E[\bbX_0]|)^2(|\E[\bbX_0]|+\mathbbm 1_{\{0\}}(\E[\bbX_0]))}\big|\bigr].
        \end{split}
    \end{equation}}
    \argument{\lref{eqc};the fact that $\P(|Y|\leq \consc)=1$}{that
    \begin{equation}\llabel{eq1}
    \begin{split}
        &\E\biggl[\biggl|\frac{Y}{\varepsilon+\sqrt{Z}}-\frac{\E[Y]}{\varepsilon+|\E[\bbX_0]|}\biggr|\biggr]\\&=\E\biggl[\biggl|\frac{Y}{\varepsilon+\sqrt{Z}}-\frac{\E[Y]}{\varepsilon+\sqrt{Z}}+\frac{\E[Y]}{\varepsilon+\sqrt{Z}}-\frac{\E[Y]}{\varepsilon+|\E[\bbX_0]|}\biggr|\biggr]\\
        &\leq \E\biggl[\biggl|\frac{Y}{\varepsilon+\sqrt{Z}}-\frac{\E[Y]}{\varepsilon+\sqrt{Z}}\bigg|\biggr]+\E\biggl[\bigg|\frac{\E[Y]}{\varepsilon+\sqrt{Z}}-\frac{\E[Y]}{\varepsilon+|\E[\bbX_0]|}\biggr|\biggr]\\
        &\leq \frac{\E[|Y-\E[Y]|]}{\varepsilon}+|\E[Y]|\E\biggl[\bigg|\frac{1}{\varepsilon+|\E[\bbX_0]|}-\frac{1}{\varepsilon+\sqrt{Z}}\biggr|\biggr]\\
        &\leq \frac{\E[|Y-\E[Y]|]}{\varepsilon}+|\E[Y]|\E\biggl[\bigg|\frac{\sqrt{Z}-\sqrt{(\E[\bbX_0])^2}}{\varepsilon^2}\biggr|\biggr]\\
        &\leq \frac{\E[|Y-\E[Y]|]}{\varepsilon}+\frac{\consc}{\varepsilon^2}\,\E\bigl[\big|\sqrt{Z}-\sqrt{(\E[\bbX_0])^2}\big|\bigr].
        \end{split}
    \end{equation}}
      \argument{\lref{def: XYZ};\lref{eqc};the fact that $\bbX_k$, $k\in \N_0$, are \iid; the fact that for all $x,y\in \R$ it holds that $(x+y)^2\leq 2x^2+2y^2$}{that
    \begin{equation}\llabel{eq2.5}
    \begin{split}
       & \E\biggl[\textstyle(1-\beta)\sum\limits_{k=0}^\infty\beta^k\frac{(\bbX_k+\E[\bbX_k])^2}{(\sqrt{Z}+|\E[\bbX_0]|)^2+\mathbbm 1_{\{(0,0)\}}(Z,\E[\bbX_0])}\biggr]\\
      & \leq \E\biggl[\textstyle(1-\beta)\sum\limits_{k=0}^\infty\beta^k\frac{2(\bbX_k)^2}{Z+1_{\{0,0\}}(Z,\E[\bbX_0])}\biggr]+\E\biggl[(1-\beta)\sum\limits_{k=0}^\infty\beta^k\frac{2(\E[\bbX_k])^2}{(\E[\bbX_0])^2+\mathbbm 1_{\{(0,0)\}}(Z,\E[\bbX_0])}\biggr]\\
      &\leq2+2=4.
        \end{split}
    \end{equation}}
    \argument{\lref{eq2.5};\lref{def: XYZ};the fact that $\bbX_k$, $k\in \N_0$, are \iid; the Cauchy-Schwarz inequality}{that
    \begin{align}
        &\E\bigl[\big|\sqrt{Z}-\sqrt{(\E[\bbX_0])^2}\big|\bigr]=\E\biggl[\frac{\bigl|Z-(\E[\bbX_0])^2\bigr|}{\sqrt{Z}+|\E[\bbX_0]|+\mathbbm 1_{\{(0,0)\}}(Z,\E[\bbX_0])}\biggr]\notag\\
        &=\E\biggl[\frac{\bigl|(1-\beta)\sum_{k=0}^\infty \beta^k[(\bbX_k)^2-(\E[\bbX_0])^2]\bigr|}{\sqrt{Z}+|\E[\bbX_0]|+\mathbbm 1_{\{(0,0)\}}(Z,\E[\bbX_0])}\biggr]\notag\\
        &= \E\biggl[\frac{\big|(1-\beta)\sum_{k=0}^\infty\beta^k (\bbX_k-\E[\bbX_k])(\bbX_k+\E[\bbX_k])\big|}{\sqrt{Z}+|\E[\bbX_0]|+\mathbbm 1_{\{(0,0)\}}(Z,\E[\bbX_0])}\biggr]\llabel{eq2}\\
       &\leq\textstyle \biggl(\E\biggl[(1-\beta)\sum\limits_{k=0}^\infty\beta^k(\bbX_k-\E[\bbX_k])^2\biggr]\biggr)^{1/2}\biggl(\E\biggl[(1-\beta)\sum\limits_{k=0}^\infty\beta^k\frac{(\bbX_k+\E[\bbX_k])^2}{(\sqrt{Z}+|\E[\bbX_0]|+\mathbbm 1_{\{(0,0)\}}(Z,\E[\bbX_0]))^2}\biggr]\biggr)^{1/2}\notag\\
        &=\textstyle \biggl(\E\biggl[(1-\beta)\sum\limits_{k=0}^\infty\beta^k(\bbX_k-\E[\bbX_k])^2\biggr]\biggr)^{1/2}\biggl(\E\biggl[(1-\beta)\sum\limits_{k=0}^\infty\beta^k\frac{(\bbX_k+\E[\bbX_k])^2}{(\sqrt{Z}+|\E[\bbX_0]|)^2+\mathbbm 1_{\{(0,0)\}}(Z,\E[\bbX_0])}\biggr]\biggr)^{1/2}\notag\\
       &\leq 2\bigl(\E\bigl[(\bbX_0-\E[\bbX_0])^2\bigr]\bigr)^{1/2}\notag.
    \end{align}}
    \argument{\lref{def: XYZ};the fact that $\bbY_k$, $k\in \N_0$, are \iid; the Cauchy-Schwarz inequality}{that
    \begin{equation}\llabel{eq2'}
        \E\bigl[|Y-\E[Y]|\bigr]\leq (1-\alpha)\sum_{k=0}^\infty\alpha^k\E[|\bbY_k-\E[\bbY_k]|]=\E[|\bbY_0-\E[\bbY_0]|]\leq \bigl(\E\bigl[(\bbY_0-\E[\bbY_0])^2\bigr]\bigr)^{1/2}.
    \end{equation}}
    \argument{\lref{eq2'};\lref{eq1};\lref{eq2}; the Cauchy-Schwarz inequality}{that
    \begin{equation}\llabel{evidence 1}
        \E\biggl[\biggl|\frac{Y}{\varepsilon+\sqrt{Z}}-\frac{\E[Y]}{\varepsilon+\sqrt{(\E[Z])^2}}\biggr|\biggr]\leq \frac{1}{\varepsilon}\bigl(\E\bigl[(\bbY_0-\E[\bbY_0])^2\bigr]\bigr)^{1/2}+\frac{2\consc}{\varepsilon^2}\bigl(\E[(\bbX_0-\E[\bbX_0])^2]\bigr)^{1/2}.
    \end{equation}}
    \argument{\lref{def: XYZ};\lref{def: T}}{that
    \begin{equation}\llabel{eq3}
        \begin{split}
           & \E\biggl[\bigg|\frac{XT}{(\varepsilon+\sqrt{Z})^2\sqrt{Z+1_{\{0\}}(Z)}}-\frac{\E[X]\bbT}{(\varepsilon+|\E[\bbX_0]|)^2\bigl(|\E[\bbX_0]|+\mathbbm 1_{\{0\}}(\E[\bbX_0])\bigr)}\bigg|\biggr]\\
           &\leq \E\biggl[\bigg|\frac{XT}{(\varepsilon+\sqrt{Z})^2\sqrt{Z+1_{\{0\}}(Z)}}-\frac{\E[X]T}{(\varepsilon+\sqrt{Z})^2\sqrt{Z+1_{\{0\}}(Z)}}\bigg|\biggr]\\
           &+\E\biggl[\bigg|\frac{\E[X]T}{(\varepsilon+\sqrt{Z})^2\sqrt{Z+1_{\{0\}}(Z)}}-\frac{\E[X]T}{(\varepsilon+|\E[\bbX_0]|)^2\bigl(|\E[\bbX_0]|+\mathbbm 1_{\{0\}}(\E[\bbX_0])\bigr)}\bigg|\biggr]\\
           &+\E\biggl[\bigg|\frac{\E[X](T-\bbT)}{(\varepsilon+|\E[\bbX_0]|)^2\bigl(|\E[\bbX_0]|+\mathbbm 1_{\{0\}}(\E[\bbX_0])\bigr)}\bigg|\biggr]\\
           &\leq\E\biggl[\bigg|\frac{\E[X]T}{(\varepsilon+\sqrt{Z})^2\sqrt{Z+1_{\{0\}}(Z)}}-\frac{\E[X]T}{(\varepsilon+|\E[\bbX_0]|)^2\bigl(|\E[\bbX_0]|+\mathbbm 1_{\{0\}}(\E[\bbX_0])\bigr)}\bigg|\biggr]\\
           &+ \frac{1}{\varepsilon^2}\E\biggl[\frac{|T(X-\E[X])|}{\sqrt{Z+1_{\{0\}}(Z)}}\biggr]+\frac{1}{\varepsilon^2} \E\biggl[\frac{|\E[X](T-\bbT)|}{|\E[\bbX_0]|+\mathbbm 1_{\{0\}}(\E[\bbX_0])}\biggr].
        \end{split}
    \end{equation}}
    \argument{\lref{def: XYZ};\lref{def: T};the fact that $\P(|\bbX_0|\leq \consc)=\P(|\bbY_0|\leq \consc)=1$; the assumption that $(\bbX_k,\bbY_k)$, $k\in \N_0$, are identically distributed;the Cauchy-Schwarz inequality}{that
    \begin{equation}\llabel{evd2: eq1}
    \begin{split}
         &\E\biggl[\frac{|\E[X](T-\bbT)|}{|\E[\bbX_0]|+\mathbbm 1_{\{0\}}(\E[\bbX_0])}\biggr]=  \frac{|\E[X]|}{|\E[\bbX_0]|+\mathbbm 1_{\{0\}}(\E[\bbX_0])}\E[|T-\bbT|]\leq\E[|T-\bbT|]\\
         &\leq (1-\beta)\sum_{k=0}^\infty\beta^k\E[|\bbX_k\bbY_k-\E[\bbX_k]\E[\bbY_k]|]=\E[|\bbX_0\bbY_0-\E[\bbX_0]\E[\bbY_0]|]\\
         &\leq \E[|\bbY_0(\bbX_0-\E[\bbX_0])|]+\E[|\E[\bbX_0](\bbY_0-\E[\bbY_0])|]\\
         &\leq \consc\E[|\bbX_0-\E[\bbX_0]|]+\consc\E[|\bbY_0-\E[\bbY_0]|]\\
         &\leq \consc\bigl(\E\bigl[(\bbX_0-\E[\bbX_0])^2\bigr]\bigr)^{1/2}+\consc\bigl(\E\bigl[(\bbY_0-\E[\bbY_0])^2\bigr]\bigr)^{1/2}.
         \end{split}
    \end{equation}}
    \argument{\lref{def: XYZ};\lref{def: T}; the fact that $\P(|\bbY_0|\leq \consc)=1$; the Cauchy-Schwarz inequality}{that it holds $\P$-a.s.\ that
    \begin{equation}\llabel{eq5}
    \begin{split}
        \textstyle T^2&=\textstyle\biggl|(1-\beta)\sum\limits_{k=0}^\infty\beta^k\bbX_k\bbY_k\biggr|^2\leq \consc^2(1-\beta) \biggl|\sum\limits_{k=0}^\infty\beta^k\bbX_k\biggr|^2=\consc^2(1-\beta) \biggl|\sum\limits_{k=0}^\infty\sqrt{\beta^k}\sqrt{\beta^k}\bbX_k\biggr|^2\\
        &\leq \textstyle \consc^2(1-\beta) \biggl(\sum\limits_{k=0}^\infty\beta^k\biggr) \biggl(\sum\limits_{k=0}^\infty\beta^k|\bbX_k|^2\biggr)=\textstyle \consc\biggl(\sum\limits_{k=0}^\infty\beta^k|\bbX_k|^2\biggr)=\consc^2Z.
        \end{split}
    \end{equation}}
    \argument{\lref{eq5};\lref{def: XYZ};the Cauchy-Schwarz inequality; the assumption that $\bbX_k$, $k\in \N$, are \iid; the fact that $0<\alpha<1$}{that
    \begin{equation}\llabel{evd2: eq2}
    \begin{split}
       \E\biggl[\frac{|T(X-\E[X])|}{\sqrt{Z+1_{\{0\}}(Z)}}\biggr]
       &\leq \consc\E\bigl[|X-\E[X]|\bigr]\leq \consc \bigl(\E\bigl[|X-\E[X]|^2\bigr]\bigr)^{1/2}\\
       &\leq \consc\biggl(\E\biggl[\textstyle\biggl((1-\alpha)\sum\limits_{k=0}^\infty\alpha^k(\bbX_k-\E[\bbX_k])\biggr)^2\biggr]\biggr)^{1/2}\\
    &= \consc\biggl(\textstyle(1-\alpha)^2\sum\limits_{k=0}^\infty\alpha^{2k}\textstyle\E\bigl[(\bbX_k-\E[\bbX_k])^2\bigr]\biggr)^2\biggr)^{1/2}\\
    &= \consc\biggl(\frac{(1-\alpha)^2}{1-\alpha^2}\E\bigl[(\bbX_0-\E[\bbX_0])^2\bigr]\biggr)^{1/2}\leq \consc\bigl(\E\bigl[(\bbX_0-\E[\bbX_0])^2\bigr]\bigr)^{1/2}.
       \end{split}
    \end{equation}}
\argument{\lref{eqc};\lref{eq5}}{that
\begin{equation}\llabel{eq7}
\begin{split}
    &\E\biggl[\bigg|\frac{\E[X]T}{(\varepsilon+\sqrt{Z})^2\sqrt{Z+\mathbbm 1_{\{0\}}(Z)}}-\frac{\E[X]T}{(\varepsilon+|\E[\bbX_0]|)^2\bigl(|\E[\bbX_0]|+\mathbbm 1_{\{0\}}(\E[\bbX_0])\bigr)}\bigg|\biggr]\\
&=\E\biggl[\bigg|\frac{\E[X]T\bigl((\varepsilon+\sqrt{Z})^2\sqrt{Z}-(\varepsilon+|\E[\bbX_0]|)^2|\E[\bbX_0]|\bigr)}{(\varepsilon+\sqrt{Z})^2\sqrt{Z+1_{\{0\}}(Z)}(\varepsilon+|\E[\bbX_0]|)^2\bigl(|\E[\bbX_0]|+\mathbbm 1_{\{0\}}(\E[\bbX_0])\bigr)}\bigg|\biggr]\\
&\leq \frac{\consc}{\varepsilon^2}\E\bigl[\big|(\varepsilon+\sqrt{Z})^2\sqrt{Z}-(\varepsilon+|\E[\bbX_0]|)^2|\E[\bbX_0]|\big|\bigr]\\
&\leq \consc\E\bigl[\big|\sqrt{Z}-\sqrt{(\E[\bbX_0])^2}\big|\bigr]+\frac{2\consc}{\varepsilon}\E\bigl[\big|Z-(\E[\bbX_0])^2\big|\bigr]\\
&+\frac{\consc}{\varepsilon^2}\E\biggl[\bigg|\frac{Z^3-|\E[\bbX_0]|^6}{\sqrt{Z^3}+\sqrt{|\E[\bbX_0]|^6}+\mathbbm 1_{\{(0,0)\}}(Z,\E[\bbX_0])}\bigg|\biggr].
\end{split}
\end{equation}}
\argument{\lref{def: XYZ};the assumption that $\bbX_k$, $k\in \N_0$, are \iid; the fact that $\P(|\bbX_0|\leq \consc)=1$}{that
\begin{equation}\llabel{eq7.1}
    \E\bigl[\big|Z-(\E[\bbX_0])^2\big|\bigr]\leq  \E\bigl[|(\bbX_0)^2-(\E[\bbX_0])^2|\bigr]=\E\bigl[(\bbX_0-\E[\bbX_0])^2\bigr]\leq 2\consc\bigl(\E\bigl[(\bbX_0-\E[\bbX_0])^2\bigr]\bigr)^{1/2}.
\end{equation}}
\argument{\lref{def: XYZ};\lref{eq7.1};the fact that for all $x,y\in \R$ it holds that $xy\leq x^2+y^2$; the fact that $\P(Z\leq \consc^2)=\P(|\bbX_0|\leq \consc)=1$}{that
\begin{equation}\llabel{eq7.2}
\begin{split}
    &\E\biggl[\bigg|\frac{Z^3-|\E[\bbX_0]|^6}{\sqrt{Z^3}+\sqrt{|\E[\bbX_0]|^6}+\mathbbm 1_{\{(0,0)\}}(Z,\E[\bbX_0])}\bigg|\biggr]\\
    &=\E\biggl[\bigg|\frac{\bigl(Z-(\E[\bbX_0])^2\bigr)\bigl((Z^2+Z(\E[\bbX_0])^2+(\E[\bbX_0])^4\bigr)}{\sqrt{Z^3}+\sqrt{|\E[\bbX_0]|^6}+\mathbbm 1_{\{(0,0)\}}(Z,\E[\bbX_0])}\bigg|\biggr]\\
    &\leq 2\E\biggl[\bigg|\frac{\bigl(Z-(\E[\bbX_0])^2\bigr)Z^2}{\sqrt{Z^3}+\sqrt{|\E[\bbX_0]|^6}+\mathbbm 1_{\{(0,0)\}}(Z,\E[\bbX_0])}\bigg|\biggr]\\
   & \quad+2\E\biggl[\bigg|\frac{\bigl(Z-(\E[\bbX_0])^2\bigr)(\E[\bbX_0])^4}{\sqrt{Z^3}+\sqrt{|\E[\bbX_0]|^6}+\mathbbm 1_{\{(0,0)\}}(Z,\E[\bbX_0])}\bigg|\biggr]\\
    &\leq 
2\E\bigl[\big|\bigl(Z-(\E[\bbX_0])^2\bigr)\sqrt{Z}\big|\bigr]+2\E\bigl[\big|\bigl(Z-(\E[\bbX_0])^2\bigr)|\E[\bbX_0]|\big|\bigr]\leq 8\consc^2\bigl(\E\bigl[(\bbX_0-\E[\bbX_0])^2\bigr]\bigr)^{1/2}.
\end{split}
\end{equation}}
\argument{\lref{eq2};\lref{eq7};\lref{eq7.1};\lref{eq7.2};}{that
\begin{equation}\llabel{evd2: eq3}
\begin{split}
   & \E\biggl[\bigg|\frac{\E[X]T}{(\varepsilon+\sqrt{Z})^2\sqrt{Z+\mathbbm 1_{\{0\}}(Z)}}-\frac{\E[X]T}{(\varepsilon+|\E[\bbX_0]|)^2\bigl(|\E[\bbX_0]|+\mathbbm 1_{\{0\}}(\E[\bbX_0])\bigr)}\bigg|\biggr]\\
   &\leq \bigg(2\consc+\frac{4\consc^2}{\varepsilon}+\frac{8\consc^3}{\varepsilon^2}\biggr)\bigl(\E\bigl[(\bbX_0-\E[\bbX_0])^2\bigr]\bigr)^{1/2}.
    \end{split}
\end{equation}}
\argument{\lref{eq3};\lref{evd2: eq1};\lref{evd2: eq2};\lref{evd2: eq3}}{that
\begin{equation}\llabel{evidence 2}
\begin{split}
    & \E\biggl[\bigg|\frac{XT}{(\varepsilon+\sqrt{Z})^2\sqrt{Z+\mathbbm 1_{\{0\}}(Z)}}-\frac{\E[X]\bbT}{(\varepsilon+|\E[\bbX_0]|)^2\bigl(|\E[\bbX_0]|+\mathbbm 1_{\{0\}}(\E[\bbX_0])\bigr)}\bigg|\biggr]\\
    &\leq
     \biggl(2\consc+\frac{2\consc^2}{\varepsilon}+\frac{8\consc^3+2\consc}{\varepsilon^2}\biggr)\big(\E\bigl[(\bbX_0-\E[\bbX_0])^2\bigr]\bigr)^{1/2}+\frac{\consc}{\varepsilon^2}\big(\E\bigl[(\bbY_0-\E[\bbY_0])^2\bigr]\bigr)^{1/2}.
     \end{split}
\end{equation}}
\argument{\lref{def: F-G};\lref{evidence 1};\lref{evidence 2};}{\lref{conclude}\dott}
\end{aproof}
In \cref{item 5: lem negative derivative} in \cref{lem: negative derivative} below we verify that the \Adam\ vector field satisfies a well-known monotonicity condition (cf., \eg, \cite[Definition 5.7.12, Definition 5.7.13, Lemma 5.7.14, and Lemma 5.7.15]{ArBePhi2024}) that we employ in our convergence rate analysis for \Adam\ in \cref{theo: convergence of Adam} in \cref{subsec: convergence with sufficient batch} below. \cref{lem: negative derivative} is inspired by Lemma 9.1 in our preliminary work \cite{DereichAdamconvergence2024}.
\newcommand{\setV}{U}
\begin{athm}{lemma}{lem: negative derivative}
    Assume \cref{setting: Adam vector field}, let
     $\setV\subseteq \R^\fd$ be compact, and let $\varrho\in \R$ satisfy for all $\theta=(\theta_1,\dots,\theta_\fd)\in \setV$, $i,j\in\{1,2,\dots,\fd\}$, $x\in \setX$ that
     \begin{equation}\llabel{def: varrho}
   \textstyle \varrho\geq \max\bigl\{|(\frac{\partial^2}{\partial\theta_i\partial\theta_j}\smalll)(\theta,x)|,|(\frac{\partial}{\partial\theta_i}\smalll)(\theta,x)|\bigr\}.
     \end{equation}
     Then
    \begin{enumerate}[label=(\roman*)]
    \item \label{item 1: lem negative derivative} there exists $\fC\in \R$ such that for all $M,p\in \N$, $n\in\N_0$, $\theta,\vartheta\in \setV$ it holds that 
    \begin{equation}
      \|\bbX_{M}^{n}(\theta)\|\leq \fC=1 \qqandqq \|\bbX_{M}^{n}(\theta)-\bbX_{M}^{n}(\vartheta)\|^p\leq \fC^p\|\theta-\vartheta\|^p,
    \end{equation}
        \item \label{item 2: lem negative derivative} it holds for all $M\in \N$ that $f_M\in C^1(\R^\fd,\R^\fd)$,
         \item  \label{item 3: lem negative derivative} it holds for all  $M\in \N$, $\theta=(\theta_1,\dots,\theta_\fd)\allowbreak\in \R^\fd$, $i,j\in \{1,2,\dots,\fd\}$ that
    \begin{equation}
    \begin{split}
        &\textstyle\frac{\partial}{\partial\theta_j}f_M^i(\theta)
        =\textstyle\E\biggl[\frac{(1-\alpha)\sum_{k=0}^\infty\alpha^k(\frac{\partial}{\partial\theta_j}\bbX_M^{k,i})(\theta)}{\varepsilon+\sqrt{(1-\beta)\sum_{k=0}^\infty\beta^k|\bbX_M^{k,i}(\theta)|^2}}\biggr]\\
        &-\E\Biggl[\textstyle\frac{\bigl[\mathbbm 1_{\{0\}}(\sum_{k=0}^\infty |\bbX_M^{k,i}(\theta)|)+(1-\alpha)\sum_{k=0}^\infty\alpha^k\bbX_{M}^{k,i}(\theta)\bigr]\bigl[(1-\beta)\sum_{k=0}^\infty\beta^k\bbX_M^{k,i}(\theta)(\frac{\partial}{\partial\theta_j}\bbX_{M}^{k,i})(\theta)\bigr]}{\bigl(\varepsilon+\sqrt{(1-\beta)\sum_{k=0}^\infty\beta^k|\bbX_M^{k,i}(\theta)|^2}\bigr)^2\sqrt{\mathbbm 1_{\{0\}}(\sum_{k=0}^\infty |\bbX_M^{k,i}(\theta)|)+(1-\beta)\sum_{k=0}^\infty\beta^k|\bbX_M^{k,i}(\theta)|^2}}\Biggr],
        \end{split}
        \end{equation}
%          \item \label{item 4: lem negative derivative} it holds for all $\theta=(\theta_1,\dots,\theta_\fd)\in \setV$, $i\in \{1,2,\dots,\fd\}$, $M\in \N$ that
%         \begin{equation}
% \biggl|f_M^i(\theta)+\biggl(\frac{\textstyle\E[(\frac{\partial}{\partial\theta_i}\smalll)(\theta,X_{1,1})]}{\textstyle\varepsilon+|\E[(\frac{\partial}{\partial\theta_i}\smalll)(\theta,X_{1,1})]|}\biggr)\biggr|\leq 2\rho\biggl(\frac{1}{\varepsilon}+\frac{\rho^2}{\varepsilon^2}\biggr) M^{-1/2},
             % \end{equation}
        \item \label{item 4: lem negative derivative} it holds for all $M\in \N$, $\theta=(\theta_1,\dots,\theta_\fd)\in \setV$, $i,j\in \{1,2,\dots,\fd\}$ that
        \begin{equation}
        \begin{split}
            \biggl|\textstyle\frac{\partial}{\partial\theta_j}f_M^i(\theta)+\frac{\varepsilon\E[(\frac{\partial^2}{\partial \theta_i\partial\theta_j}\smalll)(\theta,X_{1,1})]}{\bigl(\varepsilon+\bigl|\E[(\frac{\partial}{\partial \theta_i}\smalll)(\theta,X_{1,1})]\bigr|\bigr)^2}\biggr|\leq 4\varrho\bigl(2\varrho+\frac{2\varrho^2+1}{\varepsilon}+\frac{8\varrho^3+5\varrho}{\varepsilon^2}\bigr)M^{-1/2},
            \end{split}
        \end{equation}
        and
    \item \label{item 5: lem negative derivative} it holds for all $M\in \N$, $\theta\in \setV$, $x\in \R^\fd$ that \begin{equation}
        \spro{x, (f_M)'(\theta)x}
     \leq \bigl[\textstyle-\frac{\varepsilon \chi}{(\varepsilon+\varrho)^2}+4\fd\varrho\bigl(2\varrho+\frac{2\varrho^2+1}{\varepsilon}+\frac{8\varrho^3+5\varrho}{\varepsilon^2}\bigr)M^{-1/2}\bigr]\|x\|^2  \ifnocf.
    \end{equation}
        \end{enumerate}
        \cfout[.]
\end{athm}
\begin{aproof}
Throughout this proof assume without loss of generality that $\operatorname{Var}(X_{1,1})\neq 0$ (otherwise note that the fact that for all $\theta\in \R^d$, $M\in \N$, $i\in \{1,2,\dots,d\}$, $n\in \N_0$ it holds that
\begin{equation}
     \P\bigl( \textstyle \bbX_M^{n,i}=-(\frac{\partial}{\partial\theta_i}\smalll)(\theta,\E[X_{1,1}])\bigr)=1
\end{equation}
\begin{equation}
    \text{and}\qquad \textstyle f_M^i(\theta)=(1-\alpha)\,\E\Bigl[\frac{\sum_{n=0}^\infty\alpha^n\bbX_M^{n,i}(\theta)}{\varepsilon+\sqrt{(1-\beta)\sum_{n=0}^\infty\beta^n|\bbX_M^{n,i}(\theta)|^2}}\Bigr]=-\frac{(\frac{\partial}{\partial\theta_i}\smalll)(\theta,\E[X_{1,1}])}{\varepsilon+|(\frac{\partial}{\partial\theta_i}\smalll)(\theta,\E[X_{1,1}])|}
\end{equation}
(cf.\ \cref{def: bbX} and \cref{def: f}) establishes items \ref{item 1: lem negative derivative}, \ref{item 2: lem negative derivative}, \ref{item 3: lem negative derivative}, \ref{item 4: lem negative derivative}, and \ref{item 5: lem negative derivative}) and let $\vertiii{\cdot}\colon \R^{d\times d}\to\R$ satisfy for all $A=(a_{i,j})_{(i,j)\in \{1,2,\dots,d\}^2}\in \R^{d\times d}$ that 
\begin{equation}
    \textstyle\vertiii{A}=\bigl(\sum_{i,j=1}^d (a_{i,j})^2\bigr)^{1/2}.
\end{equation}
\argument{\cref{def: bbX};the assumption that $\smalll\in C^2(\R^\fd\times\R^\dimX)$}{that for all $\theta=(\theta_1,\dots,\theta_\fd)\in \R^\fd$, $i,j\in \{1,2,\dots,\fd\}$, $M,N\in \N$, $\omega\in\Omega$ with $\sum_{k=0}^N|\bbX_M^{k,i}(\theta,\omega)|^2>0$ it holds that
    \begin{equation}\llabel{eq: derivative calculate}
    \begin{split}
        &\displaystyle\frac{\partial}{\partial \theta_j} \textstyle\biggl(\frac{(1-\alpha)\sum_{k=0}^N\alpha^k\bbX_M^{k,i}(\theta,\omega)}{\varepsilon+\sqrt{(1-\beta)\sum_{k=0}^N\beta^k|\bbX_M^{k,i}(\theta,\omega)|^2}}\biggr)=\frac{(1-\alpha)\sum_{k=0}^N\alpha^k(\frac{\partial}{\partial\theta_j}\bbX_M^{k,i})(\theta,\omega)}{\varepsilon+\sqrt{(1-\beta)\sum_{k=0}^N\beta^k|\bbX_M^{k,i}(\theta,\omega)|^2}}\\
        &-\textstyle\frac{\bigl[(1-\alpha)\sum_{k=0}^N\alpha^k\bbX_{M}^{k,i}(\theta,\omega)\bigr]\bigl[(1-\beta)\sum_{k=0}^N\beta^k\bbX_M^{k,i}(\theta,\omega)(\frac{\partial}{\partial\theta_j}\bbX_{M}^{k,i})(\theta,\omega)\bigr]}{\bigl(\varepsilon+\sqrt{(1-\beta)\sum_{k=0}^N\beta^k|\bbX_M^{k,i}(\theta,\omega)|^2}\bigr)^2\sqrt{(1-\beta)\sum_{k=0}^N\beta^k|\bbX_M^{k,i}(\theta,\omega)|^2}}.
        \end{split}
    \end{equation}}
    \argument{\cref{def: bbX};the assumption that $\smalll\in C^2(\R^\fd\times\R^\dimX)$}{that for all $\theta=(\theta_1,\dots,\theta_\fd)\in \R^\fd$, $i,j\in \{1,2,\dots,\fd\}$, $M,N\in \N$, $\omega\in\Omega$ with $\sum_{k=0}^N|\bbX_M^{k,i}(\theta,\omega)|^2=0$ it holds that
    \begin{equation}\llabel{eq: derivative calculate2}
    \begin{split}
        &\displaystyle\frac{\partial}{\partial \theta_j} \textstyle\biggl(\frac{(1-\alpha)\sum_{k=0}^N\alpha^k\bbX_M^{k,i}(\theta,\omega)}{\varepsilon+\sqrt{(1-\beta)\sum_{k=0}^N\beta^k|\bbX_M^{k,i}(\theta,\omega)|^2}}\biggr)=\textstyle\frac{(1-\alpha)\sum_{k=0}^N\alpha^k(\frac{\partial}{\partial\theta_j}\bbX_M^{k,i})(\theta,\omega)}{\varepsilon}\\
        &=\textstyle\frac{(1-\alpha)\sum_{k=0}^N\alpha^k(\frac{\partial}{\partial\theta_j}\bbX_M^{k,i})(\theta,\omega)}{\varepsilon+\sqrt{(1-\beta)\sum_{k=0}^N\beta^k|\bbX_M^{k,i}(\theta,\omega)|^2}}-\textstyle\frac{(1-\beta)\sum_{k=0}^N\beta^k\bbX_M^{k,i}(\theta,\omega)(\frac{\partial}{\partial\theta_j}\bbX_{M}^{k,i})(\theta,\omega)}{\bigl(\varepsilon+\sqrt{(1-\beta)\sum_{k=0}^N\beta^k|\bbX_M^{k,i}(\theta,\omega)|^2}\bigr)^2}.
        \end{split}
    \end{equation}}
    \argument{\lref{def: varrho}}{that 
    \begin{equation}\llabel{eq: moment estimate}
    \begin{split}
       &\textstyle \sup\limits_{\theta=(\theta_1,\dots,\theta_\fd)\in \setV}\sup\limits_{i,j\in\{1,2,\dots,\fd\}} \bigl|(\frac{\partial}{\partial\theta_i}\smalll)(\theta,X_{1,1})\bigr|\leq \varrho\\
        &\text{and}\qquad\textstyle\sup\limits_{\theta=(\theta_1,\dots,\theta_\fd)\in \setV}\sup\limits_{i,j\in\{1,2,\dots,\fd\}} \bigl|(\frac{\partial^2}{\partial\theta_i\theta_j}\smalll)(\theta,X_{1,1})\bigr|\leq \varrho.
        \end{split}
    \end{equation}}
    \argument{\lref{eq: moment estimate};\cref{def: bbX};}{for all $M\in \N$, $n\in \N_0$, $i,j\in \{1,2,\dots,\fd\}$ that
    \begin{equation}\llabel{eq3'}
    \begin{split}
       \textstyle \sup\limits_{\theta=(\theta_1,\dots,\theta_\fd)\in \setV}|\bbX_M^{n,i}(\theta)|\leq \varrho \qqandqq\textstyle \sup\limits_{\theta=(\theta_1,\dots,\theta_\fd)\in \setV}|(\frac{\partial}{\partial \theta_j}\bbX^{n,i}_M)(\theta)|\leq \varrho.
       \end{split}
    \end{equation}}
    \argument{\lref{eq3'};}{for all $p\in \N$ that
    \begin{equation}\llabel{eq3}
    \begin{split}
     \sup_{M\in \N}\sup_{n\in \N_0}\sup_{\theta=(\theta_1,\dots,\theta_\fd)\in \setV}\sup_{i,j\in \{1,2,\dots,\fd\}}\bigl(|\bbX_{M}^{n,i}(\theta)|^p+\textstyle|(\frac{\partial}{\partial\theta_j}\bbX_{M}^{n,i})(\theta)|^p\bigr)\leq 2\varrho^p
       <\infty.
        \end{split}
    \end{equation}}
    \argument{\lref{eq3};the fundamental theorem of calculus}{\cref{item 1: lem negative derivative}\dott}
    \startnewargseq
    \argument{the Cauchy-Schwarz inequality; Jensen's inequality;the assumption that $\alpha^2<\beta$}{that for all $M\in \N$, $N\in \N\cup \{\infty\}$, $\theta\in \R^\fd$, $i\in \{1,2,\dots,\fd\}$, $\omega\in \Omega$ it holds that
    \begin{equation}\llabel{eq2}
    \begin{split}
       \textstyle \biggl|(1-\alpha)\sum\limits_{k=0}^N\alpha^k\bbX_M^{k,i}(\theta,\omega)\biggr|&=\textstyle(1-\alpha)\biggl|\sum\limits_{k=0}^N \alpha^{k}\beta ^{-k/2} \beta^{k/2} \bbX_M^{k,i}(\theta,\omega)\biggr| \\
&\le \textstyle (1-\alpha)\biggl(\sum\limits_{k=0}^N \alpha^{2k}\beta ^{-k}\biggr)^{1/2}\biggl(\sum\limits_{k=0}^N \beta ^{k} |\bbX_M^{k,i}(\theta,\omega)|^2\biggr)^{1/2}\\
&\leq \textstyle \frac{1-\alpha }{\sqrt{1-\alpha^2/\beta}\sqrt{1-\beta}} \sqrt{(1-\beta)\sum\limits_{k=0}^N\beta^k|\bbX_M^{k,i}(\theta,\omega)|^2}.
\end{split}
    \end{equation}}
    \argument{\lref{eq: derivative calculate};\lref{eq: derivative calculate2};\lref{eq2};\lref{eq3'}}{that for all $M\in \N$ and all compact $W\subseteq\R^\fd$ it holds that
    \begin{align}
            &\textstyle \sup\limits_{i,j\in \{1,2,\dots,\fd\}}\E\biggl[\sup\limits_{N\in \N}\sup\limits_{\theta=(\theta_1,\dots,\theta_\fd)\in W}\biggl|\frac{\partial}{\partial \theta_j}\biggl(\frac{\sum_{k=0}^N\alpha^k\bbX_M^{k,i}(\theta)}{\varepsilon+\sqrt{(1-\beta)\sum_{k=0}^N\beta^k|\bbX_M^{k,i}(\theta)|^2}}\biggr)\biggr|\biggr] \nonumber\\
            &\leq \textstyle \sup\limits_{i,j\in \{1,2,\dots,\fd\}}\E\biggl[\sup\limits_{N\in \N}\sup\limits_{\theta=(\theta_1,\dots,\theta_\fd)\in W}\biggl(\biggl|\displaystyle\frac{1-\alpha}{\varepsilon}\textstyle\sum\limits_{k=0}^N\alpha^k(\frac{\partial}{\partial\theta_j}\bbX_M^{k,i})(\theta)\biggr| +\textstyle\biggl|\displaystyle\frac{1-\alpha}{\varepsilon^2}\textstyle\sum\limits_{k=0}^N\alpha^k\bbX_{M}^{k,i}(\theta)(\frac{\partial}{\partial\theta_j}\bbX_{M}^{k,i})(\theta)\biggr|\biggr)\biggr] \nonumber\\
    &\leq \sup_{i,j\in \{1,2,\dots,\fd\}}\E\biggl[\sup\limits_{\theta=(\theta_1,\dots,\theta_\fd)\in W}\biggl(\biggl|\displaystyle\frac{1-\alpha}{\varepsilon}\textstyle\sum\limits_{k=0}^\infty\alpha^k(\frac{\partial}{\partial\theta_j}\bbX_M^{k,i})(\theta)\biggr|+\textstyle\biggl|\displaystyle\frac{1-\alpha}{\varepsilon^2}\textstyle\sum\limits_{k=0}^\infty\alpha^k\bbX_{M}^{k,i}(\theta)(\frac{\partial}{\partial\theta_j}\bbX_{M}^{k,i})(\theta)\biggr|\biggr)\biggr]\nonumber\\
    &<\infty. \llabel{eq4}
    \end{align}}
    \argument{\cref{def: f};\lref{eq4};\lref{eq: derivative calculate};Lebesgue's dominated convergence theorem; the fact that it holds $\P$-a.s.\ that $\liminf_{N\to\infty}\sum_{k=0}^N|\bbX_M^{k,i}(\theta)|^2>0$}{that for all $\theta=(\theta_1,\dots,\theta_\fd)\allowbreak\in \R^\fd$, $i,j\in \{1,2,\dots,\fd\}$, $M\in \N$ it holds that
    \begin{equation}\llabel{derivative}
    \begin{split}
        \displaystyle\frac{\partial}{\partial\theta_j}f_M^i(\theta)&=\textstyle\E\biggl[\frac{(1-\alpha)\sum_{k=0}^\infty\alpha^k(\frac{\partial}{\partial\theta_j}\bbX_M^{k,i})(\theta)}{\varepsilon+\sqrt{(1-\beta)\sum_{k=0}^\infty\beta^k|\bbX_M^{k,i}(\theta)|^2}}\biggr]\\
        &-\E\Biggl[\textstyle\frac{\bigl[\mathbbm 1_{\{0\}}(\sum_{k=0}^\infty |\bbX_M^{k,i}(\theta)|)+(1-\alpha)\sum_{k=0}^\infty\alpha^k\bbX_{M}^{k,i}(\theta)\bigr]\bigl[(1-\beta)\sum_{k=0}^\infty\beta^k\bbX_M^{k,i}(\theta)(\frac{\partial}{\partial\theta_j}\bbX_{M}^{k,i})(\theta)\bigr]}{\bigl(\varepsilon+\sqrt{\mathbbm 1_{\{0\}}(\sum_{k=0}^\infty |\bbX_M^{k,i}(\theta)|)+(1-\beta)\sum_{k=0}^\infty\beta^k|\bbX_M^{k,i}(\theta)|^2}\bigr)^2\sqrt{(1-\beta)\sum_{k=0}^\infty\beta^k|\bbX_M^{k,i}(\theta)|^2}}\Biggr].
        \end{split}
    \end{equation}}
\argument{\lref{eq2};\lref{eq4};\lref{derivative};Lebesgue's dominated convergence theorem;}{items \ref{item 2: lem negative derivative} and \ref{item 3: lem negative derivative}\dott}
\startnewargseq
\argument{\cref{def: bbX};\lref{def: varrho};the fact that $ \P(X_{1,1}\in\setX)=1$;the assumption that $\bbX_{n,m}$, $(n,m)\in (\N_0)^2$, are \iid;}{that for all $n\in \N_0$, $M\in \N$, $i,j\in \{1,2,\dots,\fd\}$, $\theta=(\theta_1,\dots,\theta_\fd)\in \setV$ it holds that
\begin{equation}\llabel{eqq1}
   \textstyle\E\bigl[\bigl(\bbX_M^{n,i}(\theta)-\E[\bbX_{M}^{n,i}(\theta)]\bigr)^2\bigr]=M^{-1}\E\bigl[\big|(\frac{\partial}{\partial\theta_i}\smalll)(\theta,X_{n,1})-\E[(\frac{\partial}{\partial\theta_i}\smalll)(\theta,X_{n,1})]\big|^2\bigr]\leq 4\varrho^2 M^{-1}
\end{equation}
and
\begin{equation}\llabel{eqq2}
\begin{split}
     &\textstyle\E\bigl[\bigl((\frac{\partial}{\partial\theta_j}\bbX_M^{n,i})(\theta)-\E[(\frac{\partial}{\partial\theta_j}\bbX_M^{n,i})(\theta)]\bigr)^2\bigr]\\
     &=\textstyle M^{-1}\E\bigl[\big|(\frac{\partial^2}{\partial\theta_i\partial\theta_j}\smalll)(\theta,X_{n,1})-\E[(\frac{\partial^2}{\partial\theta_i\partial\theta_j}\smalll)(\theta,X_{n,1})]\big|^2\bigr]\leq 4\varrho^2 M^{-1}.
     \end{split}
\end{equation}}
\argument{\cref{lem: explicit batch size};\lref{eq3'};\lref{derivative};\lref{eqq2}}{\cref{item 4: lem negative derivative}\dott}
\startnewargseq
In the following let $A=(A_{i,j})_{i,j\in \{1,2,\dots,\fd\}}\allowbreak\colon\R^\fd\to\R^{\fd\times\fd}$ and $B=(B_{i,j})_{i,j\in \{1,2,\dots,\fd\}}\colon\R^\fd\to\R^{\fd\times\fd}$ satisfy for all $\theta=(\theta_1,\dots,\theta_\fd)\in \R^\fd$, $i,j\in \{1,2,\dots,\fd\}$ that
\begin{equation}\llabel{def: AB}
    A_{i,j}(\theta)=\E[\textstyle(\frac{\partial^2}{\partial \theta_i\partial\theta_j}\smalll)(\theta,X_{1,1})] \displaystyle \qqandqq B_{i,j}(\theta)=\textstyle\frac{\varepsilon}{\bigl(\varepsilon+\bigl|\E[(\frac{\partial}{\partial \theta_i}\smalll)(\theta,X_{1,1})]\bigr|\bigr)^2}\mathbbm 1_{\{j\}}(i).
\end{equation}
\startnewargseq
\argument{\lref{eq: moment estimate};\lref{def: AB};Lebesgue's dominated convergence theorem}{that for all $\theta=(\theta_1,\dots,\theta_\fd)\in \R^\fd$, $i,j\in \{1,2,\dots,\fd\}$ it holds that
\begin{equation}\llabel{eqqq1}
   \textstyle \frac{\partial^2}{\partial \theta_i\partial\theta_j}\E[\smalll(\theta,X_{1,1})]=\E[\textstyle(\frac{\partial}{\partial \theta_i\theta_j}\smalll)(\theta,X_{1,1})]=A_{i,j}(\theta).
\end{equation}}
\argument{\lref{eq: moment estimate};\lref{def: AB}}{that there exists $\bfC\in (0,\infty)$ which satisfies for all $\theta\in \setV$, $x\in \R^\fd$ that
\begin{equation}\llabel{def: bfC}
   \spro{x, B(\theta)x}\geq \frac{\varepsilon}{(\varepsilon+\varrho)^2} \|x\|^2 \ifnocf.
\end{equation}}
\cfout[.]
\startnewargseq 
\argument{\lref{eqqq1};the assumption that $\R^\fd\ni \theta\mapsto \E[\smalll(\theta,X_{1,1})]\in \R$ is $\chi$-\stronglyconvex;}{that for all $\theta,x\in \R^\fd$ it holds that
\begin{equation} \llabel{arg1} 
    \spro{x,A(\theta)x}\geq \chi\|x\|^2.
\end{equation}}
\startnewargseq
\argument{\lref{def: bfC};\lref{arg1}}{that for all $\theta\in \setV$, $x\in \R^\fd$ it holds that
\begin{equation}\llabel{eq: positive matrx}
   \spro{x, A(\theta)B(\theta)x}\geq \frac{\varepsilon \chi}{(\varepsilon+\varrho)^2}\|x\|^2.
\end{equation}}
In the following we prove that there exists $\fM\in \N$ such that for all $\theta\in \setV$, $M\in \N\cap [\fM,\infty)$, $x\in \R^\fd$ it holds that
\begin{equation}\llabel{def: fM}
    \spro{x,(f_M)'(\theta)x}\leq  -\frac{\varepsilon \chi}{2(\varepsilon+\varrho)^2}\|x\|^2.
\end{equation}
\startnewargseq
\argument{\cref{item 4: lem negative derivative};\lref{def: AB}}{that for all $M\in \N$, $i,j\in \{1,2,\dots,\fd\}$, $\theta=(\theta_1,\dots,\theta_\fd)\in \R^\fd$ it holds that
\begin{equation}\llabel{eqt1}
    \textstyle \bigl|\textstyle\frac{\partial}{\partial\theta_j}f_M^i(\theta)+A_{i,j}(\theta)B_{i,i}(\theta)\bigr|\leq 4\varrho\bigl(2\varrho+\frac{2\varrho^2+1}{\varepsilon}+\frac{8\varrho^3+5\varrho}{\varepsilon^2}\bigr)M^{-1/2}.
\end{equation}}
\argument{\lref{eqt1};}{that for all $M\in \N$, $\theta\in \setV$ it holds that
\begin{equation}\llabel{eqt2}
    \textstyle\vertiii{(f_M)'(\theta)+A(\theta)B(\theta)}\leq 4\fd\varrho\bigl(2\varrho+\frac{2\varrho^2+1}{\varepsilon}+\frac{8\varrho^3+5\varrho}{\varepsilon^2}\bigr)M^{-1/2}.
\end{equation}}
\argument{\lref{eqt2};the Cauchy-Schwarz inequality}{that for all $M\in \N$, $\theta\in \setV$, $x\in \R^\fd$ it holds that
\begin{equation}\llabel{eqt3}
    \textstyle|\spro{x, (f_M)'(\theta)+A(\theta)B(\theta),x}|\leq 4\fd\varrho\bigl(2\varrho+\frac{2\varrho^2+1}{\varepsilon}+\frac{8\varrho^3+5\varrho}{\varepsilon^2}\bigr)M^{-1/2}\|x\|^2.
\end{equation}}
\argument{\lref{eqt3};\lref{eq: positive matrx}}{that for all $M\in \N$, $\theta\in \setV$, $x\in \R^\fd$ it holds that
\begin{equation}\llabel{eqt4}
\begin{split}
     \textstyle\spro{x, (f_M)'(\theta),x}&\leq | \spro{x, -A(\theta)B(\theta)x}|+|\spro{x, (f_M)'(\theta)+A(\theta)B(\theta),x}|\\
     &\leq \textstyle\bigl[-\frac{\varepsilon \chi}{(\varepsilon+\varrho)^2}+4\fd\varrho\bigl(2\varrho+\frac{2\varrho^2+1}{\varepsilon}+\frac{8\varrho^3+5\varrho}{\varepsilon^2}\bigr)M^{-1/2}\bigr]\|x\|^2.
     \end{split}
\end{equation}}
\argument{\lref{eqt4};}{\cref{item 5: lem negative derivative}\dott}
\end{aproof}
\subsection{Existence of zeros of the Adam vector field for sufficiently large batch sizes}
In the elementary result in \cref{lem: def deltaM} below we provide for every $\delta \in (0,\infty)$ a lower bound for the probability that the components of the stochastic gradient in \cref{def: bbX} deviate more than the mini-batch size $M \in \N$ multiplied by $\sqrt{ \delta }$ from 0. Our proof of \cref{lem: def deltaM} is strongly based on the arguments in the proof of Lemma 9.2 in \cite[(232)--(233)]{DereichAdamconvergence2024}.
\begin{athm}{lemma}{lem: def deltaM}
    Assume \cref{setting: Adam vector field} and let $\globalmin \in \R^d$, $p\in \N$ satisfy
    \begin{equation}\llabel{def: globmin}
        \E[\smalll(\globalmin, X_{1,1})]=\inf_{\theta\in \R^\fd}\E[\smalll(\theta, X_{1,1})].
    \end{equation}
    Then there exist $\delta\in (0,\infty)$, $\bfM\in \N$ which satisfy for all $M\in \N\cap[\bfM,\infty)$, $i\in \{1,2,\dots,\fd\}$ that
\begin{equation}\llabel{conclude}
   \textstyle \P\bigl(|\sqrt{M}\bbX^{1,i}_M(\globalmin)|\leq \sqrt{\delta}\bigr)\leq \frac{p^{-3/2}}{2}.
\end{equation}
\end{athm}
\begin{aproof}
Throughout this proof let $\varTheta=(\varTheta_1,\dots,\varTheta_\fd)\in \R^\fd$ satisfy $\varTheta=\globalmin$.
     \argument{the fact that $ \P(X_{1,1}\in\setX)=1$;the fact that $\setX$ is bounded;the assumption that  $\smalll\in C^2(\R^\fd\times\R^{\dimX},\R)$}{that for all $p\in \N$ and all compact $U\subseteq \R^\fd$ it holds that
    \begin{equation}\llabel{eq: moment estimate}
        \textstyle\E\biggl[\sup\limits_{\theta=(\theta_1,\dots,\theta_\fd)\in U}\sup\limits_{i\in\{1,2,\dots,\fd\}} \bigl(\bigl|(\frac{\partial}{\partial\theta_i}\smalll)(\theta,X_{1,1})\bigr|^p\bigr)\biggr]<\infty.
    \end{equation}}
    \argument{\lref{def: globmin};\lref{eq: moment estimate};Lebesgue's dominated convergence theorem}{that for all $i\in\{1,2,\dots,\fd\}$ it holds that
    \begin{equation}\llabel{eq2}
        \textstyle\E[(\frac{\partial}{\partial\varTheta_i}\smalll)(\varTheta,X_{1,1}]=\frac{\partial}{\partial \varTheta_i}\E[\smalll(\varTheta,X_{1,1})]=0.
    \end{equation}}
    \argument{the assumption that $X_{m,n}$, $(m,n)\in (\N_0)^2$, are \iid;}{that for all $\theta=(\theta_1,\dots,\theta_\fd)\in \R^\fd$, $n\in \N$ it \llabel{arg1} holds that $(\frac{\partial}{\partial\theta_i}\smalll)(\theta,X_{n,j})$, $j\in\N$, are \iid\ }
    \argument{\lref{arg1};\cref{def: bbX};\lref{eq: moment estimate}; \lref{eq2}; the central limit theorem}{that for all $i\in \{1,2,\dots,\fd\}$ it holds that
\begin{equation}\llabel{eq1}
\begin{split}
\textstyle \sqrt{M}\bbX_M^{1,i}(\varTheta)\rightarrow \mathcal N\bigl(0,\operatorname{Var}((\frac{\partial}{\partial\varTheta_i}\smalll)(\varTheta,X_{1,1}))\bigr)
\end{split}
\end{equation}
in distribution as $M\to\infty$\dott}
In the following for every $i\in \{1,2,\dots,d\}$ let $\bfX_i\colon \Omega\to\R$ be a normal random distributed which satisfies
\begin{equation}\llabel{def: bfX}
    \E[\bfX_i]=0\qqandqq \textstyle \operatorname{Var}(\bfX_i)=\operatorname{Var}((\frac{\partial}{\partial\varTheta_i}\smalll)(\varTheta,X_{1,1})).
\end{equation}
\startnewargseq
\argument{\lref{def: bfX};}{that there exists $\delta\in (0,\infty)$ which satisfies for all $i\in \{1,2,\dots,d\}$ that
\begin{equation}\llabel{def: delta}
    \P(|\bfX_i|\leq \delta)\leq \frac{p^{-3/2}}{4}.
\end{equation}}
\startnewargseq
\argument{\lref{eq1};\lref{def: bfX}}{that for all $i\in \{1,2,\dots,d\}$ it holds that
\begin{equation}\llabel{eq3}
    \textstyle\lim_{M\to\infty}\P(|\sqrt{M}\bbX_M^{1,i}(\varTheta)|\leq \delta)=\P(|\bfX_i|\leq \delta).
\end{equation}}
\argument{\lref{eq3};\lref{def: delta}}{that for all $i\in \{1,2,\dots,d\}$ it holds that
\begin{equation}\llabel{eq4}
   \textstyle \P(|\sqrt{M}\bbX_M^{1,i}(\varTheta)|\leq \delta)\leq\frac{p^{-3/2}}{4}.
\end{equation}}
\argument{\lref{eq4};}{that there exists $\bfM\in \N$ such that for all $M\in \N\cap[\bfM,\infty)$, $i\in \{1,2,\dots,d\}$ it holds that
\begin{equation}\llabel{eq5}
    \textstyle\P(|\sqrt{M}\bbX_M^{1,i}(\varTheta)|\leq \frac{-p^{-3/2}}{4}.
\end{equation}}
\argument{\lref{eq5};}{\lref{conclude}\dott}
\end{aproof}
In \cref{lem: estimate bbM} we provide an upper bound for fractional powers of certain terms appearing in the \Adam\ vector field. \cref{lem: estimate bbM} is inspired by Lemma 8.5 in our preliminary work \cite{DereichAdamconvergence2024}.
\begin{athm}{lemma}{lem: estimate bbM}
    Let $(\Omega,\cF,\P)$ be a probability space, let $\delta,\fc\in (0,\infty)$, $\beta\in (0,1)$, let $X_n\colon \Omega\to [0,\fc]$, $n\in \N_0$, be bounded \iid\ random variables, and let $q\in (0,\beta^{\nicefrac{3}{2}})$ satisfy $\P(X_0\leq \delta)\leq q$. Then
    \begin{equation}\llabel{conclude}
\textstyle\E\Bigl[\bigl[\mathbbm 1_{\{0\}}(\sum_{n=0}^\infty X_n)+(1-\beta)\sum_{n=0}^\infty\beta^nX_n\bigr]^{-3/2}\Bigr]\leq\displaystyle3\bigl(\delta-\beta^{-3/2}q\delta\bigr)^{-3/2}.
    \end{equation}
\end{athm}
\begin{aproof}
    Throughout this proof assume without loss of generality that \allowbreak $\allowbreak\operatorname{Var}(X_0)\allowbreak\neq 0$ (otherwise note that
the fact that $\E[ | X_0 - \E[ X_0 ] |^2 ] = \operatorname{Var}( X_0 ))$,
the fact that $X_n$, $n \in \N$, are \iid, 
the fact that \begin{equation}\P( X_0 > \delta ) = 1 - \P( X_0 \leq \delta ) \geq 1 - q > 1 - \beta^{-3/2 } > 1 - 1 = 0,\end{equation} 
the fact that 
\begin{equation}
 \textstyle [ \mathbbm 1_{ \{0\} }( \sum_{ n = 0 }^{ \infty } \E[ X_n ] ) + ( 1 - \beta ) \sum_{ n = 0 }^{ \infty } \beta^n \E[ X_n ] ]^{ - 3 / 2 }
  =
  [ \mathbbm 1_{ {0} }( \E[ X_0 ] ) + \E[ X_0 ] ]^{ - 3 / 2 },
\end{equation}
and the fact that for all $e \in ( \delta, \infty)$ it holds that
$ 
  [ \mathbbm 1_{ \{0\} }( e ) + e ]^{ - 3 / 2 }
= 
  e^{ - 3 / 2 }
<
  \delta^{ - 3 / 2 }
\leq 
  ( 1 - \beta^{ - 3 / 2 } q )^{ 3 / 2 }
= 
  ( \delta - \beta^{ - 3 / 2 } q \delta )^{ 3 / 2 }
\leq 
  3 ( \delta - \beta^{ - 3 / 2 } q \delta )^{ 3 / 2 }
$
% the fact that
%     \begin{equation}
%     \begin{split}
%     &\textstyle\E\Bigl[\bigl[\mathbbm 1_{\{0\}}(\sum_{n=0}^\infty X_n)+(1-\beta)\sum_{n=0}^\infty\beta^nX_n\bigr]^{-3/2}\Bigr]\\
%    &\textstyle=\E\Bigl[\bigl[\mathbbm 1_{\{0\}}(\E[X_0])+(1-\beta)\sum_{n=0}^\infty\beta^n\E[X_0]\bigr]^{-3/2}\Bigr]\\
%    &\textstyle=\E\Bigl[\bigl[\mathbbm 1_{\{0\}}(\E[X_0])+\E[X_0]\bigr]^{-3/2}\Bigr]\leq (\E[X_0])^{-3/2}\leq \delta^{-3/2}.
%    \end{split}  
%     \end{equation}
    establish \lref{conclude}), assume without loss of generality that $q=\P(X_0\leq \delta)$,
    let $\bbX_n\colon \Omega\to\R$, $n\in \N_0$, satisfy for all $n\in \N_0$ that 
    \begin{equation}\llabel{def: bbX}
        \bbX_n=\delta\mathbbm 1_{(\delta,\infty)}(X_n).
    \end{equation}
    \argument{the fact that $\operatorname{Var}(X_0)\neq 0$;the assumption that $X_n$, $n\in \N_0$, are \iid}{that for all $n\in \N_0$ it holds that
    \begin{equation}\llabel{eq1}
        \P(X_n=0)=\P(X_0=0)<1.
    \end{equation}}
    \argument{\lref{eq1}; the fact that for all $n\in \N_0$ it holds that $X_n\geq 0$}{that
    \begin{equation}\llabel{eq2}
         \begin{split}
         &\P\bigl(\textstyle\mathbbm 1_{\{0\}}(\sum_{n=0}^\infty X_n)=0\bigr)\\
        &=\P\bigl(\textstyle\mathbbm 1_{\{0\}}(\sum_{n=0}^\infty |X_n|)=0\bigr)=\P\bigl( \sum_{ n = 0 }^{ \infty } | X_n  | \neq 0 \bigr)=
\P\bigl( \Omega \backslash \{ \sum_{ n = 0 }^{ \infty } | X_n  | = 0 \} \bigr)\\
        &=\textstyle\P( \Omega \backslash \{ \forall \, n \in \N_0 \colon X_n = 0 \} )=1-\P(\forall\, n\in \N_0\colon X_n=0)=1-\bigl[\prod_{n=0}^\infty \P(X_n=0)\bigr]\\
        &\textstyle= 1 - \lim_{ n \to \infty } [ \P( X_0 =0) ]^n=1.
        \end{split}
    \end{equation}}
    \argument{\lref{def: bbX};the assumption that $X_n$, $n\in \N_0$, are \iid}{that $\bbX_n$, $n\in \N_0$, are \llabel{arg1} \iid\dott}
    \argument{\lref{arg1};\lref{def: bbX}; the fact $\P(X_0\leq \delta)<1$; the fact that $\delta>0$}{that for all $n\in \N_0$ it holds that
    \begin{equation}\llabel{eq3}
        \P(\bbX_n=0)=\P(\bbX_0=0)=\P(X_0\leq\delta)<1.
    \end{equation}}
    \argument{\lref{eq3};the fact that $\bbX_n$, $n\in \N_0$, are \iid; the fact that for all $n\in \N_0$ it holds that $\bbX_n\geq 0$}{that 
    \begin{equation}\llabel{eq4}
         \begin{split}
         &\P\bigl(\textstyle\mathbbm 1_{\{0\}}(\sum_{n=0}^\infty \bbX_n)=0\bigr)\\
        &=\P\bigl(\textstyle\mathbbm 1_{\{0\}}(\sum_{n=0}^\infty |\bbX_n|)=0\bigr)=\P\bigl( \sum_{ n = 0 }^{ \infty } | \bbX_n  | \neq 0 \bigr)=
\P\bigl( \Omega \backslash \{ \sum_{ n = 0 }^{ \infty } | \bbX_n  | = 0 \} \bigr)\\
        &=\textstyle\P( \Omega \backslash \{ \forall \, n \in \N_0 \colon \bbX_n = 0 \} )=1-\P(\forall\, n\in \N_0\colon \bbX_n=0)=1-\bigl[\prod_{n=0}^\infty \P(\bbX_n=0)\bigr]\\
        &\textstyle=1-\bigl[\prod_{n=0}^\infty \P(X_n=0)\bigr]=1 - \lim_{ n \to \infty } [ \P( X_0 =0) ]^n=1.
        \end{split}
    \end{equation}}
    \argument{\lref{eq4};\lref{eq2}}{that
    \begin{equation}\llabel{eq5}
    \P\bigl(\textstyle\mathbbm 1_{\{0\}}(\sum_{n=0}^\infty X_n)=\mathbbm 1_{\{0\}}(\sum_{n=0}^\infty \bbX_n)=0\bigr)=1.
    \end{equation}}
    \argument{\lref{def: bbX}; the fact that for all $n\in \N_0$ it holds that $X_n\geq 0$}{that for all $n\in \N_0$ it holds that 
    \begin{equation}\llabel{eq6}
        \bbX_n\leq X_n.
    \end{equation}}
    \argument{\lref{eq6};\lref{eq4}}{that
    \begin{equation}\llabel{eq7}
        \begin{split}
            &\textstyle\E\Bigl[\bigl[\mathbbm 1_{\{0\}}(\sum_{n=0}^\infty X_n)+(1-\beta)\sum_{n=0}^\infty\beta^nX_n\bigr]^{-3/2}\Bigr]\\
            &\leq \textstyle\E\Bigl[\bigl[\mathbbm 1_{\{0\}}(\sum_{n=0}^\infty \bbX_n)+(1-\beta)\sum_{n=0}^\infty\beta^n\bbX_n\bigr]^{-3/2}\Bigr].
        \end{split}
    \end{equation}}
    \argument{\lref{def: bbX};the assumption that $\P(\bbX_0\leq \delta)=q$}{that
    \begin{equation}\llabel{eq8}
        \P(\bbX_0=0)=\P(X_0\leq \delta)=q=1-\P(X_0>\delta)=1-\P(\bbX_0=\delta).
    \end{equation}}
    \argument{\lref{eq8}; the fact that $\bbX_n$, $n\in\N_0$, are \iid}{that for all $p\in \R$ it holds that
    \begin{equation}\llabel{eq9}
    \begin{split}
       & \textstyle\E\Bigl[\bigl[\mathbbm 1_{\{0\}}(\sum_{n=0}^\infty \bbX_n)+(1-\beta)\sum_{n=0}^\infty\beta^n\bbX_n\bigr]^{p}\Bigr]\\
       &=\textstyle\E\Bigl[\bigl[\mathbbm 1_{\{0\}}(\sum_{n=0}^\infty \bbX_n)+(1-\beta)\sum_{n=0}^\infty\beta^n\bbX_n\bigr]^{p}\mathbbm 1_{\{0\}}(\bbX_0)\Bigr]\\
       &+\textstyle\E\Bigl[\bigl[\mathbbm 1_{\{0\}}(\sum_{n=0}^\infty \bbX_n)+(1-\beta)\sum_{n=0}^\infty\beta^n\bbX_n\bigr]^{p}\mathbbm 1_{\{\delta\}}(\bbX_0)\Bigr]\\
       &=\textstyle\E\Bigl[\bigl[\mathbbm 1_{\{0\}}(\sum_{n=1}^\infty \bbX_n)+(1-\beta)\sum_{n=1}^\infty\beta^n\bbX_n\bigr]^{p}\mathbbm 1_{\{0\}}(\bbX_0)\Bigr]\\
       &+\textstyle\E\Bigl[\bigl[(1-\beta)\delta+(1-\beta)\sum_{n=1}^\infty\beta^n\bbX_n\bigr]^{p}\mathbbm 1_{\{\delta\}}(\bbX_0)\Bigr]\\
       &=\textstyle\E\Bigl[\bigl[\mathbbm 1_{\{0\}}(\sum_{n=1}^\infty \bbX_n)+(1-\beta)\sum_{n=1}^\infty\beta^n\bbX_n\bigr]^{p}\Bigr]\E\bigl[\mathbbm 1_{\{0\}}(\bbX_0)\bigr]\\
       &+\textstyle\E\Bigl[\bigl[(1-\beta)\delta+(1-\beta)\sum_{n=1}^\infty\beta^n\bbX_n\bigr]^{p}\Bigr]\E\bigl[\mathbbm 1_{\{\delta\}}(\bbX_0)\bigr]\\
       &=\textstyle\E\Bigl[\bigl[\mathbbm 1_{\{0\}}(\sum_{n=1}^\infty \bbX_n)+(1-\beta)\sum_{n=1}^\infty\beta^n\bbX_n\bigr]^{p}\Bigr]\P(\bbX_0=0)\\
       &+\textstyle\E\Bigl[\bigl[(1-\beta)\delta+(1-\beta)\sum_{n=1}^\infty\beta^n\bbX_n\bigr]^{p}\Bigr]\P(\bbX_0=\delta)\\
       &=q\,\textstyle\E\Bigl[\bigl[\mathbbm 1_{\{0\}}(\sum_{n=1}^\infty \bbX_n)+(1-\beta)\sum_{n=1}^\infty\beta^n\bbX_n\bigr]^{p}\Bigr]\\
       &+\textstyle(1-q)\,\E\Bigl[\bigl[(1-\beta)\delta+(1-\beta)\sum_{n=1}^\infty\beta^n\bbX_n\bigr]^{p}\Bigr].
        \end{split}
    \end{equation}}
    \argument{the fact that $\bbX_n$, $n\in \N_0$, are \iid;}{that $(\bbX_n)_{n\in \N}$ and $(\bbX_n)_{n\in \N_0}$ are \llabel{arg2} identically distributed\dott}
    \argument{\lref{arg2};\lref{eq9};}{for all $p\in \R$ that
    \begin{equation}\llabel{eq10}
    \begin{split}
        & \textstyle\E\Bigl[\bigl[\mathbbm 1_{\{0\}}(\sum_{n=0}^\infty \bbX_n)+(1-\beta)\sum_{n=0}^\infty\beta^n\bbX_n\bigr]^{p}\Bigr]\\
        &=q\,\textstyle\E\Bigl[\bigl[\mathbbm 1_{\{0\}}(\sum_{n=1}^\infty \bbX_{n-1})+(1-\beta)\sum_{n=1}^\infty\beta^n\bbX_{n-1}\bigr]^{p}\Bigr]\\
       &+\textstyle(1-q)\,\E\Bigl[\bigl[(1-\beta)\delta+(1-\beta)\sum_{n=1}^\infty\beta^n\bbX_{n-1}\bigr]^{p}\Bigr]\\
       &=q\beta^p\,\textstyle\E\Bigl[\bigl[\beta^{-p}\mathbbm 1_{\{0\}}(\sum_{n=1}^\infty \bbX_{n-1})+(1-\beta)\sum_{n=1}^\infty\beta^{n-1}\bbX_{n-1}\bigr]^{p}\Bigr]\\
       &+\textstyle(1-q)\beta^p\,\E\Bigl[\bigl[(1-\beta)\beta^{-1}\delta+(1-\beta)\sum_{n=1}^\infty\beta^{n-1}\bbX_{n-1}\bigr]^{p}\Bigr]\\
       &=q\beta^p\,\textstyle\E\Bigl[\bigl[\mathbbm 1_{\{0\}}(\sum_{n=0}^\infty \bbX_{n})+(1-\beta)\sum_{n=0}^\infty\beta^{n}\bbX_{n}\bigr]^{p}\Bigr]\\
       &+\textstyle(1-q)\beta^p\,\E\Bigl[\bigl[(1-\beta)\beta^{-1}\delta+(1-\beta)\sum_{n=0}^\infty\beta^{n}\bbX_{n}\bigr]^{p}\Bigr].
       \end{split}
    \end{equation}}
    \argument{\lref{eq10};}{for all $p\in \R$ with $q\beta^p<1$ that
    \begin{equation}\llabel{eq11}
    \begin{split}
         & \textstyle(1-q\beta^p)\E\Bigl[\bigl[\mathbbm 1_{\{0\}}(\sum_{n=0}^\infty \bbX_n)+(1-\beta)\sum_{n=0}^\infty\beta^n\bbX_n\bigr]^{p}\Bigr]\\
         &=\textstyle(1-q)\beta^p\,\E\Bigl[\bigl[(1-\beta)\beta^{-1}\delta+(1-\beta)\sum_{n=0}^\infty\beta^{n}\bbX_{n}\bigr]^{p}\Bigr].
         \end{split}
    \end{equation}}
    \argument{\lref{eq11};\lref{eq5}; the fact that $q\beta^{-1/2}<1$ the fact that for all $x,y,z\in [0,\infty)$ with $x+z>0$ and $y>z$ it holds that
    \begin{equation}
        (x+z)^{-1/2}-(x+y)^{-1/2}=\frac{y-z}{\bigl(\sqrt{x+z}+\sqrt{x+y}\bigr)\sqrt{(x+z)(x+y)}}\geq \frac{y-z}{2(x+y)^{3/2}}
    \end{equation}}{that
    \begin{equation}\llabel{eq12}
    \begin{split}
       &(\beta^{-1/2}-1)\, \textstyle\E\Bigl[\bigl[(1-\beta)\beta^{-1}\delta+(1-\beta)\sum_{n=0}^\infty\beta^{n}\bbX_{n}\bigr]^{-1/2}\Bigr]\\
       &= \textstyle(1-q)\beta^{-1/2}\E\Bigl[\bigl[\mathbbm 1_{\{0\}}(\sum_{n=0}^\infty \bbX_n)+(1-\beta)\sum_{n=0}^\infty\beta^n\bbX_n\bigr]^{-1/2}\Bigr]\\
       &-(1-q)\beta^{-1/2} \textstyle\E\Bigl[\bigl[(1-\beta)\beta^{-1}\delta+(1-\beta)\sum_{n=0}^\infty\beta^{n}\bbX_{n}\bigr]^{-1/2}\Bigr]\\
       &\geq\textstyle (1-q)\beta^{-1/2}\textstyle\E\Bigl[\Bigl[\frac{(1-\beta)\beta^{-1}\delta-\mathbbm 1_{\{0\}}(\sum_{n=0}^\infty \bbX_n)}{2}\Bigr]\bigl[(1-\beta)\beta^{-1}\delta+(1-\beta)\sum_{n=0}^\infty\beta^{n}\bbX_{n}\bigr]^{-3/2}\Bigr]\\
       &= \textstyle\bigl[\frac{(1-q)\beta^{-1/2}(1-\beta)\beta^{-1}\delta}{2}\bigr]\textstyle\E\Bigl[\bigl[(1-\beta)\beta^{-1}\delta+(1-\beta)\sum_{n=0}^\infty\beta^{n}\bbX_{n}\bigr]^{-3/2}\Bigr].
       \end{split}
    \end{equation}}
    \argument{\lref{eq12};Holder inequality}{that 
    \begin{equation}\llabel{eq13}
    \begin{split}
        &\E\Bigl[\bigl[(1-\beta)\beta^{-1}\delta+(1-\beta)\sum_{n=0}^\infty\beta^{n}\bbX_{n}\bigr]^{-3/2}\Bigr]\\
        &\leq \displaystyle\Bigl[\frac{2(\beta^{-1/2}-1)}{(1-q)\beta^{-1/2}(1-\beta)\beta^{-1}\delta}\Bigr]\textstyle\E\Bigl[\bigl[(1-\beta)\beta^{-1}\delta+(1-\beta)\sum_{n=0}^\infty\beta^{n}\bbX_{n}\bigr]^{-1/2}\Bigr]\\
        &= \biggl[\frac{2\beta}{(1-q)(1+\sqrt{\beta})\delta}\biggr]\textstyle\E\Bigl[\bigl[(1-\beta)\beta^{-1}\delta+(1-\beta)\sum_{n=0}^\infty\beta^{n}\bbX_{n}\bigr]^{-1/2}\Bigr]\\
        &\textstyle\leq 2\beta (1-q)^{-1}\delta^{-1}\Bigl(\E\Bigl[\bigl[(1-\beta)\beta^{-1}\delta+(1-\beta)\sum_{n=0}^\infty\beta^{n}\bbX_{n}\bigr]^{-3/2}\Bigr]\Bigr)^{1/3}.
        \end{split}
    \end{equation}}
    \argument{\lref{eq13};}{that
    \begin{equation}\llabel{eq13'}
        \begin{split}
            &\textstyle\Bigl(\E\Bigl[\bigl[(1-\beta)\beta^{-1}\delta+(1-\beta)\sum_{n=0}^\infty\beta^{n}\bbX_{n}\bigr]^{-3/2}\Bigr]\Bigr)^{2/3}\leq 2\beta (1-q)^{-1}\delta^{-1}.
        \end{split}
    \end{equation}}
    \argument{\lref{eq13'};}{that
    \begin{equation}\llabel{eq14}
    \begin{split}
        \textstyle\E\Bigl[\bigl[(1-\beta)\beta^{-1}\delta+(1-\beta)\sum_{n=0}^\infty\beta^{n}\bbX_{n}\bigr]^{-3/2}\Bigr]&\leq 2^{3/2}\beta^{3/2}(1-q)^{-3/2}\delta^{-3/2}\\
        &\leq 3\beta^{3/2}(1-q)^{-3/2}\delta^{-3/2}.
        \end{split}
    \end{equation}}
    \argument{\lref{eq14};\lref{eq11};the fact that $q\beta^{-3/2}<1$}{that
    \begin{equation}\llabel{eq15}
    \begin{split}
        &\textstyle\E\Bigl[\bigl[\mathbbm 1_{\{0\}}(\sum_{n=0}^\infty \bbX_n)+(1-\beta)\sum_{n=0}^\infty\beta^n\bbX_n\bigr]^{-3/2}\Bigr]\\
        &\leq \biggl[\frac{(1-q)\beta^{-3/2}}{1-q\beta^{-3/2}}\biggr]\textstyle\E\Bigl[\bigl[(1-\beta)\beta^{-1}\delta+(1-\beta)\sum_{n=0}^\infty\beta^{n}\bbX_{n}\bigr]^{-3/2}\Bigr]\\
        &\leq \frac{3}{(1-q)^{1/2}(1-q\beta^{-3/2})\delta^{3/2}}\leq \frac{3}{(1-q\beta^{-3/2})^{3/2}\delta^{3/2}}.
        \end{split}
    \end{equation}}
    \argument{\lref{eq15};\lref{eq7}}{\lref{conclude}\dott}
    \end{aproof}
In \cref{item 2: existence of solution} in \cref{lem: existence of solution} we provide for every sufficiently large size of the mini-batch $M \in \N$ an upper bound for the distance between a zero of the \Adam\ vector field and the global minimizer $\vartheta$
of the \SOP\ in \cref{def: globmin: existence of solution}. \cref{lem: existence of solution} is inspired by Lemma 9.2 in our preliminary work \cite{DereichAdamconvergence2024} in which a special case of \cref{lem: existence of solution} is considered.
\begin{athm}{lemma}{lem: existence of solution}
    Assume \cref{setting: Adam vector field}, let $\setTheta\subseteq \R^d$ be open and convex, let $\globalmin \in \setTheta$, $\rho\in (0,\infty)$ satisfy $\setTheta\subseteq \{x\in \R^\fd\colon \|x-\globalmin\|\leq \rho\}$ and
    \begin{equation}\label{def: globmin: existence of solution}
        \E[\smalll(\globalmin, X_{1,1})]=\inf_{\theta\in \R^\fd}\E[\smalll(\theta, X_{1,1})],
    \end{equation}
    let $\varrho\in \R$ satisfy for all $\theta=(\theta_1,\dots,\theta_\fd)\in \setTheta$, $i,j\in\{1,2,\dots,\fd\}$, $x\in \setX$ that
     \begin{equation}\llabel{def: varrho}
   \textstyle \varrho\geq \max\bigl\{|(\frac{\partial^2}{\partial\theta_i\partial\theta_j}\smalll)(\theta,x)|,|(\frac{\partial}{\partial\theta_i}\smalll)(\theta,x)|\bigr\},
   \end{equation}
    let $\bfB\in \R$ be the constant in the Burkholder-Davis-Gundy inequality for the third moment, let $\delta\in (0,\infty)$, $\bfM\in \N$ satisfy for all $M\in \N\cap[\bfM,\infty)$, $i\in \{1,2,\dots,\fd\}$ that
\begin{equation}\llabel{def: delta}
   \textstyle \P\bigl(|\sqrt{M}\bbX^{1,i}_M(\vartheta)|\leq \sqrt{\delta}\bigr)\leq \frac{\beta^{3/2}}{2}.
\end{equation}
     and let $\scrM\in \R$ satisfy 
     \begin{equation}\llabel{def: batch size}
\textstyle\scrM=18\rho^{-1}\beta^{-5/2}(1-\beta)\varepsilon^{-2}\varrho^3\delta^{-3/2}(\delta+6(\bfB+3)\varrho^2)(\varepsilon+\varrho)^2(\varepsilon \chi)^{-1}.
     \end{equation}
    Then there exists $(\roott_M)_{M\in \N}\subseteq\setTheta$ such that
    \begin{enumerate}[label=(\roman*)]
        \item 
        \label{item 1: existence of solution} 
        it holds for all $M\in \N\cap[\max\{\bfM,\scrM,64\fd^2\varrho^2(2\varrho+\varepsilon^{-1}(2\varrho+1)+\varepsilon^{-2}(8\varrho^3+5\varrho))^2(\varepsilon+\varrho)^4(\varepsilon \chi)^{-2}\},\allowbreak\infty)$ that 
        \begin{equation}
            f_M(\theta_M)=0
        \end{equation} and
        \item \label{item 2: existence of solution} it holds for all $M\in \N$ that 
            $\|\theta_M-\globalmin\|\leq \rho\scrM M^{-1}$. 
        \end{enumerate}
\end{athm}
\begin{aproof}
Throughout this proof let $h\colon [0,\infty)\to \R$ satisfy for all $x\in \R$ that
\begin{equation}\llabel{def: h}
    h(x)=\frac{1}{\varepsilon+\sqrt{x}},
\end{equation}
for every $M\in \N$, $i\in \{1,2,\dots,\fd\}$ let $\bbM_M^i\colon \R^\fd\times\Omega\to\R$ satisfy for all $\omega\in \Omega$, $\theta\in \R^\fd$ that
\begin{equation}\llabel{def: bbV}
    \bbM_M^i(\theta,\omega)=\textstyle(1-\beta)\sum\limits_{k=0}^\infty\beta^k|\bbX_M^{k,i}(\theta,\omega)|^2,
\end{equation}
let $\varTheta=(\varTheta_1,\dots,\varTheta_\fd)\in \R^\fd$ satisfy $\varTheta=\globalmin$, and assume without loss of generality that $\operatorname{Var}((\frac{\partial}{\partial\varTheta_i}\smalll)(\varTheta,\allowbreak X_{1,1}))\neq 0$ (otherwise for every $M\in \N$ let $\roott_M=\globalmin$).
In the following let $\fc\in \R$ satisfy
\begin{equation}\llabel{def: batch size 2}
\fc=3d\beta^{-5/2}(1-\beta)\varepsilon^{-2}\varrho^3(\delta^{-1/2}+6\delta^{-3/2}(\bfB+3)\varrho^2).
\end{equation}
     \argument{the fact that $ \P(X_{1,1}\in\setX)=1$;the fact that $\setX$ is bounded;the assumption that  $\smalll\in C^2(\R^\fd\times\R^{\dimX},\R)$}{that for all $p\in \N$ and all compact $U\subseteq \R^\fd$ it holds that
    \begin{equation}\llabel{eq: moment estimate}
        \textstyle\E\biggl[\sup\limits_{\theta=(\theta_1,\dots,\theta_\fd)\in U}\sup\limits_{i\in\{1,2,\dots,\fd\}} \bigl(\bigl|(\frac{\partial}{\partial\theta_i}\smalll)(\theta,X_{1,1})\bigr|^p\bigr)\biggr]<\infty.
    \end{equation}}
    \argument{\cref{def: globmin: existence of solution};\lref{eq: moment estimate};Lebesgue's dominated convergence theorem}{that for all $i\in\{1,2,\dots,\fd\}$ it holds that
    \begin{equation}\llabel{eq2''}
        \textstyle\E[(\frac{\partial}{\partial\varTheta_i}\smalll)(\varTheta,X_{1,1}]=\frac{\partial}{\partial \varTheta_i}\E[\smalll(\varTheta,X_{1,1})]=0.
    \end{equation}}
    \argument{\lref{eq2''};\cref{def: bbX};the assumption that $X_{n,m}$, $(n,m)\in\N^2$, are \iid;}{that for all $M\in \N$, $i\in \{1,2,\dots,d\}$ that
    \begin{equation}\llabel{eq2}
        \E[\bbX^{1,i}_M(\varTheta)]=0.
    \end{equation}}
\argument{\lref{def: delta};}{that for all $M\in \N\cap[\bfM,\infty)$, $i\in \{1,2,\dots,\fd\}$ it holds that
\begin{equation}\llabel{def: delta2}
    \P\bigl(|\bbX^{1,i}_M(\varTheta)|^2\leq \textstyle \frac{\delta}{M}\bigr)\leq \frac 12 \beta^{3/2}.
\end{equation}}
\argument{\cref{def: bbX};\lref{def: varrho};the fact that $ X_{1,1}\in\setX $}{that 
\begin{equation}\llabel{eq3'}
    \sup_{n\in \N_0}\sup_{M\in \N}\sup_{i\in \{1,2,\dots,\fd\}}|\bbX_M^{n,i}(\varTheta)|\leq \varrho.
\end{equation}}
\argument{\lref{eq3'};\lref{def: bbV};\lref{def: delta2};\cref{lem: estimate bbM} (applied for every $M\in \N\cap[\bfM,\infty)$, $i\in \{1,2,\dots,\fd\}$ with $\Omega\curvearrowleft\Omega$, $\cF\curvearrowleft\cF$, $\P\curvearrowleft \P$, $\delta\curvearrowleft \frac{\delta}{M}$, $\beta\curvearrowleft\beta$, $(X_n)_{n\in \N_0}\curvearrowleft (\Omega\ni\omega\mapsto |\bbX_M^{n,i}(\varTheta,\omega)|^2\in [0,\varrho^2])_{n\in \N_0}$, $\fc\curvearrowleft \varrho^2$, $q\curvearrowleft\frac{\beta^{3/2}}{2}$ in the notation of \cref{lem: estimate bbM});the fact that for all $M\in \N$, $i\in \{1,2,\dots,\fd\}$ it holds that $\bbX_{M}^{n,i}(\varTheta)$, $n\in \N_0$, are \iid;Holder inequality}{that for all $M\in \N\cap[\bfM,\infty)$ $i\in \{1,2,\dots,\fd\}$ it holds that
\begin{equation}\llabel{eq3}
\begin{split}
    &\E\bigl[|\bbM_{M}^i(\varTheta)|^{-3/2}\bigr]=\textstyle\E\biggl[\biggl((1-\beta)\sum\limits_{k=0}^\infty\beta^k|\bbX_M^{k,i}(\varTheta)|^2\biggr)^{-p}\biggr]\leq \frac{3(M/\delta)^{3/2}}{(1-\frac{1}{2}\beta^{3/2}\beta^{-3/2})^{3/2}}\leq 9\bigl(\frac{M}{\delta}\bigr)^{3/2}
    \end{split}
\end{equation}
\begin{equation}\llabel{eq3.1}
   \text{and}\qquad  \E\bigl[|\bbM_{M}^i(\varTheta)|^{-1/2}\bigr]\leq \Bigl(\E\bigl[|\bbM_{M}^i(\varTheta)|^{-3/2}\bigr]\Bigr)^{1/3}\leq \textstyle3\bigl(\frac{M}{\delta}\bigr)^{1/2}.
\end{equation}}
\argument{\cref{def: f};\lref{def: varrho};\lref{def: h};\lref{eq2};\lref{eq3'};\lref{eq3};\cite[Theorem 8.3]{DereichAdamconvergence2024};the fact that for all $M\in \N$, $i\in \{1,2,\dots,\fd\}$ it holds that $\bbX_{M}^{n,i}$, $n\in \N$, are i.i.d}{that for all $i\in \{1,2,\dots,\fd\}$, $M\in \N$ it holds that
\begin{equation}\llabel{eq4}
\begin{split}
  &  |f_M^i(\varTheta)|\leq \beta^{-3/2}\frac{(1-\alpha)(1-\beta)}{1-\alpha\beta}\varrho^3\E[|h'(\bbM_M^i(\varTheta))|]\,M^{-2}\\
    &+\beta^{-5/2}\frac{(1-\alpha)(1-\beta)^2}{1-\alpha\beta^2} \E[h''(\bbM_M^i(\varTheta))] \bigl(\textstyle2\bfB\varrho^5M^{-5/2}+5\varrho^5M^{-3}+\varrho^5M^{-4}\bigr).
    \end{split}
\end{equation}}
\startnewargseq
\argument{\lref{eq4};the fact that for all $x\in \R$ it holds that $|h'(x)|\leq \varepsilon^{-2}x^{-1/2}$ and $|h''(x)|\leq \varepsilon^{-2}x^{-3/2}$}{that for all $M\in \N\cap[\bfM,\infty)$ $i\in \{1,2,\dots,\fd\}$ it holds that
\begin{equation}\llabel{eq5}
\begin{split}
   & |f_M^i(\varTheta)|\le  \beta^{-3/2}\frac{(1-\alpha)(1-\beta)}{1-\alpha\beta}\varrho^3\varepsilon^{-2} \E[(\bbM_M^i(\varTheta))^{-1/2}]M^{-2}\\
    &+ \beta^{-5/2}\frac{(1-\alpha)(1-\beta)^2}{1-\alpha\beta^2}\varepsilon^{-2}\E[(\bbM_M^i(\varTheta))^{-3/2}]\bigl(\textstyle2\bfB\varrho^5M^{-5/2}+5\varrho^5M^{-3}+\varrho^5M^{-4}\bigr).
    \end{split}
\end{equation}}
\argument{\lref{eq5};\lref{eq3};\lref{eq3.1}; the fact that $0<\alpha<\sqrt{\beta}<1$}{that for all $M\in \N\cap[\bfM,\infty)$, $i\in \{1,2,\dots,\fd\}$ it holds that
\begin{equation}\llabel{eq5.5'}
    \begin{split}
         |f_M^i(\varTheta)|
         &\textstyle\leq \beta^{-3/2}\frac{(1-\alpha)(1-\beta)}{1-\alpha\beta}\varrho^3\varepsilon^{-2}\frac{3M^{1/2}}{\delta^{1/2}}M^{-2}\\
         &
         \textstyle+\beta^{-5/2}\frac{(1-\alpha)(1-\beta)^2}{1-\alpha\beta^2}\varepsilon^{-2}\frac{9M^{3/2}}{\delta^{3/2}}\bigl(\textstyle2\bfB\varrho^5M^{-5/2}+5\varrho^5M^{-3}+\varrho^5M^{-4}\bigr)\\
         &\leq 3\beta^{-3/2}(1-\beta)\varepsilon^{-2}\varrho^3\delta^{-1/2}M^{-1}+9\beta^{-5/2}(1-\beta)\varepsilon^{-2}\delta^{-3/2}(2\bfB+6)\varrho^5M^{-1}.
    \end{split}
\end{equation}}
\argument{\lref{eq5.5'};\lref{def: batch size 2}; the fact that for all $x=(x_1,\dots,x_d)\in \R^d$ it holds that $\|x\|\leq \sum_{i=1}^d|x_i|$}{that for all $M\in \N\cap[\bfM,\infty)$ it holds that
\begin{equation}\llabel{eq5.5}
\begin{split}
    \|f_M(\varTheta)\|&\leq 3d\beta^{-3/2}(1-\beta)\varepsilon^{-2}\varrho^3\delta^{-1/2}M^{-1}+9d\beta^{-5/2}(1-\beta)\varepsilon^{-2}\delta^{-3/2}(2\bfB+6)\varrho^5M^{-1}\\
    &\leq \fc M^{-1}.
    \end{split}
\end{equation}}
\startnewargseq
 \argument{;}{ that for all $M\in \N\cap[64\fd^2\varrho^2(2\varrho+\varepsilon^{-1}(2\varrho+1)+\varepsilon^{-2}(8\varrho^3+5\varrho))^2(\varepsilon+\varrho)^4(\varepsilon \chi)^{-2},\infty)$ it holds that
    \begin{equation}\llabel{eqM}
    \begin{split}
   & 
        -\frac{\varepsilon \chi}{(\varepsilon+\varrho)^2}+4\fd\varrho\biggl(2\varrho+\frac{2\varrho^2+1}{\varepsilon}+\frac{8\varrho^3+5\varrho}{\varepsilon^2}\biggr)M^{-1/2}\\
        &\leq  
        -\frac{\varepsilon \chi}{(\varepsilon+\varrho)^2}+4\fd\varrho\biggl(2\varrho+\frac{2\varrho^2+1}{\varepsilon}+\frac{8\varrho^3+5\varrho}{\varepsilon^2}\biggr)
        \\
        &\quad\cdot\bigl[64\fd^2\varrho^2(2\varrho+(2\varrho+1)\varepsilon^{-1}+(8\varrho^3+6\varrho)\varepsilon^{-2})^2(\varepsilon+\varrho)^4(\varepsilon \chi)^{-2}\big]^{-1/2}\\
        &=  -\frac{\varepsilon \chi}{(\varepsilon+\varrho)^2}+ \frac{\varepsilon \chi}{2(\varepsilon+\varrho)^2}=-\frac{\varepsilon \chi}{2(\varepsilon+\varrho)^2}.
        \end{split}
    \end{equation}}
    \argument{\cref{item 4: lem negative derivative} in \cref{lem: negative derivative};\cref{def: globmin: existence of solution};\lref{def: batch size};\lref{eqM}}{that for all $M\in \N\cap [64\fd^2\varrho^2(2\varrho+\varepsilon^{-1}(2\varrho+1)+\varepsilon^{-2}(8\varrho^3+5\varrho))^2(\varepsilon+\varrho)^4(\varepsilon \chi)^{-2},\infty)$, $\theta\in \R^\fd$, $x\in \R^\fd$ with $\|\theta-\varTheta\|\leq \rho$ it holds that
        \begin{equation}\llabel{def: fc}
            \spro{x, (f_M)'(\theta)x}\leq -\frac{\varepsilon \chi}{2(\varepsilon+\varrho)^2}\|x\|^2 .
        \end{equation}}
    \argument{\lref{def: fc}; the Cauchy-Schwarz inequality;}{that for all $M\in \N\cap[64\fd^2\varrho^2(2\varrho+\varepsilon^{-1}(2\varrho+1)+\varepsilon^{-2}(8\varrho^3+5\varrho))^2(\varepsilon+\varrho)^4(\varepsilon \chi)^{-2},\infty)$, $\theta\in \R^\fd$ with $\|\theta-\varTheta\|\leq \rho$ it holds that
\begin{equation}\llabel{eeq1}
\begin{split}
    \frac{\varepsilon \chi}{2(\varepsilon+\varrho)^2}\|\theta-\varTheta\|^2\leq -\spro{\theta-\varTheta, f_M(\theta)-f_M(\varTheta)}\leq \|\theta-\varTheta\|\|f_M(\theta)-f_M(\varTheta)\|.
    \end{split}
\end{equation}}
\argument{\lref{eeq1};}{that for all $M\in \N\cap[\max\{64\fd^2\varrho^2(2\varrho+\varepsilon^{-1}(2\varrho+1)+\varepsilon^{-2}(8\varrho^3+5\varrho))^2(\varepsilon+\varrho)^4(\varepsilon \chi)^{-2},6\rho^{-1}\fc(\varepsilon+\varrho)^2(\varepsilon \chi)^{-1}\},\infty)$, $\theta\in \R^\fd$ with $\|\theta-\varTheta\|=6\fc(\varepsilon+\varrho)^2(\varepsilon \chi M)^{-1}$ it holds that
\begin{equation}\llabel{eeq2}
    \|f_M(\theta)-f_M(\varTheta)\|\geq \frac{\varepsilon \chi}{2(\varepsilon+\varrho)^2}\|\theta-\varTheta\|= 3\fc M^{-1}.
\end{equation}}
\argument{\lref{eeq2};the triangle inequality}{for all $M\in \N\cap[\max\{\scrM, 64\fd^2\varrho^2(2\varrho+\varepsilon^{-1}(2\varrho+1)+\varepsilon^{-2}(8\varrho^3+5\varrho))^2(\varepsilon+\varrho)^4(\varepsilon \chi)^{-2}\},\infty)$, $\theta\in \R^\fd$ with $\|\theta-\varTheta\|= 6\fc(\varepsilon+\varrho)^2(\varepsilon cM)^{-1}$ that
\begin{equation}\llabel{def: scrM}
   \|f_M(\theta)\|\geq 3\fc M^{-1}-\|f_M(\varTheta)\|.
\end{equation}}
\argument{\cref{item 2: lem negative derivative} in \cref{lem: negative derivative};\lref{def: batch size}}{that for all $M\in \N\cap[\max\{\bfM,\allowbreak 64\fd^2\varrho^2(2\varrho+\varepsilon^{-1}(2\varrho+1)+\varepsilon^{-2}(8\varrho^3+5\varrho))^2(\varepsilon+\varrho)^4(\varepsilon \chi)^{-2}\},\infty)$ there exists $\roott_M\in \R^\fd$ which satisfies
\begin{equation}\llabel{def: varthetaM}
    \|\roott_M-\varTheta\|\leq \frac { 6\fc(\varepsilon+\varrho)^2}{\varepsilon \chi M}=\rho \scrM M^{-1} \qandq \|f_M(\roott_M)\|^2=\inf_{\xi\in \R^\fd,\, \|\xi-\varTheta\|\leq \frac { 6\fc(\varepsilon+\varrho)^2}{\varepsilon\chi M}}\|f_M(\xi)\|^2.
\end{equation}}
\startnewargseq
\argument{\lref{eq5.5};\lref{def: scrM}}{that for all $M\in \N\cap[\max\{\bfM,\scrM,64\fd^2\varrho^2(2\varrho+\varepsilon^{-1}(2\varrho+1)+\varepsilon^{-2}(8\varrho^3+5\varrho))^2(\varepsilon+\varrho)^4(\varepsilon \chi)^{-2}\},\infty)$, $\xi\in \R^\fd$ with $\|\xi-\varTheta\|=6\fc(\varepsilon+\varrho)^2(\varepsilon \chi M)^{-1}$ it holds that
\begin{equation}\llabel{eqtg1}
    |f_M(\xi)|\geq 3\fc M^{-1}-\fc M^{-1}=2\fc M^{-1}>\fc M^{-1} \geq f_M(\varTheta).
\end{equation}}
\argument{\lref{eqtg1};\lref{def: varthetaM}}{that for all $M\in \N\cap[\max\{\bfM,\scrM,64\fd^2\varrho^2(2\varrho+\varepsilon^{-1}(2\varrho+1)+\varepsilon^{-2}(8\varrho^3+5\varrho))^2(\varepsilon+\varrho)^4(\varepsilon \chi)^{-2}\},\infty)$ it holds that
\begin{equation}\llabel{eeq5}
    \|\roott_M-\varTheta\|<\frac { 6\fc(\varepsilon+\varrho)^2}{\varepsilon \chi M}\leq \rho.
\end{equation}}
\argument{\lref{eeq5};\cref{def: globmin: existence of solution}}{for all $M\in \N\cap[\max\{\bfM,\scrM,64\fd^2\varrho^2(2\varrho+\varepsilon^{-1}(2\varrho+1)+\varepsilon^{-2}(8\varrho^3+5\varrho))^2(\varepsilon+\varrho)^4(\varepsilon \chi)^{-2}\},\infty)$ that
\begin{equation}\llabel{evd1}
    \roott_M\in \setTheta.
\end{equation}}
\argument{\cref{item 2: lem negative derivative} in \cref{lem: negative derivative};\lref{def: varthetaM};\lref{eeq5}}{that for all $M\in \N\cap[\max\{\bfM,\scrM,\allowbreak64\fd^2\varrho^2(2\varrho+\varepsilon^{-1}(2\varrho+1)+\varepsilon^{-2}(8\varrho^3+5\varrho))^2(\varepsilon+\varrho)^4(\varepsilon \chi)^{-2},\infty)$ it holds that
\begin{equation}\llabel{eeq6}
    \nabla_{\roott_M} \|f_M(\roott_M)\|^2=0.
\end{equation}}
\argument{\lref{eeq6};}{that for all $M\in \N\cap[\max\{\bfM,\scrM,64\fd^2\varrho^2(2\varrho+\varepsilon^{-1}(2\varrho+1)+\varepsilon^{-2}(8\varrho^3+5\varrho))^2(\varepsilon+\varrho)^4(\varepsilon \chi)^{-2},\infty)$ it holds that
\begin{equation}\llabel{eeq7}
     ((f_M)'(\roott_M))f_M(\roott_M)=0.
\end{equation}}
\argument{\lref{eeq7};\lref{def: fc};\lref{def: varthetaM}}{that for all $M\in \N\cap[\max\{\bfM,\scrM,64\fd^2\varrho^2(2\varrho+\varepsilon^{-1}(2\varrho+1)+\varepsilon^{-2}(8\varrho^3+5\varrho))^2(\varepsilon+\varrho)^4(\varepsilon \chi)^{-2}\},\allowbreak\infty)$ it holds that
\begin{equation}\llabel{evidence 1}
    f_M(\roott_M)=0.
\end{equation}}
\argument{\lref{evidence 1};\lref{def: batch size 2};\lref{def: varthetaM};\lref{evd1}}{\cref{item 1: existence of solution,item 2: existence of solution}\dott}
\end{aproof}
\renewcommand{\bfn}{N}
\subsection{Convergence rates for locally monotone stochastic optimization problems}

In \cref{prop: convergence general} below we employ the main result of our preliminary work \cite[Theorem 2.5]{DereichAdamconvergence2024} to establish convergence with convergence rates of the \Adam\ optimization process to a zero of the \Adam\ vector field under local monotonicity assumptions on the \Adam\ vector field (see \cref{item 4: convergence general} in \cref{prop: convergence general}). Our arguments in \cref{prop: convergence general} are inspired by the proof of Theorem 1.1 in our preliminary work \cite{DereichAdamconvergence2024} in which convergence rates for \Adam\ are established for the special class of the simple quadratic \SOPs\ in \cref{main theorem 2} above.
    \begin{athm}{prop}{prop: convergence general}
        Let $(\Omega,\cF,\P)$ be a probability space, let $\fd,\dimX\in \N$, let $\setTheta\subseteq \R^\fd$ be bounded and measurable, let $\setX\subseteq \R^\dimX$ be compact, let $X_{n,m}\colon \Omega\allowbreak\to \setX$, $(n,m) \in (\N_0)^2$, be \iid\ random variables, let $\smalll=(\smalll(\theta,x))_{(\theta,x)\in \R^\fd\times\R^{\dimX}}\allowbreak\in C^2(\R^\fd\times\R^{\dimX},\R)$, $\alpha\in (0,1)$, $\beta\in (\alpha^2,1)$, $\fc,\varepsilon,p\in (0,\infty)$,  let $(\gamma_n)_{n\in \N}\subseteq (0,\infty)$ be non-increasing, assume $\limsup_{n\to\infty}\gamma_n=0$ and $\limsup_{n\to\infty}((\gamma_n)^{-2}(\gamma_n-\gamma_{n+1}))< \fc$, and for every $M\in \N$ let $f_M=(f_M^1,\dots,f_M^\fd)\colon \R^\fd\to\R^\fd$ satisfy for all $i\in \{1,2,\dots,\fd\}$, $\theta=(\theta_1,\dots,\theta_\fd)\in \R^\fd$ that
     \begin{equation}\llabel{def: f}
         \textstyle f_M^i(\theta)=-(1-\alpha)\,\E\Biggl[\frac{\textstyle\sum_{n= 0}^\infty\alpha^n \frac 1M \textstyle\sum_{m=1}^M(\nabla_{\theta_i} \smalll)(\theta,X_{n,m})}{\textstyle\varepsilon+\bigl[(1-\beta)\sum_{n=0}^\infty\beta^n|\textstyle \frac 1M \sum_{m=1}^M(\nabla_{\theta_i} \smalll)(\theta,X_{n,m})|^2\bigr]^{1/2}}\Biggr].
     \end{equation}
       Then there exists $\fC\in \N$ such that for all $\vartheta\in  \{ \theta \in \setTheta \colon f_M( \theta ) = 0 \}$, $\bfn\in \N\cap[\fC,\infty)$, $M\in \N$ and all stochastic processes $\bfm=(\bfm_n)_{n\in \N_0}=((\bfm_{n}^{1},\dots,\bfm_{n}^{\fd}))_{n\in \N_0}\colon \N_0\allowbreak\times\Omega\to \R^\fd$, $\democrator=(\democrator_n)_{n\in \N_0}=((\democrator_{n}^{1},\dots,\democrator_{n}^{\fd}))_{n\in \N_0}\allowbreak\colon \N_0 \times\Omega\to\R^\fd$, and $\Theta=(\Theta_n)_{n\in \N_0}=((\Theta_{n}^{1},\dots,\Theta_{n}^{\fd}))_{n\in \N_0}\allowbreak\colon\allowbreak \N_0\times\Omega\to\R^\fd$ which satisfy that
     \begin{enumerate}[label=(\roman*)]
         \item it holds for all $n\in \N$, $i\in \{1,2,\dots,\fd\}$ that 
    \begin{equation}\llabel{setting: eq1: conjecture}
     \bfm_0=0,\qquad\bfm_n^i=\alpha \bfm_{n-1}^i+(1-\alpha)\bigl[\textstyle \frac 1M \sum_{m=1}^M(\nabla_{\theta_i} \smalll)(\Theta_{n-1},X_{n,m})\bigr],
    \end{equation}
    \begin{equation}   \llabel{setting: eq2: conjecture}
        \democrator_n=0,\qquad \democrator_{n}^{i}=\beta \democrator_{n-1}^{i}+(1-\beta)\bigl|\textstyle \frac 1M \sum_{m=1}^M(\nabla_{\theta_i} \smalll)(\Theta_{n-1},X_{n,m})\bigr|^2,
    \end{equation}
    \begin{equation}\llabel{setting: eq3: conjecture}
     \textstyle\text{and}\qquad \Theta_{n}^{i}=\Theta_{n-1}^{i} -\gamma_n\bigl[\frac{\bfm_{n}^{i}}{1-\alpha^n}\bigr] \Bigl[\varepsilon+\bigl[\frac{\democrator_{n}^{i}}{1-\beta^n}\bigr]^{\nicefrac{1}{2}}\Bigr]^{-1}, 
 \end{equation}
     \item it holds that $\Theta_0$ and $(X_{n,m})_{(n,m)\in \N^2}$ are independent,
     \item it holds that $\inf_{n\in \N}\P(\Theta_n\in \setTheta)=1$, and
     \item \label{item 4: convergence general} it holds for all $\theta\in \bar{V}$ that $\spro{f_M(\theta)-f_M(\globalmin),\theta-\globalmin}\leq -\fc\|\theta-\globalmin\|^2$
     \end{enumerate}
   we have that
    \begin{equation}\llabel{conclude}
        \bigl(\E\bigl[\|\Theta_\bfn-\globalmin \|^p\bigr]\bigr)^{\nicefrac{1}{p}}\leq \fC\sqrt{\gamma_{\bfn+1}}.
    \end{equation}
    \end{athm}
    \begin{aproof}
        Throughout this proof assume without loss of generality that $p>2$, for every $M\in \N$, $\theta\in \R^\fd$ let $\Psi_M^{\theta}\colon [0,\infty)\to\R^\fd$ satisfy for all $t\in [0,\infty)$ that
\begin{equation}\llabel{def: Psi}
    \Psi_M^{\theta}(t)=\theta,
    \end{equation}
    for every $k\in \N_0$ let $\varrho_k\in \R$ satisfy 
\begin{equation}\llabel{def:varrho}
\varrho_{k}=(1-\alpha) \varepsilon^{-1}\Bigl(\alpha^{k}+ 
\frac{1}{\sqrt{1-\alpha^{2}/\beta}}\beta^{k/2}\Bigr),
\end{equation}
 let $\ell_{\varrho}\subseteq (\R^\fd)^{\N_0}$ and $\|\cdot\|_{\ell_{\varrho}}\colon \ell_{\varrho}\to[0,\infty]$ satisfy for all $\mathbf x=(x_{k})_{k\in\N_{0}}\in \ell_{\varrho}$ that
\begin{equation}\llabel{def: norm1}
\|\mathbf x\|_{\ell_{\varrho}}=\sum_{k=0}^\infty \varrho_{k}\, \|x_{k}\|,
\end{equation}
 for every $m\in \R^\fd$, $v=(v_1,\dots,v_d)\in [0,\infty)^\fd$ let $(m,v)_{\ell_{\varrho}}\in \R\cup\{\infty\}$ satisfy 
\begin{multline}\llabel{def: norm2}
\textstyle(m,v)_{\ell_{\varrho}}=\inf\Bigl(\{\infty\}\\
\textstyle\cup\Bigl \{ \|\bfx\|_{\ell_{\varrho}}\colon \bfx=(x_k)_{k\in \N_0}=((x_{k}^{(i)})_{i\in\{1,2,\dots,\fd})_{k\in \N_0}\in \ell_{\varrho},\ m=(1-\alpha) \sum\limits_{k=0}^\infty \alpha^{k} x_{k}\\
\textstyle \text{ \ and } \qquad \Bigl(\forall\, i\in\{1,2,\dots,\fd\}\colon
v_i=(1-\beta) \sum\limits_{k=0}^\infty \beta^{k} |x_{k} ^{(i)}|^{2}\Bigr)\Bigr\}\Bigr),
\end{multline}
and for every $M\in \N$, $n\in \N_0$ let $\bbX_M^n=(\bbX_M^{n,1},\dots,\bbX_M^{n,\fd})\colon \allowbreak\R^\fd\times\Omega\to \R^\fd$ satisfy for all $\theta=(\theta_1,\dots,\theta_\fd)\in \R^\fd$ that
    \begin{equation}\llabel{def: bbX}
        \bbX_M^{n,i}(\theta)=-\frac 1M \biggl[\textstyle \sum\limits_{j=1}^M(\frac{\partial}{\partial\theta_i}\smalll)(\theta,X_{n,j})\biggr].
    \end{equation}
    \argument{\lref{def: Psi}}{that for all $\vartheta\in \setTheta$, $M\in \N$, $t\in [0,\infty)$ with $ f_M(\vartheta)=0$ it holds that
    \begin{equation}\llabel{evidence 11}
\Psi_M^\vartheta(t)=\Psi_M^\vartheta(0)+\int_{0}^tf_M(\Psi_M^\vartheta(s))\,\d s.
    \end{equation}}
    \argument{\lref{def: Psi};}{that for all $\vartheta\in \setTheta$, $M\in \N$, $x\in \R^\fd$, $t\in [0,\infty)$ it holds that
    \begin{equation}\llabel{evidence 1'}
        \spro{f_M(x)-f_M(\vartheta),x-\vartheta}=\spro{f_M(x)-f_M(\Psi_M^\vartheta(t)),x-\Psi_M^\vartheta(t)}.
    \end{equation}}
 \argument{\lref{def: bbX};\cref{item 1: lem negative derivative} in \cref{lem: negative derivative}; the fact that for all $n,m\in \N$ it holds that $\P(\|X_{n,m}\|\leq \bfx)$}{that there exists $\consc\in (0,\infty)$ which satisfies that for all $M\in \N$, $n\in \N_0$, $\theta_1,\theta_2\in \overline{\setTheta}$ it holds that
    \begin{equation}\llabel{def: C2}
       \|\bbX_{M}^{n}(\theta_1)\|\leq \consc \qqandqq \bigl(\E\bigl[\|\bbX_{M}^{n}(\theta_1)-\bbX_{M}^{n}(\theta_2)\|^p\bigr]\bigr)^{1/p}\leq \consc\|\theta_1-\theta_2\|.
    \end{equation}}
In the following let $\cK\in \R$ satisfy 
\begin{equation}\llabel{def: K}
    \cK\geq \textstyle2\bigl[\sup_{x\in \overline{V}}\|x\|\bigr]+\consc+\consc\textstyle \bigl[\sum_{k=0}^\infty\varrho_k\bigr]+1.
\end{equation} 
\startnewargseq
\argument{\lref{def: C2};\lref{def: K}}{that for all $M\in \N$, $n\in\N_0$, $\theta_1,\theta_2\in \overline{\setTheta}$ it holds that
\begin{equation}\llabel{evidence 1}
\bigl(\E[\|\bbX_M^n(\theta_1)\|^p]\bigr)^{1/p}\le\cK \qqandqq \bigl(\E[\|\bbX_M^n(\theta_1)-\bbX_M^n(\theta_2)\|^p]\bigr)^{1/p}\leq \cK \|\theta_1-\theta_2\|.
\end{equation}}
\argument{\lref{def: norm1};\lref{def: norm2};\lref{def: C2};the assumption that $(\gamma_n)_{n\in \N}$ is non-increasing;the fact that $\cK\geq \sqrt{\gamma_1}\consc$}{that for all $M\in \N$, $\bfn\in \N_0$, $(\theta_n)_{n\in \N_0}\subseteq \R^\fd$ and all stochastic processes $\bfm_n=(\bfm_n)_{n\in \N_0}=((\bfm_{n}^{1},\dots,\bfm_{n}^{\fd})\allowbreak)_{n\in \N_0}\allowbreak\colon \N_0\allowbreak\times\Omega\to \R^\fd$, $\democrator=(\democrator_n)_{n\in \N_0}=((\democrator_{n}^{1},\dots,\democrator_{n}^{\fd}))_{n\in \N_0}\colon \N_0 \times\Omega\to\R^\fd$ which satisfy for all $ n\in \N$, $i\in \{1,2,\dots,\fd\}$ that  $\|\theta_n\|\leq \fC$,
    \begin{equation}
     \bfm_0=0,\qquad\bfm_n=\alpha \bfm_{n-1}+(1-\alpha)\bbX_M^{n}(\theta_n),
    \end{equation}
    and
    \begin{equation} 
        \democrator_0=0,\qquad \democrator_{n}^{i}=\beta \democrator_{n-1}^{i}+(1-\beta)|\bbX_M^{n,i}(\theta_n)|^2
    \end{equation}
    we have that it holds $\P$-a.s.\ that
\begin{equation}\llabel{evidence 4}
\begin{split}
 \sqrt{\gamma_{n+1}} (\bfm_\bfn,\democrator_\bfn)_{\ell_{\varrho}} &\leq\sqrt{\gamma_1} \|(\bbX_M^{\bfn}(\theta_\bfn),\bbX_M^{\bfn-1}(\theta_{\bfn-1}),\dots,\bbX_M^{0}(\theta_0),0,0,\dots)\|_{\ell_\varrho}\\
   &\leq\sqrt{\gamma_1} \consc\bigl[\textstyle \textstyle\sum_{k=0}^\infty\varrho_k \bigr]\leq \cK.
    \end{split}
\end{equation}}
\argument{\lref{def: K}}{that for all $\theta,\vartheta\in \bar{\setTheta}$ it holds that
\begin{equation}\llabel{evidence 5}
    \|\theta-\vartheta\|\leq  \textstyle2\bigl[\sup_{x\in \bar{\setTheta}}\|x\|\bigr]\leq \cK.
\end{equation}}
    \argument{\lref{def: norm1};\lref{def: norm2};\lref{def: bbX};\lref{evidence 11};\lref{evidence 1'};\lref{evidence 1};\lref{evidence 4};\lref{evidence 5};the assumption that $\limsup_{n\to\infty}\allowbreak \gamma_n=0$;the assumption that $(\gamma_n)_{n\in \N}$ is non-increasing; the assumption $\limsup_{n\to\infty}\gamma_n=0$; the assumption that $\limsup_{n\to\infty}((\gamma_n)^{-2}(\gamma_n-\gamma_{n+1}))< 2\fc$; \cite[Theorem 2.5]{DereichAdamconvergence2024} (applied with $\alpha \curvearrowleft \alpha$, $\beta\curvearrowleft\beta$, $\varepsilon\curvearrowleft\varepsilon$, $(\gamma_n)_{n\in \N}\curvearrowleft (\gamma_n)_{n\in \N}$, $\cK\curvearrowleft \cK$, $c_1\curvearrowleft \fc$, $c_2\curvearrowleft \frac 12 \fc$, $c\curvearrowleft \frac12 \fc$, $p\curvearrowleft p$, $d\curvearrowleft \fd$ in the notation of \cite[Theorem 2.5]{DereichAdamconvergence2024})}{that there exist $\eta_1,\eta_2\in (0,\infty)$, $\fN\in \N$ which satisfy for all $\vartheta\in  \{ \theta \in \setTheta \colon f_M( \theta ) = 0 \}$, $\bfn\in \N\cap[\fN,\infty)$, $M\in \N$ and all stochastic processes $\bfm_n=(\bfm_n)_{n\in \N_0}=((\bfm_{n}^{1},\dots,\bfm_{n}^{\fd}))_{n\in \N_0}\colon \N_0\allowbreak\times\Omega\to \R^\fd$, $\democrator=(\democrator_n)_{n\in \N_0}=((\democrator_{n}^{1},\dots,\democrator_{n}^{\fd}))_{n\in \N_0}\colon \N_0 \times\Omega\to\R^\fd$, and $\Theta=(\Theta_n)_{n\in \N_0}=((\Theta_{n}^{1},\dots,\Theta_{n}^{\fd}))_{n\in \N_0}\colon \N_0\times\Omega\to\R^\fd$ which satisfy that
     \begin{enumerate}[label=(\Roman*)]
         \item it holds for all $n\in \N$, $i\in \{1,2,\dots,\fd\}$ that 
    \begin{equation}
     \bfm_0=0,\qquad\bfm_{n}^{i}=\alpha \bfm_{n-1}^{i}+(1-\alpha)\bigl[\textstyle \frac 1M \sum_{m=1}^M(\nabla_{\theta_i} \smalll)(\Theta_{n-1},X_{n,m})\bigr],
    \end{equation}
    \begin{equation}
        \democrator_0=0,\qquad \democrator_{n}^{i}=\beta \democrator_{n-1}^{i}+(1-\beta)\bigl|\textstyle \frac 1M \sum_{m=1}^M(\nabla_{\theta_i} \smalll)(\Theta_{n-1},X_{n,m})\bigr|^2,
    \end{equation}
    \begin{equation}
     \textstyle\text{and}\qquad \Theta_{n}^{i}=\Theta_{n-1}^{i} -\gamma_n\bigl[\frac{\bfm_{n}^{i}}{1-\alpha^n}\bigr] \Bigl[\varepsilon+\bigl[\frac{\democrator_{n}^{i}}{1-\beta^n}\bigr]^{\nicefrac{1}{2}}\Bigr]^{-1}, 
 \end{equation}
     \item it holds that $\Theta_0$ and $(X_{n,m})_{(n,m)\in \N^2}$ are independent,
     \item it holds that $\inf_{n\in \N}\P(\Theta_n\in \setTheta)$, and
     \item it holds for all $\theta\in \setTheta$ that $\spro{f_M(\theta)-f_M(\globalmin),\theta-\globalmin}\leq -\fc\|\theta-\globalmin\|^2$
     \end{enumerate}
    we have that
\begin{equation}\llabel{eqq1}
\begin{split}
&\E[|\Theta_\bfn-\globalmin|^p]^{1/p}= \E\bigl[\E\bigl[|\Theta_\bfn-\globalmin|^p|\sigma((\Theta_n)_{n\in \{1,2,\dots,\fN\}},(X_{n,m})_{(n,m)\in \{1,2,\dots,\fN\}^2})\bigr]\bigr]\\
&\le \E\Bigl[\Bigl\|\eta_1 +(1+\varepsilon)\frac{\|\Theta_{\fN}-\globalmin\|}{\sqrt{\gamma_{n+1}}}+\eta_1 \sqrt{(\bfm_{\fN},\democrator_{\fN})_{\ell_{\varrho}}}\Bigr\|^p
\Bigr]^{1/p} \sqrt{\gamma _{\bfn+1}}\\
&\leq  \E\Bigl[\Bigl\|\eta_1 +(1+\varepsilon)\frac{\cK}{\sqrt{\gamma _{n+1}}}+\eta_1 \sqrt{\consc} \Bigr\|^p\Bigr]^{1/p} \sqrt{\gamma _{\bfn+1}}\\
&\leq \eta_2  \sqrt{\gamma _{\bfn+1}}.
\end{split}
\end{equation}}
\argument{\lref{eqq1};}{\lref{conclude}\dott}
 \end{aproof}
 \renewcommand{\bfm}{\mathbf{m}}
\renewcommand{\democrator}{\mathbb{V}}
 \renewcommand{\bfn}{\mathbf{n}}
  \subsection{Convergence rates for strongly convex stochastic optimization problems}\label{subsec: convergence with sufficient batch}
\begin{athm}{prop}{theo: convergence of Adam}
    Assume \cref{setting: Adam vector field}, for every $M,n\in \N_0$ let  $\bfm_n^M=(\bfm^M_{n,1},\dots,\bfm^M_{n,\fd})\colon \Omega\allowbreak\to \R^\fd$, $\democrator_n^M=(\democrator^M_{n,1},\dots,\democrator^M_{n,\fd})\colon  \Omega\to\R^\fd$, and $\Theta_n^M=(\Theta^M_{n,1},\dots,\Theta^M_{n,\fd})\colon \Omega\to\R^\fd$ be random variables, assume for all $M,n\in \N$, $i\in \{1,2,\dots,\fd\}$ that
    \begin{equation}\label{setting: eq1: conjecture}
     \bfm_0^M=0,\qquad\bfm_n^M=\alpha \bfm_{n-1}^M+(1-\alpha)\bbX_{M}^n(\Theta_{n-1}^M),
    \end{equation}
    \begin{equation}   \label{setting: eq2: conjecture}
        \democrator_0^M=0,\qquad \democrator_{n,i}^M=\beta \democrator_{n-1,i}^M+(1-\beta)|\bbX_M^{n,i}(\Theta_{n-1}^M)|^2,
    \end{equation}
    \begin{equation}\label{setting: eq3: conjecture}
     \textstyle\text{and}\qquad \Theta_{n,i}^M=\Theta_{n-1,i}^M +\gamma_n\bigl[\frac{\bfm_{n,i}^M}{1-\alpha^n}\bigr] \Bigl[\varepsilon+\bigl[\frac{\democrator_{n,i}^M}{1-\beta^n}\bigr]^{\nicefrac{1}{2}}\Bigr]^{-1}, 
 \end{equation}
     assume for all $M\in \N$ that $\Theta_0^M$ and $(X_{n,m})_{(n,m)\in (\N_0)^2}$ are independent, let $\setTheta\subseteq \R^d$ be open and convex, assume for all $M,n\in \N$ that
     \begin{equation}\llabel{eq: boundTheta}
        \P(\Theta_n^M\in \setTheta)=1 ,
     \end{equation}
     let $\globalmin \in \setTheta$, $\rho\in (0,\infty)$ satisfy $\setTheta\subseteq \{x\in \R^\fd\colon \|x-\globalmin\|\leq \rho\} $ and $\E[\smalll(\globalmin, X_{1,1})]=\inf_{\theta\in \R^\fd}\allowbreak\E[\smalll(\theta, X_{1,1})]$, let $\varrho\in 
    \R$ satisfy for all $\theta=(\theta_1,\dots,\theta_\fd)\in \setTheta$, $i,j\in\{1,2,\dots,\fd\}$, $x\in \setX$ that
     \begin{equation}\llabel{def: varrho}
   \textstyle \varrho\geq \max\{|(\frac{\partial^2}{\partial\theta_i\partial\theta_j}\smalll)(\theta,x)|,|(\frac{\partial}{\partial\theta_i}\smalll)(\theta,x)|\},
     \end{equation}
      let $p\in (2,\infty)$ satisfy 
        $\limsup_{n\to\infty}(2(\gamma_n)^{-2}(\gamma_n-\gamma_{n+1}))<\varepsilon \chi(\varepsilon+\varrho)^{-2}$ and $\sum_{ n = 1 }^{ \infty } ( \gamma_n )^p < \infty$,
      let $\bfB\in \R$ be the constant in the Burkholder-Davis-Gundy inequality for the third moment, and let $\cM,\scrM\in \R$, $\delta\in (0,\infty)$, $\bfM\in \N$ satisfy $\cM=64\fd^2\varrho^2(2\varrho+\varepsilon^{-1}(2\varrho+1)+\varepsilon^{-2}(8\varrho^3+5\varrho))^2(\varepsilon+\varrho)^4(\varepsilon \chi)^{-2}$,
\begin{equation}\llabel{def: scrM}
   \textstyle \scrM= 18\rho^{-1}\beta^{-5/2}\allowbreak(1-\beta)\allowbreak\varepsilon^{-3}\allowbreak\varrho^3\allowbreak
     \delta^{-3/2}(\delta+6(\bfB+3)\varrho^2)(\varepsilon+\varrho)^2\chi^{-1},
\end{equation}
\begin{equation}\llabel{def: delta}
    \text{and}\qquad \textstyle\sup_{M\in \N\cap[\bfM,\infty)}\sup_{i\in \{1,2,\dots,\fd\}}\P\bigl(|\sqrt{M}\bbX^{1,i}_M(\globalmin)|\leq \sqrt{\delta}\bigr)\leq \frac {\beta^{3/2}}{2}.
\end{equation}
Then
     there exist $\fC\in \R$ and $(\roott_M)_{M\in \N}\subseteq \R^\fd$ which satisfy that
     \begin{enumerate}[label=(\roman*)]
         \item \label{item 1: theo: convergence of Adam} it holds for all $M\in \N\cap[\max\{\bfM,\cM,\scrM\},\infty) $ that $\P\bigl(\textstyle\limsup_{n\to\infty}\|\Theta^M_{n}-\roott_M\|=0\bigr)=1$,
         \item \label{item 2: theo: convergence of Adam} it holds for all $M\in \N\cap[\max\{\bfM,\cM,\scrM\},\infty) $, $n\in \N_0$ that $\bigl(\E[\|\Theta^M_{n}-\roott_M\|^p]\bigr)^{\nicefrac{1}{p}}\leq \fC\sqrt{\gamma_{n+1}}$, and
         \item \label{item 3: theo: convergence of Adam} it holds for all $M\in \N$ that $\|\roott_M-\globalmin\|\leq \rho\scrM M^{-1}$.
     \end{enumerate}
\end{athm}
\begin{aproof}
Throughout this proof let $\fM\in \R$ satisfy
     \begin{equation}\llabel{def: fM}
         \textstyle \fM=\max\{\bfM,\cM,\scrM \}.
     \end{equation}
\argument{\cref{lem: negative derivative}; \cref{lem: existence of solution};\lref{def: fM}}{that there exists $(\roott_M)_{M\in\N}\subseteq \setTheta$ which satisfies for all $N\in \N$, $M\in \N\cap[\fM,\allowbreak\infty)$, $\varTheta\in \setTheta$, $x\in \R^\fd$ that 
\begin{equation}\llabel{def: c12}
  \|\roott_N-\globalmin\|\leq \rho\scrM M^{-1},\qquad f_M(\roott_M)=0,\qqandqq \spro{x,( f_M)'(\varTheta)x}\leq -\frac{\varepsilon \chi}{2(\varepsilon+\varrho)^2} \|x\|^2.
\end{equation}}
\startnewargseq
\argument{\lref{def: c12}; the assumption that $\setTheta$ is convex;the fundamental theorem of calculus}{that for all $\varTheta\in\R^\fd$, $M\in \N\cap[\fM,\infty)$ with $\|\theta\|\leq \max\{\bfc,\|\vartheta\|+1\}$ it holds that
\begin{equation}\llabel{eq:934657}
\spro{f_M(\varTheta)-f_M(\roott_M),\varTheta-\roott_M}\leq -\frac{\varepsilon \chi}{2(\varepsilon+\varrho)^2}\|\varTheta-\roott_M\|^2.
\end{equation}}
\argument{\lref{eq: boundTheta};\lref{eq:934657}; the fact that $\limsup_{n\to\infty}\gamma_n=0$; the assumption that $\limsup_{n\to\infty}\allowbreak((\gamma_n)^{-2}(\gamma_n-\gamma_{n+1}))<\frac 12\varepsilon \chi(\varepsilon+\varrho)^{-2}$;\cref{prop: convergence general} (applied with $\Omega\curvearrowleft\Omega$, $\cF\curvearrowleft\cF$, $\P\curvearrowleft\P$, $\fd\curvearrowleft\fd$, $\dimX\curvearrowleft\dimX$, $(X_{n,m})_{(n,m)\in (\N_0)^2}\curvearrowleft (X_{n,m})_{(n,m)\in (\N_0)^2}$, $\smalll\curvearrowleft\smalll$, $\alpha\curvearrowleft\alpha$, $\beta\curvearrowleft\beta$, $\fC\curvearrowleft \max\{ \bfc,\|\globalmin\|+1\}$, $\fc\curvearrowleft \frac{\varepsilon \chi}{4(\varepsilon+\varrho)^2}$, $\varepsilon\curvearrowleft\varepsilon$, $(\gamma_n)_{n\in \N}\curvearrowleft (\gamma_n)_{n\in \N}$, $(f_M)_{M\in \N}\curvearrowleft(f_M)_{M\in \N}$, $(\globalmin_M)_{M\in \N}\curvearrowleft(\roott_M)_{M\in \N}$, $p\curvearrowleft 2p$ in the notation of \cref{prop: convergence general})}{that there exist $\eta\in (0,\infty)$ and $\fn\in \N$ which satisfy for all $M\in \N\cap[\fM,\infty)$, $n\in \N\cap[\fn,\infty)$ that
\begin{equation}\llabel{eqq1}
\begin{split}
 \bigl(\E[\|\Theta^M_{n}-\roott_M\|^{2p}]\bigr)^{\nicefrac{1}{2p}}\leq \eta  \sqrt{\gamma _{n+1}}.
\end{split}
\end{equation}}
\startnewargseq
\argument{the fact that $\setTheta$ is bounded}{that there exists $\bfc\in \R^\fd$ which satisfies
\begin{equation}\llabel{def: bfc}
    \setTheta\subseteq\{x\in \R^\fd\colon \|x\|\leq \bfc\}.
\end{equation}}
\argument{\cref{setting: eq1: conjecture};\cref{setting: eq2: conjecture};\cref{setting: eq3: conjecture}}{that for all $i\in \{1,2,\dots,\fd\}$ it holds that
\begin{equation}\llabel{eqp1'}
\begin{split}
 &\textstyle\Theta_{0,i}^M=\Theta_{1,i}^M -\gamma_1\bigl[\frac{\bfm_{1,i}^M}{1-\alpha}\bigr] \Bigl[\varepsilon+\bigl[\frac{\democrator_{1,i}^M}{1-\beta}\bigr]^{\nicefrac{1}{2}}\Bigr]^{-1}\\
&=\textstyle\Theta_{1,i}^M-\gamma_1\biggl[ \displaystyle\frac 1M \textstyle\sum\limits_{m=1}^M(\nabla_{\theta_i} \smalll)(\Theta^{M}_{0},X_{1,m})\biggr]\biggl[\varepsilon+\biggl[\displaystyle \frac 1M\textstyle \sum\limits_{m=1}^M(\nabla_{\theta_i} \smalll)(\Theta^{M}_{0},X_{1,m})\biggr]\biggr]^{-1}.
\end{split}
\end{equation}}
\argument{\lref{eq: boundTheta};\lref{eqp1'};\lref{def: bfc}}{that for all $M\in \N$, $i\in \{1,2,\dots,\fd\}$ it holds $\P$-a.s.\ that
\begin{equation}\llabel{eqp2'}
    |\Theta^{M}_{0,i}|\leq |\Theta^{M}_{1,i}|+\gamma_1\leq \bfc+\gamma_1.
\end{equation}}
\startnewargseq
\argument{\lref{eqp2'};}{that for all $M\in \N$ it holds that
\begin{equation}\llabel{eqp2''}
   \P\bigl( \|\Theta^{M}_{0}\|\leq \fd(\bfc+\gamma_1)\bigr)=1.
\end{equation}}
\argument{\lref{eqp2''};\lref{eq: boundTheta}}{for all $M\in \N$, $n\in \N_0$ that
\begin{equation}\llabel{eqp3'}
    \P\bigl(\|\Theta^{M}_{n}\|\leq \fd(\bfc+\gamma_1)\bigr)=1.
\end{equation}}
\argument{\lref{eqp3'};\lref{def: bfc};the assumption that for all $M\in \N$ it holds that $\roott_M\in \setTheta$; the assumption that $(\gamma_n)_{n\in\N}$ is non-increasing}{for all $M\in \N\cap[\fM,\infty)$, $n\in\N_0\cap[0,\fn]$ that
\begin{equation}\llabel{eqp4}
    \bigl(\E[\|\Theta^M_{n}-\roott_M\|^{2p}]\bigr)^{\nicefrac{1}{2p}}\leq d(\bfc+\gamma_1)+\bfc\leq [d(\bfc+\gamma_1)+\bfc](\gamma_{\fn+1})^{-1/2}\sqrt{\gamma_{n+1}}.
\end{equation}}
\argument{\lref{eqp4};\lref{eqq1}}{for all $M\in \N\cap[\fM,\infty)$, $n\in\N_0$ that
\begin{equation}\llabel{eqq1'}
     \bigl(\E[\|\Theta^M_{n}-\roott_M\|^{2p}]\bigr)^{\nicefrac{1}{2p}}\leq \bigl[\eta+ (d(\bfc+\gamma_1)+\bfc)^{-1/2}\bigr] \sqrt{\gamma _{n+1}}.
\end{equation}}
\argument{\lref{eqq1'};the assumption that $\sum_{n=1}^\infty(\gamma_n)^{p}  <\infty$}{that for all $M\in \N\cap[\fM,\infty)$ it holds that
\begin{equation}\llabel{eqq2}
\textstyle\E\biggl[\sum\limits_{n=0}^\infty\|\Theta_n^M-\roott_M\|^{2p}\biggr]\leq\bigl[\eta+ (d(\bfc+\gamma_1)+\bfc)^{-1/2}\bigr]^{2p} \biggl[\sum\limits_{n=0}^\infty(\gamma_{n+1})^p\biggr] <\infty.
\end{equation}}
\argument{\lref{eqq2};}{for all $M\in \N\cap[\fM,\infty)$ that it holds $\P$-a.s.\ that \llabel{arg1}
$\sum_{n=0}^\infty\|\Theta_n^M-\roott_M\|^{2p}<\infty$\dott}
\argument{\lref{arg1};}{that for all $M\in \N\cap[\fM,\infty)$ it holds that
\begin{equation}\llabel{eqq3}
    \limsup_{n\to\infty}\|\Theta^M_n-\roott_M\|=0.
\end{equation}}
\argument{\lref{def: c12};\lref{eqq3};\lref{eqq1'}}{items \ref{item 1: theo: convergence of Adam}, \ref{item 2: theo: convergence of Adam}, and \ref{item 3: theo: convergence of Adam}\dott}
\end{aproof}
\begin{athm}{theorem}{cor: convegence of Adam}
    Let $(\Omega,\cF,\P)$ be a probability space, let $\fd,\dimX\in \N$, let $\setX\subseteq \R^\dimX$ be compact, let $X_{n,m}\colon \Omega\allowbreak\to K$, $(n,m) \in \N^2$, be \iid\ random variables, let $\smalll=(\smalll(\theta,x))_{(\theta,x)\in \R^\fd\times\R^{\dimX}}\allowbreak\in C^2(\R^\fd\times\R^{\dimX},\R)$, $\varepsilon\in (0,\infty)$,
     let $(\gamma_n)_{n\in \N}\subseteq (0,\infty)$ be non-increasing, for every $M,n\in \N_0$, $\para\in \R^2$ let 
    $
\mandV{k}{M}{\para}{n}=(\mandVcom{k}{M}{\para}{n}{1},\dots,\allowbreak\mandVcom{k}{M}{\para}{n}{\fd})\colon\allowbreak \Omega\to\R^\fd$, $k\in\N$, and $\Adamprocess{M}{\para}{n}=(\Adamprocesscom{M}{\para}{n}{1},\dots,\Adamprocesscom{M}{\para}{n}{\fd})\colon\allowbreak \Omega\to\R^\fd$ be random variables, assume for all $k,M,n\in \N$, $i\in \{1,2,\dots,\fd\}$, $\para=(\para_1,\para_2)\in (0,1)^2$ that
    \begin{equation}\llabel{def: bbV}
    \mandV{k}{M}{\para}{0}=0,\qquad \mandVcom{k}{M}{\para}{n}{i}= \beta_k \mandVcom{k}{M}{\para}{n-1}{i}+(1-\beta_k)\bigl[\textstyle \frac 1M \sum_{m=1}^M(\nabla_{\theta_i} \smalll)(\Adamprocess{M}{\para}{n-1},X_{n,m})\bigr]^k,
        \end{equation}
        \begin{equation}\llabel{def: Theta}
     \textstyle\text{and}\qquad \Adamprocesscom{M}{\para}{n}{i}=\Adamprocesscom{M}{\para}{n-1}{i} -\gamma_n\Bigl[\frac{\mandVcom{1}{M}{\para}{n}{i}}{1-(\para_1)^n}\Bigr] \Bigl[\varepsilon+\Bigl[\frac{\mandVcom{2}{M}{\para}{n}{i}}{1-(\para_2)^n}\Bigr]^{\nicefrac{1}{2}}\Bigr]^{-1}, 
     \end{equation}
      let $\loss\colon \R^\fd\to\R$ satisfy for all $\theta\in \R^\fd$ that $\loss(\theta)=\E[\smalll(\theta,X_{1,1})]$,
     let $p,\chi\in (0,\infty)$, assume that $\cL$ is $\chi$-\stronglyconvex, let $\setTheta\subseteq\R^\fd$ be open, bounded,  and convex, let $\globalmin \in \setTheta$ satisfy $\loss(\globalmin)=\inf_{\theta\in \R^\fd}\loss(\theta)$,
    assume for all $\theta=(\theta_1,\dots,\theta_\fd)\in \setTheta$, $i,j\in\{1,2,\dots,\fd\}$, $x\in\setX$ that
     \begin{equation}\llabel{def: varrho}
     \begin{split}
   \textstyle 2\bigl[\varepsilon+ \max\{|(\frac{\partial^2}{\partial\theta_i\partial\theta_j}\smalll)(\theta,x)|,|(\frac{\partial}{\partial\theta_i}\smalll)(\theta,x)|\}\bigr]^2\bigl[\limsup_{n\to\infty}  ( \gamma_n )^{ - 2 } ( \gamma_n - \gamma_{ n + 1 } ) \bigr]\leq \varepsilon \chi,
   \end{split}
     \end{equation}
     assume 
       $\sum_{ n = 1 }^{ \infty } ( \gamma_n )^p < \infty$, and assume for all $M\in \N$, $\para\in \R^2$ that $\Adamprocess{M}{\para}{0}$ and $(X_{n,m})_{(n,m)\in \N^2}$ are independent \cfload. Then there exist $(\scrc_{\para})_{\para\in \R^2}\subseteq \N$ and $(\roott^M_{\para})_{(M,\para)\in \N\times\R^2}\subseteq \R^\fd$ such that for all $\para = ( \para_1, \para_2 ) \in (0,1) \times [ p^{ - 1 }, 1 )$ with $( \para_1 )^2 < \para_2$ and $\inf_{M,n\in \N}\P(\Adamprocess{M}{\para}{n}\in \setTheta)=1$ we have that
     \begin{enumerate}[label=(\roman*)]
        \item \label{item 1: cor} it holds that $\sum_{M=\scrc_0}^\infty\E\bigl[\limsup_{n\to\infty}\bigl(\|\Adamprocess{M}{\para}{n}-\roott^M_{\para}\|+|\loss(\Adamprocess{M}{\para}{n})-\loss(\roott^{M}_{\para})|\bigr)\bigr]=0$,
         \item \label{item 2: cor} it holds for all $M\in \N\cap[\scrc_0,\infty)$, $n\in \N_0$ that $\bigl(\E\bigl[\|\Adamprocess{M}{\para}{n}-\roott^M_{\para}\|^p+|\loss(\Adamprocess{M}{\para}{n})-\loss(\roott^M_{\para})|^p\bigr]\bigr)^{1/p}\leq \scrc_{\para}\sqrt{\gamma_{n+1}}$, 
         \item \label{item 3: cor} it holds for all $M\in \N$ that $\|\roott^M_{\para}-\globalmin\|\leq \scrc_0M^{-1}(1-\beta_2) $,
          \item \label{item 3.1: cor} it holds for all $M\in \N$ that $|\loss(\roott^M_{\para})-\loss(\globalmin)|\leq \scrc_0 M^{-2} (1-\beta_2)^2$,
         \item \label{item 4: cor} it holds for all $M \in \N$, $n \in \N_0$ that
$\bigl( \E\bigl[ \| \Adamprocess{M}{\para}{n} - \globalmin\|^p \bigr] \bigr)^{ 1/p } 
\leq \scrc_0\mathbbm 1_{[1,\scrc_0]}(M)+\scrc_0 M^{ - 1 }( 1 - \beta_2) +
\scrc_{\para}\sqrt{ \gamma_{n+1} } $, and
\item \label{item 4.1: cor} it holds for all $M \in \N$, $n \in \N_0$ that
$\bigl( \E\bigl[  | \loss( \Adamprocess{M}{\para}{n})- \loss( \globalmin ) |^p
\bigr] \bigr)^{ 1/p } 
\leq
\scrc_0\mathbbm 1_{[1,\scrc_0]}(M)+\scrc_0M^{ - 2 } (1-\beta_2)^2 $+$\scrc_{\para}\gamma_{n+1}$.
     \end{enumerate}
\end{athm}
\begin{aproof}
Throughout this proof assume without loss of generality that $p> 2$ and let $\varrho\in 
    \R$ satisfy for all $\theta=(\theta_1,\dots,\theta_\fd)\in \setTheta$, $i,j\in\{1,2,\dots,\fd\}$, $x\in \setX$ that
     \begin{equation}\llabel{def: varrho2}
   \textstyle \varrho\geq \max\{|(\frac{\partial^2}{\partial\theta_i\partial\theta_j}\smalll)(\theta,x)|,|(\frac{\partial}{\partial\theta_i}\smalll)(\theta,x)|\}.
     \end{equation}
     \argument{the assumption that $\setTheta$ is bounded}{that there exists $\fc\in \R$ which satisfies
    \begin{equation}\llabel{def: fc}
        \setTheta\subseteq\{x\in \R^d\colon \|x\|\leq \fc\}.
    \end{equation}}
    \startnewargseq
    \argument{\cref{lem: def deltaM};}{that there exists $\delta,\bfM\in (0,\infty)$ which satisfies \begin{equation}\llabel{def: delta}
    \textstyle\sup_{M\in \N\cap[\bfM,\infty)}\sup_{i\in \{1,2,\dots,\fd\}}\P\bigl(|\sqrt{M}\bbX^{1,i}_M(\globalmin)|\leq \delta\bigr)\leq \frac {p^{-3/2}}{2}.
\end{equation}
     }
\startnewargseq
     In the following let $\bfB\in \R$ be the constant in the Burkholder-Davis-Gundy inequality for the third moment and let $\cM,\scrM\in \R$ satisfy $\cM=64\fd^2\varrho^2(2\varrho+\varepsilon^{-1}(2\varrho+1)+\varepsilon^{-2}(8\varrho^3+5\varrho))^2(\varepsilon+\varrho)^4(\varepsilon \chi)^{-2}$ and \begin{equation}\llabel{def: scrM}
   \textstyle \scrM= 18\rho^{-1}p^{5/2}\allowbreak\allowbreak\varepsilon^{-3}\allowbreak\varrho^3\allowbreak
     \delta^{-3/2}(\delta+6(\bfB+3)\varrho^2)(\varepsilon+\varrho)^2\chi^{-1},
\end{equation} 
\argument{\lref{def: scrM};}{that for all $\beta_2\in [p^{-1},1)$ it holds that
\begin{equation}\llabel{evd1}
    \scrM\geq 18\rho^{-1}(\beta)^{-5/2}\allowbreak(1-\beta_2)^{1/2}\allowbreak\varepsilon^{-3}\allowbreak\varrho^3\allowbreak
   \delta^{-3/2}(3\delta+6(\bfB+3)\varrho^2)(\varepsilon+\varrho)^2\chi^{-1}.
\end{equation}}
\argument{\lref{evd2};}{for all  $\beta_2\in [p^{-1},1)$ that
\begin{equation}\llabel{evd2}
    \textstyle\sup_{M\in \N\cap[\bfM,\infty)}\sup_{i\in \{1,2,\dots,\fd\}}\P\bigl(|\sqrt{M}\bbX^{1,i}_M(\globalmin)|\leq \delta\bigr)\leq \frac {(\beta_2)^{3/2}}{2}.
\end{equation}}
\argument{the assumption that $\setTheta$ is open; the assumption that $\globalmin\in \setTheta$}{that there exists $\rho\in (0,1)$ such that
\begin{equation}\llabel{def: rho}
    \{x\in \R^\fd\colon \|x-\globalmin\|\leq \rho\}\subseteq \setTheta.
\end{equation}}
    \argument{\lref{def: rho};\lref{evd1};\lref{evd2};the assumption that $\loss$ is $\chi$-\stronglyconvex;the assumption that for all $n,m\in \N$ it holds that $X_{n,m}\in \setX$;\cref{theo: convergence of Adam}}{that there exist $\consc\in \R$, $(\scrc_{\para})_{\para\in \R^2}\subseteq\N$, and $(\roott^{M}_{\para})_{(M,\para)\in \N\times\R^2}\subseteq\setTheta$ which satisfy for all $\para = ( \para_1, \para_2 ) \in (0,1) \times [ p^{ - 1 }, 1 )$ with $( \para_1 )^2 < \para_2$ and $\inf_{M,n\in \N}\P(\Adamprocess{M}{\para}{n}\in \setTheta)=1$ that
    \begin{enumerate}[label=(\Roman*)]
        \item \llabel{item 1} it holds for all $M\in \N\cap[\consc,\infty)$ that $\P\bigl(\limsup_{n\to\infty}\|\Adamprocess{M}{\para}{n}-\roott^M_{\para}\|=1\bigr)=1$,
        \item \llabel{item 2} it holds for all $M\in \N\cap[\consc,\infty)$, $n\in \N_0$ that $\bigl(\E[\|\Adamprocess{M}{\para}{n}-\roott^M_{\para}\|^{2p}]\bigr)^{2/p}\leq \scrc_{\para}\sqrt{\gamma_{n+1}}$, and 
        \item \llabel{item 3} it holds for all $M\in \N$ that $\|\roott^M_{\para}-\globalmin\|\leq \consc (1-\beta_2)M^{-1}$.
    \end{enumerate}}  
    \startnewargseq
      \argument{\unskip, \eg, \cite[Lemma 5.8.6]{ArBePhi2024};the fact that $\loss\in C^1(\R^\fd,\R)$; the assumption that $\setTheta$ is bounded; the assumption that $\loss(\globalmin)=\inf_{\theta\in \R^\fd}\loss(\theta)$}{that there exists $\fC\in\R$ which satisfies for all $\theta,\varTheta\in \setTheta$ that
    \begin{equation}\llabel{def: fC}
        |\loss(\theta)-\loss(\varTheta)|\leq \fC\|\theta-\varTheta\|\qqandqq \loss(\theta)-\loss(\globalmin)\leq \fC\|\theta-\globalmin\|^2.
    \end{equation}
    }
    \argument{\lref{item 1};the fact that $\loss$ is continuous}{that for all $M\in \N\cap[\consc,\infty)$, $\para = ( \para_1, \para_2 ) \in (0,1) \times [ p^{ - 1 }, 1 )$ with $( \para_1 )^2 < \para_2$ and $\inf_{M,n\in \N}\P(\Adamprocess{M}{\para}{n}\in \setTheta)=1$ it holds that
    \begin{equation}\llabel{evidence 1}
\P\biggl(\limsup\limits_{n\to\infty}\bigl[\|\Adamprocess{M}{\para}{n}-\roott^M_\para\|+|\loss(\Adamprocess{M}{\para}{n})-\loss(\theta_M^\para)|\bigr]=0\biggr)=1.
    \end{equation}}
    \argument{\lref{item 2};\lref{item 3};\lref{def: fC}}{that for all $M\in \N\cap [\consc,\infty)$, $n\in \N_0$, $\para = ( \para_1, \para_2 ) \in (0,1) \times [ p^{ - 1 }, 1 )$ with $( \para_1 )^2 < \para_2$ and $\inf_{M,n\in \N}\P(\Adamprocess{M}{\para}{n}\in \setTheta)=1$ it holds that
    \begin{equation}\llabel{evidence 2}
    \begin{split}
     &\bigl(\E\bigl[\|\Adamprocess{M}{\para}{n}-\roott^M_{\para}\|^p+|\loss(\Adamprocess{M}{\para}{n}-\loss(\roott^M_{\para})|^p\bigr]\bigr)^{1/p}\\
     &\leq \bigl(\E\bigl[\|\Adamprocess{M}{\para}{n}-\roott^M_{\para}\|^p+\fC^p\|\Adamprocess{M}{\para}{n}-\roott^M_{\para}\|^p\bigr]\bigr)^{1/p}\\
     &\leq (1+\fC^p)^{1/p} \bigl(\E[\|\Adamprocess{M}{\para}{n}-\roott^M_{\para}\|^p]\bigr)^{1/p}\leq (1+\fC^p)^{1/p} \scrc_{\para_1}\sqrt{\gamma_{n+1}}.
        \end{split}
    \end{equation}}
    \argument{\lref{item 3};\lref{def: fC}}{that for all $M\in \N\cap[\consc,\infty)$, $\para = ( \para_1, \para_2 ) \in (0,1) \times [ p^{ - 1 }, 1 )$ with $( \para_1 )^2 < \para_2$ and $\inf_{M,n\in \N}\P(\Adamprocess{M}{\para}{n}\in \setTheta)=1$ it holds that
    \begin{equation}\llabel{evidence 3}
        |\loss(\roott^M_{\para})-\loss(\globalmin)|\leq \fC\|\roott^M_{\para}-\globalmin\|^2\leq \fC \consc^2(1-\para_2)^2M^{-2}.
    \end{equation}}
    \argument{\lref{item 2};\lref{item 3};\lref{def: fC}}{that for all $M\in \N\cap[\consc,\infty)$, $n\in \N_0$, $\para = ( \para_1, \para_2 ) \in (0,1) \times [ p^{ - 1 }, 1 )$ with $( \para_1 )^2 < \para_2$ and $\inf_{M,n\in \N}\P(\Adamprocess{M}{\para}{n}\in \setTheta)=1$ it holds that
    \begin{equation}\llabel{evidence 4'}
    \begin{split}
         \bigl(\E\bigl[\|\Adamprocess{M}{\para}{n}-\globalmin\|^{2p}\bigr)^{1/2p}
     &\leq  \bigl(\E[\|\Adamprocess{M}{\para}{n}-\roott^M_{\para}\|^{2p}]\bigr)^{1/2p}+  \bigl(\E[\|\roott^M_{\para}-\globalmin\|^{2p}]\bigr)^{1/2p}\\
     &\leq \scrc_{\para}\sqrt{\gamma_{n+1}}+ \consc(1-\para_2)M^{-1}.
     \end{split}
    \end{equation}}
    \argument{\lref{def: fc};the assumption that $\globalmin\in \setTheta$}{for all $M\in \N\cap[1,\consc]$, $n\in \N_0$, $\para = ( \para_1, \para_2 ) \in (0,1) \times [ p^{ - 1 }, 1 )$ with $( \para_1 )^2 < \para_2$ and $\inf_{M,n\in \N}\P(\Adamprocess{M}{\para}{n}\in \setTheta)=1$ that
    \begin{equation}\llabel{evidence 4''}
    \begin{split}
         \bigl(\E\bigl[\|\Adamprocess{M}{\para}{n}-\globalmin\|^{2p}\bigr)^{1/2p}
     &\leq  2\fc.
     \end{split}
    \end{equation}}
    \argument{\lref{evidence 4'};\lref{evidence 4''}}{for all $M\in \N$, $n\in \N_0$, $\para = ( \para_1, \para_2 ) \in (0,1) \times [ p^{ - 1 }, 1 )$ with $( \para_1 )^2 < \para_2$ and $\inf_{M,n\in \N}\P(\Adamprocess{M}{\para}{n}\in \setTheta)=1$ that
    \begin{equation}\llabel{evidence 4}
    \begin{split}
         \bigl(\E\bigl[\|\Adamprocess{M}{\para}{n}-\globalmin\|^{2p}\bigr)^{1/2p}
     \leq 2\fc\mathbbm 1_{[1,c]}(M)+\scrc_{\para}\sqrt{\gamma_{n+1}}+ \consc(1-\para_2)M^{-1}.
     \end{split}
    \end{equation}}
    \argument{\lref{def: fC};\lref{evidence 4}}{that for all $M\in \N$, $n\in \N_0$, $\para = ( \para_1, \para_2 ) \in (0,1) \times [ p^{ - 1 }, 1 )$ with $( \para_1 )^2 < \para_2$ and $\inf_{M,n\in \N}\P(\Adamprocess{M}{\para}{n}\in \setTheta)=1$ it holds that
    \begin{equation}\llabel{evidence 4.1}
    \begin{split}
        \bigl(\E\bigl[|\loss(\Adamprocess{M}{\para}{n}-\loss(\globalmin)|^p\bigr]\bigr)^{1/p}
     &\leq \fC \bigl(\E[\|\Adamprocess{M}{\para}{n}-\globalmin\|^{2p}]\bigr)^{1/p}\\
     &\leq \fC \bigl(2\fc\mathbbm 1_{[1,c]}(M)+\scrc_{\para}\sqrt{\gamma_{n+1}}+ \consc(1-\para_2)M^{-1}\bigr)^2\\
     &\leq 12\fC\fc\mathbbm 1_{[1,c]}(M)+3\fC(\scrc_\beta)\gamma_{n+1}+3\fC c^2(1-\beta_2)^2M^{-2}.
     \end{split}
    \end{equation}}
    \argument{\lref{item 3};\lref{evidence 1};\lref{evidence 2};\lref{evidence 3};\lref{evidence 4};\lref{evidence 4.1}; the fact that  $\limsup_{n\to\infty}\gamma_n=0$}{items \ref{item 1: cor}, \ref{item 2: cor}, \ref{item 3: cor}, \ref{item 3.1: cor}, \ref{item 4: cor}, and \ref{item 4.1: cor}\dott}
\end{aproof}
\begin{athm}{cor}{cor: convegence of Adam2}
    Let $(\Omega,\cF,\P)$ be a probability space, let $\fd,\dimX,p\in \N$, let $X_{n,m}\colon \Omega\allowbreak\to \R^{\dimX}$, $(n,m) \in \N^2$, be bounded \iid\ random variables, let $\smalll=(\smalll(\theta,x))_{(\theta,x)\in \R^\fd\times\R^{\dimX}}\allowbreak\in C^2(\R^\fd\times\R^{\dimX},\R)$, $\varepsilon\in (0,\infty)$, let $(\gamma_n)_{n\in \N}\subseteq (0,\infty)$ be non-increasing, for every $M,n\in \N_0$, $\para\in \R^2$ let 
    $
\mandV{k}{M}{\para}{n}=(\mandVcom{k}{M}{\para}{n}{1},\dots,\mandVcom{k}{M}{\para}{n}{\fd})\colon\allowbreak \Omega\to\R^\fd$, $k\in\N$, and $\Adamprocess{M}{\para}{n}=(\Adamprocesscom{M}{\para}{n}{1},\dots,\Adamprocesscom{M}{\para}{n}{\fd})\colon\allowbreak \Omega\to\R^\fd$ be random variables, assume for all $k,M,n\in \N$, $i\in \{1,2,\dots,\fd\}$, $\para=(\para_1,\para_2)\in (0,1)^2$ that
    \begin{equation}\llabel{def: bbV}
    \mandV{k}{M}{\para}{0}=0,\qquad \mandVcom{k}{M}{\para}{n}{i}= \beta_k \mandVcom{k}{M}{\para}{n-1}{i}+(1-\beta_k)\bigl[\textstyle \frac 1M \sum_{m=1}^M(\nabla_{\theta_i} \smalll)(\Adamprocess{M}{\para}{n-1},X_{n,m})\bigr]^k,
        \end{equation}
        \begin{equation}\llabel{def: Theta}
     \textstyle\text{and}\qquad \Adamprocesscom{M}{\para}{n}{i}=\Adamprocesscom{M}{\para}{n-1}{i} -\gamma_n\Bigl[\frac{\mandVcom{1}{M}{\para}{n}{i}}{1-(\para_1)^n}\Bigr] \Bigl[\varepsilon+\Bigl[\frac{\mandVcom{2}{M}{\para}{n}{i}}{1-(\para_2)^n}\Bigr]^{\nicefrac{1}{2}}\Bigr]^{-1}, 
     \end{equation}
     assume 
       $\limsup_{n\to\infty}( ( \gamma_n )^{ - 2 } ( \gamma_n - \gamma_{ n + 1 } )+\sum_{m=n}^\infty(\gamma_m)^p ) =0$, assume for all $M\in \N$, $\para\in \R^2$ that $\Adamprocess{M}{\para}{0}$ and $(X_{n,m})_{(n,m)\in \N^2}$ are independent, let $\loss\colon \R^\fd\to\R$ satisfy for all $\theta\in \R^\fd$ that $\loss(\theta)=\E[\smalll(\theta,X_{1,1})]$, assume that $\cL$ is \stronglyconvexx, and let $\globalmin \in \R^\fd$ satisfy $\loss(\globalmin)=\inf_{\theta\in \R^\fd}\loss(\theta)$ \cfload. Then there exist $(\scrc_{\para})_{\para\in \R^2}\subseteq \R$ such that for all $M \in \N$, $n \in \N_0$, $\para = ( \para_1, \para_2 ) \in (0,1) \times [ p^{ - 1 }, 1 )$ with $( \para_1 )^2 < \para_2$ and $\inf_{N,m\in \N}\P\bigl(\|\Adamprocess{N}{\para}{m}\|\leq p\bigr)=1$ it holds that
\begin{equation} \llabel{conclude}
   \bigl( \E\bigl[ \| \Adamprocess{M}{\para}{n} - \globalmin\|^p \bigr] \bigr)^{ 1/p } \leq  \scrc_0 M^{ - 1 } ( 1 - \beta_2 \mathbbm 1_{ [\scrc_0,\infty)} ( M ) ) + \scrc_{ \beta } ( \gamma_{ n + 1 } )^{ 1 / 2 }.
\end{equation}
\end{athm}
\begin{aproof}
Throughout this proof let $\fc\in \R$ satisfy
\begin{equation}\llabel{def: fc}
\P(X_{1,1}\leq [-\fc,\fc]^\dimX)=1,
\end{equation}
 for every $n,m\in\N$ let $\bfX_{n,m}\colon\Omega\to \R^\dimX$ satisfy
 \begin{equation}\llabel{def: bfX}
 \bfX_{n,m}=X_{n,m}\mathbbm 1_{[-\fc,\fc]^\dimX}(X_{n,m}), 
 \end{equation}
 and let  $
\mandVV{k}{M}{\para}{n}=(\mandVcomV{k}{M}{\para}{n}{1},\dots,\mandVcomV{k}{M}{\para}{n}{\fd})\colon\allowbreak \Omega\to\R^\fd$, $(M,n,k,\beta)\in (\N_0)^2\times\N\times\R^2$, and $\AdamprocessV{M}{\para}{n}=(\AdamprocesscomV{M}{\para}{n}{1},\dots,\AdamprocesscomV{M}{\para}{n}{\fd})\colon\allowbreak \Omega\to\R^\fd$, $(M,n,\beta)\in (\N_0)^2\times\R^2$, satisfy for all $k,M,n\in \N$, $i\in \{1,2,\dots,\fd\}$, $\para=(\para_1,\para_2)\in (0,1)^2$ that
    \begin{equation}\llabel{def: bbV2}
    \mandVV{k}{M}{\para}{0}=0,\qquad \mandVcomV{k}{M}{\para}{n}{i}= \beta_k \mandVcomV{k}{M}{\para}{n-1}{i}+(1-\beta_k)\bigl[\textstyle \frac 1M \sum_{m=1}^M(\nabla_{\theta_i} \smalll)(\AdamprocessV{M}{\para}{n-1},X_{n,m})\bigr]^k,
        \end{equation}
        \begin{equation}\llabel{def: Theta2}
     \AdamprocesscomV{M}{\para}{0}{i}=\Adamprocesscom{M}{\para}{0}{i},\qquad\textstyle\text{and}\qquad \AdamprocesscomV{M}{\para}{n}{i}=\AdamprocesscomV{M}{\para}{n-1}{i} -\gamma_n\Bigl[\frac{\mandVcomV{1}{M}{\para}{n}{i}}{1-(\para_1)^n}\Bigr] \Bigl[\varepsilon+\Bigl[\frac{\mandVcomV{2}{M}{\para}{n}{i}}{1-(\para_2)^n}\Bigr]^{\nicefrac{1}{2}}\Bigr]^{-1}, 
     \end{equation}
\argument{the assumption that $\cL$ is \stronglyconvexx}{that \llabel{argg1} there exists $\chi\in (0,\infty)$ such that $\loss$ is $\chi$-\stronglyconvex\dott}
    \argument{\lref{argg1};the fact that $\bfX_{n,m}\in [-\fc,\fc]^\dimX$;\lref{def: bbV2};\lref{def: Theta2};\cref{item 4: cor} in \cref{cor: convegence of Adam}}{that there exists $(\scrc_{\para})_{\para\in \R^2}\subseteq \R$ which satisfies for all $M \in \N$, $n \in \N_0$, $\para = ( \para_1, \para_2 ) \in (0,1) \times [ p^{ - 1 }, 1 )$ with $( \para_1 )^2 < \para_2$ and $\inf_{N,m\in \N}\P\bigl(\|\AdamprocessV{N}{\para}{m}\|\leq p\bigr)=1$ that
\begin{equation} \llabel{eq1}
   \bigl( \E\bigl[ \| \AdamprocessV{M}{\para}{n} - \globalmin\|^p \bigr] \bigr)^{ 1/p } \leq  \scrc_0 M^{ - 1 } ( 1 - \beta_2 \mathbbm 1_{ [\scrc_0,\infty)} ( M ) ) + \scrc_{ \beta } ( \gamma_{ n + 1 } )^{ 1 / 2 }.
\end{equation}
}
\argument{ the assumption that $X_{n,m}$, $(n,m)\in \N^2$, are \iid;\lref{def: fc};\lref{def: bfX};}{that \llabel{eqq1.1} $(X_{n,m})_{(n,m)\in \N^2}$ and $(\bfX_{n,m})_{(n,m)\in (\N_0)^2}$ are identically distributed\dott}
\argument{\lref{eqq1.1};\lref{def: bbV};\lref{def: Theta};\lref{def: bbV2};\lref{def: Theta2}}{that for all $M\in \N$, $n\in \N_0$, $\beta\in (0,1)$ it holds that \llabel{eqq1.2} $\Adamprocess{M}{\para}{n}$  and $\AdamprocessV{M}{\para}{n}$  are  \llabel{arg2} identiacally distributed\dott}
\argument{\lref{arg2};\lref{eq1}}{\lref{conclude}\dott}
\end{aproof}
\begin{athm}{cor}{cor: convegence of Adam2 intro}
     Let $(\Omega,\cF,\P)$ be a probability space, let $\fd,\dimX,p\in \N$, let $X_{n,m}\colon \Omega\allowbreak\to \R^{\dimX}$, $(n,m) \in \N^2$, be bounded \iid\ random variables, let $\smalll=(\smalll(\theta,x))_{(\theta,x)\in \R^\fd\times\R^{\dimX}}\allowbreak\in C^2(\R^\fd\times\R^{\dimX},\R)$, $\para_1\in (0,1)$, $\para_2\in ((\para_1)^2,1)$, $\varepsilon\in (0,\infty)$, let $(\gamma_n)_{n\in \N}\subseteq (0,\infty)$ be non-increasing, for every $M,n\in \N_0$ let  $
\democrator^{k,M}_{n}=(\democrator^{k,M}_{n,1},\dots,\democrator^{k,M}_{n,\fd})\colon\allowbreak \Omega\to\R^\fd$, $k\in\{1,2\}$, and $\Theta_n^M=(\Theta^M_{n,1},\dots,\Theta^M_{n,\fd})\colon\allowbreak \Omega\to\R^\fd$ be random variables, assume for all $M,n\in \N$, $i\in \{1,2,\dots,\fd\}$, $k\in \{1,2\}$ that
    \begin{equation}\label{eq1: main theorem'}
    \democrator_0^{k,M}=0,\qquad \democrator_{n,i}^{k,M}= \beta_k \democrator_{n-1,i}^{k,M}+(1-\beta_k)\bigl[\textstyle \frac 1M \sum_{m=1}^M(\nabla_{\theta_i} \smalll)(\Theta_{n-1}^M,X_{n,m})\bigr]^k,
        \end{equation}
        \begin{equation}\label{eq2: main theorem'}
   \Theta_{n,i}^M=\Theta_{n-1,i}^M -\gamma_n\Bigl[\frac{\democrator_{n,i}^{1,M}}{1-(\para_1)^n}\Bigr] \Bigl[\varepsilon+\Bigl[\frac{\democrator_{n,i}^{2,M}}{1-(\para_2)^n}\Bigr]^{\nicefrac{1}{2}}\Bigr]^{-1},
     \end{equation}
     and $\P\bigl(\sup_{M,n\in \N}\|\Theta_n^M\|\leq p\bigr)=1$,
     assume 
        $\limsup_{n\to\infty}(  ( \gamma_n )^{ - 2 } ( \gamma_n - \gamma_{ n + 1 } ) +\sum_{m=n}^\infty(\gamma_m)^p) =0$, assume that $(\Theta_0^M)_{M\in \N}$ and $(X_{n,m})_{(n,m)\in \N^2}$ are independent, let $\loss\colon \R^\fd\to\R$ satisfy for all $\theta\in \R^\fd$ that $\loss(\theta)=\E[\smalll(\theta,X_{1,1})]$, assume that $\cL$ is \stronglyconvexx, and let $\globalmin \in \R^\fd$ satisfy $\loss(\globalmin)=\inf_{\theta\in \R^\fd}\loss(\theta)$ \cfload. Then there exists $\scrc\in \R$ such that for all $M,n \in \N$ it holds that
\begin{equation}\label{conclude: main theorem'} 
\bigl( \E\bigl[ \| \Theta^M_n - \globalmin\|^p \bigr] \bigr)^{ \nicefrac{1}{p} } 
\leq
\scrc( M^{ - 1 }+\sqrt{\gamma_{n} }   ).
     \end{equation}
\end{athm}
\begin{aproof}
    \argument{the assumption that $\P\bigl(\sup_{M,n\in \N}\|\Theta_n^M\|\leq p\bigr)=1$;the fact that for all $n\in \N$ it holds that $\gamma_{n+1}\leq \gamma_n$;\cref{cor: convegence of Adam2}}{that there exists $\scrc\in \R$ such that for all $M,n\in \N$ it holds that
    \begin{equation}
\begin{split}
    \bigl( \E\bigl[ \| \Adamprocess{M}{\para}{n} - \globalmin\|^p \bigr] \bigr)^{ 1/p } 
    &\leq  \scrc M^{ - 1 } ( 1 - \beta_2 \mathbbm 1_{ [\scrc_0,\infty)} ( M ) ) + \scrc ( \gamma_{ n + 1 } )^{ 1 / 2 }\\
    &\leq \scrc M^{ - 1 } +\scrc ( \gamma_{ n  } )^{ 1 / 2 }.
\end{split}
    \end{equation}}
\end{aproof}

In \cref{lem: unsummable} and \cref{lem: verify 3} below we recall two few elementary facts regarding the assumptions on the learning rates that we employ in \cref{cor: convegence of Adam} and \cref{cor: convegence of Adam2} above. Only for completeness we include here detailed proofs for \cref{lem: unsummable} and \cref{lem: verify 3}.
\begin{athm}{lemma}{lem: unsummable}
    Let $(\gamma_n)_{n\in \N}\subseteq (0,\infty)$ be non-increasing and assume $\limsup_{n\to\infty}( ( \gamma_n )^{ - 2 } ( \gamma_n - \gamma_{ n + 1 } ) ) <\infty$. Then 
   \begin{equation}\llabel{item 2}
       \sum_{n=1}^\infty \gamma_n=\infty.
   \end{equation} 
\end{athm}
\begin{aproof}
% \argument{the assumption that $\limsup_{n\to\infty}( \gamma_n + ( \gamma_n )^{ - 2 } ( \gamma_n - \gamma_{ n + 1 } ) ) =0$}{that
% \begin{equation}\llabel{eq1}
%     \limsup_{n\to\infty} \gamma_n =0.
% \end{equation}} 
% \argument{\lref{eq1};the assumption that for all $n\in \N$ it holds that $\gamma_n>0$}{that 
% \begin{equation}\llabel{eq2}
%     \sup_{n\in \N} \gamma_n<\infty.
% \end{equation}}
\argument{the assumption that $\limsup_{n\to\infty}( ( \gamma_n )^{ - 2 } ( \gamma_n - \gamma_{ n + 1 } ) ) <\infty $;}{that there exists $c\in (0,\infty)$ which satisfies
\begin{equation}\llabel{def: c}
  \textstyle  c>\sup_{n\in \N}(( \gamma_n )^{ - 2 } ( \gamma_n - \gamma_{ n + 1 } )).
\end{equation}}
\startnewargseq
% \argument{\lref{def: c};}{that for all $k \in \N$ it holds that
% \begin{equation}\llabel{eq1}
% \begin{split}
% \textstyle \gamma_k = \sum\limits_{ n = k }^{ \infty } ( \gamma_n - \gamma_{ n + 1 } ) \leq \sum\limits_{ n = k }^{ \infty } \bigl[ c ( \gamma_n )^2 \bigr] 
% \leq \sum\limits_{ n = k }^{ \infty } \bigl[ c \gamma_k \gamma_n \bigr] 
% = c \gamma_k \Bigl[ \sum\limits_{ n = k }^{ \infty } \gamma_n \Bigr] .
% \end{split}
% \end{equation}}
 We prove \lref{item 2} by contradiction. In the following we thus assume that 
\begin{equation}\llabel{assume}
\textstyle \sum_{n=1}^\infty \gamma_n<\infty.
\end{equation}
\startnewargseq
\argument{\lref{assume};}{that \llabel{arg1} $\limsup_{n\to\infty}\gamma_n=0$\dott}
\argument{\lref{arg1};\lref{assume}}{that for all $k \in \N$ it holds that
\begin{equation}\llabel{eq1}
\textstyle \gamma_k = \sum\limits_{ n = k }^{ \infty } ( \gamma_n - \gamma_{ n + 1 } ) \leq \sum\limits_{ n = k }^{ \infty } \bigl[ c ( \gamma_n )^2 \bigr] 
\leq \sum\limits_{ n = k }^{ \infty } \bigl[ c \gamma_k \gamma_n \bigr] 
= c \gamma_k \biggl[ \sum\limits_{ n = k }^{ \infty } \gamma_n \biggr] .
\end{equation}}
\argument{\lref{eq1};}{that \llabel{arg2} $1 \leq c ( \lim_{ k \to \infty } [ \sum_{ n = k }^{ \infty } \gamma_n ] )$\dott} 
\argument{\lref{arg2};}{that \llabel{arg3} $0 < \lim_{ k \to \infty } [ \sum_{ n = k }^{ \infty } \gamma_n ]$\dott}
\argument{\lref{arg3};}{that
\begin{equation}\llabel{arg4}
    \textstyle \sum_{n=1}^\infty\gamma_n=\infty.
\end{equation}}
% \argument{\lref{assume};}{that there exists $k\in \N$ which satisfies
% \begin{equation}\llabel{def: k}
%     \sum_{n=k}^\infty \gamma_n < c^{-1}.
% \end{equation}}
% \startnewargseq
% \argument{\lref{def: c};the assumption that for all $n\in \N$ it holds that $\gamma_n>0$; the fact that $\limsup_{n\to\infty}\gamma_n=0$}{that
% \begin{equation}\llabel{eq4}
%     \sum_{n=k}^\infty (\gamma_n)^2\geq \sum_{n=k}^\infty c^{-1}(\gamma_n-\gamma_{n+1})= c^{-1}\gamma_k.
% \end{equation}}
% \argument{\lref{assume};\lref{eq4}; the fact that $(\gamma_n)_{n\in \N}$ is non-increasing}{that
% \begin{equation}
%     c^{-1}=(\gamma_{k})^{-1}c^{-1}\gamma_k\leq (\gamma_k)^{-1}\textstyle\biggl[\sum\limits_{n=k}^\infty (\gamma_n)^2\biggr]=\sum\limits_{n=k}^\infty (\gamma_n)^2(\gamma_k)^{-1}\leq \sum\limits_{n=k}^\infty\gamma_n<c^{-1}.
% \end{equation}}
This contradiction proves \lref{item 2}.
\end{aproof}
\begin{athm}{lemma}{lem: verify 3}
    Let $r\in (0,\infty)$ and for every $n\in \N$ let $\gamma_n\in \R$ satisfy $\gamma_n=n^{-r}$. Then the following two statements are equivalent:
    \begin{enumerate}[label=(\roman*)]
        \item \label{item 1: verify 3} It holds that $r< 1$.
        \item \label{item 2: verify 3} It holds that $\limsup_{n\to\infty}( ( \gamma_n )^{ - 2 } ( \gamma_n - \gamma_{ n + 1 } ) ) =0$.
    \end{enumerate}
\end{athm}
\begin{aproof}
In the following we prove that (\ref{item 2: verify 3}$\rightarrow$\ref{item 1: verify 3}). In the following we thus assume that
\begin{equation}\llabel{assume 0}
   \textstyle \limsup_{n\to\infty}( ( \gamma_n )^{ - 2 } ( \gamma_n - \gamma_{ n + 1 } ) ) =0.
\end{equation}
\argument{\lref{assume 0};the fact that $(\gamma_n)_{n\in \N}$ is non-increasing;\cref{lem: unsummable}}{that
\begin{equation}\llabel{eqt2}
   \textstyle \sum_{n=1}^\infty\gamma_n=\infty.
\end{equation}}
\argument{\lref{eqt2};the fact that for all $s\in (1,\infty)$ it holds that $\sum_{n=1}^\infty n^{-s}<\infty$}{that
\begin{equation}\llabel{eqt3}
    r\leq 1.
\end{equation}}
\argument{\lref{eqt3};the fact that 
  $\textstyle  \limsup_{n\to\infty}(n^{-2}(n^{-1}-(n+1)^{-1})=1>0$;}
{that \llabel{argt1} $r<1$\dott}
\argument{\lref{argt1};}{(\ref{item 2: verify 3}$\rightarrow$\ref{item 1: verify 3})\dott}
\startnewargseq
In the following we prove that (\ref{item 1: verify 3}$\rightarrow$\ref{item 2: verify 3}). In the following we thus assume that
\begin{equation}\llabel{assume}
    r<1.
\end{equation}
% \argument{\lref{assume}; the assumption that for all $n\in \N$ it holds that $\gamma_n=n^{-r}$}{that
% \begin{equation}\llabel{eq1}
%     \textstyle\limsup_{n\to\infty} \gamma_n=0.
% \end{equation}}
\argument{\lref{assume};the assumption that for all $n\in \N$ it holds that $\gamma_n=n^{-r}$; the mean value theorem}{that
\begin{equation}\llabel{eq3}
\begin{split}
   &\limsup_{n\to\infty}(( \gamma_n )^{ - 2 } ( \gamma_n - \gamma_{ n + 1 } ))=\limsup_{n\to\infty}( n^{2r}[n^{-r}-(n+1)^{-r}])\\
   &\leq \limsup_{n\to\infty}( n^{2r}(-r)n^{-r-1}(n-(n+1)))=\limsup_{n\to\infty}( n^{2r}rn^{-r-1})\\
   &=\limsup_{n\to\infty}( n^{r-1}r)=0.
   \end{split}
\end{equation}}
\argument{\lref{eq3};}{(\ref{item 1: verify 3}$\rightarrow$\ref{item 2: verify 3})\dott}
\end{aproof}

\section{Non-convergence of Adam for non-symmetrically distributed data}\label{sec: converge to wrong point}
\newcommand{\Advf}[2]{Z_{#1,#2}}
\newcommand{\Advfcom}[3]{Z_{#1,#2}^{#3}}
In this section we reveal in \cref{lem: non-convergence of Adam2} below in the situation of a class of simple quadratic \SOPs\ that for \emph{every} $\beta_1 \in (0,1)$, $\beta_2 \in ( (\beta_1)^2, 1)$, $\varepsilon \in (0,\infty)$ we have that the \Adam\ optimizer with parameters $\beta_1, \beta_2$, and $\varepsilon$ \emph{fails} to converge in probability to the solution of the considered \SOP\ (to the unique global minimum point of the considered \SOP) if the data in the considered \SOP\ is not symmetrically distributed. Our proof of \cref{lem: non-convergence of Adam2} employs an application of Corollary 1.10 in our preliminary work \cite{DereichJentzenKassing2025}.

 For convenience of the reader a slightly generalized variant of Corollary 1.10 in  \cite{DereichJentzenKassing2025} is formulated here in \cref{theo: convergence point of Adam copy 2} below. \cref{theo: convergence point of Adam copy 2} demonstrates in a very general setting that if the \Adam\ optimizer converges in probability to an arbitrary random variable $\vartheta \colon \Omega \to \R^d$, then this random variable must $\P$-almost surely be a zero of the \Adam\ vector field. The \Adam\ vector field has been introduced in \cite[Definition 2.4]{DereichAdamconvergence2024} and is formulated here as the function $f\colon \R^d \to \R^d$ in \cref{item 2: theo: convergence point of Adam copy 3 pre} in \cref{theo: convergence point of Adam copy 3 pre} below.
\subsection{Adam limit theorem -- limit points of the Adam optimizer}
\begin{athm}{prop}{theo: convergence point of Adam copy 3 pre}[\Adam\ vector field formulation]
    Let $( \Omega, \cF, \P )$ be a probability space, 
let $ \fd, \mathscr{d} \in \N $, 
let $X_{ n, m } \colon \Omega\to \R^{ \mathscr{d} }$, $ (n,m) \in (\N_0)^2 $, be \iid\ random variables, assume for all $p\in\N$ that $\E[\|X_{1,1}\|^p ] <\infty$, let $\smalll \in C^1( \R^\fd \times \R^{ \scrd }, \R )$ have at most polynomially growing derivatives, let $M \in \N$, $\beta_1,\beta_2 \in (0,1)$, $\varepsilon \in (0,\infty)$,
and for every $\theta=(\theta_1,\dots,\theta_d)\in \R^d$, $N\in \N$ let $\Advf{\theta}{N} =(\Advfcom{\theta}{N}{1},\dots,\Advfcom{\theta}{N}{d})\colon \Omega \to \R^\fd$ satisfy for all $i\in \{1,2,\dots,\fd\}$ that 
    \begin{equation}\llabel{def: Adam vector field}
         \Advfcom{\theta}{N}{i} =\frac{(1-\beta_1)\textstyle\sum_{n= 0}^N(\beta_1)^n M^{-1} \textstyle\sum_{m=1}^M(\nabla_{\theta_i} \smalll)(\theta,X_{n,m})}{\textstyle\varepsilon+\bigl[(1-\para_2)\sum_{n=0}^N(\beta_2)^n\bigl|\textstyle M^{-1} \sum_{m=1}^M(\nabla_{\theta_i} \smalll)(\theta,X_{n,m})\bigr|^2\bigr]^{1/2}}.
    \end{equation}
    Then 
\begin{enumerate}[label=(\roman*)]
 \item \llabel{item 1} it holds for all $\theta\in \R^d$, $p\in (0,\infty)$ that $\E\bigl[ \sup_{ N \in \N } \| \Advf{\theta}{N}\|^p \bigr] < \infty$,
    \item \label{item 2: theo: convergence point of Adam copy 3 pre} there exists a unique $f \colon \R^d \to \R^d$ which satisfies for all $\theta \in \R^d$ that
$\limsup_{ N \to \infty } \| f( \theta ) - \E[ \Advf{\theta}{N} ] \| = 0$,
\item \llabel{item 3} it holds that $f \in C( \R^d, \R^d )$, 
\item \llabel{item 4} it holds that $\P\bigl(\limsup_{N\to\infty}\sup_{n,m\in \N\cap[N,\infty)}\|\Advf{\theta}{n}-\Advf{\theta}{m}\|=0\bigr)=1$,
\item \llabel{item 5} it holds for every $\theta\in\R^d$, $p\in (0,\infty)$ and every random variable $\cZ\colon \Omega\to\R^d$ with $\limsup_{N\to\infty} \E\bigl[\min\{1,\allowbreak\|\cZ-\Advf{\theta}{N}\|\}\bigr]=0$ that $\E\bigl[ \| \cZ \|^p+\|\cZ\| \bigr] < \infty$ and $\f(\theta)=\E[\cZ]$, and
\item \llabel{item 6} there exists a continuous random field $\cZ\colon \R^d\times \Omega \to \R^d$ which satisfies for all $\theta \in \R^d$, $p\in (0,\infty)$ that
$\E\bigl[ \| \cZ( \theta ) \|^p \bigr] < \infty$, $\limsup_{ N \to \infty } \E\bigl[ \| \cZ(\theta)- \Advf{\theta}{N} \|^p \bigr] = 0$, and $\P( \limsup_{ N \to \infty } \| \cZ(\theta) \allowbreak- \Advf{\theta}{N} \| = 0 ) = 1$.
\end{enumerate}
\end{athm}
\begin{aproof}
\argument{the assumption that $\smalll$ have at most polynomially growing derivaties}{that there exist $c\in\R$, $\fp\in \N$ which satisfy for all $\theta\in \R^d$, $x\in \R^\dimX$, $i\in \{1,2,\dots,d\}$ that
\begin{equation}\llabel{eq1}
   | (\nabla_{\theta_i} \smalll)(\theta,x)|\leq c+c\|\theta\|^\fp+c\|x\|^\fp.
\end{equation}}
\startnewargseq 
\argument{\lref{def: Adam vector field};}{that for all $\theta\in \R^d$, $N\in \N$, $i\in \{1,2,\dots,d\}$ it holds that
\begin{equation}\llabel{eq2}
\begin{split}
 |\Advfcom{\theta}{N}{i}|&\leq \varepsilon^{-1}\bigl|(1-\beta_1)\textstyle\sum_{n= 0}^N(\beta_1)^n M^{-1} \textstyle\sum_{m=1}^M(\nabla_{\theta_i} \smalll)(\theta,X_{n,m})\bigr|\\
 &\textstyle\leq \varepsilon^{-1}(1-\beta_1)\sum_{n= 0}^N(\beta_1)^n \bigl|M^{-1} \textstyle\sum_{m=1}^M(\nabla_{\theta_i} \smalll)(\theta,X_{n,m})\bigr|\\
 &\textstyle\leq \varepsilon^{-1}(1-\beta_1)\sum_{n= 0}^\infty(\beta_1)^n\bigl| M^{-1} \textstyle\sum_{m=1}^M (\nabla_{\theta_i} \smalll)(\theta,X_{n,m})\bigr|.
    \end{split}
\end{equation}}
\argument{\lref{eq2};}{that for all $\theta\in \R^d$, $i\in \{1,2,\dots,d\}$ it holds that
\begin{equation}\llabel{eq2'}
\begin{split}
\textstyle\sup_{N\in \N} |\Advfcom{\theta}{N}{i}|
 \textstyle\leq \varepsilon^{-1}(1-\beta_1)\sum_{n= 0}^\infty(\beta_1)^n\bigl| M^{-1} \textstyle\sum_{m=1}^M (\nabla_{\theta_i} \smalll)(\theta,X_{n,m})\bigr|.
    \end{split}
\end{equation}}
\argument{\lref{eq1};the fact that for all $n,p\in \N$, $x_1,\dots ,x_n\in \R$ it holds that $\bigl|\sum_{i=1}^nx_m\bigr|^p\leq n^{p-1}\sum_{i=1}^n |x_i|^p$}{that for all $\alpha\in (0,1)$, $\theta\in \R^d$, $i\in \{1,2,\dots,d\}$, $n\in \N_0$, $p\in\N$ it holds that
\begin{equation}\llabel{eq2.5}
\begin{split}
&\bigl|M^{-1} \textstyle\sum_{m=1}^M (\nabla_{\theta_i} \smalll)(\theta,X_{n,m})\bigr|^p\leq M^{-1} \textstyle\sum_{m=1}^M |(\nabla_{\theta_i} \smalll)(\theta,X_{n,m})|^p\\
    &\leq M^{-1} \textstyle\sum_{m=1}^M (c+c\|\theta\|^\fp+c \|X_{n,m}\|^\fp)^p\\
    &\leq \textstyle 3^{p-1} M^{-1} \textstyle\sum_{m=1}^M (c^p+c^p\|\theta\|^{p\fp}+c^p \|X_{n,m}\|^{p\fp}).
    \end{split}
\end{equation}}
\argument{\lref{eq2.5};\lref{eq2'}; the triangle inequality;the assumption that $X_{n,m}$, $n,m\in (\N_0)^2$, are \iid; the fact that for all $p\in \N$ it holds that $\E\bigl[\|X_{0,0}\|^{p\fp}\bigr]<\infty$}{that for all $\theta\in \R^d$, $i\in\{1,2,\dots,d\}$, $p\in \N$ it holds that
\begin{equation}\llabel{eq2.5'}
\begin{split}
    &\textstyle\bigl(\E\bigl[\sup_{N\in \N} |\Advfcom{\theta}{N}{i}|^p\bigr]\bigr)^{1/p}\\
    &\textstyle=\varepsilon^{-1}(1-\beta_1)\bigl(\E\bigl[\bigl(\sum_{n= 0}^\infty(\beta_1)^n| M^{-1} \textstyle\sum_{m=1}^M (\nabla_{\theta_i} \smalll)(\theta,X_{n,m})|\bigr)^p\bigr]\bigr)^{\nicefrac{1}{p}}\\
    &\leq \varepsilon^{-1}(1-\beta_1)\textstyle\sum_{n=0}^\infty \bigl(\E\bigl[\bigl((\beta_1)^n| M^{-1} \textstyle\sum_{m=1}^M (\nabla_{\theta_i} \smalll)(\theta,X_{n,m})|\bigr)^p\bigr]\bigr)^{1/p}\\
    &=\varepsilon^{-1}(1-\beta_1)\textstyle\sum_{n=0}^\infty (\beta_1)^n\bigl(\E\bigl[| M^{-1} \textstyle\sum_{m=1}^M (\nabla_{\theta_i} \smalll)(\theta,X_{n,m})|^p\bigr]\bigr)^{1/p}\\
    &=\varepsilon^{-1}\bigl(\E\bigl[| M^{-1} \textstyle\sum_{m=1}^M (\nabla_{\theta_i} \smalll)(\theta,X_{0,m})|^p\bigr]\bigr)^{1/p}\\
    &\leq\textstyle \varepsilon^{-1}\bigl(3^{p-1} M^{-1} \textstyle\sum_{m=1}^M \bigl(c^p+c^p\|\theta\|^{p\fp}+c^p \E\bigl[\|X_{0,m}\|^{p\fp}\bigr]\bigr)\bigr)^{1/p}\\
    &\leq \textstyle 3\varepsilon^{-1}\bigl(c^p+c^p\|\theta\|^{p\fp}+c^p \E\bigl[\|X_{0,0}\|^{p\fp}\bigr]\bigr)^{1/p}<\infty.
    \end{split}
\end{equation}}
\argument{\lref{eq2.5'};the triangle inequality; the fact that for all $x=(x_1,\dots,x_d)\in \R^d$ it holds that $\|x\|\leq \sum_{i=1}^d|x_i|$}{that for all $\theta\in \R^d$, $p\in \N$ it holds that
\begin{equation}\llabel{eq2.5''}
    \begin{split}
        &\textstyle\bigl(\E\bigl[\sup_{N\in \N} \|\Advf{\theta}{N}\|^p\bigr]\bigr)^{1/p}\leq \textstyle\bigl(\E\bigl[\sup_{N\in \N} (\sum_{i=1}^d|\Advfcom{\theta}{N}{i}|)^p\bigr]\bigr)^{1/p}\\
        &\textstyle \leq \bigl(\E\bigl[\bigl(\sum_{i=1}^d\sup_{N\in \N} |\Advfcom{\theta}{N}{i}|\bigr)^p\bigr]\bigr)^{1/p}\leq \sum_{i=1}^d\bigl(\E\bigl[\sup_{N\in \N} |\Advfcom{\theta}{N}{i}|^p\bigr]\bigr)^{1/p}<\infty.
    \end{split}
\end{equation}}
\argument{\lref{eq2.5''};}{\lref{item 1}\dott}
\startnewargseq
\argument{\lref{eq2.5};the assumption that $X_{n,m}$, $n,m\in (\N_0)^2$, are \iid;the fact that for all $p\in \N$ it holds that $\E[\|X_{0,0}\|^{p\fp}<\infty$}{that for all $\fC\in [0,\infty)$, $\alpha\in (0,1)$, $i\in \{1,2,\dots,d\}$, $p\in \N$ it holds that
\begin{equation}\llabel{eq3''}
     \begin{split}
       &\textstyle \E\biggl[\sup\limits_{\theta\in \R^d,\,\|\theta\|\leq \fC}\textstyle(1-\alpha)\sum_{n= 0}^\infty\alpha^n \bigl|M^{-1} \textstyle\sum_{m=1}^M (\nabla_{\theta_i} \smalll)(\theta,X_{n,m})\bigr|^p\biggr]\\
       &\leq \textstyle 3^{p-1}\E\bigl[(1-\alpha)\sum_{n= 0}^\infty\alpha^n M^{-1} \textstyle\sum_{m=1}^M \bigl(c^p+c^p\fC^{p\fp}+c^p \|X_{n,m}\|^{p\fp}\bigr)\bigr]\\
       &= \textstyle 3^{p-1}\E\bigl[(1-\alpha)\sum_{n= 0}^\infty\alpha^n (c^p+c^p\fC^{p\fp}+c^p \|X_{n,m}\|^{p\fp})\bigr]\\
       &=3^{p-1}c^p+3^{p-1}c^p\fC^{p\fp}+3^{p-1}c^p\,\E\bigl[ \|X_{0,0}\|^{p\fp}\bigr]<\infty.
    \end{split}
\end{equation}}
\argument{\lref{eq3''};}{for every compact $V\subseteq\R^d$ and every $\alpha\in (0,1)$, $i\in \{1,2,\dots,d\}$, $p\in\N$ that
\begin{equation}\llabel{eq3'}
   \textstyle \E\bigl[\sup_{\theta\in V}\textstyle(1-\alpha)\sum_{n= 0}^\infty\alpha^n \bigl|M^{-1} \textstyle\sum_{m=1}^M (\nabla_{\theta_i} \smalll)(\theta,X_{n,m})\bigr|^p\bigr]<\infty.
\end{equation}}
% \argument{\lref{eq3'};}
% {for all $\alpha\in (0,1)$, $\theta\in \R^d$, $i\in\{1,2,\dots,d\}$, $p\in \N$ that
% \begin{equation}\llabel{eq3}
%     \begin{split}
%        &\textstyle \E\bigl[\textstyle(1-\alpha)\sum_{n= 0}^\infty\alpha^n \bigl|M^{-1} \textstyle\sum_{m=1}^M (\nabla_{\theta_i} \smalll)(\theta,X_{n,m})\bigr|^p\bigr]<\infty.
%     \end{split}
% \end{equation}}
% \argument{\lref{eq3};\lref{eq2'};}{that for all $\theta\in \R^d$, $i\in \{1,2,\dots,d\}$ it holds that
% \begin{equation}\llabel{eq4}
%     \begin{split}
%        & \textstyle\E\bigl[\sup_{N\in \N}\|\Advfcom{\theta}{N}{i}\|\big]\leq \textstyle \bigl(\E\bigl[\varepsilon^{-1}(1-\beta_1)\sum_{n= 0}^\infty(\beta_1)^n\bigl| M^{-1} \textstyle\sum_{m=1}^M (\nabla_{\theta_i} \smalll)(\theta,X_{n,m})\bigr|\bigr]\bigr)^{1/p}\\
% &\textstyle=\varepsilon^{-1}\sum_{i=1}^d\E\bigl[(1-\beta_1)\sum_{n= 0}^\infty(\beta_1)^n\bigl| M^{-1} \textstyle\sum_{m=1}^M (\nabla_{\theta_i} \smalll)(\theta,X_{n,m})\bigr|\bigr]<\infty.
%     \end{split}
% \end{equation}}
% \argument{\lref{eq4};}{\lref{item 1}\dott}
\startnewargseq
In the following let $\cZ=(\cZ_1,\dots,\cZ_d)\colon\R^d\times\Omega\to\R^d $ satisfy for all  $i\in \{1,2,\dots,d\}$, $\theta\in \R^d$, $\omega\in \Omega$ with $\sum_{n= 0}^\infty(\beta_1)^n\bigl| M^{-1} \textstyle\sum_{m=1}^M (\nabla_{\theta_i} \smalll)(\theta,X_{n,m}(\omega))\bigr|<\infty$ and
 \begin{equation}
 \textstyle \sum_{n= 0}^\infty(\beta_2)^n\bigl| M^{-1} \textstyle\sum_{m=1}^M (\nabla_{\theta_i} \smalll)(\theta,\allowbreak X_{n,m}(\omega))\bigr|^2<\infty
\end{equation}
 that
\begin{equation}\llabel{def: fZ}
    \cZ_i(\theta,\omega)=\frac{(1-\beta_1)\textstyle\sum_{n= 0}^\infty(\beta_1)^n M^{-1} \textstyle\sum_{m=1}^M(\nabla_{\theta_i} \smalll)(\theta,X_{n,m}(\omega))}{\textstyle\varepsilon+\bigl[(1-\para_2)\sum_{n=0}^\infty(\beta_2)^n\bigl|\textstyle M^{-1} \sum_{m=1}^M(\nabla_{\theta_i} \smalll)(\theta,X_{n,m}(\omega))\bigr|^2\bigr]^{1/2}}.
\end{equation}
\argument{\lref{eq3'};}{that for every compact $V\subseteq \R^d$  it holds that
\begin{equation}\llabel{eq5}
\begin{split}
  & \textstyle \P\bigl(\forall\, i\in \{1,2,\dots,d\}\colon\sup_{\theta\in V}\sum_{n= 0}^\infty(\beta_1)^n\bigl| M^{-1} \textstyle\sum_{m=1}^M (\nabla_{\theta_i} \smalll)(\theta,X_{n,m}(\omega))\bigr|<\infty \bigr)\\
  &=\textstyle \P\bigl(\forall\, i\in \{1,2,\dots,d\}\colon\sup_{\theta\in V}\sum_{n= 0}^\infty(\beta_2)^n\bigl| M^{-1} \textstyle\sum_{m=1}^M (\nabla_{\theta_i} \smalll)(\theta,X_{n,m}(\omega))\bigr|^2<\infty \bigr)\\
  &=1.
   \end{split}
\end{equation}}
\argument{\lref{eq5};\lref{def: Adam vector field};\lref{def: fZ}}{for all $\theta\in \R^d$, $i\in \{1,2,\dots,d\}$ that
\begin{equation}\llabel{eq6}
   \textstyle \P\bigl(\limsup_{N\to\infty}|\Advfcom{\theta}{N}{i}-\cZ_i(\theta)|=0\bigr)=1.
\end{equation}}
\argument{\lref{eq6};}{for all $\theta\in \R^d$ that
\begin{equation}\llabel{eq7}
    \textstyle \P\bigl(\limsup_{N\to\infty}\|\Advf{\theta}{N}-\cZ(\theta)\|=0\bigr)=1.
\end{equation}}
\argument{\lref{eq7};\lref{item 1};the dominated convergence theorem}{that for all $\theta\in \R^d$, $p\in (0,\infty)$ it holds that $\E\bigl[\|\cZ(\theta)\|^p\bigr]<\infty$ and
\begin{equation}\llabel{eq8'}
    \limsup_{N\to\infty} \E\bigl[\|\Advf{\theta}{N}-\cZ(\theta)\|^p\bigr]=0.
\end{equation}}
\argument{\lref{eq8'}}{for all $\theta\in \R^d$, $p\in (0,\infty)$ that
\begin{equation}\llabel{eq8}
    \limsup_{N\to\infty} \|\E[\Advf{\theta}{N}]-\E[\cZ(\theta)]\|=0.
\end{equation}}
In the following let $f\colon \R^d\to\R^d$ satisfy for all $\theta\in \R^d$ that
\begin{equation}\llabel{def: f}
    f(\theta)=\E[\cZ(\theta)].
\end{equation}
\startnewargseq
\argument{\lref{eq8};\lref{def: f}}{that for all $\theta\in \R^d$ it holds that
\begin{equation}\llabel{eq9}
     \limsup_{N\to\infty} \|\E[\Advf{\theta}{N}]-f(\theta)\|=0.
\end{equation}}
\argument{\lref{eq9};}{\cref{item 2: theo: convergence point of Adam copy 3 pre}\dott}
\startnewargseq
\argument{\lref{eq5}; the fact that for all $x\in \R^\dimX$ it holds that $\R^d\ni\theta=(\theta_1,\dots,\theta_d)\mapsto (\nabla_{\theta_i} \smalll)(\theta,x)\in \R$ is continuous}{that for all $(\theta_n)_{n\in \N_0}$, $i\in \{1,2,\dots,d\}$ with $\limsup_{n\to\infty}\|\theta_n-\theta_0\|=0$ it holds $\P$-a.s.\ that
\begin{equation}\llabel{eq10'}
    \limsup_{n\to\infty}(1-\beta_1)\textstyle\sum_{n= 0}^\infty(\beta_1)^n M^{-1} \textstyle\sum_{m=1}^M|(\nabla_{\theta_i} \smalll)(\theta_n,X_{n,m})-(\nabla_{\theta_i} \smalll)(\theta_0,X_{n,m})|=0
\end{equation}
\begin{multline}\llabel{eq10}
   \text{and}\qquad  \limsup_{n\to\infty}(1-\beta_1)\textstyle\sum_{n= 0}^\infty(\beta_1)^n \bigl|(M^{-1} \textstyle\sum_{m=1}^M(\nabla_{\theta_i} \smalll)(\theta_n,X_{n,m}))^2\\-(M^{-1} \textstyle\sum_{m=1}^M(\nabla_{\theta_i} \smalll)(\theta_0,X_{n,m}))^2\bigl|=0.
\end{multline}}
\argument{\lref{eq10};}{that for all $(\theta_n)_{n\in \N_0}$, $i\in \{1,2,\dots,d\}$ with $\limsup_{n\to\infty}\|\theta_n-\theta_0\|=0$ it holds $\P$-a.s.\ that
\begin{equation}\llabel{eq11}
    \limsup_{n\to\infty}\|\cZ(\theta_n)-\cZ(\theta_0)\|=0.
\end{equation}}
\argument{\lref{eq3'};\lref{def: fZ};\lref{eq5};}{that for all compact $V\subseteq \R^d$ it holds that
\begin{equation}\llabel{eq12}
\begin{split}
   &\textstyle \E\bigl[\sup_{\theta\in V} \|\cZ(\theta)\|\bigr]\leq \textstyle\sum_{i=1}^d\E\bigl[\sup_{\theta\in V} |\cZ_i(\theta)|\bigr]\\
   &\textstyle\leq\sum_{i=1}^d \textstyle\E\bigl[\sup_{\theta\in V}\textstyle(1-\beta_1)\sum_{n= 0}^\infty(\beta_1)^n \bigl|M^{-1} \textstyle\sum_{m=1}^M (\nabla_{\theta_i} \smalll)(\theta,X_{n,m})\bigr|^p\bigr]<\infty.
   \end{split}
\end{equation}}
\argument{\lref{eq11};\lref{eq12};the dominated convergence theorem}{for all $(\theta_n)_{n\in \N_0}$, $i\in \{1,2,\dots,d\}$ with $\limsup_{n\to\infty}\|\theta_n-\theta_0\|=0$ that
\begin{equation}\llabel{eq13}
    \limsup_{n\to\infty}\|\E[\cZ(\theta_n)]-\E[\cZ(\theta_0)]\|=0.
\end{equation}}
\argument{\lref{eq13};\lref{def: f}}{\lref{item 3}\dott}
\startnewargseq
\argument{\lref{eq7};}{\lref{item 4}\dott}
\startnewargseq
\argument{\lref{eq7};}{that for every $\theta\in \R^d$ and every random variable $\fZ\colon \Omega\to\R^d$ with $\limsup_{N\to\infty} \E\bigl[\min\{1,\allowbreak\|\fZ-\Advf{\theta}{N}\|\bigr]=0$ it holds $\P$-a.s.\ that
\begin{equation}\llabel{eq14}
    \cZ(\theta)=\fZ.
\end{equation}}
\argument{\lref{eq14};\lref{eq8'};\lref{def: f};}{for every $\theta\in \R^d$, $p\in (0,\infty)$ and every random variable $\fZ\colon \Omega\to\R^d$ with $\limsup_{N\to\infty} \E\bigl[\min\{1,\allowbreak\|\fZ-\Advf{\theta}{N}\|\bigr]=0$ that
\begin{equation}\llabel{eq15}
    \E\bigl[\|\fZ\|^p\bigr]=\E\bigl[\|\cZ(\theta)\|^p\bigr]<\infty\qqandqq \E[\fZ]=\E[\cZ(\theta)]=f(\theta).
\end{equation}}
\argument{\lref{eq15};}{\lref{item 5}\dott}
\startnewargseq
\argument{\lref{eq7};\lref{eq8'};\lref{eq11};}{\lref{item 6}\dott}
\end{aproof}
\renewcommand{\mandV}[2]{\Theta^{#1}_{#2}}
\renewcommand{\mandVcom}[3]{\Theta^{#1}_{#2,#3}}
\begin{athm}{theorem}{theo: convergence point of Adam copy 2}[\textcolor{red}{\Adam\ limit theorem}]
    Let $( \Omega, \cF, \P )$ be a probability space, 
let $ \fd, \mathscr{d} \in \N $, 
let $X_{ n, m } \colon \Omega\to \R^{ \mathscr{d} }$, $ (n,m) \in (\N_0)^2 $, be \iid\ random variables, let $\smalll \in C^2( \R^\fd \times \R^{ \scrd }, \R )$ have at most polynomially growing derivatives, let $M \in \N$, $\beta_1 \in (0,1)$, $\beta_2 \in ((\beta_1)^2,1)$, $\varepsilon \in (0,\infty)$, let $( \gamma_n )_{ n\in \N }\subseteq (0,\infty)$ be non-increasing, for every $k,n\in \N_0$ let 
    $
\mandV{k}{n}=(\mandVcom{k}{n}{1},\dots,\mandVcom{k}{n}{\fd})\colon\allowbreak \Omega\to\R^\fd$ be a random variable, assume for all $k\in \{1,2\}$, $i\in \{1,2,\dots,\fd\}$ that
    \begin{equation}
    \mandV{k}{0}=0,\qquad \mandVcom{k}{n}{i}= \beta_k \mandVcom{k}{n-1}{i}+(1-\beta_k)\bigl[\textstyle \frac 1M \sum_{m=1}^M(\nabla_{\theta_i} \smalll)(\mandV{0}{n-1},X_{n,m})\bigr]^k
        \end{equation}
        \begin{equation}
 \text{and}\qquad\mandVcom{0}{n}{i}=\mandVcom{0}{n-1}{i} -\gamma_n [1-(\para_1)^n]^{-1}\bigl[\varepsilon+\bigl[[1-(\para_2)^n]^{-1}\mandVcom{2}{n}{i}\bigr]^{\nicefrac{1}{2}}\bigr]^{-1}\mandVcom{1}{n}{i}, 
     \end{equation}
 assume\footnote{Note that for all $x \in \R$ it holds that $\lceil x \rceil = \min( \{ z \in \Z \colon z \geq x \} )$.}  $\inf_{n\in \N}\gamma_n=0<\sup_{r\in (0,1)}\liminf_{n\to\infty}((\gamma_n)^{-1}\gamma_{n+\lceil(\gamma_n)^{-1}r\rceil})$, assume for all $m\in \N$ that $\E[\|X_{1,1}\|^m ] <\infty= \sum_{ n = 1 }^{ \infty } \gamma_n $, assume that $\mandV{0}{0}$ and $(X_{n,m}\allowbreak)_{(n,m)\in \N^2}$ are independent, 
let $ \vartheta \colon\Omega\to \R^\fd $ be a random variable which satisfies 
$
  \limsup_{ n \to \infty }\E[ \min\{ 1, \| \mandV{0}{n} - \vartheta \|\} ] = 0
$, 
and let $f =(f_1,\dots,f_\fd)\colon \R^\fd \to \R^\fd$ satisfy for all $\theta=(\theta_1,\dots,\theta_\fd)\in \R^\fd$, $i\in \{1,2,\dots,\fd\}$ that 
    \begin{equation}\llabel{def: Adam vector field}
          f_i(\theta)=\lim_{N\to\infty}\E\Biggl[\frac{(1-\para_1)\textstyle\sum_{n= 0}^N(\para_1)^n M^{-1} \textstyle\sum_{m=1}^M(\nabla_{\theta_i} \smalll)(\theta,X_{n,m})}{\textstyle\varepsilon+\sqrt{(1-\para_2)\sum_{n=0}^N(\para_2)^n\bigl|\textstyle M^{-1} \sum_{m=1}^M(\nabla_{\theta_i} \smalll)(\theta,X_{n,m})\bigr|^2}}\Biggr]
    \end{equation}
    (cf.\ \cref{theo: convergence point of Adam copy 3 pre}).
    Then $\P( f( \vartheta )= 0 ) = 1$.
\end{athm}
\begin{aproof}
\argument{\cite[Corollary 1.10]{DereichJentzenKassing2025};}{ that $\E\bigl[ \| f( \vartheta )\| \bigr] = 0$\dott}
\end{aproof}
\renewcommand{\mandV}[4]{\mathbb{V}^{#1,#2,#3}_{#4}}
\renewcommand{\mandVcom}[5]{\mathbb{V}^{#1,#2,#3}_{#4,#5}}
\subsection{Adam converges to wrong points}
\begin{athm}{lemma}{prove order measure}
    Let $(\Omega,\cF,\P)$ be a probability space, let $p\in (-\infty,0)$, $q\in (0,\infty)$, $\varrho\in (0,1)$ satisfy $-p>q$ and $p\varrho+(1-\varrho)q=0$, let $X \colon \Omega\to\{0,q\}$ and $Y\colon \Omega\to\{0,-p\}$ be random variables, and assume 
    \begin{equation}\llabel{def: XY}
        \P(X=0)=\varrho\qqandqq \P(Y=0)=1-\varrho.
    \end{equation}
    Then there exists a random variable $Z\colon\Omega\to \{0,-p\}$ which satisfies
    \begin{enumerate}[label=(\roman*)]
        \item \label{item 1: prove order measure} it holds that $\P(Z=0)=\P(Y=0)$, $\P(Z=-p)=\P(Y=p)$,
        \item \label{item 2: prove order measure} it holds $\P$-a.s.\ that $\E[Z|X]= X$, and
        \item \label{item 3: prove order measure} it holds that $\P(X=Z)<1$.
    \end{enumerate}
\end{athm}
\begin{aproof}
\argument{the assumption that $-p>q$;the assumption that $p\varrho+(1-\varrho)q=0$}{that \llabel{arg1} $\varrho<\frac 12$\dott}
    \argument{\lref{arg1};}{that there exists $Z\colon \Omega\to\{0,-p\}$ which satisfies
        \begin{equation}\llabel{def: Z1}
            \P(Z=0,X=0)=\varrho,\qquad \P(Z=-p,X=0)=0,        \end{equation}
        \begin{equation}\llabel{def: Z2}
          \P(Z=0,X=q)=1-2\varrho,\qquad  \text{and}\qquad \P(Z=-p,X=q)=\varrho.
        \end{equation}}
        \startnewargseq
        \argument{\lref{def: XY};\lref{def: Z1};\lref{def: Z2};the assumption that for all $\omega\in \Omega$ it holds that $X(\omega)\in \{0,q\}$}{that 
        \begin{equation}\llabel{evidence 1}
            \P(Z=0)=\P(Z=0,X=0)+\P(Z=0,X=q)=1-2\varrho+\varrho=1-\varrho=\P(Y=0)
        \end{equation}
        and
        \begin{equation}\llabel{evidence 2}
            \P(Z=-p)=\P(Z=-p, X=0)+\P(Z=-p,X=q)=0+\varrho=\varrho=\P(Y=q).
        \end{equation}}
        \argument{\lref{def: Z1};\lref{def: Z2};the assumption that $p\varrho+(1-\varrho)q=0$}{that it holds that
        \begin{equation}\llabel{eq3}
            \E[Z|X=0]=\frac{0\P(Z=0,X=0)+(-p)\P(Z=-p,X=0)}{\P(Z=0,X=0)+\P(Z=-p,X=0)}=\frac{0\varrho+(-p)0}{\varrho+0}=0.
        \end{equation}
        and
        \begin{equation}\llabel{eq4}
            \E[Z|X=q]=\frac{0\P(Z=0,X=q)+(-p)\P(Z=-p,X=q)}{\P(Z=0,X=q)+\P(Z=-p,X=q)}=\frac{0(1-2\varrho)+(-p)\varrho}{1-2\varrho+\varrho}=q.
        \end{equation}}
        \argument{\lref{eq4};the assumption that for all $\omega\in \Omega$ it holds that $X(\omega)\in \{0,q\}$}{that it holds $\P$-a.s.\ that
        \begin{equation}\llabel{evidence 3}
            \E[Z|X]=X.
        \end{equation}}
        \argument{\lref{def: Z2}}{that
        \begin{equation}\llabel{evidence 4}
            \P(X\neq Z)\geq \rho>0.
        \end{equation}}
        \argument{\lref{evidence 1};\lref{evidence 2};\lref{evidence 3};\lref{evidence 4}}{items \ref{item 1: prove order measure} and \ref{item 2: prove order measure}\dott}
\end{aproof}
In the next notion, \cref{definition: two point distributed} below, we briefly recall the concept of a two-point distributed random variable. 
\begin{definition}\label{definition: two point distributed}
    Let $( \Omega, \cF, \P )$ be a probability space, let $\fd \in \N$, and let $X \colon \Omega \to \R^\fd$ be a random variable. Then we say that $X$ is $\P$-\tpd\ (we say that $X$ is \tpd) if and only if there exist $v, w \in \R^\fd$ such that $\P( X \in \{ v, w \} ) = 1$.
\end{definition}
\renewcommand{\cL}{\mathbb{L}}
\begin{athm}{prop}{lem: non-convergence of Adam2'}
    Let $(\Omega,\cF,\P)$ be a probability space, let $X_{n} \colon \Omega \to \R$, $n \in \N$, be \iid\ \tpd\ random variables, assume $\E[X_1]=0$, assume that $X_1$ is not symmetrically distributed, let $\cL=(\cL(\theta,x))_{(\theta,x)\in \R^2} \colon \R^2\allowbreak\to\R$ satisfy for all $\theta,x\in \R$ that $\cL(\theta,x)=|\theta-x|^2$, let $\para_1\in (0,1)$, $\para_2\in ((\para_1)^2,1)$, $\varepsilon\in (0,\infty)$, let $(\gamma_n)_{n\in \N}\subseteq (0,\infty)$ be non-increasing, for every $n\in \N_0$ let  $
\democrator_{n}^{k}\colon\allowbreak \Omega\to\R$, $k\in\{1,2\}$, and $\Theta_n\colon\allowbreak \Omega\to\R$ be random variables, assume for all $n\in \N$, $k\in \{1,2\}$ that
    \begin{equation}\llabel{def: V}
    \democrator_0^{k}=0,\qquad \democrator_{n}^{k}= \beta_k \democrator_{n-1}^{k}+(1-\beta_k)[(\nabla_{\theta} \cL)(\Theta_{n-1},X_{n})]^k,
        \end{equation}
        \begin{equation}\llabel{def: Theta}
    \textstyle\text{and}\qquad \Theta_{n}=\Theta_{n-1} -\gamma_n\Bigl[\frac{\democrator_{n}^{1}}{1-(\para_1)^n}\Bigr] \Bigl[\varepsilon+\Bigl[\frac{\democrator_{n}^{2}}{1-(\para_2)^n}\Bigr]^{\nicefrac{1}{2}}\Bigr]^{-1}, 
     \end{equation}
     assume 
        $\inf_{n\in \N}\gamma_n=0<\sup_{r\in (0,1)}\liminf_{n\to\infty}((\gamma_n)^{-1}\gamma_{n+\lceil(\gamma_n)^{-1}r\rceil})$, and assume that $\Theta_0$ and $(X_{n})_{n\in \N}$ are independent \cfload. Then 
        \begin{equation}\llabel{conclude}
            \limsup_{ n \to \infty }\E\bigl[ \min\{ 1, | \Theta_n |\} \bigr] > 0.
            \end{equation}
\end{athm}
\begin{aproof}
    Throughout this proof let $\varrho\colon (\N_0)^2\to \N$ be a bijection, for every $n,m\in \N_0$ let $\bfX_{n,m}\colon \Omega\to\R$ satisfy $\bfX_{n,m}=X_{\varrho(n,m)}$, let $\enpointp,\enpointq\in \R$ satisfy $\P(X_1\in \{\enpointp,\enpointq\})=1$, assume without loss of generality that  $-\enpointp > \enpointq$, let $\varrho\in \R$ satisfy $\varrho= \P( X_{ 1} = \enpointp )$,
 let $f \colon \R \to \R$ satisfy for all $\theta \in \R$ that
\begin{multline}\llabel{def: Adam}
 \textstyle f(\theta)=(1-\beta_1) \,\E\biggl[ \textstyle \bigl[\sum_{n=0}^\infty (\beta_1)^{n} (\nabla_{\theta} \cL)(\theta,\bfX_{n,1})\bigr]\\
  \textstyle\cdot\Bigl[\bigl[(1-\beta_2) \sum_{n=0}^\infty (\beta_2)^{n} |(\nabla_{\theta} \cL)(\Theta_{n-1},\bfX_{n,1})|^{2}\bigr]^{\nicefrac{1}{2}}+\varepsilon\Bigr]^{-1}\biggr],
\end{multline}
 for every $n \in \N_0$ let $\mathbb{F}_n \subseteq \cF$ satisfy
\begin{equation}\llabel{def: bbF}
  \mathbb F_n = \sigma\bigl(
    ( \bfX_{ v,1} )_{ v \in ( \N_0 \backslash \{n\} )}
  \bigr),
  \end{equation}
  and let $\bfV_n^1\colon\Omega\to\R$, $\bfV_n^2\colon\Omega\to\R$, $\varTheta_n\colon \Omega\to\R$, $n\in \N_0$, satisfy for all $n\in \N$, $k\in \{1,2\}$ that 
 \begin{equation}\llabel{def: V2}
    \bfV_0^{k}=0,\qquad \democrator_{n}^{k}= \beta_k \bfV_{n-1}^{k}+(1-\beta_k)[(\nabla_{\theta} \cL)(\varTheta_{n-1},\bfX_{n,1})]^k,
        \end{equation}
        \begin{equation}\llabel{def: Theta2}
  \varTheta_0=\Theta_0,\qquad  \textstyle\text{and}\qquad \varTheta_{n}=\varTheta_{n-1} -\gamma_n\Bigl[\frac{\bfV_{n}^{1}}{1-(\para_1)^n}\Bigr] \Bigl[\varepsilon+\Bigl[\frac{\bfV_{n}^{2}}{1-(\para_2)^n}\Bigr]^{\nicefrac{1}{2}}\Bigr]^{-1}.
     \end{equation}
\argument{the fact that $\bfX_{ 1,1 }$ is centered}{that
\begin{equation}\llabel{eqq1}
    0 = \E[ \bfX_{ 1,1 } ] = \enpointp \P( X_{ 1, 1 } = \enpointp ) + \enpointq \P( \bfX_{ 1, 1 } = \enpointq )
= \enpointp \varrho + \enpointq ( 1 - \varrho ) = \enpointq + \varrho ( \enpointp - \enpointq ).
\end{equation}}
\argument{\lref{def: Adam};the assumption that for all $\theta,x\in \R$ it holds that $\cL(\theta,x)=|\theta-x|^2$}{that for all $\theta\in \R$ it holds that
\begin{equation}\llabel{eq1}
\begin{split}
        &f( \theta )\\
       & =
  \sum_{ n = 0 }^{ \infty }
  (\beta_1)^n ( 1 - \beta_1 ) \E\biggl[ 2(\theta-\bfX_{n,1}) \Bigl[\Bigl[(1-\beta_2) \textstyle\sum\limits_{k=0}^\infty (\beta_2)^{k} |2(\theta-\bfX_{n,1})|^{2}\Bigr]^{\nicefrac{1}{2}}+\varepsilon\Bigr]^{-1}\biggr].
  \end{split}
\end{equation}}
\argument{the fact that $\varrho=\P(\bfX_{1,1}=p)$; the fact that $\E[\bfX_{1,1}]=0$; the assumption that $\enpointp<0$}{that \llabel{arg1} $\varrho>0$\dott}
\argument{\lref{eqq1};\lref{arg1};the fact that $\enpointp^2>\enpointq^2$; the fact that $\bfX_{n,1}$, $n\in \N_0$, are \iid;}{that for all $n \in \N_0$, 
$c \in [0,\infty)$, $\kappa \in (0,\infty)$ it holds that
\begin{equation}\llabel{eq3}
\begin{split}
\E\biggl[ 
  \frac{ -2\bfX_{n,1}}{
    \varepsilon + [ c + 
    \kappa |-2\bfX_{n,1}|^2 ]^{ 1 / 2 }
  } 
\biggr]
&= \E\biggl[ \frac{ - 2 \bfX_{ 1,1 } }{ \varepsilon + [ c + 4 \kappa | X_{ 1} |^2 ]^{ 1 / 2 }} \biggr] \\
&= 
\frac{ - 2 \enpointp \varrho }{ \varepsilon + [ c + 4 \kappa \enpointp^2 ]^{ 1 / 2 } }
+
\frac{ - 2 \enpointq ( 1 - \varrho ) }{ \varepsilon + [ c + 4 \kappa \enpointq^2 ]^{ 1 / 2 } }\\
&< \frac{ - 2 \enpointp \varrho }{ \varepsilon + [ c + 4 \kappa \enpointq^2 ]^{ 1 / 2 } }
+
\frac{ - 2 \enpointq ( 1 - \varrho ) }{ \varepsilon + [ c + 4 \kappa \enpointq^2 ]^{ 1 / 2 } }\\
&=\frac{ - 2 \enpointp \varrho- 2 \enpointq ( 1 - \varrho ) }{ \varepsilon + [ c + 4 \kappa \enpointq^2 ]^{ 1 / 2 } }\\
&=0.
\end{split}
\end{equation}}
\argument{\lref{def: bbF};the fact that $\bfX_{n,1}$, $n\in \N_0$, are independent}{that \llabel{arg2} for all $n\in \N_0$ it holds that $\sigma(\bfX_{n,1})$ and $\mathbb F_n$ are independent\dott}
\argument{\lref{eq3};\lref{arg2};the factorization lemma for conditional expectations (see, \eg, \cite[Proposition 3.15]{DeRoAr2024nonconvergence})}{that for all $n\in \N_0$ it holds $\P$-a.s.\ that
\begin{equation}\llabel{eq4}
\begin{split}
  \E\biggl[ -2\bfX_{n,1} \Bigl[\Bigl[(1-\beta_2) \textstyle\sum\limits_{k=0}^\infty (\beta_2)^{k} |-2\bfX_{n,1}|^{2}\Bigr]^{\nicefrac{1}{2}}+\varepsilon\Bigr]^{-1}\Big| \mathbb F_n\biggr]< 0.
   \end{split}
\end{equation}}
\argument{\lref{eq4};}{for all $n\in \N_0$ that
\begin{equation}\llabel{eq5}
    \E\biggl[ -2\bfX_{n,1} \Bigl[\Bigl[(1-\beta_2) \textstyle\sum\limits_{k=0}^\infty (\beta_2)^{k} |-2\bfX_{n,1}|^{2}\Bigr]^{\nicefrac{1}{2}}+\varepsilon\Bigr]^{-1}\biggr]< 0.
\end{equation}}
\argument{\lref{eq5};\lref{eq1}}{that
\begin{equation}\llabel{eq6}
    f(0)=\sum_{ n = 0 }^{ \infty }
  (\beta_1)^n ( 1 - \beta_1 ) \E\biggl[ -2\bfX_{n,1} \Bigl[\Bigl[(1-\beta_2) \textstyle\sum\limits_{k=0}^\infty (\beta_2)^{k} |-2\bfX_{n,1}|^{2}\Bigr]^{\nicefrac{1}{2}}+\varepsilon\Bigr]^{-1}\biggr]<0.
\end{equation}}
\argument{\lref{eq6};\lref{def: V2};\lref{def: Theta2};\cref{theo: convergence point of Adam copy 2} (applied with $\Omega \curvearrowleft\Omega$, $\cF\curvearrowleft\cF$, $\P\curvearrowleft\P$, $d\curvearrowleft$, $\scrd\curvearrowleft d$, $(X_{n,m})_{(n,m)\in (\N_0)^2}\curvearrowleft (\bfX_{n,m})_{(n,m)\in(\N_0)^2}$, $\bbL\curvearrowleft (\R^2\ni (\theta,x)\mapsto \|\theta-x\|^2\in \R)$, $M\curvearrowleft 1$, $\beta_1\curvearrowleft\beta_1$, $\beta_2\curvearrowleft\beta_2$, $\varepsilon \curvearrowleft\varepsilon$, $(\gamma_n)_{n\in \N}\curvearrowleft(\gamma_n)_{n\in \N}$, $(\Theta_n^0)_{n\in \N_0} \curvearrowleft (\varTheta_n)_{n\in \N_0}$, $(\Theta_n^k)_{(k,n)\in \{1,2\}\times \N_0}\curvearrowleft (\bfV_{n}^k)_{(k,n)\in \{1,2\}\times \N_0}$ in the notation of \cref{theo: convergence point of Adam copy 2})}{that
\begin{equation}\llabel{eq7}
 \limsup_{ n \to \infty }\E\bigl[ \min\{ 1, | \varTheta_n |\} \bigr] > 0.
\end{equation}}
\argument{the assumption that $X_n$, $n\in \N$, are \iid;the fact that $\varrho$ is a bijection;the fact that for all $(n,m)\in (\N_0)^2$ it holds that $\bfX_{n,m}=X_{\varrho(n,m)}$}{that \llabel{argg1} $(\bfX_{n,1})_{n\in \N_0}$ and $(X_n)_{n\in \N_0}$ are identically distributed\dott}
\argument{\lref{argg1};\lref{def: V};\lref{def: Theta};\lref{def: V2};\lref{def: Theta2}}{that \llabel{argg2} $(\Theta_n)_{n\in \N_0}$ and $(\varTheta_n)_{n\in \N_0}$ are identically distributed\dott}
\argument{\lref{argg2};\lref{eq7}}{\lref{conclude}\dott}
\end{aproof}
\begin{athm}{cor}{lem: non-convergence of Adam2}
    Let $(\Omega,\cF,\P)$ be a probability space, let $\fd\in \N$, let $X_{ n } \colon \Omega \to \R^d$, $n \in \N$, be \iid\ \tpd\ random variables, assume $X_1$ is not symmetrically distributed, let $\cL=(\cL(\theta,x))_{(\theta,x)\in \R^\fd\times\R^\fd} \colon \R^\fd\times\R^\fd\to\R$ satisfy for all $\theta,x\in \R^\fd$ that $\cL(\theta,x)=\|\theta-x\|^2$, let $\para_1\in (0,1)$, $\para_2\in ((\para_1)^2,1)$, $\varepsilon\in (0,\infty)$, let $(\gamma_n)_{n\in \N}\subseteq (0,\infty)$ be non-increasing, for every $n\in \N_0$ let  $
\democrator_{n}^{k}=(\democrator_{n,1}^{k},\dots,\democrator_{n,\fd}^{k})\colon\allowbreak \Omega\to\R^\fd$, $k\in\{1,2\}$, and $\Theta_n=(\Theta_{n,1},\dots,\Theta_{n,\fd})\colon\allowbreak \Omega\to\R^\fd$ be random variables, assume for all $n\in \N$, $k\in \{1,2\}$, $i\in \{1,2,\dots,\fd\}$ that
    \begin{equation}\llabel{def: bbV}
    \democrator_0^{k}=0,\qquad \democrator_{n,i}^{k}= \beta_k \democrator_{n-1,i}^{k}+(1-\beta_k)\bigl[\textstyle (\nabla_{\theta_i} \cL)(\Theta_{n-1},X_{n})\bigr]^k,
        \end{equation}
        \begin{equation}
    \textstyle\text{and}\qquad \Theta_{n,i}=\Theta_{n-1,i} -\gamma_n\Bigl[\frac{\democrator_{n,i}^{1}}{1-(\para_1)^n}\Bigr] \Bigl[\varepsilon+\Bigl[\frac{\democrator_{n,i}^{2}}{1-(\para_2)^n}\Bigr]^{\nicefrac{1}{2}}\Bigr]^{-1}, 
     \end{equation}
     assume
       $\inf_{n\in \N}\gamma_n=0<\sup_{r\in (0,1)}\liminf_{n\to\infty}((\gamma_n)^{-1}\gamma_{n+\lceil(\gamma_n)^{-1}r\rceil})$, and assume that $\Theta_0$ and $(X_{n})_{n\in \N}$ are independent \cfload. Then
        \begin{equation}\llabel{conclude}
           \limsup_{ n \to \infty }\E\bigl[ \min\{ 1, \| \Theta_n - \E[X_1]\|\} \bigr] > 0.
           \end{equation}
\end{athm}
\begin{aproof}
Throughout this proof assume without loss of generality that $\E[X_1]=0$, for every $n\in \N_0$ let $\bfX_{n}=(\bfX_{n}^1,\dots,\bfX_{n}^\fd)\allowbreak\colon\allowbreak \Omega\to\R^d$ satisfy $\bfX_{n}=X_{n}$ and let $\fp=(\fp_1,\dots,\fp_\fd)$, $\fq=(\fq_1,\dots,\fq_\fd)\in \R^\fd$ satisfy $\P(\bfX_1\in \{\fp,\fq\})=1$.
\argument{the assumption that $X_1$ is not symmetrically distributed; the assumption that $\E[X_1]=0$; the fact that $\P(X_1\in \{\fp,\fq\})=1$}{that \llabel{arg1} $-\fp\neq \fq$\dott}
\argument{\lref{arg1};}{that there exists $i\in \{1,2,\dots,\fd\}$ which satisfies
\begin{equation}\llabel{def: i}
    -\fp_i\neq \fq_i.
\end{equation}
}
\startnewargseq
\argument{\cref{lem: non-convergence of Adam2'};\lref{def: i}; the fact that $\P(\bfX^i_{1}\in \{\fp_i,\fq_i\})$}{that 
\begin{equation}\llabel{eq1}
     \limsup_{n\to\infty}\E\bigl[ | \Theta_{n,i}\allowbreak- \E[ \bfX_{ 1 }^i ] | \bigr]>0.
\end{equation}}
\argument{\lref{eq1};}{\lref{conclude}\dott}
\end{aproof}

\renewcommand{\cL}{\mathcal{L}}
\section{Convergence rates for Adam for symmetrically distributed data}\label{sec: convergence symmetric}
In this section we establish in \cref{prop: convergence in symmetric case} in \cref{subsec: convegence symmetry} below for the class of \SOPs\ considered in \cref{main theorem 2} that \Adam\ converges to a critical point of the \SOP\ if the data of the \SOP\ is
\emph{symmetrically distributed} and we combine the conclusion of \cref{prop: convergence in symmetric case} with \cref{lem: non-convergence of Adam2} from \cref{sec: converge to wrong point} to establish in \cref{main theorem 2'} in \cref{subsec: chracterization} below a slightly stronger variant of the \Adam\ symmetry theorem in \cref{main theorem 2} in the introduction. \cref{main theorem 2} then follows directly from  \cref{main theorem 2'}.
\subsection{Half-concavity property for the Adam vector field}
\begin{athm}{lemma}{Adam concave}[Half-concavity property for the Adam vector field]
    Let $ \cons \in [0,\infty) $, $ \kappa, \varepsilon \in (0,\infty) $ and let $\adamvector \colon \R \to \R$ satisfy for all $x \in \R$ that
\begin{equation}\llabel{def: adam}
  \adamvector( x ) = x \bigl[ [ \cons + \kappa x^2 ]^{ 1 / 2 } + \varepsilon \bigr]^{ - 1 } .
\end{equation}
Then 
\begin{enumerate}[label=(\roman*)]
    \item \label{item 1: half convex} it holds for all $x \in \R$ that $\adamvector(x) = - \adamvector( - x )$, 
     \item \label{item 3: half convex} it holds that $f|_{ \R\backslash\{0\}} \in C^{ \infty }( \R\backslash\{0\}, \R )$, 
    \item \label{item 2: half convex} it holds that $f\in C^1(\R,\R)$, 
    \item \label{item 4: half convex} it holds for all $x\in \R$ that $f'(x)>0$,
    and
    \item \label{item 5: half convex} it holds for all $x \in (0,\infty)$ that $\adamvector''(x) < 0<f(x)$.
\end{enumerate}
\end{athm}
\begin{aproof}
\argument{\lref{def: adam}; the fact that for all $x\in \R$ it holds that $x^2=(-x)^2$}{\cref{item 1: half convex}\dott}
\startnewargseq
\argument{\cref{item 1: half convex};\lref{def: adam};the fact that $((0,\infty)\ni x\mapsto x\in \R)\in C^{\infty}((0,\infty),\R)$;the fact that $((0,\infty)\ni x\mapsto x^2\in (0,\infty))\in C^{\infty}((0,\infty),(0,\infty))$;the fact that $((0,\infty)\ni x\mapsto [\cons+\kappa x]^{1/2}\in (0,\infty))\in C^{\infty}((0,\infty),(0,\infty))$; the fact that $((0,\infty)\ni x\mapsto (x+\varepsilon)^{-1}\in \R)\in C^{\infty}((0,\infty),\R)$; the chain rule}{\cref{item 3: half convex}\dott}
\startnewargseq
\argument{\cref{item 3: half convex};\lref{def: adam}}{that for all $x\in \R\backslash\{0\}$ it holds that
\begin{equation}\llabel{eq1}
\begin{split}
    \adamvector'(x)&=\frac{1}{\varepsilon+\sqrt{\cons+\kappa x^2}}-\frac{\kappa x^2}{(\varepsilon+\sqrt{\cons+\kappa x^2})^2\sqrt{\cons+\kappa x^2}}.
    \end{split}
\end{equation}}
\argument{\lref{def: adam};}{that for all $x\in \R\backslash\{0\}$ it holds that
\begin{equation}\llabel{eq1.1}
\frac{f(x)-f(0)}{x}= \frac{f(x)}{x}=\bigl[ [ \cons + \kappa x^2 ]^{ 1 / 2 } + \varepsilon \bigr]^{ - 1 }   
\end{equation}}
\argument{\lref{eq1.1};}{that
\begin{equation}\llabel{eq1.2}
   \textstyle \limsup_{\R\backslash \{0\} \ni x\to 0} \bigl|\frac{f(x)-f(0)}{x}-\frac{1}{\varepsilon+\sqrt{\cons}}\bigr|=0.
\end{equation}}
\argument{\lref{eq1}}{that 
\begin{equation}\llabel{eq1.3}
   \frac{1}{\varepsilon+\sqrt{\cons}}=\lim_{\R\backslash\{0\}\ni x\to 0} f'(x).
\end{equation}}
\argument{\lref{eq1.3};\lref{eq1.2};\cref{item 3: half convex}}{\cref{item 2: half convex}\dott}
\startnewargseq
\argument{\lref{eq1};the fact that for all $x\in \R\backslash\{0\}$ it holds that $x^2> 0$; the assumption that $\varepsilon> 0$; the assumption that $c\geq 0$; the assumption that $\kappa>0$}{that for all $x\in \R\backslash\{0\}$ it holds that
\begin{equation}\llabel{eqq1}
\begin{split}
    \adamvector'(x) &>\frac{1}{\varepsilon+\sqrt{\cons+\kappa x^2}}-\frac{\kappa x^2}{(\varepsilon+\sqrt{\cons+\kappa x^2})\sqrt{\cons+\kappa x^2}\sqrt{\cons+\kappa x^2}}\\
    &=\frac{1}{\varepsilon+\sqrt{\cons+\kappa x^2}}-\frac{\kappa x^2}{(\varepsilon+\sqrt{\cons+\kappa x^2})(\cons+\kappa x^2)}\\
    &\geq \frac{1}{\varepsilon+\sqrt{\cons+\kappa x^2}}-\frac{1}{\varepsilon+\sqrt{\cons+\kappa x^2}}=0.
    \end{split}
\end{equation}}
\argument{\lref{eqq1};the fact that $f'(0)=[ \varepsilon+\cons ^{ 1 / 2 }  ]^{ - 1 } >0$}{\cref{item 4: half convex}\dott}
\startnewargseq
\argument{\lref{eq1};\cref{item 3: half convex}}{
that for all $x\in \R\backslash \{0\}$ it holds that
\begin{equation}\llabel{eq2}
\begin{split}
    \adamvector''(x)&=\frac{-\kappa x}{(\varepsilon+\sqrt{\cons+\kappa x^2})^2\sqrt{\cons+\kappa x^2}}-\frac{2\kappa x}{(\varepsilon+\sqrt{\cons+\kappa x^2})^2\sqrt{\cons+\kappa x^2}}\\
    &+\frac{\kappa x^2\bigl[2\kappa x(\varepsilon+\sqrt{\cons+\kappa x^2})+(\varepsilon+\sqrt{\cons+\kappa x^2})^2\frac{\kappa x}{\sqrt{\cons+\kappa x^2}}\bigr]}{(\varepsilon+\sqrt{\cons+\kappa x^2})^4(\cons+\kappa x^2)}\\
    &=\frac{-3\kappa x}{(\varepsilon+\sqrt{\cons+\kappa x^2})^2\sqrt{\cons+\kappa x^2}}+\frac{\kappa x^2\bigl[2\kappa x+(\varepsilon+\sqrt{\cons+\kappa x^2})\frac{\kappa x}{\sqrt{\cons+\kappa x^2}}\bigr]}{(\varepsilon+\sqrt{\cons+\kappa x^2})^3(\cons+\kappa x^2)}.
    \end{split}
\end{equation}}
\argument{\lref{eq2};}{for all $x\in \R\backslash\{0\}$ that
\begin{equation}\llabel{eq3}
\begin{split}
    \adamvector''(x)
    &\textstyle=\frac{\kappa x}{(\varepsilon+\sqrt{\cons+\kappa x^2})^3(\cons+\kappa x^2)} \Bigl[-3\bigl(\varepsilon+\sqrt{\cons+\kappa x^2}\bigr)\sqrt{\cons+\kappa x^2}+\kappa x^2\Bigl[2+\frac{\varepsilon+\sqrt{\cons+\kappa x^2}}{\sqrt{\cons+\kappa x^2}}\Bigr]\Bigr]\\
    &\textstyle=\frac{\kappa x}{(\varepsilon+\sqrt{\cons+\kappa x^2})^3(\cons+\kappa x^2)} \Bigl[-3\varepsilon\sqrt{\cons+\kappa x^2}-3(\cons+\kappa x^2)+2\kappa x^2+\frac{\kappa x^2(\varepsilon+\sqrt{\cons+\kappa x^2})}{\sqrt{\cons+\kappa x^2}}\Bigr]\\
    &=\textstyle\frac{\kappa x}{(\varepsilon+\sqrt{\cons+\kappa x^2})^3(\cons+\kappa x^2)} \Bigl[-3\varepsilon\sqrt{\cons+\kappa x^2}-3(\cons+\kappa x^2)+2\kappa x^2+\frac{\varepsilon\kappa x^2}{\sqrt{\cons+\kappa x^2}}+\kappa x^2\Bigr]\\
    &=\textstyle\frac{\kappa x}{(\varepsilon+\sqrt{\cons+\kappa x^2})^3(\cons+\kappa x^2)} \Bigl[-3\varepsilon\sqrt{\cons+\kappa x^2}-3\cons+\frac{\varepsilon\kappa x^2}{\sqrt{\cons+\kappa x^2}}\Bigr].
    \end{split}
\end{equation}}
\argument{the assumption that $\cons\geq 0$; the assumption that $\kappa>0$; the assumption that $\varepsilon>0$}{that for all $x\in \R\backslash\{0\}$ it holds that
\begin{equation}\llabel{eq4}
    \frac{\varepsilon\kappa x^2}{\sqrt{\cons+\kappa x^2}}-3\varepsilon \sqrt{\cons+\kappa x^2}=\frac{\varepsilon\kappa x^2-3\varepsilon(\cons+\kappa x^2)}{\sqrt{\cons+\kappa x^2}}=\frac{-2\varepsilon\kappa x^2-3\varepsilon\cons}{\sqrt{\cons+\kappa x^2}}<0.
\end{equation}}
\argument{\lref{eq4};the assumption that $\cons\geq 0$}{that for all $x\in (0,\infty)$ it holds that
\begin{equation}\llabel{eq4.5}
    -3\varepsilon\sqrt{\cons+\kappa x^2}-3\cons+\frac{\varepsilon\kappa x^2}{\sqrt{\cons+\kappa x^2}}<0.
\end{equation}}
\argument{\lref{eq4.5};\lref{def: adam};\lref{eq3};the assumption that $\cons\geq 0$; the assumption that $\varepsilon>0$; the assumption that $\kappa>0$}{that for all $x\in (0,\infty)$ it holds that
\begin{equation}\llabel{eq5}
    \adamvector''(x)<0<f(x).
\end{equation}}
\argument{\lref{eq5};}{\cref{item 5: half convex}\dott}
\end{aproof}
\subsection{Monotonicity of the Adam vector field for symmetrically distributed data}
In \cref{Lemma X} below we recall a well-known elementary fact on identically distributed random variables and in \cref{lem: tg1}, \cref{Lemma Y}, \cref{lem: sum of independent}, and \cref{cor: sum of independent} below we collect a few elementary properties for odd function evaluations, scalar multiplications, sums, and scaled averages of (independent) symmetrically distributed random variables. Only for completeness we include here detailed proofs for \cref{Lemma X}, \cref{lem: tg1}, \cref{Lemma Y}, \cref{lem: sum of independent}, and \cref{cor: sum of independent}.

\begin{athm}{lemma}{Lemma X}
     Let $(\Omega, \cF,\P)$ be a probability space, let $d \in 
     \N$, let $X \colon\Omega\to\R^d$ and $Y \colon\Omega\to \R^d$ be identically distributed random variables, and let $\varphi \colon \R^d \to \R^d$ be measurable. Then it holds that $\varphi( X )$ and $\varphi( Y )$ are identically distributed.
\end{athm}
\begin{aproof}
    \argument{the assumption that $X$ and $Y$ are identically distributed;}{that \llabel{arg0} for all $A\in \cB(\R^d)$ it holds that $ \P(X\in A)=\P(Y\in A)$\dott}
   \argument{\lref{arg0}}{for all $A \in \cB(\R^d)$ that
   \begin{equation}\llabel{arg0.1}
       \P( \varphi( X ) \in A ) = \P( X \in \varphi^{ -1 }( A ) ) = \P( Y \in \varphi^{ - 1 }( A ) ) = \P( \varphi(Y ) \in A ).
   \end{equation}}
   \argument{\lref{arg0.1};}{that $\varphi(X)$ and $\varphi(Y)$ are identically distributed\dott}
\end{aproof}
\begin{athm}{lemma}{lem: tg1}
    Let $( \Omega, \cF, \P )$ be a probability space, let $d \in \N$, let $X \colon \Omega \to \R^d$  be a symmetrically distributed random variable, let $f \colon \R^d\to \R^d$ be measurable, assume $\E\bigl[ \| X \| + \| f(X) \| \bigr] < \infty$ and $\E[ X ] = 0$, and assume for all $x \in \R^d$ that $f(x) = - f( - x )$. Then $\E[ f( X ) ] = 0$.
\end{athm}
\begin{aproof}
    \argument{the assumption that $X$ is symmetrically distributed;the fact that $\E\bigl[\|X\|\bigr]<\infty$; the assumption that $\E[X]=0$;}{that \llabel{arg1} $X$ and $- X$ are identically distributed\dott}
    \argument{\lref{arg1};the fact that $\E\bigl[\|f(X)\|\bigr]<\infty$; the assumption that for all $x\in \R^d$ it holds that $f(x)=-f(-x)$}{ that
\begin{equation}\llabel{eq1}
  2 \E[ f( X ) ] = \E[ f( X ) ] + \E[ f( X ) ] = \E[ f( X ) ] - \E[ f( - X ) ] = \E[ f(X) ] - \E[ f( X ) ] = 0.
\end{equation}}
\argument{\lref{eq1};}{that $\E[ f( X ) ] = 0$\dott}
\end{aproof}
\begin{athm}{lemma}{lem: monotonicity fix other}
    Let $ \cons \in [0,\infty) $, $ \kappa, \varepsilon \in (0,\infty) $, let $(\Omega,\cF,\P)$ be a probability space, let $X\colon \Omega\allowbreak\to \R$ be a symmetrically distributed random variable,  assume $\E\bigl[|X|\bigr]<\infty$ and $\E[X]=0$, and let $f\colon \R\to\R$ satisfy for all $\theta\in \R$ that
    \begin{equation}\llabel{def: f}
         f(\theta)=\E\biggl[\frac{X-\theta}{\varepsilon+\sqrt{\cons + \kappa |X-\theta|^2 } }\biggr].
    \end{equation}
    Then 
    \begin{enumerate}[label=(\roman*)]
        \item \label{item 1: fix other} it holds that $f(0)=0$ and
        \item  \label{item 2: fix other} it holds that f is strictly decreasing.
    \end{enumerate}
\end{athm}
\begin{aproof}
Throughout this proof for every $k\in \N$ let $Y_k\colon \Omega\to\R$ satisfy $Y_k=\max\{(-1)^kX,0\}$.
    \argument{\cref{Adam concave};}{that for all $\omega\in \Omega$, $\theta_1,\theta_2\in \R$ with $\theta_1<\theta_2$ it holds that
    \begin{equation}\llabel{eq1}
        \frac{X(\omega)-\theta_1}{\varepsilon+\sqrt{\cons + \kappa |X(\omega)-\theta_1|^2 } }>\frac{X(\omega)-\theta_2}{\varepsilon+\sqrt{\cons + \kappa |X(\omega)-\theta_2|^2 } }.
    \end{equation}}
    \argument{\lref{eq1};\lref{def: f}}{for all $\theta_1,\theta_2\in \R$ with $\theta_1<\theta_2$ that
    \begin{equation}\llabel{eq2}
        f(\theta_1)=\E\biggl[\frac{X-\theta_1}{\sqrt{\cons + \kappa |X-\theta_1|^2 } + \varepsilon}\biggr]>\E\biggl[\frac{X-\theta_2}{\sqrt{\cons + \kappa |X-\theta_2|^2 } + \varepsilon}\biggr]=f(\theta_2).
    \end{equation}}
    \argument{the assumption that $\E\bigl[|X|\bigr]<\infty$}{that
    \begin{equation}\llabel{arg1}
        \E\Biggl[|X|+ \biggl|\frac{X}{\varepsilon+\sqrt{\cons + \kappa X^2 }}\biggr|\Biggr]\leq  \E\biggl[ |X| + \frac{ \kappa^{ - 1/2 } \kappa^{ 1 / 2 } | X | }{ \varepsilon + \kappa^{ 1 / 2 } | X | } \biggr] 
\leq 
  \E\bigl[ | X | + \kappa^{ - 1/2 } \bigr]<\infty.
    \end{equation}}
    \argument{\lref{arg1};\cref{lem: tg1};the assumption that $X$ is symmetrically distributed; the assumpion that $\E[X]=0$;the fact that for all $x\in \R$ it holds that \begin{equation}
        \frac{x}{\varepsilon+\sqrt{\cons+\kappa x^2}}=-\biggl[\frac{-x}{\varepsilon+\sqrt{\cons+\kappa (-x)^2}}\biggr]
    \end{equation};}{that
    \begin{equation}\llabel{eq3}
        \begin{split}
           \E\biggl[\frac{X}{\varepsilon+\sqrt{\cons + \kappa X^2 } }\biggr]=
         0.
        \end{split}
    \end{equation}}
    \argument{\lref{eq3};\lref{def: f}}{that
    \begin{equation}\llabel{eq4}
        f(0)=\E\biggl[\frac{X}{\varepsilon+\sqrt{\cons + \kappa X^2 } }\biggr]=0.
    \end{equation}}
    \argument{\lref{eq4};\lref{eq2}}{\cref{item 1: fix other,item 2: fix other}\dott}
\end{aproof}
\begin{athm}{lemma}{Lemma Y}
    Let $( \Omega, \cF,\P )$ be a probability space, let $d \in \N$, let $X \colon\Omega \to \R^d$ and $Y \colon \Omega\to \R^d$ be independent symmetrically distributed random variables, and assume $\E\bigl[\|X\|+\|Y\| \bigr] < \infty$ and $\E[ X ] = \E[ Y ] = 0$. Then it holds that $X + Y$ and $X - Y$ are identically distributed.
\end{athm}
\begin{aproof}
    Throughout this proof for every $A\in \cB(\R^d)$ let $\phi_A\colon \R^d\to \R$ and $\psi_A\colon \R^d\to\R$ satisfy for all $x\in \R^d$ that
    \begin{equation}\llabel{def: phi1}
        \phi_A(x)=\E\bigl[\mathbbm 1_{\{x+Y\in A\}}\bigr]\qqandqq \psi_A(x)=\E\bigl[\mathbbm 1_{\{x-Y\in A\}}\bigr].
    \end{equation}
    \startnewargseq
     \argument{the assumption that $Y$ is symmetrically distributed;the fact that $\E[Y]=0$}{that \llabel{eq1} $Y$ and $-Y$ are identically distributed\dott}
    \argument{\lref{eq1};\cref{Lemma X}}{that for all $x\in \R^d$ it holds that \llabel{arg1} $x+Y$ and $x-Y$ are identically distributed\dott}
    \argument{\lref{arg1};}{for all $A\in \cB(\R^d)$, $x\in \R^d$ that \llabel{arg2} $\P( x + Y \in A ) = \P( x - Y \in A )$\dott} 
    \argument{\lref{arg2};\lref{def: phi1}}{for all $A\in \cB(\R^d)$, $x\in \R^d$ that
    \begin{equation}\llabel{eq3'}
    \begin{split}
        &\phi_A(x)=\psi_A(x).
        \end{split}
    \end{equation}}
    \argument{\lref{def: phi1};the assumption that $X$ and $Y$ are independent;the factorization lemma for conditional expectations (see, \eg, \cite[Proposition 3.15]{DeRoAr2024nonconvergence})}{that for all $A\in \cB(\R^d)$ it holds $\P$-a.s.\ that
    \begin{equation}\llabel{eq4'}
    \begin{split}
       \E\bigl[\mathbbm 1_{\{X+Y\in A\}}|X\bigr]=\phi_A(X)\qqandqq \E\bigl[\mathbbm 1_{\{X-Y\in A\}}|X\bigr]=\psi_A(X).
        \end{split}
        \end{equation}}
        \argument{\lref{eq4'};\lref{def: phi1};\lref{eq3'}}{for all $A\in \cB(\R^d)$ that
        \begin{equation}\llabel{evidence 2}
        \begin{split}
           & \P(X+Y\in A)=\E\bigl[\mathbbm 1_{\{X+Y\in A\}}\bigr]=\E\Bigl[\E\bigl[\mathbbm 1_{\{X+Y\in A\}}|X\bigr]\Bigr]=\E[\phi_A(X)]\\
           &=\E[\psi_A(X)]=\E\Bigl[\E\bigl[\mathbbm 1_{\{X-Y\in A\}}|X\bigr]\Bigr]=\E\bigl[\mathbbm 1_{\{X-Y\in A\}}\bigr]=\P(X-Y\in A).
            \end{split}
        \end{equation}}
        \argument{\lref{evidence 2};}{that $X+Y$ and $X-Y$ are identically distributed\dott}
\end{aproof}
\begin{athm}{lemma}{lem: sum of independent}
    Let $(\Omega,\cF,\P)$ be a probability space, let $d\in \N$, let $X\colon \Omega\to \R^d$ and $Y\colon \Omega\to\R^d$ be independent symmetrically distributed random variables, assume $\E\bigl[\|X\|+\|Y\|\bigr]<\infty$ and $\E[X]=\E[Y]=0$, and let $\alpha\in \R$. Then
    \begin{enumerate}[label=(\roman*)]
        \item \label{item 1: sum of independent} it holds that $\alpha X$ is symmetrically distributed and
        \item \label{item 2: sum of independent} it holds that $X+Y$ is symmetrically distributed.
    \end{enumerate}
\end{athm}
\begin{aproof}
    \argument{\cref{Lemma X} (applied with $\Omega\curvearrowleft\Omega$, $\cF\curvearrowleft\cF$, $\P\curvearrowleft\P$, $d\curvearrowleft d$, $X\curvearrowleft X$, $Y\curvearrowleft -X$, $\varphi\curvearrowleft (\R^d\ni x\mapsto \alpha x\in \R^d)$ in the notation of \cref{Lemma X})}
    {that \llabel{arg1} $\alpha X$ and $-\alpha X$ are identically distributed\dott}
    \argument{\lref{arg1};the fact that $\E\bigl[\|\alpha X\|\bigr]=|\alpha|\,\E\bigl[\|X\|\bigr]<\infty$; the fact that \begin{equation}
        \E[\alpha X]=\alpha \E[X]=0
    \end{equation}}
    {\cref{item 1: sum of independent}\dott}
    \startnewargseq
    \argument{\cref{Lemma Y} (applied with $\Omega\curvearrowleft\Omega$, $\cF\curvearrowleft\cF$, $\P\curvearrowleft\P$, $d\curvearrowleft d$, $X\curvearrowleft X$, $Y\curvearrowleft Y$ in the notation of \cref{Lemma Y})}{that $X+Y$ and $X-Y$ are identically distributed\dott}
      \argument{\cref{Lemma Y} (applied with $\Omega\curvearrowleft\Omega$, $\cF\curvearrowleft\cF$, $\P\curvearrowleft\P$, $d\curvearrowleft d$, $X\curvearrowleft -Y$, $Y\curvearrowleft X$ in the notation of \cref{Lemma Y})}{that \llabel{argg1} $-Y+X$ and $-Y-X$ are identically distributed\dott}
      \argument{\lref{argg1};the fact that $X+Y$ and $X-Y$ are identically distributed}{that \llabel{argg2} $X+Y$ and $-X-Y$ are identically distributed\dott}
      \argument{\lref{argg2};the fact that $\E\bigl[\|X+Y\|\bigr]\leq \E\bigl[\|X\|+\|Y\|\bigr]<\infty$; the fact that \begin{equation}
          \E[X+Y]=\E[X]+\E[Y]=0
      \end{equation}}{\cref{item 2: sum of independent}\dott}
\end{aproof}
\begin{athm}{cor}{cor: sum of independent}
     Let $(\Omega,\cF,\P)$ be a probability space, let $M,d\in \N$, let $X_m\colon \Omega\to\R^d$, $m\in \{1,2,\dots,M\}$, be independent symmetrically distributed random variables, and assume for all $m\in \{1,2,\dots,M\}$ that $\E\bigl[\|X_m\|\bigr]<\infty$ and $\E[X_m]=0$. Then $\frac{2}{M} \sum_{m=1}^MX_m$ is symmetrically distributed.
\end{athm}
\begin{aproof}
\argument{the assumption that $X_m$, $m\in \{1,2,\dots,M\}$, are independent; the assumption that for all $m\in \N$ it holds that $\E\bigl[\|X_m\|\bigr]<\infty$; the assumption that for all $m\in \{1,2,\dots,M\}$ it holds that $\E[X_m]=0$; \cref{item 2: sum of independent}\ in \cref{lem: sum of independent};induction}{that \llabel{arg1} $\sum_{m=1}^MX_m$ is symmetrically distributed\dott}
    \argument{\lref{arg1};\cref{item 1: sum of independent} in \cref{lem: sum of independent} (applied with $\Omega\curvearrowleft \Omega$, $\cF\curvearrowleft\cF$, $\P\curvearrowleft\P$, $d\curvearrowleft d$, $X\curvearrowleft \sum_{m=1}^MX_m$, $Y\curvearrowleft 0$, $\alpha\curvearrowleft\frac 2M$ in the notation of \cref{lem: sum of independent})}{that \llabel{arg2} $\frac 2M\sum_{m=1}^MX_m$ is symmetrically distributed\dott}
\end{aproof}
\begin{athm}{prop}{lem: monotonicity Adam vector field at origin pre}
     Let $(\Omega,\cF,\P)$ be a probability space, let $\fc\in [0,\infty)$, let $X_{n,m}\colon \Omega\allowbreak\to [-\fc,\fc]$, $(n,m) \in (\N_0)^2$, be \iid\ random variables, assume that $X_{1,1}$ is symmetrically distributed, assume $\E[X_{1,1}]=0$, let $\alpha\in (0,1)$, $\beta\in (\alpha^2,1)$, $\varepsilon\in (0,\infty)$, $M\in \N$, let $f\colon\R\to\R$ satisfy for all $\theta\in \R$ that
     \begin{equation}\llabel{eq1}
    \begin{split}
        f(\theta)&=(1-\alpha)\,\E\Biggl[\frac{\textstyle\sum_{n= 0}^\infty \alpha^nM^{-1}\textstyle \sum_{m=1}^M2(X_{n,m}-\theta)}{\textstyle\varepsilon+\sqrt{(1-\beta)\sum_{n=0}^\infty\beta^n
          |M^{-1}\sum_{m=1}^M2(X_{n,m}-\theta)|^2}}\Biggr],
          \end{split}
    \end{equation}
    and let $K\in (0,\infty)$.
    Then
    \begin{enumerate}[label=(\roman*)]
        \item \label{item 1: monotinicity symmetric pre} it holds that $f\in C^1(\R,\R)$,
        \item \label{item 3: monotonicity symmetric pre} it holds for all  $\theta\in (0,\infty)$ that $ f'( 0 ) < 0 = f( 0 ) < - f( \theta ) = f( - \theta )$, and
        \item \label{item 4: monotonicity symmetric pre} there exists $\fC\in (0,\infty)$ such that for all $\theta\in [-K,K]$ it holds that $\theta f(\theta)\leq -\fC\theta^2$.
    \end{enumerate}
\end{athm}
\begin{aproof}
    Throughout this proof assume without loss of generality that $\operatorname{Var}(X_{1,1})\allowbreak \neq0$ (otherwise note that \lref{eq1}, the fact that for all $\theta \in [-K,K]$ it holds that $\theta[ - 2 \theta ( \varepsilon + 2 | \theta | )^{ - 1 } ] 
= - 2 \theta^2 ( \varepsilon + 2 | \theta | )^{ - 1 } 
\leq - 2 ( \varepsilon + 2 K )^{ - 1 } \theta^2 
$, the fact that $\E\bigl[ | X_{ 1, 1 } |^2 \bigr] = \E\bigl[ | X_{ 1, 1 } - \E[ X_{ 1, 1 }] |^2 \bigr] = \operatorname{Var}( X_{ 1, 1 } )$, and \cref{Adam concave} (applied with $c \curvearrowleft 0$, $\kappa \curvearrowleft 4$, $\varepsilon \curvearrowleft \varepsilon$ in the notation of \cref{Adam concave}) establish items \ref{item 1: monotinicity symmetric pre}, \ref{item 3: monotonicity symmetric pre}, and \ref{item 4: monotonicity symmetric pre}) and for every $n\in \N_0$ let $\bbX_n\colon \Omega\to\R$ satisfy
\begin{equation}\llabel{def: bbX}
    \bbX_n=\frac 2M\sum_{m=1}^MX_{n,m}.
\end{equation}
\argument{\lref{eq1};\lref{def: bbX}}{that for all $\theta\in \R$ it holds that
\begin{equation}\llabel{eq2}
      f(\theta)=(1-\alpha)\,\E\Biggl[\frac{\textstyle\sum_{n= 0}^\infty\alpha^n(\bbX_n-2\theta)}{\textstyle\varepsilon+\sqrt{(1-\beta)\sum_{n=0}^\infty\beta^n|\bbX_n-2\theta|^2}}\Biggr].
\end{equation}}
\argument{the assumption that $\E[ X_{ 1, 1 } ] = 0$}{that for all $\chi \in (0,2)$, $\theta \in \R$ it holds that
\begin{equation}\llabel{eq2.1}
\begin{split}
   & \E\bigl[ | \theta - X_{ 1, 1 } |^2 \bigr] - \frac{ \chi \theta^2 }{ 2 }
= \E\bigl[ \theta^2 - 2 \theta X_{ 1, 1 } + ( X_{ 1, 1 } )^2 \bigr] - \frac{ \chi \theta^2 }{ 2 }\\
&
= \theta^2 + \operatorname{Var}( X_{ 1, 1 } ) - \frac{ \chi \theta^2 }{ 2 }
= \Bigl( 1 - \frac{ \chi }{ 2 } \Bigr) \theta^2 + \operatorname{Var}( X_{ 1, 1 } ).
\end{split}
\end{equation}}
\argument{\lref{eq2.1};}{for all $\chi \in (0,2)$ that
\begin{equation}\llabel{eq2.2}
\R \ni \theta \mapsto \E\bigl[ | \theta - X_{ 1,1 } |^2 \bigr] \in \R 
\end{equation}
is $\chi$-\stronglyconvex \dott}
    \argument{\lref{eq2};\lref{eq2.2};\cref{item 2: lem negative derivative} in \cref{lem: negative derivative} (applied with $d \curvearrowleft 1$, $\dimX \curvearrowleft 1$, $\smalll \curvearrowleft (\R \times \R \ni (\theta,x) \mapsto | \theta - x |^2 \in \R)$ in the notation of \cref{lem: negative derivative})}{\cref{item 1: monotinicity symmetric pre}\dott}
    \startnewargseq
    \argument{\cref{item 3: lem negative derivative} in \cref{lem: negative derivative};\lref{eq2}; \lref{eq2.2};the assumption that $\operatorname{Var}(X_{1,1})\neq 0$}
    {that for all $\theta\in \R$ it holds that
    \begin{equation}\llabel{eq3'}
    \begin{split}
       \textstyle f'(\theta)
       &=\textstyle-\E\biggl[\frac{2}{\varepsilon+\sqrt{(1-\beta)\sum_{n=0}^\infty\beta^n|\bbX_n-2\theta|^2}}\biggr]\\
        &\textstyle+\E\biggl[\frac{2\bigl[\mathbbm 1_{\{0\}}(\sum_{n=0}^\infty |\bbX_n-2\theta|)+(1-\alpha)\sum_{n=0}^\infty\alpha^n(\bbX_{n}-2\theta)\bigr]\bigl[(1-\beta)\sum_{n=0}^\infty\beta^n (\bbX_{n}-2\theta)\bigr]}{\bigl(\varepsilon+\sqrt{(1-\beta)\sum_{n=0}^\infty\beta^n|\bbX_{n}-2\theta|^2}\bigr)^2\sqrt{\mathbbm 1_{\{0\}}(\sum_{n=0}^\infty |\bbX_n-2\theta|)+(1-\beta)\sum_{n=0}^\infty\beta^n|\bbX_{n}-2\theta|^2}}\biggr].
        \end{split}
    \end{equation}}
    \argument{the assumption that $\operatorname{Var}(X_{1,1})\neq 0$;the assumption that $X_{n,m}$, $(n,m)\in (\N_0)^2$, are \iid;\lref{def: bbX}}{that for all $n\in \N_0$ it holds that
    \begin{equation}\llabel{eq3.1}
    \begin{split}
\operatorname{Var}(\bbX_n)&=\operatorname{Var}\biggl( \displaystyle\frac{ 2 }{ M } \textstyle\sum\limits_{ m=1 }^{M  } X_{ n,m } \biggr) = \displaystyle\frac{ 4 }{ M^2 } \operatorname{Var}\biggl(\textstyle\sum\limits_{ m=1}^{ M} X_{ n, m } \biggr) \\
&\textstyle
= \displaystyle\frac{ 4 }{ M^2 } \textstyle\sum\limits_{ m=1}^{ M } \operatorname{Var}( X_{ n, m } ) = \displaystyle\frac{4\operatorname{Var}(X_{1,1})}{M}\neq 0.  
\end{split}
    \end{equation}}
    \argument{\lref{eq3.1};}{that for all $n\in \N_0$, $\theta\in \R$ it holds that
    \begin{equation}\llabel{eq3.1.1'}
        \P(\bbX_n=2\theta)< 1.
    \end{equation}}
    \argument{\lref{eq3.1.1'};the fact that $\bbX_n$, $n\in \N_0$, are \iid}{that for all $n\in \N_0$, $\theta\in \R$ it holds that
    \begin{equation}\llabel{eq3.1.1}
        \P(\bbX_n=2\theta)=\P( \bbX_0 = 2 \theta ) < 1.
    \end{equation}}
    \argument{\lref{eq3.1.1};the fact that $\bbX_n$, $n\in \N_0$, are \iid}{that for all $\theta\in \R$ it holds that
    \begin{equation}\llabel{eq3.2}
    \begin{split}
        &\P\bigl(\textstyle\mathbbm 1_{\{0\}}(\sum_{n=0}^\infty |\bbX_n-2\theta|)=0\bigr)=\P\bigl( \sum_{ n = 0 }^{ \infty } | \bbX_n - 2 \theta | \neq 0 \bigr)=
\P\bigl( \Omega \backslash \{ \sum_{ n = 0 }^{ \infty } | \bbX_n - 2 \theta | = 0 \} \bigr)\\
        &=\textstyle\P( \Omega \backslash \{ \forall \, n \in \N_0 \colon \bbX_n = 2 \theta \} )=1-\P(\forall\, n\in \N_0\colon \bbX_n=2\theta)=1-\bigl[\prod_{n=0}^\infty \P(\bbX_n=2\theta)\bigr]\\
        &\textstyle= 1 - \lim_{ n \to \infty } [ \P( \bbX_0 = 2 \theta ) ]^n=1.
        \end{split}
    \end{equation}}
    \argument{\lref{eq3.2};\lref{eq3'}}{for all $\theta\in \R$ that
    \begin{equation}\llabel{eq3}
    \begin{split}
       &\textstyle f'(\theta)
       =\textstyle-\E\biggl[\frac{2}{\varepsilon+\sqrt{(1-\beta)\sum_{n=0}^\infty\beta^n|\bbX_n-2\theta|^2}}\biggr]\\
        &\textstyle+\E\biggl[\frac{2\bigl[(1-\alpha)\sum_{n=0}^\infty\alpha^n(\bbX_{n}-2\theta)\bigr]\bigl[(1-\beta)\sum_{n=0}^\infty\beta^n (\bbX_{n}-2\theta)\bigr]}{\bigl(\varepsilon+\sqrt{(1-\beta)\sum_{n=0}^\infty\beta^n|\bbX_{n}-2\theta|^2}\bigr)^2\sqrt{\mathbbm 1_{\{0\}}(\sum_{n=0}^\infty |\bbX_n-2\theta|)+(1-\beta)\sum_{n=0}^\infty\beta^n|\bbX_{n}-2\theta|^2}}\biggr].
        \end{split}
\end{equation}}
    \argument{\lref{eq3};}{that 
    \begin{align}
      &  f'(0)
        =-\E\Biggl[\frac{2}{\varepsilon+\sqrt{(1-\beta)\sum_{n=0}^\infty\beta^n|\bbX_n|^2}}\Biggr]\llabel{eq4}\\
        &+\E\Biggl[\frac{2\bigl[(1-\alpha)\sum_{n=0}^\infty\alpha^n\bbX_n\bigr]\bigl[(1-\beta)\sum_{n=0}^\infty\beta^n \bbX_n\bigr]}{\bigl(\varepsilon+\sqrt{(1-\beta)\sum_{n=0}^\infty\beta^n|\bbX_n|^2}\,\bigr)^{\! 2}\sqrt{\mathbbm 1_{\{0\}}(\sum_{n=0}^\infty |\bbX_n|)+(1-\beta)\sum_{n=0}^\infty\beta^n|\bbX_n|^2}}\Biggr]\notag.
    \end{align}}
    \argument{\lref{eq4};}{that
    \begin{align}
            &f'(0)
            =-\E\Biggl[\frac{2}{\varepsilon+\sqrt{(1-\beta)\sum_{n=0}^\infty\beta^n|\bbX_n|^2}}\Biggr]\notag\\
        &\textstyle+(1-\alpha)(1-\beta)\sum\limits_{n=0}^\infty\E\Biggl[\textstyle\frac{2(\alpha\beta)^n|\bbX_n|^2}{\bigl(\varepsilon+\sqrt{(1-\beta)\sum_{k=0}^\infty\beta^k|\bbX_k|^2}\bigr)^{\!2}\sqrt{\mathbbm 1_{\{0\}}(\sum_{k=0}^\infty |\bbX_k|)+(1-\beta)\sum_{k=0}^\infty\beta^k|\bbX_k|^2}}\Biggr]\llabel{eq5}\\
        &\textstyle+(1-\alpha)(1-\beta)\textstyle\sum\limits_{n,m\in \N_0,\,n\neq m}\textstyle\E\biggl[\textstyle\frac{2\alpha^n\beta^m\bbX_n\bbX_m}{\bigl(\varepsilon+\sqrt{(1-\beta)\sum_{k=0}^\infty\beta^k|\bbX_k|^2}\bigr)^{\!2}\sqrt{\mathbbm 1_{\{0\}}(\sum_{k=0}^\infty |\bbX_k|)+(1-\beta)\sum_{k=0}^\infty\beta^k|\bbX_k|^2}}\biggr]\notag.
    \end{align}}
    \argument{\lref{def: bbX};the assumption that $X_{n,m}$, $(n,m)\in (\N_0)^2$, are \iid; the assumption that $X_{1,1}$ is symmetrically distributed;the assumption that $\E[X_{1,1}]=0$;\cref{cor: sum of independent}}{that for all $n\in \N_0$ it holds that $\bbX_n$ is \llabel{arg1} symmetrically distributed\dott}
    \argument{\lref{arg1};\cref{lem: tg1}}{that for all $n\in \N_0$, $\kappa\in (0,\infty)$, $x\in \R$, $y\in [0,\infty)$ it holds that
    \begin{equation}\llabel{eq6}
\E\Biggl[\frac{x\bbX_n}{\bigl(\varepsilon+\sqrt{(1-\beta)y+\kappa|\bbX_n|^2}\,\bigr)^{\! 2}\sqrt{(1-\beta)y+\mathbbm 1_{\{0\}}(y)+\mathbbm 1_{\{0\}}(\bbX_n)+\kappa|\bbX_n|^2}}\Biggr]=0.
    \end{equation}}
    \argument{\lref{eq6};}{for all $n, m \in \N_0$, $x \in \R$, $y \in [0,\infty)$ with $n \neq m$ that
    \begin{equation}\llabel{eqq6}
    \begin{split}
      \textstyle \E\biggl[ 
    \frac{ 
      \alpha^n \beta^m \bbX_n x 
    }{ 
      \bigl( 
        \varepsilon + \sqrt{ ( 1 - \beta ) y + ( 1 - \beta ) \beta^n | \bbX_n |^2 }  
      \bigr)^{\!2}
      \sqrt{
        \mathbbm 1_{ \{0\} }( y ) + \mathbbm 1_{ \{0\} }( \bbX_n ) + ( 1 - \beta ) y + ( 1 - \beta ) \beta^n | \bbX_n |^2
      }
    }\biggr]=0.
    \end{split}
    \end{equation}}
    \argument{\unskip, \eg, \cite[Theorem 2.26]{klenkeprobability};the fact that $\bbX_n$, $n\in \N_0$, are independent}{that for all $n\in \N_0$ it holds that \begin{equation}\llabel{eq6'} 
       \bbX_n\qqandqq (\bbX_k)_{ k\in \N_0\backslash\{n\}}
    \end{equation} 
    are independent\dott}
    \argument{\lref{eq6'};}{that for all $n,m\in \N_0$ with $n\neq m$ it holds that
    \begin{equation}\llabel{eq6''}
        \bbX_n\qqandqq \bigl(\bbX_m,\textstyle\sum_{k\in \N_0\backslash\{n\}}\beta^k|\bbX_k|^2\bigr)
    \end{equation}
    are independent\dott}
    \argument{\lref{eqq6};\lref{eq6''};the factorization lemma for conditional expectations (see, \eg, \cite[Proposition 3.15]{DeRoAr2024nonconvergence})}{that for all $n,m\in \N_0$ with $n\neq m$ it holds $\P$-a.s.\ that
    \begin{samepage}
    \begin{multline}\llabel{eq6.1}
        \textstyle\E\biggl[\frac{\alpha^n\beta^m\bbX_n\bbX_m}{\bigl(\varepsilon+\sqrt{(1-\beta)(\sum_{k\in \N_0\backslash\{n\}}\beta^k|\bbX_k|^2)+(1-\beta)\beta^n|\bbX_n|^2}\bigr)^{\!2}\sqrt{\mathbbm 1_{\{0\}}(\sum_{k\in \N_0\backslash\{n\}} \beta^k|\bbX_k|^2)+\mathbbm 1_{\{0\}}(\bbX_n)+(1-\beta)\sum_{k=0}^\infty\beta^k|\bbX_k|^2}}\\
      \Big| \bigl(\bbX_m,\textstyle\sum_{k\in \N_0\backslash\{n\}}\beta^k|\bbX_k|^2\bigr)\biggr]=0.
    \end{multline}
    \end{samepage}}
    \argument{\lref{eq6.1};}{that for all $n,m\in \N_0$ with $n\neq m$ it holds that
    \begin{equation}\llabel{eq7.1}
\textstyle\E\biggl[\frac{\alpha^n\beta^m\bbX_n\bbX_m}{\bigl(\varepsilon+\sqrt{(1-\beta)\sum_{k=0}^\infty\beta^k|\bbX_k|^2}\bigr)^{\!2}\sqrt{\mathbbm 1_{\{0\}}(\sum_{k\in \N_0\backslash\{n\}} \beta^k|\bbX_k|^2)+\mathbbm 1_{\{0\}}(\bbX_n)+(1-\beta)\sum_{k=0}^\infty\beta^k|\bbX_k|^2}}\biggr]=0.
    \end{equation}}
    \argument{\lref{eq3.1.1};the fact that $\bbX_n$, $n \in \N_0$, are \iid}{that for all $n\in \N_0$ it holds that 
     \begin{equation}\llabel{eq7.1'}
    \begin{split}
       &\textstyle \P\bigl(\bigr\{\mathbbm 1_{\{0\}}(\sum_{k\in \N_0\backslash\{n\}} \beta^k|\bbX_k|^2)+\mathbbm 1_{\{0\}}(\bbX_n)=\mathbbm 1_{\{0\}}(\sum_{k=0}^\infty |\bbX_k|)\bigr\}\cup \{\bbX_n= 0\}\bigr)\\
       &\textstyle=1- \P\bigl(\bigr\{\mathbbm 1_{\{0\}}(\sum_{k\in \N_0\backslash\{n\}} \beta^k|\bbX_k|^2)+\mathbbm 1_{\{0\}}(\bbX_n)\neq \mathbbm 1_{\{0\}}(\sum_{k=0}^\infty |\bbX_k|)\bigr\}\cap \{\bbX_n\neq 0\}\bigr)
       \\
       &\textstyle =1- \P\bigl(\bigr\{\mathbbm 1_{\{0\}}(\sum_{k\in \N_0\backslash\{n\}} \beta^k|\bbX_k|^2)\neq \mathbbm 1_{\{0\}}(\sum_{k=0}^\infty |\bbX_k|)\bigr\}\cap \{\bbX_n\neq 0\}\bigr)\\
       &=\textstyle1- \P\bigl(\bigr\{\mathbbm 1_{\{0\}}(\sum_{k\in \N_0\backslash\{n\}} \beta^k|\bbX_k|^2)\neq 0\bigr\}\cap \{\bbX_n\neq 0\}\bigr)\\
       &=\textstyle1- \P\bigl(\bigr\{\sum_{k\in \N_0\backslash\{n\}} \beta^k|\bbX_k|^2= 0\bigr\}\cap \{\bbX_n\neq 0\}\bigr)\\
       &=\textstyle 1- \P\bigl(\bigr\{\forall\, k\in \N_0\backslash\{n\}\colon \bbX_k=0\bigr\}\cap \{\bbX_n\neq 0\}\bigr)\geq 1 - \P\bigl( \forall \, k \in \N_0\backslash\{n\} \colon \bbX_k = 0 \bigr) \\
       &\textstyle = 1 - \bigl[ \prod_{ k\in \N_0\backslash\{n\}} \P( \bbX_k=0) \bigr] = 1 - \lim_{ k\to\infty } [ \P( \bbX_0=0 ) ]^{k-1} = 1 - 0 = 1.
       \end{split}
    \end{equation}}
    \argument{\lref{eq7.1'}}{that for all $n\in \N_0$ it holds that
    \begin{align}
       &\textstyle \P\biggl(\frac{\bbX_n}{\sqrt{\mathbbm 1_{\{0\}}(\sum_{k\in \N_0\backslash\{n\}} \beta^k|\bbX_k|^2)+\mathbbm 1_{\{0\}}(\bbX_n)+(1-\beta)\sum_{k=0}^\infty\beta^k|\bbX_k|^2}}=\frac{\bbX_n}{\sqrt{\mathbbm 1_{\{0\}}(\sum_{k=0}^\infty |\bbX_k|)+(1-\beta)\sum_{k=0}^\infty\beta^k|\bbX_k|^2}}\biggr)\notag\\
       &=1\llabel{eq7.2}.
    \end{align}}
    \argument{\lref{eq7.1};\lref{eq7.2};}{that for all $n,m\in \N_0$ with $n\neq m$ it holds that
    \begin{equation}\llabel{eq7}
        \textstyle\E\biggl[\textstyle\frac{\alpha^n\beta^m\bbX_n\bbX_m}{\bigl(\varepsilon+\sqrt{(1-\beta)\sum_{k=0}^\infty\beta^k|\bbX_k|^2}\bigr)^2\sqrt{\mathbbm 1_{\{0\}}(\sum_{k=0}^\infty |\bbX_k|)+(1-\beta)\sum_{k=0}^\infty\beta^k|\bbX_k|^2}}\biggr]=0.
    \end{equation}}
    \argument{\lref{eq7};\lref{eq5}}{that
    \begin{equation}\llabel{eq8}
        \begin{split}
            &f'(0)=-\E\Biggl[\frac{2}{\varepsilon+\sqrt{(1-\beta)\sum_{n=0}^\infty\beta^n|\bbX_n|^2}}\Biggr]\\
        &+(1-\alpha)(1-\beta)\sum_{n=0}^\infty\E\Biggl[\textstyle\frac{2(\alpha\beta)^n|\bbX_n|^2}{\bigl(\varepsilon+\sqrt{(1-\beta)\sum_{k=0}^\infty\beta^k|\bbX_k|^2}\bigr)^2\sqrt{\mathbbm 1_{\{0\}}(\sum_{k=0}^\infty |\bbX_k|)+(1-\beta)\sum_{k=0}^\infty\beta^k|\bbX_k|^2}}\Biggr].
        \end{split}
    \end{equation}}
    \argument{\lref{eq8};the fact that $0<\alpha<1$; the fact that $0<\beta<1$; the fact that $\varepsilon>0$}{that
    \begin{equation}\llabel{eq9}
        \begin{split}
            f'(0)&= -\E\Biggl[\frac{2}{\varepsilon+\sqrt{(1-\beta)\sum_{n=0}^\infty\beta^n|\bbX_n|^2}}\Biggr]\\
        &+(1-\alpha)\sum_{n=0}^\infty\alpha^n\,\E\Biggl[\textstyle\frac{2(1-\beta)\beta^n|\bbX_n|^2}{\bigl(\varepsilon+\sqrt{(1-\beta)\sum_{k=0}^\infty\beta^k|\bbX_k|^2}\bigr)^2\sqrt{\mathbbm 1_{\{0\}}(\sum_{k=0}^\infty |\bbX_k|)+(1-\beta)\sum_{k=0}^\infty\beta^k|\bbX_k|^2}}\Biggr]\\
        &<-\E\Biggl[\frac{2}{\varepsilon+\sqrt{(1-\beta)\sum_{n=0}^\infty\beta^n|\bbX_n|^2}}\Biggr]\\
        &+(1-\alpha)\sum_{n=0}^\infty\alpha^n\,\E\Biggl[\frac{2}{\varepsilon+\sqrt{(1-\beta)\sum_{k=0}^\infty\beta^k|\bbX_k|^2}}\Biggr]
        \\
        &=-\E\Biggl[\frac{2}{\varepsilon+\sqrt{(1-\beta)\sum_{n=0}^\infty\beta^n|\bbX_n|^2}}\Biggr]+\E\Biggl[\frac{2}{\varepsilon+\sqrt{(1-\beta)\sum_{k=0}^\infty\beta^k|\bbX_k|^2}}\Biggr]=0.
        \end{split}
    \end{equation}}
    \argument{the fact that for all $n\in \N_0$ it holds that $\E[\bbX_n]=0$; the fact that for all $n\in \N_0$ it holds that $\bbX_n$ is symmetrically distributed;\cref{lem: monotonicity fix other}}{that for all  $n\in \N_0$, $c\in [0,\infty)$, $\kappa\in (0,\infty)$, $\theta\in (-\infty,0)$ it holds that
    \begin{equation}\llabel{need to prove'}
         \E\Biggl[\frac{\bbX_{n}-\theta}{\varepsilon+\sqrt{c+\kappa |\bbX_n-\theta|^2}}\Biggr]>0.
    \end{equation}}
    \argument{\lref{need to prove'}}{for all $n \in \N_0$, $y \in [0,\infty)$ that
\begin{equation}\llabel{need to prove}
  \E\biggl[ 
    \frac{ \bbX_n - \theta 
    }{ 
      \varepsilon
      + 
      \sqrt{ 
        ( 1 - \beta ) y + ( 1 - \beta ) \beta^n | \bbX_n - \theta |^2 
      } 
    } 
  \biggr]
  > 0.
\end{equation} }
    \argument{\lref{eq6'};}{that for all $n\in \N_0$, $\theta\in \R$ it holds that \llabel{eq10''} $\bbX_n$ and $\sum_{k\in \N_0\backslash\{n\}}\beta^k|\bbX_k\allowbreak-\theta|^2$ are independent\dott}
    \argument{\lref{need to prove};\lref{eq10''}; the factorization lemma for conditional expectations (see, \eg, \cite[Proposition 3.15]{DeRoAr2024nonconvergence})}{that for all $n\in \N_0$, $\theta\in (-\infty,0)$ it holds $\P$-a.s.\ that 
    \begin{equation}\llabel{eq10'}
    \begin{split}
       & \E\Biggl[\frac{\textstyle\bbX_n-\theta}{\textstyle\varepsilon+\sqrt{(1-\beta)\bigl[\sum_{k\in \N_0\backslash\{n\}}\beta^k|\bbX_k-\theta|^2\bigr]+(1-\beta)\beta^n|\bbX_n-\theta|^2}}\biggl|\textstyle\sum\limits_{k\in \N_0\backslash\{n\}}\beta^k|\bbX_k-\theta|^2 \Biggr]\\
       &>0.
        \end{split}
    \end{equation}}
    \argument{\lref{eq10'}}{that for all $n\in \N_0$, $\theta\in (-\infty,0)$ it holds that
    \begin{equation}\llabel{eq10}
    \begin{split}
       & \E\Biggl[\frac{\textstyle(1-\alpha)\alpha^n(\bbX_n-\theta)}{\textstyle\varepsilon+\sqrt{(1-\beta)\sum_{k=0}^\infty\beta^k|\bbX_k-\theta|^2}}\Biggr]\\
       &=\E\Biggl[\frac{\textstyle(1-\alpha)\alpha^n(\bbX_n-\theta)}{\textstyle\varepsilon+\sqrt{(1-\beta)\bigl[\sum_{k\in \N_0\backslash\{n\}}\beta^k|\bbX_k-\theta|^2\bigr]+(1-\beta)\beta^n|\bbX_n-\theta|^2}}\Biggr]>0.
        \end{split}
    \end{equation}}
    \argument{\lref{eq10};\lref{eq2}}{that for all $\theta\in (-\infty,0)$ it holds that
    \begin{equation}\llabel{evd1}
    \begin{split}
         f(\theta)&=(1-\alpha)\,\E\Biggl[\frac{\textstyle\sum_{n= 0}^\infty\alpha^n(\bbX_n-\theta)}{\textstyle\varepsilon+\sqrt{(1-\beta)\sum_{n=0}^\infty\beta^n|\bbX_n-\theta|^2}}\Biggr]\\
         &=\sum_{n=0}^\infty\E\Biggl[\frac{\textstyle(1-\alpha)\alpha^n(\bbX_n-\theta)}{\textstyle\varepsilon+\sqrt{(1-\beta)\sum_{k=0}^\infty\beta^k|\bbX_k-\theta|^2}}\Biggr]> 0.
         \end{split}
    \end{equation}}
    \argument{the fact that for all $n\in \N_0$ it holds that $\bbX_n$ is symmetrically distributed;the fact that for all $n\in \N_0$ it holds that $\E[\bbX_n]=0$;}{that \llabel{arg11} for all $n\in \N_0$ it holds that $\bbX_n$ and $-\bbX_n$ are identically distributed\dott}
    \argument{\lref{arg11}; the fact that $\bbX_n$, $n\in \N_0$, are \iid}{that \llabel{arg12} $(\bbX_n)_{n\in \N_0}$ and $(-\bbX_n)_{n\in \N_0}$ are identically distributed\dott}
    \argument{\lref{eq2}; \lref{arg12}}{that for all $\theta \in \R$ it holds that
    \begin{equation}\llabel{evd2}
    \begin{split}
         f(\theta)+f(-\theta)&=(1-\alpha)\,\E\Biggl[\frac{\textstyle\sum_{n= 0}^\infty\alpha^n(\bbX_n-2\theta)}{\textstyle\varepsilon+\sqrt{(1-\beta)\sum_{n=0}^\infty\beta^n|\bbX_n-2\theta|^2}}\Biggr]\\
         &+(1-\alpha)\,\E\Biggl[\frac{\textstyle\sum_{n= 0}^\infty\alpha^n(\bbX_n+2\theta)}{\textstyle\varepsilon+\sqrt{(1-\beta)\sum_{n=0}^\infty\beta^n|\bbX_n+2\theta|^2}}\Biggr]\\
         &=(1-\alpha)\,\E\Biggl[\frac{\textstyle\sum_{n= 0}^\infty\alpha^n(\bbX_n-2\theta)}{\textstyle\varepsilon+\sqrt{(1-\beta)\sum_{n=0}^\infty\beta^n|\bbX_n-2\theta|^2}}\Biggr]\\
         &+(1-\alpha)
         \,\E\Biggl[\frac{\textstyle\sum_{n= 0}^\infty\alpha^n(-\bbX_n+2\theta)}{\textstyle\varepsilon+\sqrt{(1-\beta)\sum_{n=0}^\infty\beta^n|-\bbX_n+2\theta|^2}}\Biggr]=0.
    \end{split}
    \end{equation}}
    \argument{\lref{evd2};\lref{evd1};}{for all $\theta\in (0,\infty)$ that
    \begin{equation}\llabel{eqdb1}
        f(\theta)=-f(-\theta)<0.
    \end{equation}}
    \argument{\lref{eqdb1};\lref{evd1};\cref{item 1: monotinicity symmetric pre}}{that
    \begin{equation}\llabel{eqdb}
        f(0)=0.
    \end{equation}}
    \argument{\lref{eq9};\lref{evd2};\lref{eqdb1};\lref{eqdb}}{\cref{item 3: monotonicity symmetric pre}\dott}
    \startnewargseq
    \argument{\lref{eq9};\cref{item 1: monotinicity symmetric pre}}{that there exists $\delta \in (0,K)$ which satisfies
\begin{equation}\llabel{eq9.1}
  \textstyle-\bigl[\inf_{ \theta \in [-\delta,\delta] } | f'( \theta ) |\bigr] =\sup_{ \theta \in [ - \delta, \delta ] } f'( \theta ) < 0.
\end{equation}}
\startnewargseq
\argument{\lref{eq9.1};the fundamental theorem of calculus}{that for all $\theta \in [ - \delta, \delta ]$ it holds that
\begin{equation}\llabel{eq*}
\begin{split}
    & \textstyle | f( \theta ) - f(0) | 
=
  \bigl| \int_0^{ \theta } f'(s) \, \d s \bigr|=\big| \int_{ \min\{0,\theta\} }^{ \max\{ 0, \theta \} } f'(s)\, \d s \big| 
=
  \int_{ \min\{0, \theta \} }^{ \max\{0, \theta\} }
  | f'(s) | \,\d s 
\\
&\textstyle\geq
  | \max\{0,\theta\} - \min\{ 0, \theta \} |
  \bigl[
    \inf_{ s \in [ - \delta, \delta ] }
    | f'(s) |
  \bigr]
=\textstyle
  | \theta |
  \bigl[
    \inf_{ s \in [ - \delta, \delta ] }
    | f'(s) |
  \bigr].
\end{split}
\end{equation}}
\argument{\lref{evd1};\lref{eqdb1};\lref{eqdb}}{that for all $\theta \in \R$ it holds that \llabel{arg13}
$\theta ( f( \theta ) - f( 0 ) ) = \theta f( \theta ) \leq 0$\dott}
\argument{\lref{arg13};\lref{eq*}}{that for all $\theta \in [ - \delta, \delta ]$ it holds that
\begin{equation}\llabel{eq**}
\textstyle
  \theta ( f( \theta ) - f( 0 ) )
=
  - | \theta | | f( \theta ) - f( 0 ) |
\leq 
  - | \theta |^2 \bigl[ \inf_{ s \in[-\delta,\delta] } | f'(s) | \bigr]= -\bigl[ \inf_{ s \in[-\delta,\delta] } | f'(s) | \bigr]\theta^2.  
\end{equation}}
    \argument{\cref{item 1: monotinicity symmetric pre,item 3: monotonicity symmetric pre};}{that
    \begin{equation}\llabel{def: bfC}
      \textstyle - \bigl[ \inf_{ \theta \in [\delta,K] } | f( \theta ) | \bigr]=\sup_{\theta\in[\delta,K]} f(\theta)<0.
    \end{equation}}
    \argument{\lref{def: bfC};}{that for all $\theta \in [\delta,K]$ it holds that
    \begin{equation}\llabel{eq15}
      \textstyle \theta f(\theta) \leq \bigl[ \sup_{s\in [\delta,K] } f(s) \bigr] \theta\leq -\bigl[\inf_{s\in[\delta,K]} |f(s)|\bigr]\theta\leq -K^{-1}\bigl[\inf_{s\in[\delta,K]} |f(s)|\bigr]\theta^2.
    \end{equation}}
    \argument{\lref{eq15};\lref{evd2};}{for all $\theta\in [-K,-\delta]$ that
    \begin{equation}\llabel{eq15'}
    \begin{split}
        \textstyle\theta f(\theta)&=\theta ( - f( - \theta ) ) =(-\theta)f(-\theta)\\
        &\leq \textstyle-K^{-1}\bigl[\inf_{s\in[\delta,K]} |f(s)|\bigr](-\theta)^2\\
        &\textstyle= -K^{-1}\bigl[\inf_{s\in[\delta,K]} |f(s)|\bigr]\theta^2.
        \end{split}
    \end{equation}}
    \argument{\lref{eq15};\lref{eq15'}}{for all $\theta \in [ -K, - \delta ] \cup [\delta, K ]$ that
\begin{equation}\llabel{eq15main}
  \theta f( \theta ) \leq \textstyle -K^{-1}\bigl[\inf_{s\in[\delta,K]} |f(s)|\bigr]\theta^2.
\end{equation}}
    \argument{\lref{eqdb};\lref{eq**};\lref{eq15main}}{that for all $\theta\in [-K,K]$ it holds that
    \begin{equation}\llabel{eq16}
      \textstyle  \theta(f(\theta)-f(0))=\theta f(\theta)\leq \Bigl[\max\Bigl\{-\bigl[ \inf_{ s \in[-\delta,\delta] } | f'(s) | \bigr],-K^{-1}\bigl[\inf_{s\in[\delta,K]} |f(s)|\bigr] \Bigr\}\Bigr]\theta^2.
    \end{equation}}
    \argument{\lref{eq16};\lref{eq9.1};\lref{def: bfC}}{\cref{item 4: monotonicity symmetric pre}\dott}
\end{aproof}
\newcommand{\xx}{\theta}
\begin{athm}{cor}{lem: monotonicity Adam vector field at origin}
     Let $(\Omega,\cF,\P)$ be a probability space, let $d\in \N$, $\fc\in [0,\infty)$, let $X_{n,m}=(X_{n,m}^1,\dots,X_{n,m}^d)\colon \Omega\allowbreak\to [-\fc,\fc]^d$, $(n,m) \in (\N_0)^2$, be \iid\ random variables, assume that $X_{1,1}$ is symmetrically distributed, let $\alpha\in (0,1)$, $\beta\in (\alpha^2,1)$, $\varepsilon\in (0,\infty)$, $M\in \N$, let $f=(f_1,\dots,f_d)\colon\R^d\allowbreak\to\R^d$ satisfy for all $i\in \{1,2,\dots,d\}$, $\theta=(\theta_1,\dots,\theta_d)\in \R^d$ that
     \begin{equation}\llabel{eq1}
    \begin{split}
        f_i(\theta)&=(1-\alpha)\,\E\Biggl[\frac{\textstyle\sum_{n= 0}^\infty \alpha^nM^{-1}\textstyle \sum_{m=1}^M2(X_{n,m}^i-\theta_i)}{\textstyle\varepsilon+\sqrt{(1-\beta)\sum_{n=0}^\infty\beta^n
          |M^{-1}\sum_{m=1}^M2(X_{n,m}^i-\theta_i)|^2}}\Biggr],
          \end{split}
    \end{equation}
    and let $V\subseteq \R^d$ be compact.
    Then
    \begin{enumerate}[label=(\roman*)]
        \item \label{item 1: monotinicity symmetric} it holds that $f\in C^1(\R^d,\R^d)$,
        \item \label{item 3: monotonicity symmetric} it holds that $f(\E[X_{1,1}])=0$, and
        \item \label{item 4: monotonicity symmetric} there exists $\fC\in (0,\infty)$ such that for all $\theta\in V$ it holds that $\spro{\theta-\E[X_{1,1}], f(\theta)}\leq -\fC\|\theta-\E[X_{1,1}]\|^2$.
    \end{enumerate}
\end{athm}
\begin{aproof}
Throughout this proof for every $i\in \{1,2,\dots,d\}$ let $g_i\colon \R\to\R$ satisfy for all $\theta\in \R$ that
\begin{equation}\llabel{def: Adam}
    g_i(\theta)= (1-\alpha)\,\E\Biggl[\frac{\textstyle\sum_{n= 0}^\infty \alpha^nM^{-1}\textstyle \sum_{m=1}^M2(X_{n,m}^i-\E[X_{n,m}^i]-\theta)}{\textstyle\varepsilon+\sqrt{(1-\beta)\sum_{n=0}^\infty\beta^n
          |M^{-1}\sum_{m=1}^M2(X_{n,m}^i-\E[X_{n,m}^i]-\theta)|^2}}\Biggr]
\end{equation}
and for every $i\in \{1,2,\dots,d\}$, $n,m\in \N$ let $Y^i_{n,m}\colon \Omega\to\R$ satisfy for all $\omega\in \Omega$ that
\begin{equation}\llabel{def: Y}
    Y^i_{n,m}(\omega)=X_{n,m}^i(\omega)-\E[X_{n,m}^i].
    \end{equation}
      \argument{\lref{def: Adam};\lref{def: Y}}{for all $i\in \{1,2,\dots,d\}$, $\theta\in \R$ that
    \begin{equation}\llabel{def: Adam 2}
        g_i(\theta)=(1-\alpha)\,\E\Biggl[\frac{\textstyle\sum_{n= 0}^\infty \alpha^nM^{-1}\textstyle \sum_{m=1}^M2(Y_{n,m}^i-\theta)}{\textstyle\varepsilon+\sqrt{(1-\beta)\sum_{n=0}^\infty\beta^n
          |M^{-1}\sum_{m=1}^M2(Y_{n,m}^i-\theta)|^2}}\Biggr].
    \end{equation}}
    \argument{the assumption that $V$ is compact;}{that there exists $K\in (0,\infty)$ which satisfies
    \begin{equation}\llabel{def: K}
        V\subseteq -[K,K]^d.
    \end{equation}}
    \startnewargseq
    \argument{\lref{def: Y}; the fact that for all $i\in \{1,2,\dots,d\}$, $n,m\in \N$ it holds that $-\fc\leq |X_{n,m}^i|\leq \fc$}{that for all $i\in \{1,2,\dots,d\}$, $n,m\in \N$ it holds that
    \begin{equation}\llabel{argg0}
        -2\fc\leq |Y^i_{n,m}|\leq 2\fc.
    \end{equation}}
    \argument{\lref{argg0};the fact that $X_{n,m}$, $(n,m)\in \N^2$, are \iid;\lref{def: Y};}{that for all $i\in \{1,2,\dots,d\}$ \llabel{argg1} it holds that $\Omega\ni\omega\mapsto Y^i_{n,m}(\omega)\in [-2\fc,2\fc]$, $(n,m)\in (\N_0)^2$, are \iid\ }
    \argument{\lref{argg1};\lref{def: Adam 2}; the fact that for all $i\in \{1,2,\dots,d\}$ it holds that $Y_{1,1}^i$ is symmetrically distributed;the fact that for all $i\in \{1,2,\dots,d\}$ it holds that $\E[Y_{1,1}^i]=0$;\cref{lem: monotonicity Adam vector field at origin pre} (applied for every $i\in \{1,2,\dots,d\}$ with $\Omega\curvearrowleft\Omega$, $\cF\curvearrowleft\cF$, $\P\curvearrowleft\P$, $\fc\curvearrowleft 2\fc$, $(X_{n,m})_{(n,m)\in \N^2}\curvearrowleft (\Omega\ni\omega\mapsto Y^i_{n,m}(\omega)\in [-2\fc,2\fc])_{(n,m)\in \N^2}$, $\alpha\curvearrowleft\alpha$, $\beta\curvearrowleft\beta$, $\varepsilon\curvearrowleft\varepsilon$, $M\curvearrowleft M$, $f\curvearrowleft g_i$ in the notation of \cref{lem: monotonicity Adam vector field at origin pre})}{that there exists $\fC\in (0,\infty)$ which satisfies for all $i\in\{1,2,\dots,d\}$ that
    \begin{enumerate}[label=(\Roman*)]
        \item \llabel{item 1} it holds that $g_i\in C^1(\R,\R)$,
        \item \llabel{item 2} it holds for all  $\theta\in (0,\infty)$ that $ (g_i)'( 0 ) < 0 = g_i( 0 ) < - g_i( \theta ) = g_i( - \theta )$, and
        \item \llabel{item 3} it holds for all $\theta\in \bigl[-K-\|\E[X_{1,1}]\|,K+\|\E[X_{1,1}]\|\bigr]$ that $\theta g_i(\theta)\leq -\fC\theta^2$.
    \end{enumerate}}
    \startnewargseq
    \argument{\lref{eq1};\lref{def: Adam}}{that for all $\theta=(\theta_1,\dots,\theta_d)\in \R^d$, $i\in \{1,2,\dots,d\}$ it holds that
    \begin{equation}\llabel{eq2}
        f_i(\theta)=g_i(\theta_i-\E[X_{n,m}^i]).
    \end{equation}}
    \argument{\lref{eq2};\lref{item 1}}{for all $i\in \{1,2,\dots,d\}$ that \llabel{arg1}  $f_i\in C^1(\R^d,\R)$\dott}
    \argument{\lref{arg1};}{\cref{item 1: monotinicity symmetric}\dott}
    \startnewargseq
    \argument{\lref{eq2};\lref{item 2};}{that for all $i\in \{1,2,\dots,d\}$ it holds that
    \begin{equation}\llabel{eq3}
        f_i(\E[X_{1,1}])=g_i(0)=0.
    \end{equation}}
    \argument{\lref{eq3};}{\cref{item 3: monotonicity symmetric}\dott}
    \startnewargseq
    \argument{\lref{eq2};\lref{item 3}; the fact that for all $\theta=(\theta_1,\dots,\theta_d)\in [-K,K]^d$, $i\in \{1,2,\dots,d\}$ it holds that
    \begin{equation}
    \begin{split}
      & -K-\|\E[X_{1,1}]\|\leq -K-|\E[X_{1,1}^i]| \leq \theta_i-\E[X_{1,1}^i]\\
      &\leq \theta_i+|\E[X_{1,1}^i]|\leq \theta_i+\|\E[X_{1,1}]\|\leq K+\|\E[X_{1,1}]\| 
       \end{split}
    \end{equation}}{that for all $\theta=(\theta_1,\dots,\theta_d)\in [-K,K]^d$, $i\in\{1,2,\dots,d\}$ it holds that
    \begin{equation}\llabel{eq4}
        (\theta_i-\E[X_{1,1}^i])f_i(\theta)= (\theta_i-\E[X_{1,1}^i])g_i(\theta_i-\E[X_{1,1}^i])\leq -\fC(\theta_i-\E[X_{1,1}^i])^2.
    \end{equation}}
    \argument{\lref{eq4};}{for all $\theta\in [-K,K]^d$ that
    \begin{equation}\llabel{eq5}
        \spro{\theta-\E[X_{1,1}], f(\theta)}\leq -\fC \big\|\theta-\E[X_{1,1}]\big\|^2.
    \end{equation}}
    \argument{\lref{eq5};\lref{def: K}}{\cref{item 4: monotonicity symmetric}\dott}
\end{aproof}
\subsection{Convergence of Adam to global minima for symmetrically distributed data}\label{subsec: convegence symmetry}
\renewcommand{\consc}{\mathfrak{C}}
\begin{athm}{prop}{prop: convergence in symmetric case}
     Let $(\Omega,\cF,\P)$ be a probability space, let $\fd\in \N$, let $X_{n,m}\colon \Omega\allowbreak\to \R^{\fd}$, $(n,m) \in \N^2$, be bounded \iid\ random variables, let $\smalll=(\smalll(\theta,x))_{(\theta,x)\in \R^\fd\times\R^{\fd}}\allowbreak\colon\R^\fd\times\R^{\fd}\to\R$ satisfy for all $\theta,x\in \R^\fd$ that $\smalll(\theta,x)=\|\theta-x\|^2$, let $\alpha\in (0,1)$, $\beta\in (\alpha^2,1)$, $\varepsilon\in (0,\infty)$, let $(\gamma_n)_{n\in \N}\subseteq (0,\infty)$ be non-increasing, for every $M\in \N$, $n\in \N_0$ let  $\bfm_n^M=(\bfm^M_{n,1},\dots,\bfm^M_{n,\fd})\colon \Omega\to \R^\fd$, $\democrator_n^M=(\democrator^M_{n,1},\dots,\democrator^M_{n,\fd})\colon  \Omega\to\R^\fd$, and $\Theta_n^M=(\Theta^M_{n,1},\dots,\Theta^M_{n,\fd})\colon\allowbreak \Omega\to\R^\fd$ be random variables, assume for all $M,n\in \N$, $i\in \{1,2,\dots,\fd\}$ that
    \begin{equation}\llabel{def: bbV}
     \bfm_n^M=0,\qquad\bfm_n^M=\alpha \bfm_{n-1}^M+(1-\alpha)\bigl[\textstyle \frac 1M \sum_{m=1}^M(\nabla_{\theta_i} \smalll)(\Theta_{n-1}^M,X_{n,m})\bigr],
    \end{equation}
    \begin{equation}  \llabel{def: bbV2}
        \democrator_n^M=0,\qquad \democrator_{n,i}^M=\beta \democrator_{n-1,i}^M+(1-\beta)\bigl|\textstyle \frac 1M \sum_{m=1}^M(\nabla_{\theta_i} \smalll)(\Theta_{n-1}^M,X_{n,m})\bigr|^2,
    \end{equation}
    \begin{equation}\llabel{def: Theta}
     \textstyle\text{and}\qquad \Theta_{n,i}^M=\Theta_{n-1,i}^M -\gamma_n\bigl[\frac{\bfm_{n,i}^M}{1-\alpha^n}\bigr] \Bigl[\varepsilon+\bigl[\frac{\democrator_{n,i}^M}{1-\beta^n}\bigr]^{\nicefrac{1}{2}}\Bigr]^{-1}, 
 \end{equation}
     assume $(\Theta_0^M)_{M\in \N}$ and $(X_{n,m})_{(n,m)\in \N^2}$ are independent, assume $\limsup_{n\to\infty}(( \gamma_n )^{ - 2 } ( \gamma_n - \gamma_{ n + 1 } ) )\allowbreak =0$, 
    let $\cL \colon \R^\fd \to 
     \R$ satisfy for all $\theta\in \R^\fd$ that 
     $\cL( \theta ) = \E[\smalll(\theta,X_{1,1})]$, let $\vartheta \in \R^{\fd}$ satisfy $\cL( \vartheta ) = \inf_{ \theta \in \R^\fd } \cL( \theta )$, assume that $X_{ 1, 1 }$ is symmetrically distributed, and let $ p \in 
         \N$ satisfy $\sum_{n=1}^\infty(\gamma_n)^p<\infty$ and $\P\bigl(\sup_{M\in \N}\|\Theta_0^M\|\leq p\bigr)=1$ \cfload. 
      Then there exists 
         $\scrc \in \R$ such that for all $M,n \in \N$ it holds that \begin{equation}\llabel{conclude}
             ( \E[ \| \Theta_n^M - \vartheta \|^p ] )^{ 1 / p } \leq \scrc \sqrt{ \gamma_{n+1} }.
         \end{equation}
\end{athm}
\begin{aproof}
    Throughout this proof let $\fc\in \R$ satisfy $\P(\|X_{1,1}\|\leq \fc)=1$, for every $n,m\in \N_0$ let $\bfX_{n,m}=(\bfX_{n,m}^1,\dots,\allowbreak\bfX_{n,m}^d)\colon\Omega\to \R^d$ satisfy 
\begin{equation}\llabel{def: bfX}
    \bfX_{n,m}=X_{n+1,m+1}\mathbbm 1_{\{\|X_{n+1,m+1}\|\leq \fc\}},
\end{equation}
let  $\bbM_n^M=(\bbM^M_{n,1},\dots,\bbM^M_{n,\fd})\colon \Omega\to \R^\fd$, $(M,n)\in \N\times\N_0$, $\bfV_n^M=(\bfV^M_{n,1},\dots,\bfV^M_{n,\fd})\colon  \Omega\to\R^\fd$, $(M,n)\in \N\times\N_0$, and $\varTheta_n^M=(\varTheta^M_{n,1},\dots,\varTheta^M_{n,\fd})\allowbreak\colon \Omega\to\R^\fd$, $(M,n)\in \N\times\N_0$, satisfy for all $M,n\in \N$, $i\in \{1,2,\dots,\fd\}$ that
    \begin{equation}\llabel{def: bbV'}
     \bbM_n^M=0,\qquad\bbM_n^M=\alpha \bfm_{n-1}^M+(1-\alpha)\bigl[\textstyle \frac 1M \sum_{m=1}^M(\nabla_{\theta_i} \smalll)(\varTheta_{n-1}^M,\bfX_{n,m})\bigr],
    \end{equation}
    \begin{equation}  \llabel{def: bbV2'}
        \bfV_n^M=0,\qquad \bfV_{n,i}^M=\beta \bfV_{n-1,i}^M+(1-\beta)\bigl|\textstyle \frac 1M \sum_{m=1}^M(\nabla_{\theta_i} \smalll)(\varTheta_{n-1}^M,\bfX_{n,m})\bigr|^2,
    \end{equation}
    \begin{equation}\llabel{def: Theta'}
     \varTheta^M_0=\Theta_0^M,\qquad\textstyle\text{and}\qquad \varTheta_{n,i}^M=\varTheta_{n-1,i}^M -\gamma_n\bigl[\frac{\bbM_{n,i}^M}{1-\alpha^n}\bigr] \Bigl[\varepsilon+\bigl[\frac{\bfV_{n,i}^M}{1-\beta^n}\bigr]^{\nicefrac{1}{2}}\Bigr]^{-1}, 
 \end{equation}
and for every $M\in \N$ let $f_M=(f_M^1,\dots,f_M^d)\colon \R^d\to \R^d$ satisfy for all $i\in \{1,2,\dots,d\}$, $\theta=(\theta_1,\dots,\theta_d)\in \R^d$ that
\begin{equation}\llabel{def: Adam1}
          f_M^i(\theta)=(1-\alpha)\,\E\Biggl[\frac{\textstyle-\sum_{n= 0}^\infty\alpha^n\textstyle M^{-1}\sum_{m=1}^M(\nabla_{\theta_i} \smalll)(\theta,\bfX_{n,m})}{\textstyle\varepsilon+\sqrt{(1-\beta)\sum_{n=0}^\infty\beta^n\bigl|\textstyle M^{-1}\sum_{m=1}^M(\nabla_{\theta_i} \smalll)(\theta,\bfX_{n,m})\bigr|^2}}\Biggr].
\end{equation}
\argument{the fact that for all $\theta,x\in \R^d$ it holds that $\smalll(\theta,x)=\|\theta-x\|^2$;}{that for all $i\in \{1,2,\dots,d\}$, $\theta=(\theta_1,\dots,\theta_d)\in \R^d$, $x=(x_1,\dots,x_d)\in \R^d$ it holds that
\begin{equation}\llabel{eq1.1}
   \textstyle  (\nabla_{\theta_i} \smalll)(\theta,x)=2(\theta-x_i)-2(x_i-\theta_i).
\end{equation}}
\argument{\lref{eq1.1};\lref{def: Adam1}}{for all $i\in \{1,2,\dots,d\}$, $\theta=(\theta_1,\dots,\theta_d)\in \R^d$ that
\begin{equation}\llabel{def: Adam}
    f_M^i(\theta)=(1-\alpha)\,\E\Biggl[\frac{\textstyle\sum_{n= 0}^\infty\alpha^n\textstyle M^{-1}\sum_{m=1}^M2(\bfX_{n,m}^i-\theta_i)}{\textstyle\varepsilon+\sqrt{(1-\beta)\sum_{n=0}^\infty\beta^n\bigl|\textstyle M^{-1}\sum_{m=1}^M2(\bfX_{n,m}^i-\theta_i)\bigr|^2}}\Biggr].
\end{equation}}
\argument{\lref{def: bfX};}{that for all $n,m\in \N_0$ it holds that
\begin{equation}\llabel{def: fc}
    \|\bfX_{n,m}\|\leq \fc.
\end{equation}}
\argument{\lref{def: fc}}{for all $M,n\in \N$,  $i\in \{1,2,\dots,d\}$ it holds that
\begin{equation}\llabel{eq1.2}
   \bigl|\textstyle\frac{1}{M} \sum_{m=1}^M \bfX_{n,m}^i\bigr|\leq \fc.
\end{equation}}
\argument{\lref{eq1.1};\lref{eq1.2};}{that for all $\theta=(\theta_1,\dots,\theta_d)\in\R^d$, $M,n\in\N$, $i\in \{1,2,\dots,d\}$ it holds  that 
    \begin{equation}\llabel{eqp1}
    \begin{split}
        &2(\theta_i-\fc)=2\theta_i-2\fc \leq 2\theta_i-\displaystyle\frac 2M \textstyle\sum\limits_{m=1}^M\bfX_{n,m}^i=\displaystyle\frac 1M \textstyle\sum\limits_{m=1}^M2(\theta_i-\bfX_{n,m}^i)\\
        &\textstyle = \displaystyle\frac 1M \textstyle\sum\limits_{m=1}^M(\nabla_{\theta_i} \smalll)(\theta,\bfX_{n,m}) \leq 2\theta_i+2\fc=2(\theta_i+\fc).
        \end{split}
    \end{equation}}
    \argument{\lref{def: bbV'};\lref{def: bbV2'};\unskip, \eg, \cite[item (i) in Lemma 6.3.17]{ArBePhi2024}}{that for all $M,n\in \N$, $i\in \{1,2,\dots,d\}$ it holds that
    \begin{equation}\llabel{eq1.3}
        \bbM_{n,i}^M=\textstyle(1-\alpha)\biggl[\sum\limits_{k=1}^n\alpha^{n-k}\textstyle \frac 1M \bigl[\textstyle\sum_{m=1}^M(\nabla_{\theta_i} \smalll)(\varTheta_{k-1}^M,\bfX_{k,m})\bigr]\biggr]
        \end{equation}
        \begin{equation}\llabel{eq1.4}
            \text{and}\qquad  \bfV_{n,i}^M=(1-\beta)\biggl[\textstyle\sum\limits_{k=1}^n\beta^{n-k}\bigl| \frac 1M \bigl[\textstyle\sum_{m=1}^M(\nabla_{\theta_i} \smalll)(\varTheta_{k-1}^M,\bfX_{k,m})\bigr]\bigr|^2 \textstyle \biggr].
    \end{equation}}
    \argument{\lref{eq1.4};\lref{def: Theta'};}{that for all $M,n\in \N$, $i\in \{1,2,\dots,d\}$ it holds $\P$-a.s that
    \begin{equation}\llabel{eq1.5}
         \varTheta_{n,i}^M=\varTheta_{n-1,i}^M -\frac{\gamma_n}{1-\alpha^n}\frac{\sum_{k=1}^n(1-\alpha)\alpha^{n-k}\textstyle \frac 1M [\textstyle\sum_{m=1}^M(\nabla_{\theta_i} \smalll)(\varTheta_{k-1}^M,\bfX_{k,m})]}{ \varepsilon+\bigl[\sum_{k=1}^n(\frac{1-\beta}{1-\beta^n})\beta^{n-k}\bigl| \frac 1M [\textstyle\sum_{m=1}^M(\nabla_{\theta_i} \smalll)(\varTheta_{k-1}^M,\bfX_{k,m})]\bigr|^2 \textstyle\bigr]^{\nicefrac{1}{2}}}.
    \end{equation}}
    In the following let $A\subseteq\Omega$ satisfy for all $\omega\in A$, $M,n\in \N$, $i\in \{1,2,\dots,d\}$ that
    \begin{equation}\llabel{def: A}
    \begin{split}
           &\varTheta_{n,i}^M(\omega)=\varTheta_{n-1,i}^M(\omega) \\
           &-\frac{\gamma_n}{1-\alpha^n}\frac{\sum_{k=1}^n(1-\alpha)\alpha^{n-k}\textstyle \frac 1M [\textstyle\sum_{m=1}^M(\nabla_{\theta_i} \smalll)(\varTheta_{k-1}^M(\omega),\bfX_{k,m}(\omega))]}{ \varepsilon+\bigl[\sum_{k=1}^n(\frac{1-\beta}{1-\beta^n})\beta^{n-k}\bigl| \frac 1M [\textstyle\sum_{m=1}^M(\nabla_{\theta_i} \smalll)(\varTheta_{k-1}^M(\omega),\bfX_{k,m}(\omega))]\bigr|^2 \textstyle\bigr]^{\nicefrac{1}{2}}}.
           \end{split}
    \end{equation}
    \startnewargseq
\argument{\lref{eqp1};\lref{def: A};the fact that $\inf_{n\in \N}\frac{1-\beta}{1-\beta^n}\geq 1-\beta>0$; the fact that $\sup_{n\in \N}\frac{\gamma_n}{1-\alpha^n}\leq \frac{\gamma_1}{1-\alpha}$;\cite[Corollary 2.5]{DeRoAr2024nonconvergence} (applied for every $M\in \N$, $i\in \{1,2,\dots,d\}$, $\omega\in A$ with $\mathfrak{d}\curvearrowleft d$, $i\curvearrowleft i$, $\varepsilon\curvearrowleft\varepsilon$, $\eta \curvearrowleft 2$, $\rho\curvearrowleft2$, $\alpha\curvearrowleft\alpha$, $\beta\curvearrowleft\beta$, $\fc\curvearrowleft \fc$, $\mathbb M\curvearrowleft 0$, $\mathcal M\curvearrowleft 0$, $(G_n)_{n\in \N}\curvearrowleft \bigl(\frac 1M \bigl[\textstyle\sum_{m=1}^M(\nabla_{\theta_i} \smalll)(\cdot,\bfX_{n,m}(\omega))\bigr]\bigr)_{n\in \N}$, $\kappa\curvearrowleft \bigl(\frac{1-\beta}{1-\beta^n}\bigr)_{n\in \N}$, $(\gamma_n)_{n\in \N}\curvearrowleft \bigl(\frac{1-\alpha}{1-\alpha^n}\bigr)_{n\in \N}$, $\Theta\curvearrowleft (\varTheta_n^M(\omega))_{n\in \N_0}$ in the notation of \cite[Corollary 2.5]{DeRoAr2024nonconvergence})}{that for all $M\in \N$, $i\in \{1,2,\dots,d\}$, $\omega\in A$, $n\in \N_0$ it holds that
\begin{equation}\llabel{eqp1'}
   \textstyle |\varTheta_n^i(\omega)|\leq \fc+3\max\Bigl\{|\varTheta_{0,i}^M(\omega)|,\fc,\frac{\gamma_1(2+\alpha)\beta^{1/2}}{(1-\alpha)(1-\beta)^{1/2}(\beta^{1/2}-\alpha)}\Bigr\}.
\end{equation}}
\argument{\lref{eq1.5};\lref{def: A}}{that \llabel{argg1} $\P(A)=1$\dott}
\argument{\lref{eqp1'};\lref{argg1}; the fact that $\P\bigl(\sup_{M\in \N}\|\varTheta_0^M\|\leq p\bigr)=1$}{that there exists $\consc\in (0,\infty)$ such that 
\begin{equation}\llabel{eq11}
      \textstyle \P\bigg(\sup\limits_{(M,n)\in \N\times \N_0}\|\varTheta_n^M\|\leq \consc\biggr)=1.
     \end{equation}}
     \argument{\lref{eq11}; the fact that $\|\bfX_{1,1}\|\leq \fc$ }{that there exists $\consc\in (0,\infty)$ such that
     \begin{equation}\llabel{eq1}
          \textstyle \P\biggl(\|\bfX_{1,1}\|+\sup\limits_{(M,n)\in \N\times\N_0}\|\varTheta_n^M\|\leq \consc\biggr)=1.
     \end{equation}}
\argument{\cref{item 3: monotonicity symmetric,item 4: monotonicity symmetric} in \cref{lem: monotonicity Adam vector field at origin};\lref{def: Adam}}{that for all $i\in \{1,2,\dots,d\}$ there exists $c\in (0,\infty)$ such that for all $M\in \N$, $\theta=(\theta_1,\dots,\theta_d)\in\R^d$ with $\|\theta\|\leq \consc$ that
\begin{equation}\llabel{eq:934657}
\spro{\theta-\E[\bfX_{1,1}],f_M(\theta)-f_M(\E[\bfX_{1,1}])}=\spro{\theta-\E[\bfX_{1,1}],f_M(\theta)}\leq -c\bigl\|\theta-\E[\bfX_{1,1}]\bigr\|^2.
\end{equation}}
\argument{\cref{item 3: monotonicity symmetric} in \cref{lem: monotonicity Adam vector field at origin};\lref{def: Adam}}{that \llabel{arg1} for all $M\in \N$ it holds that $f_M(\E[\bfX_{1,1}])=0$\dott}
\argument{\lref{arg1};\lref{eq1};\lref{eq:934657};the assumption that $\Theta_0$ and $(X_{n,m})_{(n,m)\in \N^2}$ are independent; the assumption that $\sum_{n=1}^\infty(\gamma_n)^p<\infty$;the assumption that 
\begin{equation}
\textstyle\limsup_{n\to\infty}\bigl((\gamma_n)^{-2}(\gamma_n-\gamma_{n+1})\bigr)=0;
\end{equation}
the assumption that $(\gamma_n)_{n\in \N}$ is non-increasing; \cref{prop: convergence general} (applied with $\Omega\curvearrowleft\Omega$, $\cF\curvearrowleft\cF$, $\P\curvearrowleft\P$, $\fd\curvearrowleft\fd$, $\dimX\curvearrowleft\fd$, $(X_{n,m})_{(n,m)\in (\N_0)^2}\curvearrowleft (X_{n,m})_{(n,m)\in (\N_0)^2}$, $\smalll\curvearrowleft\smalll$, $\alpha\curvearrowleft\alpha$, $\beta\curvearrowleft\beta$, $\fC\curvearrowleft \consc$, $\fc\curvearrowleft c$, $\varepsilon\curvearrowleft\varepsilon$, $(\gamma_n)_{n\in \N}\curvearrowleft (\gamma_n)_{n\in \N}$, $(f_M)_{M\in \N}\curvearrowleft(f_M)_{M\in \N}$, $(\globalmin_M)_{M\in \N}\curvearrowleft(\E[\bfX_{1,1}])_{M\in \N}$, $p\curvearrowleft p$ in the notation of \cref{prop: convergence general})}{that there exist $\eta\in (0,\infty)$ and $\fn\in \N$ which satisfy for all $M\in \N$, $n\in \N\cap[\fn,\infty)$ that
\begin{equation}\llabel{eqq1}
\begin{split}
&\bigl(\E\bigl[\bigl\|\varTheta_n^M-\E[\bfX_{1,1}]\bigr\|^p\bigl]\bigr)^{1/p} \leq \eta  \sqrt{\gamma _{n+1}}.
\end{split}
\end{equation}}
\startnewargseq
\argument{\lref{def: bfX}; the assumption that $X_{n,m}$, $(n,m)\in \N^2$, are \iid; the fact that $\P(\|X_{1,1}\|\leq \fc)=1$;}{that \llabel{eqq1.1} $(X_{n,m})_{(n,m)\in \N^2}$ and $(\bfX_{n,m})_{(n,m)\in (\N_0)^2}$ are identically distributed\dott}
\argument{\lref{eqq1.1};\lref{def: bbV};\lref{def: bbV2};\lref{def: Theta};\lref{def: bbV'};\lref{def: bbV2'};\lref{def: Theta'}}{that for all $M,n\in \N$ it holds that \llabel{eqq1.2} $\varTheta_n^M$ and $\Theta_n^M$ are identiacally distributed\dott}
\argument{\lref{eqq1.2};\lref{eqq1};}{for all $M\in \N$, $n\in \N\cap[\fn,\infty)$ that
\begin{equation}\llabel{eqq1'}
\bigl(\E\bigl[\bigl\|\Theta_n^M-\E[\bfX_{1,1}]\bigr\|^p\bigl]\bigr)^{1/p} \leq \eta  \sqrt{\gamma _{n+1}}.
\end{equation}}
\argument{\lref{eqq1'};\lref{eq11};the fact that $\globalmin=\E[\bfX_{1,1}]$}{\lref{conclude}\dott}
\end{aproof}

\subsection{Characterization of the convergence of Adam to global minima}\label{subsec: chracterization}
\begin{athm}{lemma}{lem: tg3}
    Let $(\gamma_n)_{n\in \N}\subseteq (0,\infty)$ be non-increasing, assume $\lim_{n\to\infty}\gamma_n=0$, and assume $\limsup_{n\to\infty}( ( \gamma_n )^{ - 2 } ( \gamma_n - \gamma_{ n + 1 } ) )<\infty$. Then
    \begin{equation}\llabel{conclude}
        \textstyle \sup_{r\in (0,1)}\liminf_{n\to\infty}((\gamma_n)^{-1}\gamma_{n+\lceil(\gamma_n)^{-1}r\rceil})>0.
    \end{equation}
\end{athm}
\begin{aproof}
    Throughout this proof let $c\in \R$ satisfy
    \begin{equation}\llabel{def: c}
       \textstyle c= 1+\limsup_{n\to\infty}( ( \gamma_n )^{ - 2 } ( \gamma_n - \gamma_{ n + 1 } ) ).
    \end{equation}
    \argument{\lref{def: c};the assumption that $\lim_{n\to\infty}\gamma_n=0$}{that there exists $N\in \N$ which satisfies for all $n\in \N\cap [N,\infty)$ that
    \begin{equation}\llabel{eq1}
        ( \gamma_n )^{ - 2 } ( \gamma_n - \gamma_{ n + 1 } )\leq 2c\qqandqq \gamma_n\leq (8c)^{-1} .
    \end{equation}}
    \startnewargseq
    \argument{\lref{eq1};}{that for all $n\in [N,\infty)$ it holds that
    \begin{equation}\llabel{eq2}
        \gamma_n - \gamma_{ n + 1 }\leq 2c(\gamma_n)^2.
    \end{equation}}
    \argument{\lref{eq2};induction;the fact that for all $m,n\in \N$ with $m\geq n$ it holds that $\gamma_m\leq \gamma_n$}{that for all $n\in \N\cap[N,\infty)$ it holds that
    \begin{equation}\llabel{eq3}
        \begin{split}
            &\gamma_n-\gamma_{n+\lceil(8c\gamma_n)^{-1}\rceil}=\sum_{k=n+1}^{n+\lceil(8c\gamma_n)^{-1}\rceil} (\gamma_n-\gamma_{n-1})\leq \sum_{k=n+1}^{n+\lceil(8c\gamma_n)^{-1}\rceil} 2c(\gamma_k)^{2}\\
            &\leq \sum_{k=n+1}^{n+\lceil(8c\gamma_n)^{-1}\rceil} 2c(\gamma_n)^{2}= 2c\lceil(8c\gamma_n)^{-1}\rceil(\gamma_n)^2.
        \end{split}
    \end{equation}}
    \argument{\lref{eq1};}{for all $n\in \N\cap[N,\infty)$ that \llabel{arg1} $(8c\gamma_n)^{-1}\geq 1$\dott}
    \argument{\lref{arg1};}{for all $n\in \N\cap[N,\infty)$ that \llabel{arg2} $\lceil(8c\gamma_n)^{-1}\rceil\leq 2(8c\gamma_n)^{-1}=(4c\gamma_n)^{-1}$\dott}
    \argument{\lref{arg2};\lref{eq3};}{for all $n\in \N\cap[N,\infty)$ that
    \begin{equation}\llabel{eq4}
        \gamma_n-\gamma_{n+\lceil(8c\gamma_n)^{-1}\rceil}\leq 2c (4c\gamma_n)^{-1}(\gamma_n)^2=\frac{\gamma_n}{2}.
    \end{equation}}
    \argument{\lref{eq4};}{for all $n\in \N\cap[N,\infty)$ that
    \begin{equation}\llabel{eq5}
        \gamma_{n+\lceil(8c\gamma_n)^{-1}\rceil}\geq \frac{\gamma_n}{2}.
    \end{equation}}
    \argument{\lref{eq5};}{that
    \begin{equation}\llabel{eq6}
        \textstyle \liminf_{n\to\infty}((\gamma_n)^{-1}\gamma_{n+\lceil(\gamma_n)^{-1}(8c)^{-1}\rceil})\geq \frac 12>0.
    \end{equation}}
    \argument{\lref{eq6};the fact that $0<(8c)^{-1}<1$}{\lref{conclude}\dott}
    \end{aproof}
    \renewcommand{\mandV}[3]{\Theta^{#1,#2}_{#3}}
\renewcommand{\mandVcom}[4]{\Theta^{#1,#2}_{#3,#4}}
\newpage
\begin{samepage}
\begin{athm}{theorem}{main theorem 2'}[\Adam\ symmetry theorem]
     Let $(\Omega,\cF,\P)$ be a probability space, let $\fc,\fd\in \N$, let $X_{n,m}\colon \Omega\allowbreak\to \R^{\fd}$, $(n,m) \in \N^2$, be \iid\ \tpd\ random variables, let $\smalll=(\smalll(\theta,x))_{(\theta,x)\in \R^\fd\times\R^{\fd}}\allowbreak\colon\R^\fd\times\R^{\fd}\to\R$ satisfy for all $\theta,x\in \R^\fd$ that $\smalll(\theta,x)=\|\theta-x\|^2$, let $\para_1\in (0,1)$, $\para_2\in ((\para_1)^2,1)$, $\varepsilon\in (0,\infty)$, let $(\gamma_n)_{n\in \N}\subseteq (0,\infty)$ be non-increasing, for every $k,M,n\in \N_0$ let 
    $
\mandV{k}{M}{n}=(\mandVcom{k}{M}{n}{1},\dots,\mandVcom{k}{M}{n}{\fd})\colon\allowbreak \Omega\to\R^\fd$ be a random variable, assume for all $k\in \{1,2\}$, $M,n\in \N$, $i\in \{1,2,\dots,\fd\}$ that
    \begin{equation}\label{eq1: main theorem 2'}
    \mandV{k}{M}{0}=0,\qquad \mandVcom{k}{M}{n}{i}= \beta_k \mandVcom{k}{M}{n-1}{i}+(1-\beta_k)\bigl[\textstyle \frac 1M \sum_{m=1}^M(\nabla_{\theta_i} \smalll)(\mandV{0}{M}{n-1},X_{n,m})\bigr]^k,
        \end{equation}
        \begin{equation}\label{eq2: main theorem 2'}
 \mandVcom{0}{M}{n}{i}=\mandVcom{0}{M}{n-1}{i} -\gamma_n [1-(\para_1)^n]^{-1}\bigl[\varepsilon+\bigl[[1-(\para_2)^n]^{-1}\mandVcom{2}{M}{n}{i}\bigr]^{\nicefrac{1}{2}}\bigr]^{-1}\mandVcom{1}{M}{n}{i}, 
     \end{equation}
     and $\P(\|\mandV{0}{M}{0}\|\leq \fc)=1$,
     assume $\limsup_{n\to\infty}( ( \gamma_n )^{ - 2 } ( \gamma_n - \gamma_{ n + 1 } ) +\sum_{m=n}^\infty(\gamma_m)^\fc) =0$, assume that $(\Theta_0^{0,M})_{M\in \N}$ and $(X_{n,m})_{(n,m)\in \N^2}$ are independent, let $\cL \colon \R^\fd \to 
     \R$ satisfy for all $\theta\in \R^\fd$ that 
     $\cL( \theta ) = \E[\smalll(\theta,X_{1,1})]$, and let $\vartheta \in \R^{\fd}$ satisfy $\cL( \vartheta ) = \inf_{ \theta \in \R^\fd } \cL( \theta )$ \cfload.
      Then the following four statements are equivalent:
     \begin{enumerate}[label=(\roman*)]
      \item \label{item 4: main theorem 2'} It holds for all $M \in \N$ that $\limsup_{ n \to \infty } \E\bigl[ \min\bigl\{ 1, \cL(  \mandV{0}{M}{n}) - \inf_{ \theta \in \R^d } \cL( \theta ) \bigr\} \bigr] = 0$.
      \item \label{item 1: main theorem 2'} It holds for all $M\in \N$ that $\limsup_{ n \to \infty }\E[ \min\{ 1, \| \Theta_n^{0,M} - \vartheta\|\} ] = 0$.
     % \item \label{item 1: main theorem 2} It holds for all $M\in \N$, $\epsilon\in (0,\infty)$ that $\limsup_{n\to\infty}\P(\|\Theta_n^M-\vartheta\|\geq \epsilon)=0$.
     % \item \label{item 1: main theorem 2} It holds for all $M \in \N$, $\theta \in \R^\fd$ with $\limsup_{ n \to \infty } \E[ \min\{ 1, \|\Theta_n^M-\theta\| \} ] = 0$ that $\theta = \vartheta$.
     % \item \label{item 1: main theorem 2'} For every $M \in \N$ and every random variable $\chi \colon \Omega \to \R^{ \fd }$ with $\limsup_{ n \to \infty } \E[ \min\{ 1, \| \Theta^M_n - \chi \| \} ]$ = 0 it holds that $\P( \chi = \vartheta ) = 1$.
         % \item \label{item 3: main theorem 2} It holds for all $M \in \N$ that $\limsup_{ n \to \infty } \E[ \min\{ 1, \| \Theta^M_n - \vartheta \| \} ]  = 0$.
         \item \label{item 2: main theorem 2'} For all $ p \in 
        (0,\infty)$ there exists 
         $c \in \R$ such that for all $M\in \N$, $n \in \N_0$ it holds that $( \E[ \| \Theta_n^{0,M} - \vartheta \|^p ] )^{ 1 / p } \leq c \sqrt{ \gamma_{n+1} }$.
         % \item \label{item 3: main theorem 2} For all $ p \in 
         % \N$ there exists 
         % $c \in \R$ such that for all $M,n \in \N$ it holds that $( \E[ \| \Theta_n^M - \vartheta \|^p ] )^{ 1 / p } \leq c \sqrt{ \gamma_n }$.
         \item \label{item 3: main theorem 2'} It holds that $X_{ 1, 1 }$ is symmetrically distributed.
     \end{enumerate}
\end{athm}
\end{samepage}
\begin{aproof}
Throughout this proof let $\spro{\cdot\,,\cdot}\colon \R^d\times\R^d\to\R$ satisfy for all $v = ( v_1, \dots, v_d )$, $w = ( w_1, \dots, w_d ) \in \R^d$ that 
\begin{equation}
    \spro{v,w} = \sum_{ i = 1 }^{ d } v_iw_i.
\end{equation}
    \startnewargseq
  \argument{the fact that $  \sum_{ n = 1 }^{ \infty } ( \gamma_n )^\fc  < \infty$}{that \llabel{arg1} $\limsup_{ n \to \infty } \gamma_n = 0$\dott}
    \argument{\lref{arg1};}{that (\cref{item 2: main theorem 2'}$\rightarrow$\cref{item 1: main theorem 2'})\dott}
    \startnewargseq
   \argument{the fact that for all $\theta, x \in \R^d$ it holds that $\smalll( \theta, x ) = \| \theta - x \|^2$}{that for all $\theta \in \R^d$ it holds that
\begin{equation}\llabel{eq1}
\begin{split}
  \cL( \theta ) &= \E\bigl[ \| \theta - X_{ 1, 1 } \|^2 \bigr] 
  = \E\bigl[ \| ( \theta - \E[ X_{ 1, 1 } ] ) + ( \E[ X_{ 1, 1 } ] - X_{ 1, 1 } ) \|^2 \bigr]\\
  &
  = \E\bigl[\| \theta - \E[ X_{ 1, 1 } \|^2 + 2 \spro {\theta - \E[ X_{1,1} ], \E[ X_{1,1} ] - X_{ 1,1 } } + \| X_{ 1, 1 } - \E[ X_{ 1, 1 } ] \|^2 \bigr]\\
  &
  = \| \theta - \E[ X_{ 1, 1 } ] \|^2 + \E[ \| X_{ 1, 1 } - \E[ X_{1,1} ] \|^2 ] .
  \end{split}
\end{equation}}
\argument{\lref{eq1}; the fact that $\cL( \vartheta ) = \inf_{ \theta\in \R^d } \cL( \theta )$}{that 
\begin{equation}\llabel{argg1}
    \vartheta = \E[ X_{1,1}].
\end{equation}}
\argument{the fact that  $\limsup_{n\to\infty}( ( \gamma_n )^{ - 2 } ( \gamma_n - \gamma_{ n + 1 } ) ) =0$; the assumption that $(\gamma_n)_{n\in \N}\subseteq (0,\infty)$ is non-increasing; the fact that $\limsup_{n\to\infty}\gamma_n=0$;\cref{lem: tg3}}{that
\begin{equation}\llabel{argg1'}
    \textstyle \sup_{r\in (0,1)}\liminf_{n\to\infty}((\gamma_n)^{-1}\gamma_{n+\lceil(\gamma_n)^{-1}r\rceil})>0.
\end{equation}}
\argument{\lref{argg1};\lref{argg1'};\cref{lem: non-convergence of Adam2}}{that (\cref{item 1: main theorem 2'}$\rightarrow$\cref{item 3: main theorem 2'})\dott}
    \startnewargseq
    \argument{\cref{prop: convergence in symmetric case};}{that (\cref{item 3: main theorem 2'}$\allowbreak\rightarrow$\cref{item 2: main theorem 2'})\dott}
          \startnewargseq
         \argument{\lref{eq1};\lref{argg1}; the assumption that $\cL( \vartheta ) = \inf_{ \theta \in \R^\fd } \cL( \theta )$ }{that for all $M,n\in \N$ it holds that
          \begin{equation}\llabel{eq3}
          \begin{split}
           \textstyle\cL(  \mandV{0}{M}{n}) - \inf_{ \theta \in \R^d } \cL( \theta )= \cL(  \mandV{0}{M}{n}) -\cL(\vartheta)=\| \mandV{0}{M}{n}-\vartheta\|^2.
           \end{split}
          \end{equation}}
          \argument{\lref{eq3};the fact that $\lim_{n\to\infty}\gamma_n=0$}{that (\cref{item 2: main theorem 2'}$\allowbreak\rightarrow$\cref{item 4: main theorem 2'})\dott}
              \startnewargseq
              \argument{\lref{eq3};the Cauchy-Schwarz inequality;}{that for all $M,n\in \N$ it holds that
              \begin{equation}\llabel{eq4}
              \begin{split}
              \textstyle \E\bigl[ \min\bigl\{ 1, \cL(  \mandV{0}{M}{n}) - \inf_{ \theta \in \R^d } \cL( \theta ) \bigr\} \bigr]&= \E[ \min\{ 1, \|\mandV{0}{M}{n} -\vartheta\|^2 \} ]\\
              &\geq \bigl(\E[ \min\{ 1, \|\mandV{0}{M}{n} -\vartheta\| \} ]\bigr)^{2}.
              \end{split}
              \end{equation}}
              \argument{\lref{eq4};}{that (\cref{item 4: main theorem 2'}$\allowbreak\rightarrow$\cref{item 1: main theorem 2'})\dott}
             \end{aproof}
\section{Numerical simulations}\label{sec: numerial simulation}
In this section we complement the analytical findings of this work by means of several numerical simulations that indicate the sharpness of the convergence rates established in \cref{main theorem} and \cref{cor: convegence of Adam}, respectively, and that illustrate the practical appearance of the revealed theoretical phenomena in the \Adam\ symmetry theorem in \cref{main theorem 2} and \cref{main theorem 2'}, respectively.

\subsection{Description of the Adam optimizer and the Adam vector field}

In the following framework, \cref{setting: numerical simulation} below, we describe a setup that we employ in our numerical simulations in this section.
\begingroup
\renewcommand{\mandV}[3]{\Theta^{#1}_{#3}}
\renewcommand{\mandVcom}[4]{\Theta^{#1}_{#3}}
\newcommand{\numGDsteps}{600000}
\begin{setting}
\label{setting: numerical simulation}
Let $(\Omega,\cF,\P)$ be a probability space,
let $\para_1,\para_2\in (0,1)$, $\varepsilon\in (0,\infty)$, $M \in \N$,
for every $n,m \in \N$ let $X_{n,m}\colon \Omega\allowbreak\to \R$ be a random variable with $\E\bigl[|X_{n,m}|\bigr]<\infty$ and $\E[X_{n,m}]=0$,
let $\smalll=(\smalll(\theta,x))_{(\theta,x)\in \R\times\R}\allowbreak\colon\R\times\R\to\R$ satisfy for all $\theta,x\in \R$ that $\smalll(\theta,x)=|\theta-x|^2$,  
let $(\gamma_n)_{n\in \N}\subseteq (0,\infty)$ satisfy for all $n \in \N$ that 
$\gamma_n = n^{-0.99}$, 
for every $k,n\in \N_0$ let 
    $
\mandV{k}{M}{n}\colon\allowbreak \Omega\to\R$ be a random variable, 
assume for all $k\in \{1,2\}$, $n\in \N$ that
    \begin{equation}\label{eq1: main theorem 2 simulation}
    \mandV{k}{M}{0}=0,\qquad \mandVcom{k}{M}{n}{i}= \beta_k \mandVcom{k}{M}{n-1}{i}+(1-\beta_k)\bigl[\textstyle \frac 1M \sum_{m=1}^M(\nabla_{\theta_i} \smalll)(\mandV{0}{M}{n-1},X_{n,m})\bigr]^k,
        \end{equation}
        \begin{equation}\label{eq2: main theorem 2 simulation}
        \andq
 \mandVcom{0}{M}{n}{i}=\mandVcom{0}{M}{n-1}{i} -\gamma_n [1-(\para_1)^n]^{-1}\bigl[\varepsilon+\bigl[[1-(\para_2)^n]^{-1}\mandVcom{2}{M}{n}{i}\bigr]^{\nicefrac{1}{2}}\bigr]^{-1}\mandVcom{1}{M}{n}{i},
     \end{equation}
     assume that $\Theta^0_0$ is standard normal,
and let $f \colon \R \to \R$ satisfy for all $\theta\in \R$ that 
    \begin{equation}\label{def: Adam vector field numerical}
          f(\theta)=\lim_{N\to\infty}\E\Biggl[\frac{(1-\para_1)\textstyle\sum_{n= 0}^N(\para_1)^n M^{-1} \textstyle\sum_{m=1}^M(\nabla_{\theta} \smalll)(\theta,X_{n,m})}{\textstyle\varepsilon+\sqrt{(1-\para_2)\sum_{n=0}^N(\para_2)^n\bigl|\textstyle M^{-1} \sum_{m=1}^M(\nabla_{\theta} \smalll)(\theta,X_{n,m})\bigr|^2}}\Biggr]
    \end{equation}
(cf.\ \cref{theo: convergence point of Adam copy 3 pre}).
\end{setting}

In the following result, \cref{numerical simulation: consequences 1} below, we briefly note two elementary consequences of \cref{setting: numerical simulation} and, in particular, we observe that the sequence of learning rates $(\gamma_n )_{ n \in \N} \subseteq (0, \infty)$ in  \cref{setting: numerical simulation} satisfies the assumption on the learning rates in \cref{main theorem 2} (cf.\ also \cref{assumption: gamma1}).
\begin{athm}{lemma}{numerical simulation: consequences 1}
Assume \cref{setting: numerical simulation}.
Then 
\begin{enumerate}[label=(\roman*)]
\item \label{item 1: numerical simulation: consequences 1} it holds for all $p \in (\frac{100}{99}, \infty)$ that $\limsup_{n\to\infty}( ( \gamma_n )^{ - 2 } ( \gamma_n - \gamma_{ n + 1 } )+ \sum_{m=n}^\infty(\gamma_m)^p) =0$ and
\item \label{item 2: numerical simulation: consequences 1} it holds that 
$
  \big\{ \vartheta \in \R \colon \E\bigl[\smalll(\vartheta,X_{1,1})\bigr]= \inf_{ \theta \in \R } \E\bigl[\smalll(\theta,X_{1,1})\bigr] 
  \big\}
  = \{ 0 \}
$.
\end{enumerate}
\end{athm}

\begin{aproof}
\argument{\cref{lem: verify 3};the assumption that for all $n\in \N$ it holds that $\gamma_n=n^{-0.99}$} {that
\begin{equation}\llabel{eq1}
\limsup_{n\to\infty}( ( \gamma_n )^{ - 2 } ( \gamma_n - \gamma_{ n + 1 } )=0.
\end{equation}}
\argument{the assumption that for all $n\in \N$ it holds that $\gamma_n=n^{-0.99}$; the fact that for all $r\in (1,\infty)$ it holds that $\sum_{n=1}^\infty n^{-r}<\infty$;}{that for all $p \in (\frac{100}{99}, \infty)$ it holds that 
\begin{equation}\llabel{eq2}
\sum_{n=1}^\infty (\gamma_n)^p=\sum_{n=1}^\infty n^{-\frac{99p}{100}}<\infty.
\end{equation}}
\argument{\lref{eq1};\lref{eq2}}{\cref{item 1: numerical simulation: consequences 1}\dott}
 \startnewargseq
    \argument{the assumption that for all $\theta,x\in \R$ it holds that $\smalll(\theta,x)=|\theta-x|^2$; the fact that $\E[X_{1,1}]=0$}{that for all $\theta\in \R$ it holds that
    \begin{equation}\llabel{eq5}
    \E[\smalll(\theta,X_{1,1})]=\E\bigl[|\theta-X_{1,1}|^2\bigr]=\E[(X_{1,1})^2-2\theta X_{1,1}+\theta^2]=\theta^2+\E[(X_{1,1})^2].
    \end{equation}}
    \argument{\lref{eq5};}{\cref{item 2: numerical simulation: consequences 1}\dott}
\end{aproof}

\subsection{Two-point distributions}

In the following framework, \cref{setting: tpd} below, we describe the setup on which our numerical simulations in \cref{fig:Adam training trajectory,fig:cv_rate_beta2,fig:cv_rate_batch_size,fig:symmetry theorem} below are based.
\begin{setting}
\label{setting: tpd}

Assume \cref{setting: numerical simulation},
let $v \in (-\infty,0)$, $w\in (0,\infty)$,  
and assume for all $n,m\in\N$ that $ X_{ n, m } ( \Omega ) 
\subseteq \{ v, w \}$.
\end{setting}

In the following result, \cref{numerical simulation: consequences}  below, we briefly note two elementary consequences of \cref{setting: tpd}.
\begin{athm}{lemma}{numerical simulation: consequences}
Assume \cref{setting: tpd}.
Then
\begin{enumerate}[label=(\roman*)]
\item \label{numerical simulation: consequences: item 2} it holds for all $n,m\in \N$ that $\P(X_{n,m}=v) = \frac{w}{w-v}$ and
\item \label{numerical simulation: consequences: item 3} it holds for all $n,m\in \N$ that $\P(X_{n,m}=w) = \frac{v}{v-w}$ 
\end{enumerate}
\end{athm}
\begin{aproof}
\argument{the asumption that $\E[X_{1,1}]=0$;the fact that $X_{n,m}\in \{v,w\}$}{that
\begin{equation}\llabel{eq3}
v\P(X_{1,1}=v)+w\P(X_{1,1}=w)=0\qqandqq \P(X_{1,1}=v)+ \P(X_{1,1}=w)=1.
\end{equation}}
\argument{\lref{eq3};}{that
\begin{equation}\llabel{eq4}
\P(X_{1,1}=v)=\frac{w}{w-v}\qqandqq \P(X_{1,1}=w)=\frac{v}{v-w}.
\end{equation}}
\argument{\lref{eq4};}{items \ref{numerical simulation: consequences: item 2} and \ref{numerical simulation: consequences: item 3}\dott}
\end{aproof}

\subsubsection{Simulations for the convergence rate with respect to $n$}

In \cref{fig:Adam training trajectory} we consider \cref{setting: tpd} and approximately plot
$
  \E\big[ | \mandVcom{0}{M}{n}{i} | \big]
$
and 
$
  \E\big[ | \mandVcom{0}{M}{n}{i} - \vartheta| \big]
$
against 
$n \in \{1, 2, \ldots, \numGDsteps\}$
where 
$M = 1$, $\beta_1= 0.9$, $\beta_2 \in (\cup_{i = 0}^{3} \{ 1-2^{(-4 - 5i/3)} \})$, $\varepsilon=10^{-8}$, $v=-1$, and $w \in \{0.1,  1\}$ and where $\vartheta$ is a zero of the \Adam\ vector field $f$ (cf.\ \cref{def: Adam vector field numerical}).

\begin{figure}[htbp]
    \centering
    \begin{subfigure}[b]{0.4\textwidth}
        \centering
        \includegraphics[width=\textwidth]{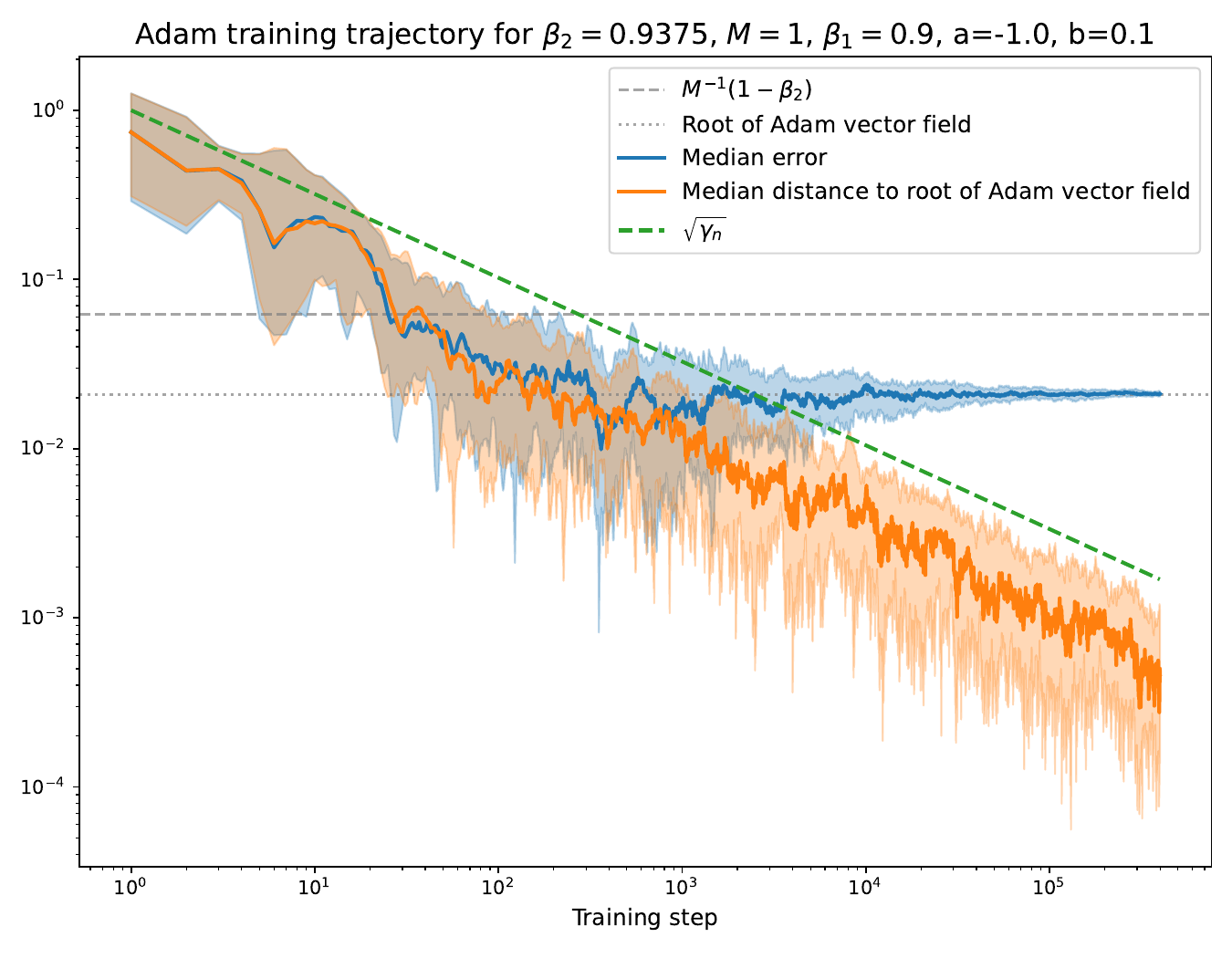}
        % \caption{Subfigure 1}
    \end{subfigure}
    % \hfill
    % \hspace{-1cm}
    \begin{subfigure}[b]{0.4\textwidth}
        \centering
        \includegraphics[width=\textwidth]{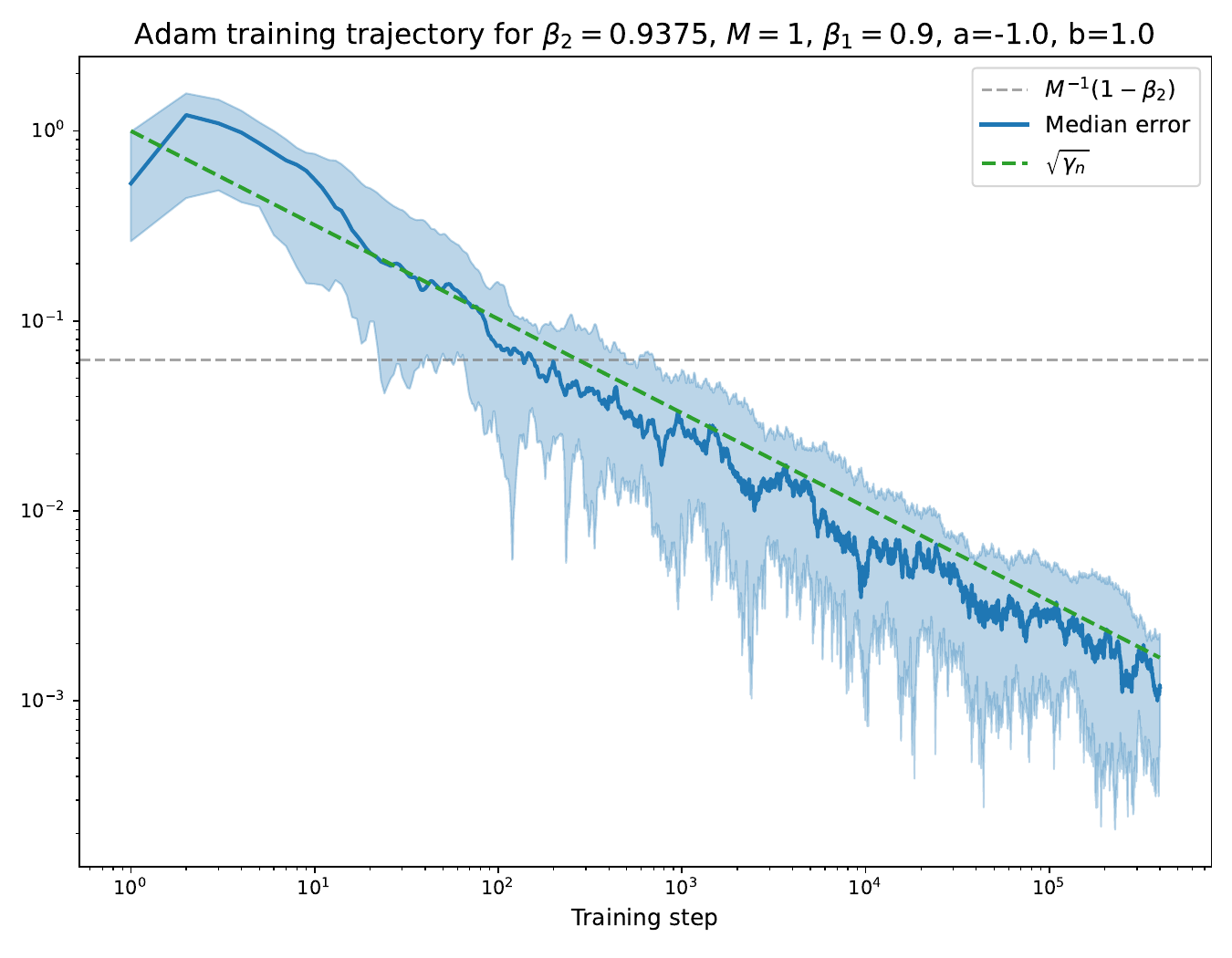}
        % \caption{Subfigure 2}
    \end{subfigure}
    
    % \vspace{-0.2cm}
    
    \begin{subfigure}[b]{0.4\textwidth}
        \centering
        \includegraphics[width=\textwidth]{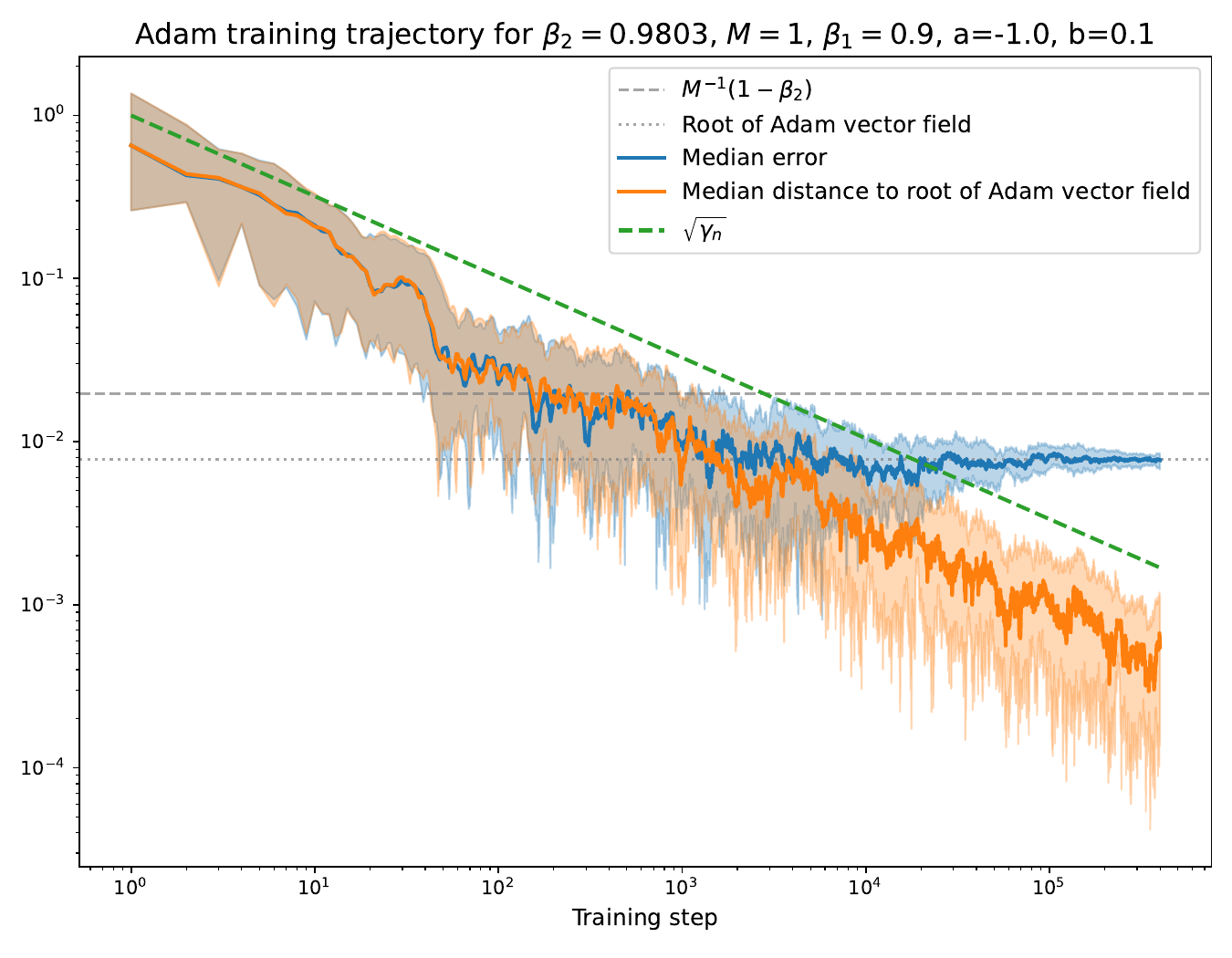}
        % \caption{Subfigure 1}
    \end{subfigure}
    % \hfill
    % \hspace{-1cm}
    \begin{subfigure}[b]{0.4\textwidth}
        \centering
        \includegraphics[width=\textwidth]{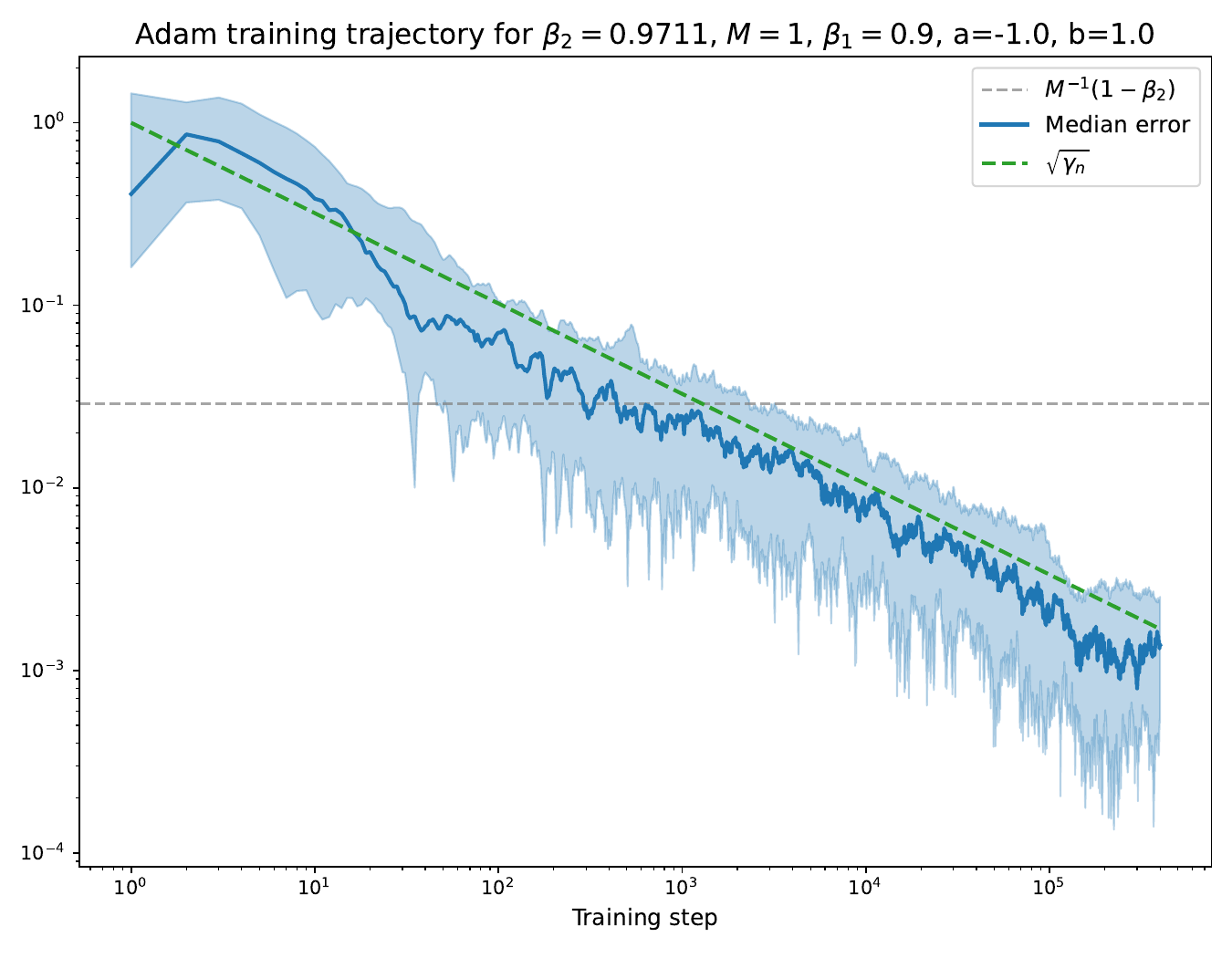}
        % \caption{Subfigure 2}
    \end{subfigure}
    
    % \vspace{-0.2cm}

    \begin{subfigure}[b]{0.4\textwidth}
        \centering
        \includegraphics[width=\textwidth]{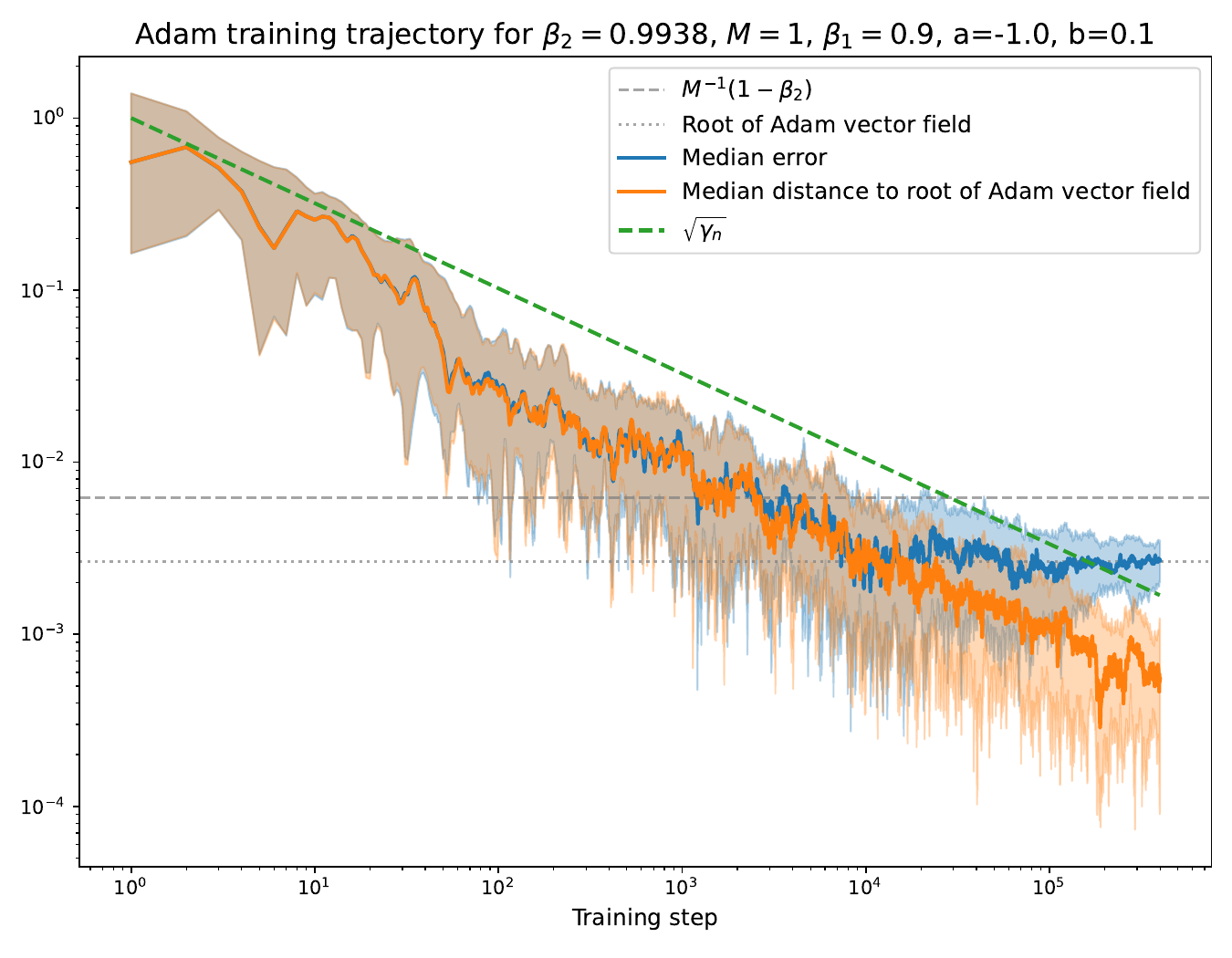}
        % \caption{Subfigure 1}
    \end{subfigure}
    % \hfill
    % \hspace{-1cm}
    \begin{subfigure}[b]{0.4\textwidth}
        \centering
        \includegraphics[width=\textwidth]{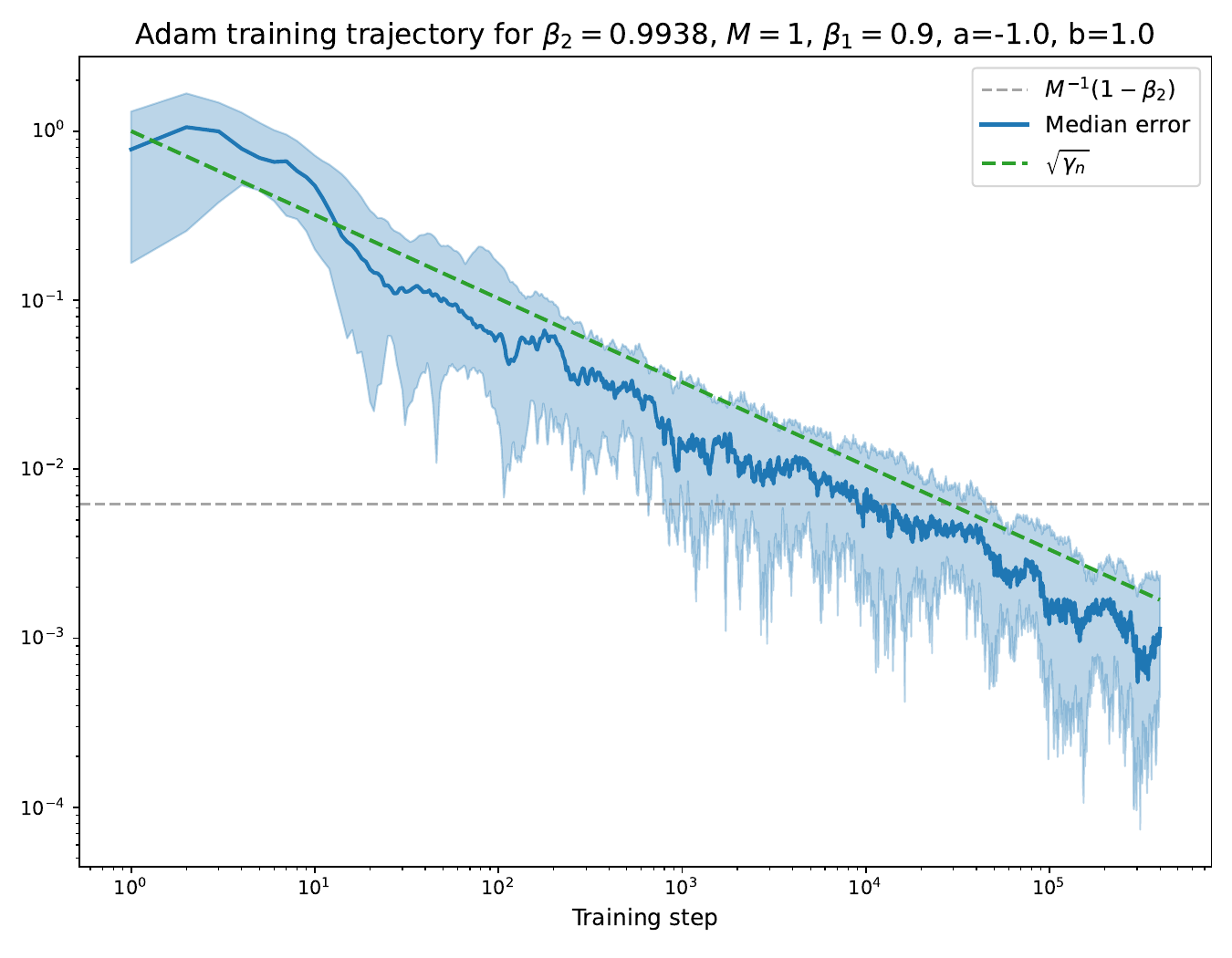}
        % \caption{Subfigure 2}
    \end{subfigure}
    
    % \vspace{-0.2cm}

    \begin{subfigure}[b]{0.4\textwidth}
        \centering
        \includegraphics[width=\textwidth]{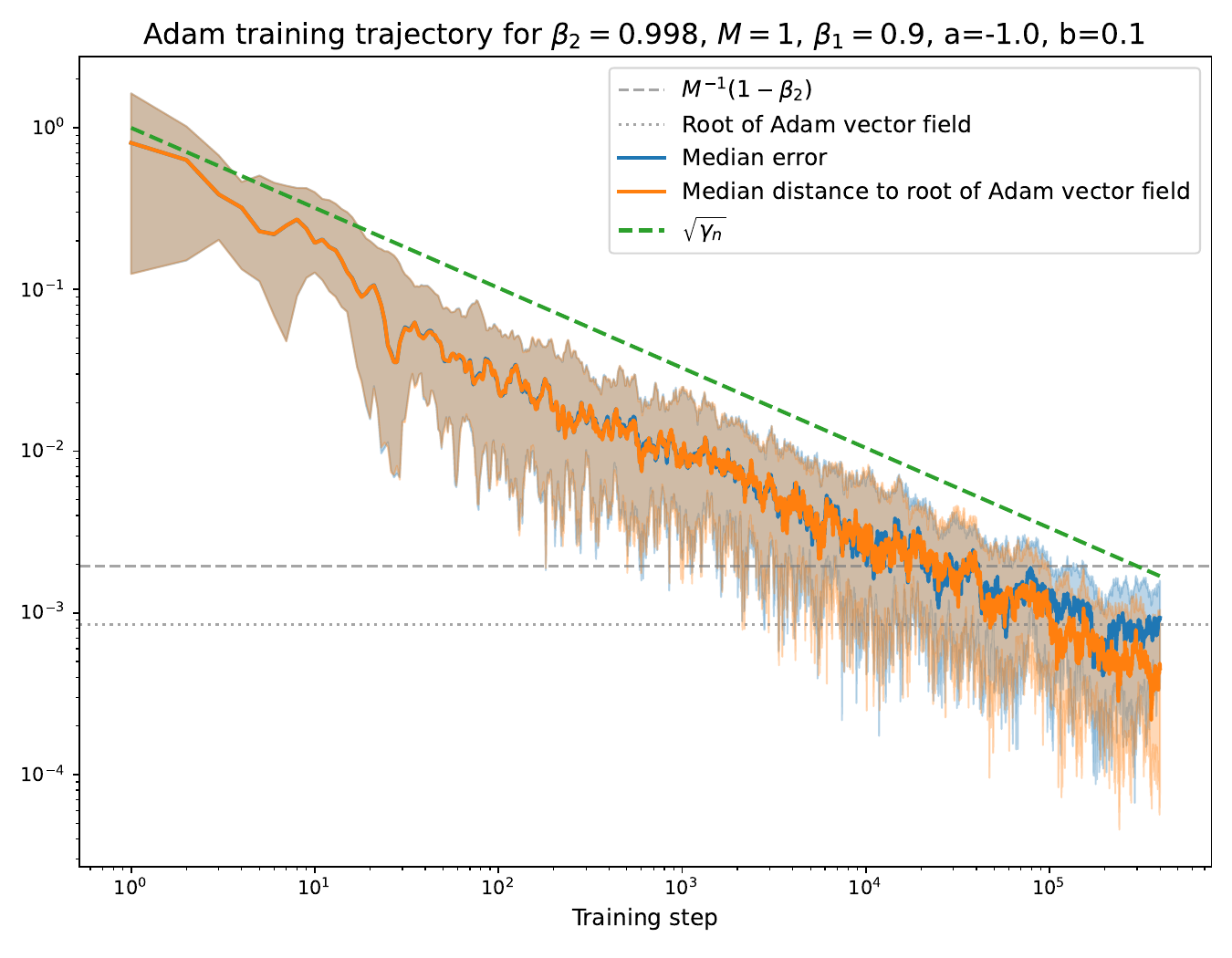}
        % \caption{Subfigure 1}
    \end{subfigure}
    % \hfill
    % \hspace{-1cm}
    \begin{subfigure}[b]{0.4\textwidth}
        \centering
        \includegraphics[width=\textwidth]{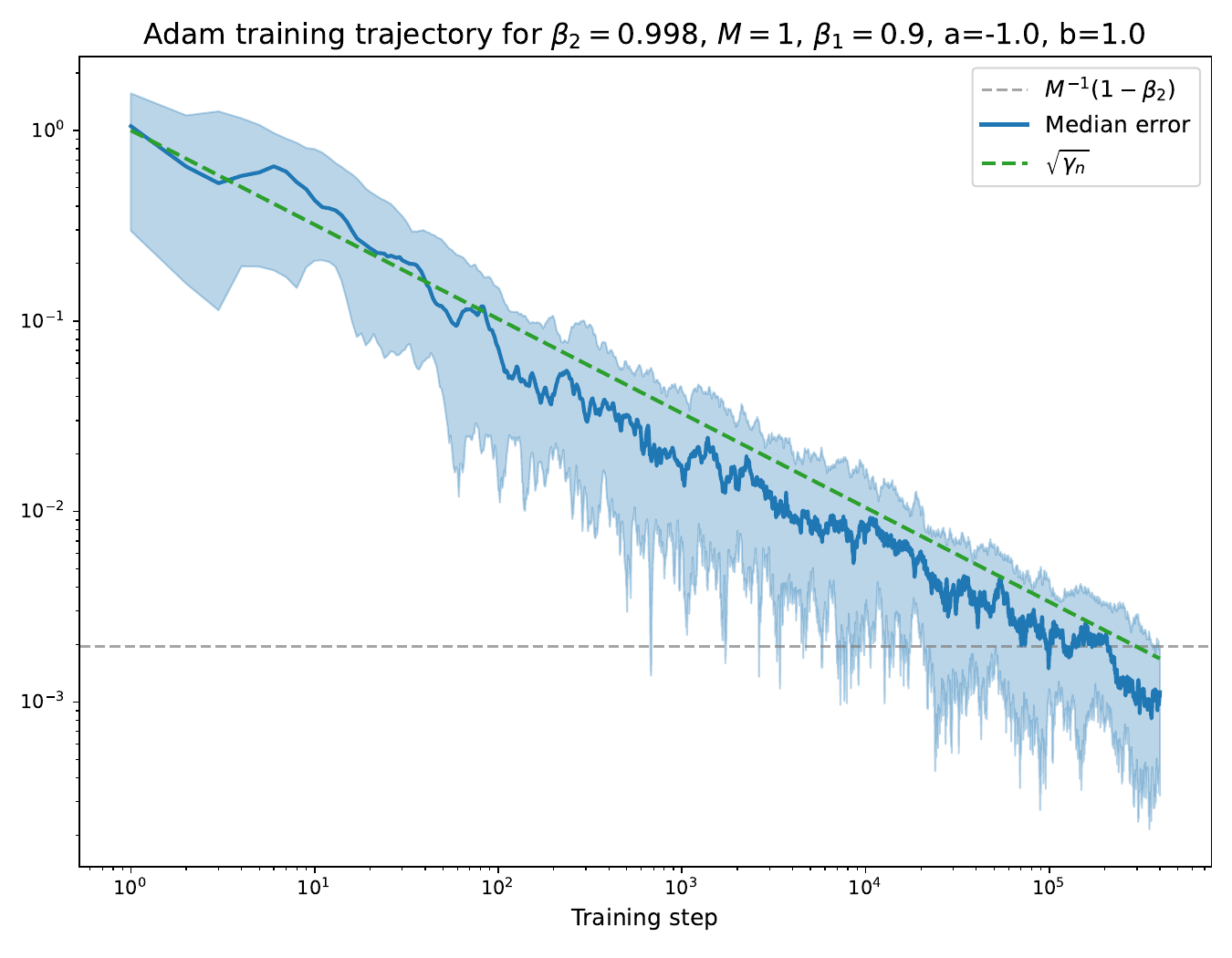}
        % \caption{Subfigure 2}
    \end{subfigure}

    \caption{
    Approximative plot of
    $
      \E\big[ | \mandVcom{0}{M}{n}{i} | \big]
    $
    and 
    $
      \E\big[ | \mandVcom{0}{M}{n}{i} - \vartheta| \big]
    $
    against  
    $n \in \{1, 2, \ldots, \numGDsteps\}$
    where  
    $M = 1$, $\beta_1= 0.9$, $\beta_2 \in (\cup_{i = 0}^{3} \{ 1-2^{(-4 - 5i/3)} \})$, $\varepsilon=10^{-8}$, $v=-1$, and $w \in \{0.1,  1\}$ 
    and where $\vartheta$ is a zero of the Adam vector field $f$ (cf.\ \cref{def: Adam vector field numerical,setting: tpd}).
    \emph{Left:} We consider the case of asymmetrically distributed data ($w=0.1$). \emph{Right:} We consider the case of symmetrically distributed data ($w=1$).
    }
    \label{fig:Adam training trajectory}
\end{figure}

\subsubsection{Simulations for the convergence rates with respect to $ M $ and $ \beta_2 $}

In \cref{fig:cv_rate_beta2} we consider \cref{setting: tpd} and approximately plot
$
  \E\big[ | \mandVcom{0}{M}{\numGDsteps}{i} | \big]
$
and
a zero of the Adam vector field $f$ (cf.\ \cref{def: Adam vector field numerical})
against
$\beta_2 \in ( \cup_{i = 0}^{9} \{ 1-2^{(-4 - 5i/9)} \} )$
where $M\in \{1, 2\}$, $\beta_1\in \{0.6,0.9\}$, $\varepsilon=10^{-8}$, $v=-1$, and $w \in \{0.1, 1\}$.
Furthermore, in \cref{fig:cv_rate_batch_size} we consider \cref{setting: tpd} and approximately plot
$
  \E\big[ | \mandVcom{0}{M}{\numGDsteps}{i} | \big]
$
and
a zero of the Adam vector field $f$ (cf.\ \cref{def: Adam vector field numerical})
against
$M \in \{1, 2, \ldots, 6\}$
where $\beta_1 \in \{0.6, 0.9\}$, $\beta_2 \in \{0.95, 0.99\}$, $\varepsilon=10^{-8}$, $v=-1$, and $w = 0.1$.

\begin{figure}[htbp]
    \centering
    \begin{subfigure}[b]{0.49\textwidth}
        \centering
        \includegraphics[width=\textwidth]{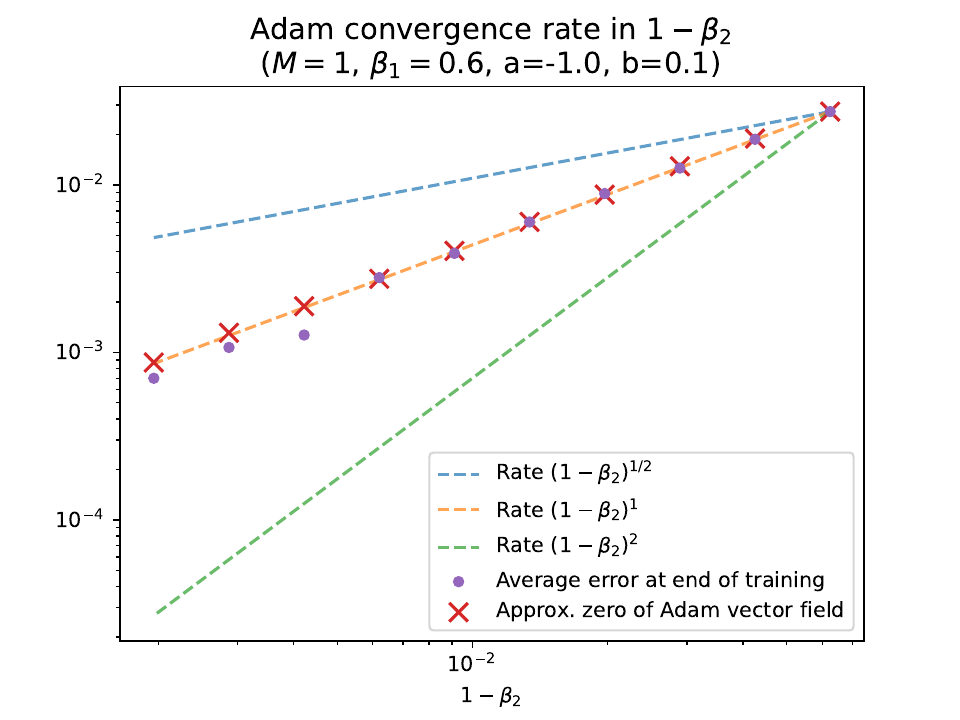}
        % \caption{Subfigure 1}
        \label{fig:sub1}
    \end{subfigure}
    % \hfill
    \hspace{-1cm}
    \begin{subfigure}[b]{0.49\textwidth}
        \centering
        \includegraphics[width=\textwidth]{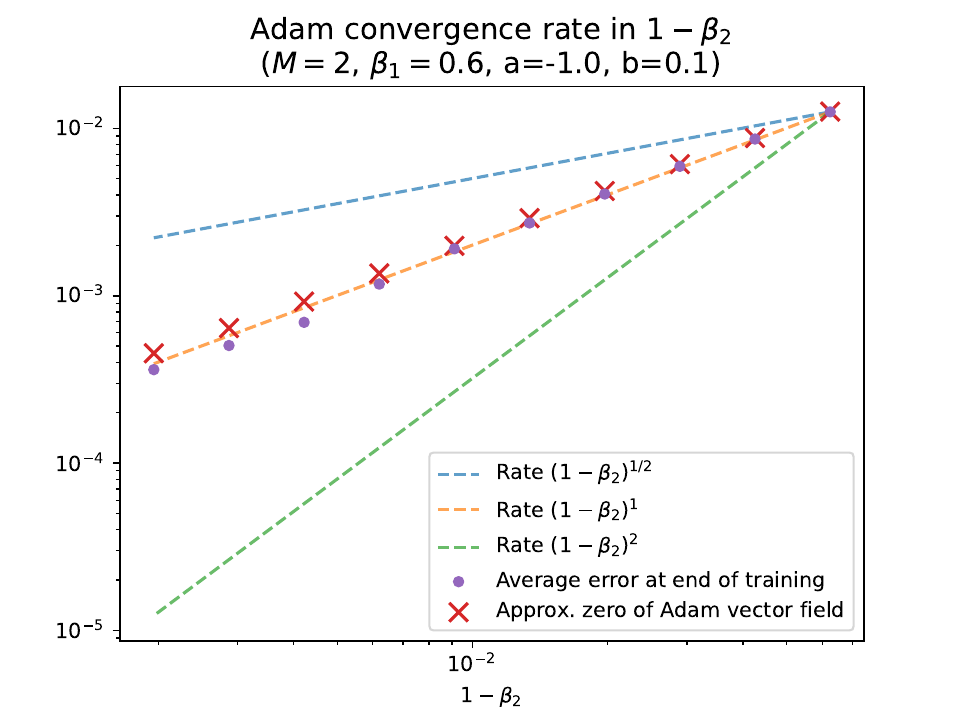}
        % \caption{Subfigure 2}
        \label{fig:sub2}
    \end{subfigure}
    
    \vspace{-0.5cm}
    
    \begin{subfigure}[b]{0.49\textwidth}
        \centering
        \includegraphics[width=\textwidth]{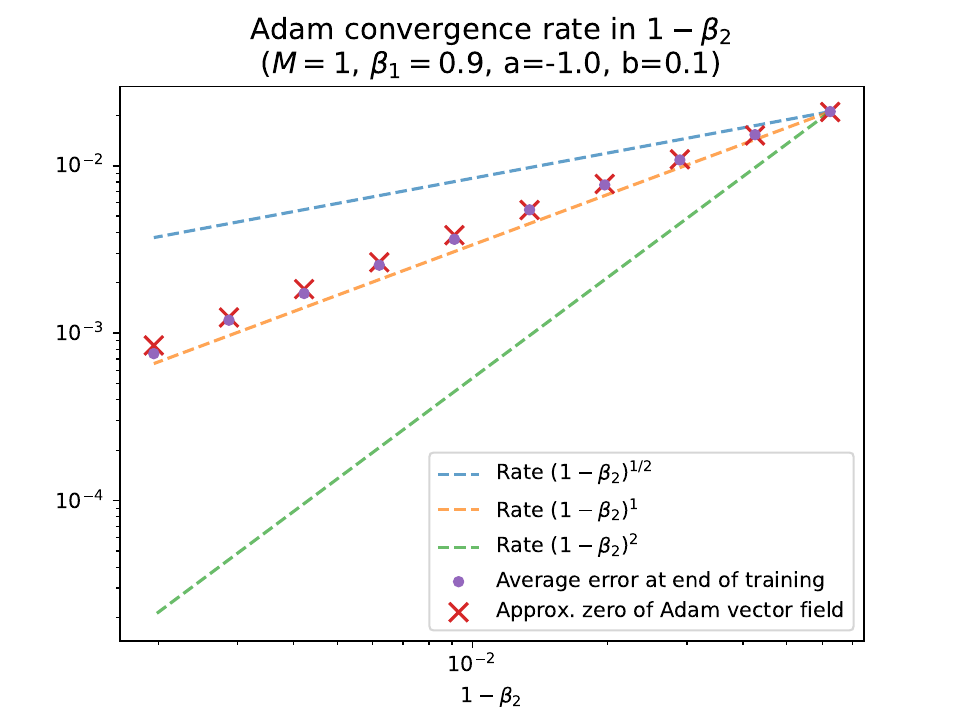}
        % \caption{Subfigure 3}
        \label{fig:sub3}
    \end{subfigure}
    % \hfill
    \hspace{-1cm}
    \begin{subfigure}[b]{0.49\textwidth}
        \centering
        \includegraphics[width=\textwidth]{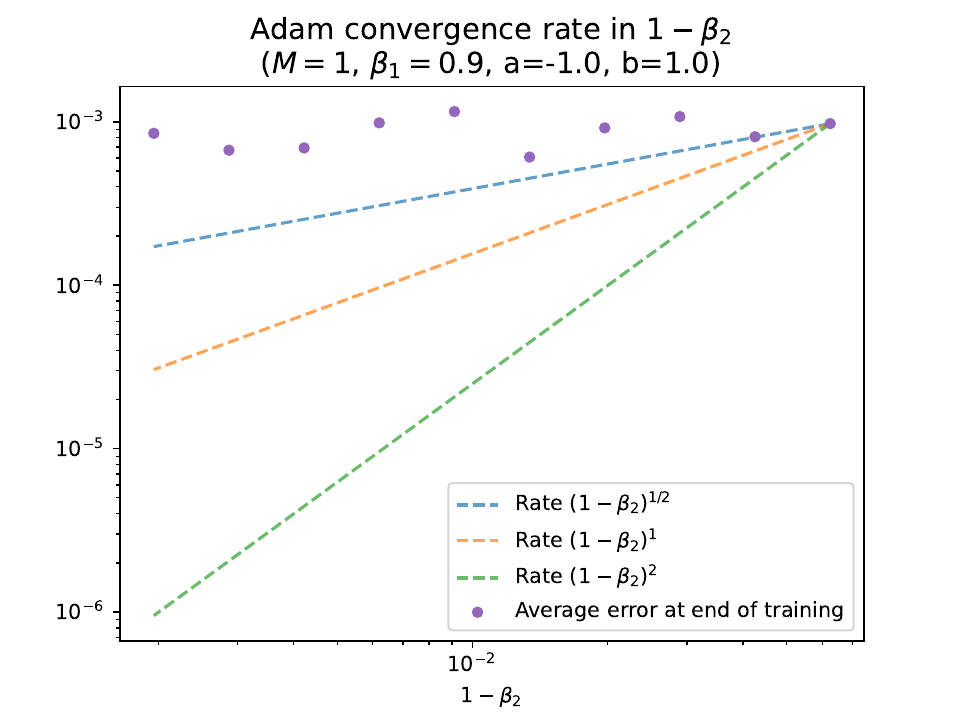}
        % \caption{Subfigure 4}
        \label{fig:sub4}
    \end{subfigure}
    
    \caption{
    Approximative plot 
    of
    $
      \E\big[ | \mandVcom{0}{M}{\numGDsteps}{i} | \big]
    $
    and
    a zero of the Adam vector field $f$ (cf.\ \cref{def: Adam vector field numerical})
    against
    $\beta_2 \in ( \cup_{i = 0}^{9} \{ 1-2^{(-4 - 5i/9)} \} )$
    where $M\in \{1, 2\}$, $\beta_1\in \{0.6,0.9\}$, $\varepsilon=10^{-8}$, $v=-1$, and $w \in \{0.1, 1\}$
    (cf.\ \cref{setting: tpd}).
    }
    \label{fig:cv_rate_beta2}
\end{figure}

\begin{figure}[htbp]
    \centering
    \begin{subfigure}[b]{0.49\textwidth}
        \centering
        \includegraphics[width=\textwidth]{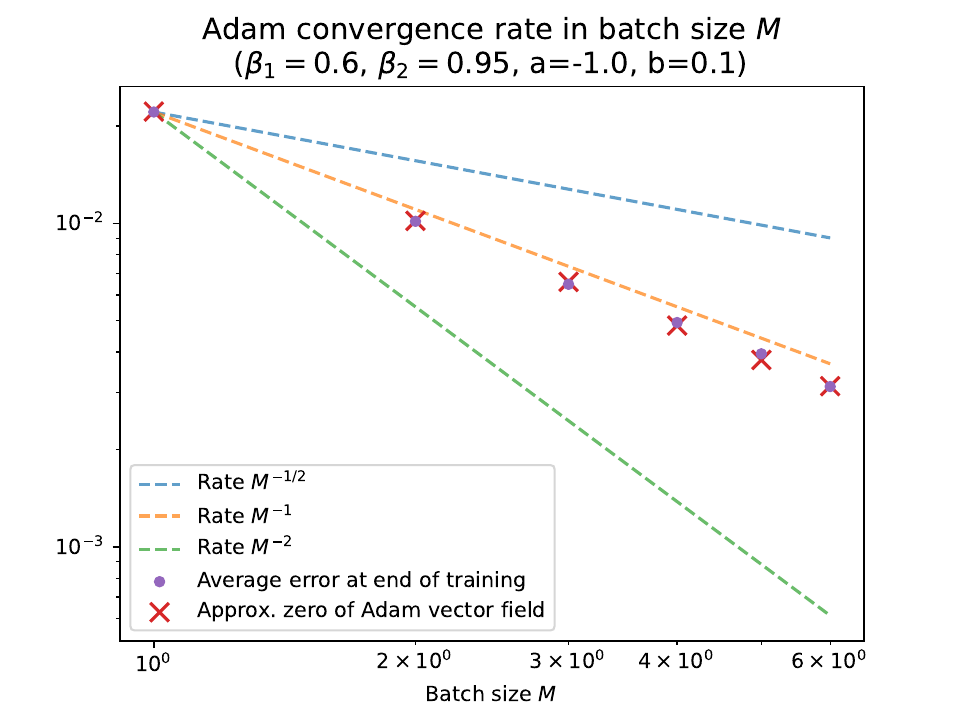}
        % \caption{Subfigure 1}
        % \label{fig:sub1}
    \end{subfigure}
    % \hfill
    \hspace{-1cm}
    \begin{subfigure}[b]{0.49\textwidth}
        \centering
        \includegraphics[width=\textwidth]{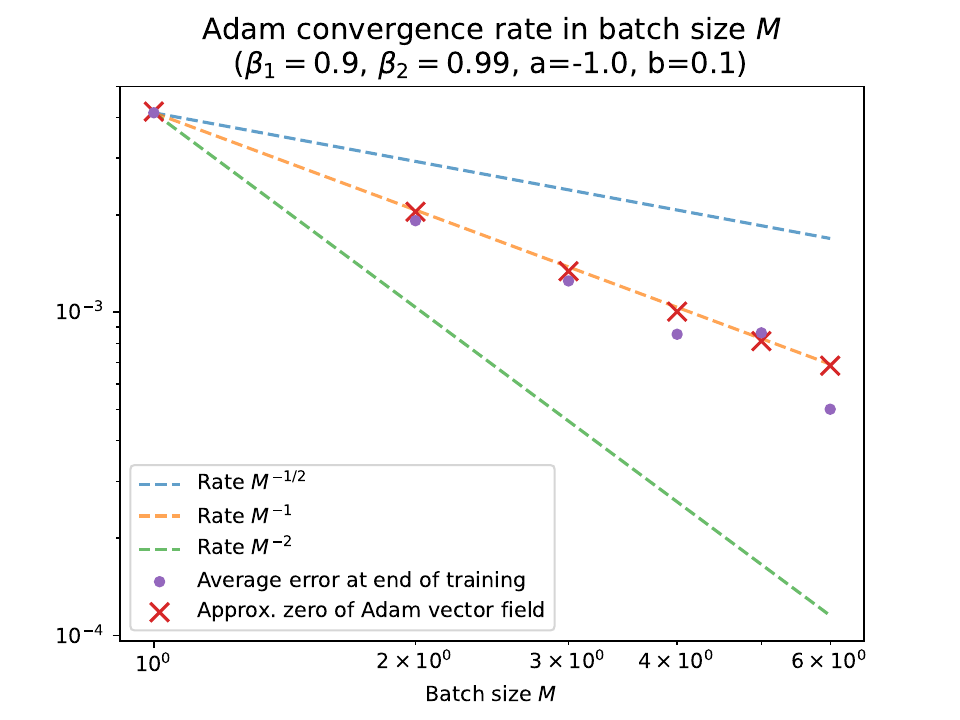}
        % \caption{Subfigure 2}
        % \label{fig:sub2}
    \end{subfigure}
    \caption{
    Approximative plot 
    of
    $
      \E\big[ | \mandVcom{0}{M}{\numGDsteps}{i} | \big]
    $
    and
    a zero of the Adam vector field $f$ (cf.\ \cref{def: Adam vector field numerical})
    against
    $M \in \{1, 2, \ldots, 6\}$
    where $\beta_1 \in \{0.6, 0.9\}$, $\beta_2 \in \{0.95, 0.99\}$, $\varepsilon=10^{-8}$, $v=-1$, and $w = 0.1$
    (cf.\ \cref{setting: tpd}).
    }
    \label{fig:cv_rate_batch_size}
\end{figure}

\subsubsection{Simulations for the Adam symmetry theorem}

In \cref{fig:symmetry theorem} we consider \cref{setting: tpd} and approximately plot
$
  \E\big[ \mandVcom{0}{M}{\numGDsteps}{i} \big]
$
and
a zero of the Adam vector field $f$ (cf.\ \cref{def: Adam vector field numerical})
against
$w \in (\cup_{i = -8}^{7} \{ 2^{i} \})$
where $M=1$, $\beta_1=0.9$, $\beta_2 = 0.9$, $\varepsilon=10^{-8}$, and $v=-1$.

\begin{figure}[htbp]
    \centering
    \begin{subfigure}[b]{0.49\textwidth}
        \centering
        \includegraphics[width=\textwidth]{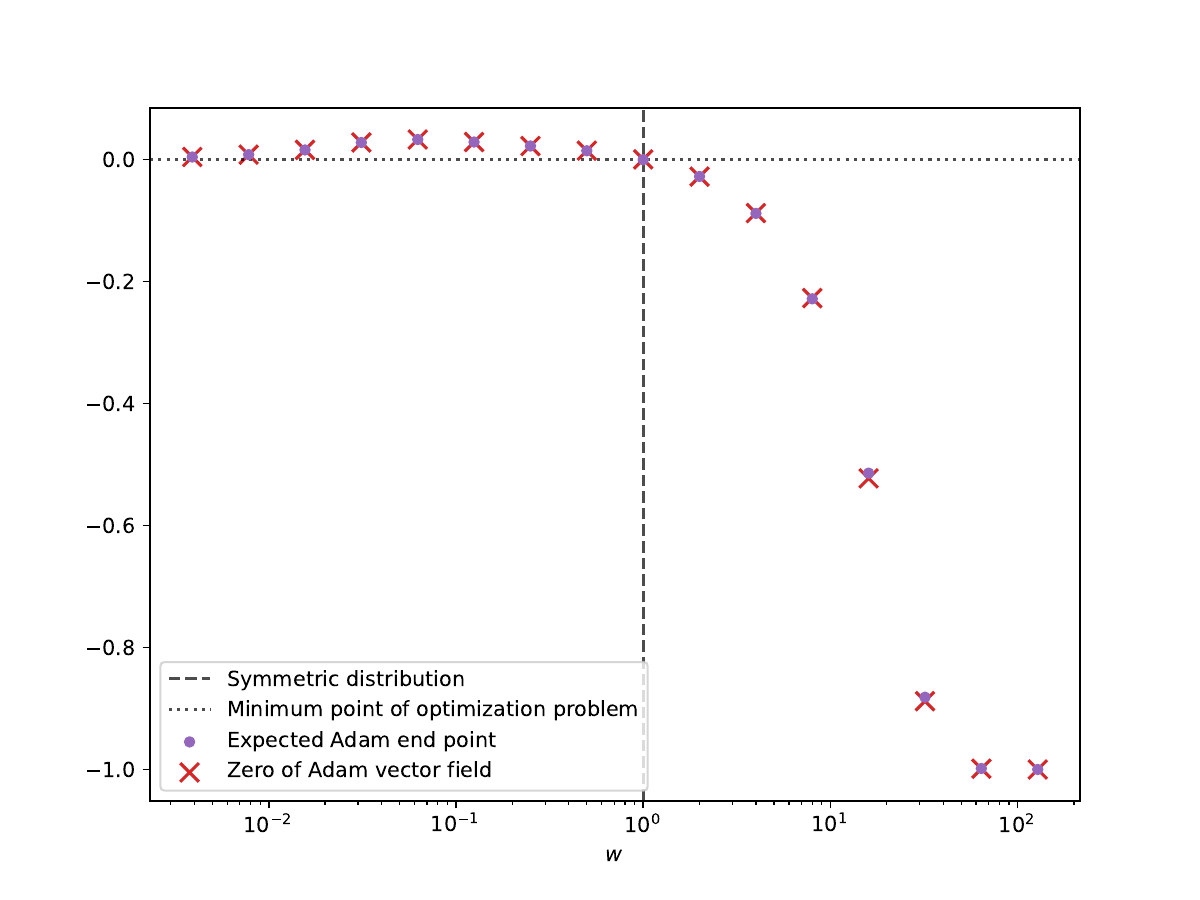}
        % \caption{Subfigure 1}
      
    \end{subfigure}
    % \hfill
    \hspace{-1cm}
    \begin{subfigure}[b]{0.49\textwidth}
        \centering
        \includegraphics[width=\textwidth]{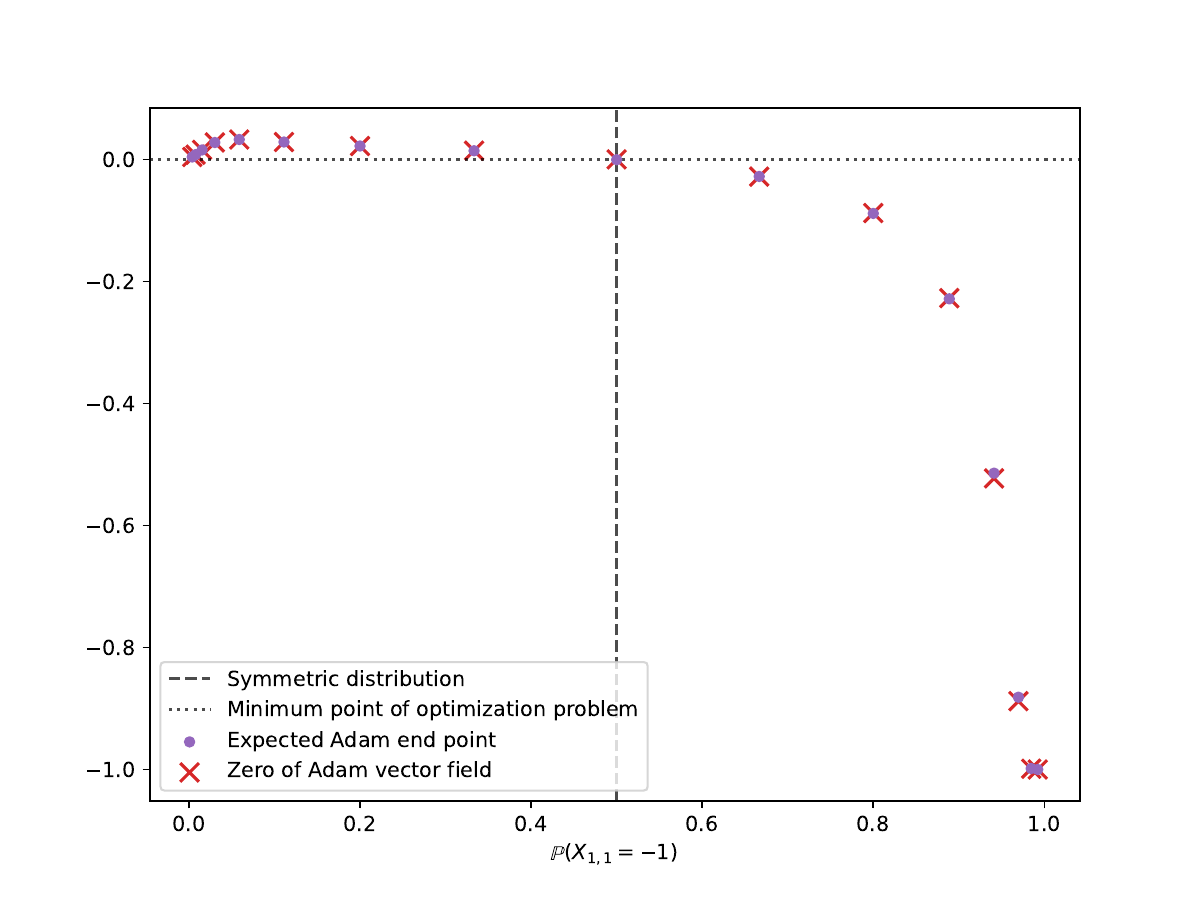}
        % \caption{Subfigure 2}

    \end{subfigure}
    \caption{
    Approximative plot 
    of
    $
      \E\big[ \mandVcom{0}{M}{\numGDsteps}{i} \big]
    $
    and
    a zero of the Adam vector field $f$ (cf.\ \cref{def: Adam vector field numerical})
    against
    $w \in (\cup_{i = -8}^{7} \{ 2^{i} \})$
    where $M=1$, $\beta_1=0.9$, $\beta_2 = 0.9$, $\varepsilon=10^{-8}$, and $v=-1$
    (cf.\ \cref{setting: tpd}).
    \emph{Left:} We measure the asymmetry of the distribution on the $x$-axis using $w \in (0,\infty)$. 
    \emph{Right:} We measure the asymmetry of the distribution on the $x$-axis using $\P(X_{1,1}=-1) = \frac{w}{1+w} \in (0,1)$.
    }
    \label{fig:symmetry theorem}
\end{figure}

\subsection{Skew normal distributions}

\newcommand{\normalpdf}{\phi}
\newcommand{\normalcdf}{\Phi}
\newcommand{\skewnormalpdf}{\varphi}
\begin{minipage}{\linewidth}
\begin{athm}{lemma}{skew normal distribution}
Let $\alpha, \xi \in \R$, $\omega \in (0, \infty)$
and
let $\normalpdf \colon \R \to [0,\infty)$, $\normalcdf \colon \R \to [0,1]$, and $\skewnormalpdf \colon \R \to [0,\infty)$ satisfy for all 
$x\in \R$ that
\begin{equation}\llabel{def: normalpdf}
  \normalpdf(x) = \frac{1}{\sqrt{2\pi}} \exp\!\left(-\frac{x^2}{2}\right),
\qquad
  \normalcdf(x) = \int_{-\infty}^x \normalpdf(t) \, \d t,
\end{equation}
\begin{equation}
\label{T_B_D}
\begin{split} 
\andq
  \skewnormalpdf(x) 
= 
  \frac{2}{\omega} \normalpdf \!\left( \frac{x-\xi}{\omega} \right) \normalcdf \!\left( \frac{\alpha (x-\xi)}{\omega} \right).
\end{split}
\end{equation}
Then
\begin{enumerate}[label=(\roman*)] 
\item \label{item 1: skew normal distribution} it holds that
\begin{equation}\llabel{eq1}
  \int_{-\infty}^\infty \skewnormalpdf(x) \, \d x = 1
\end{equation}
and
\item \label{item 2: skew normal distribution} it holds that
\begin{equation}\llabel{eq2}
  \int_{-\infty}^\infty x \skewnormalpdf(x) \, \d x = \xi + \omega \alpha \sqrt{\frac{2}{\pi(1+\alpha^2)}}.
\end{equation}
\end{enumerate}
\end{athm}
\end{minipage}

\begin{aproof}
\argument{\cite[Proposition 1.1]{Azzalini_2013}}{\cref{item 1: skew normal distribution}\dott}
\argument{\cite[(2.5) and (2.6)]{Azzalini_2013}}{\cref{item 2: skew normal distribution}\dott}
\end{aproof}

In the following framework, \cref{setting: skew normal} below, we describe the setup on which our numerical simulations in \cref{fig:skew 3d plot 1.1,fig:skew 3d plot 1.2,fig:skew 3d plot 1.3,fig:skew 3d plot 2.1,fig:skew 3d plot 2.2,fig:skew 3d plot 2.3,fig:Skew Adam training trajectory,fig:skew_cv_rate_beta2} below are based.

\begin{samepage}
\begin{setting}
\label{setting: skew normal}
Assume \cref{setting: numerical simulation}, 
let $\alpha, \xi \in \R$, $\omega \in (0, \infty)$ satisfy
$\xi = -\omega \alpha \sqrt{\frac{2}{\pi(1+\alpha^2)}}$,
let $\skewnormalpdf \colon \R \to [0,\infty)$ satisfy for all $x\in \R$ that
\begin{equation}\llabel{def: skew normal pdf}
  \skewnormalpdf(x) 
= 
  \frac{2}{\omega} \normalpdf \!\left( \frac{x-\xi}{\omega} \right) \normalcdf \!\left( \frac{\alpha (x-\xi)}{\omega} \right),
\end{equation}
and assume for all $n,m\in\N$, $B \in \mathcal{B}(\R)$ that
\begin{equation}
\label{T_B_D 2}
\begin{split} 
  \P(X_{n,m} \in B) = \int_B \skewnormalpdf(x) \, \d x
\end{split}
\end{equation}
(cf.\ \cref{skew normal distribution}).
\end{setting}
\end{samepage}

In \cref{fig:skew normal densities} we consider \cref{setting: skew normal} and plot $\skewnormalpdf$ where $\alpha \in \{-6, 4, \ldots, 6\}$ and $\omega = 3$.

\begin{figure}[htbp]
    \centering
    \begin{subfigure}[b]{0.49\textwidth}
        \centering
        \includegraphics[width=\textwidth]{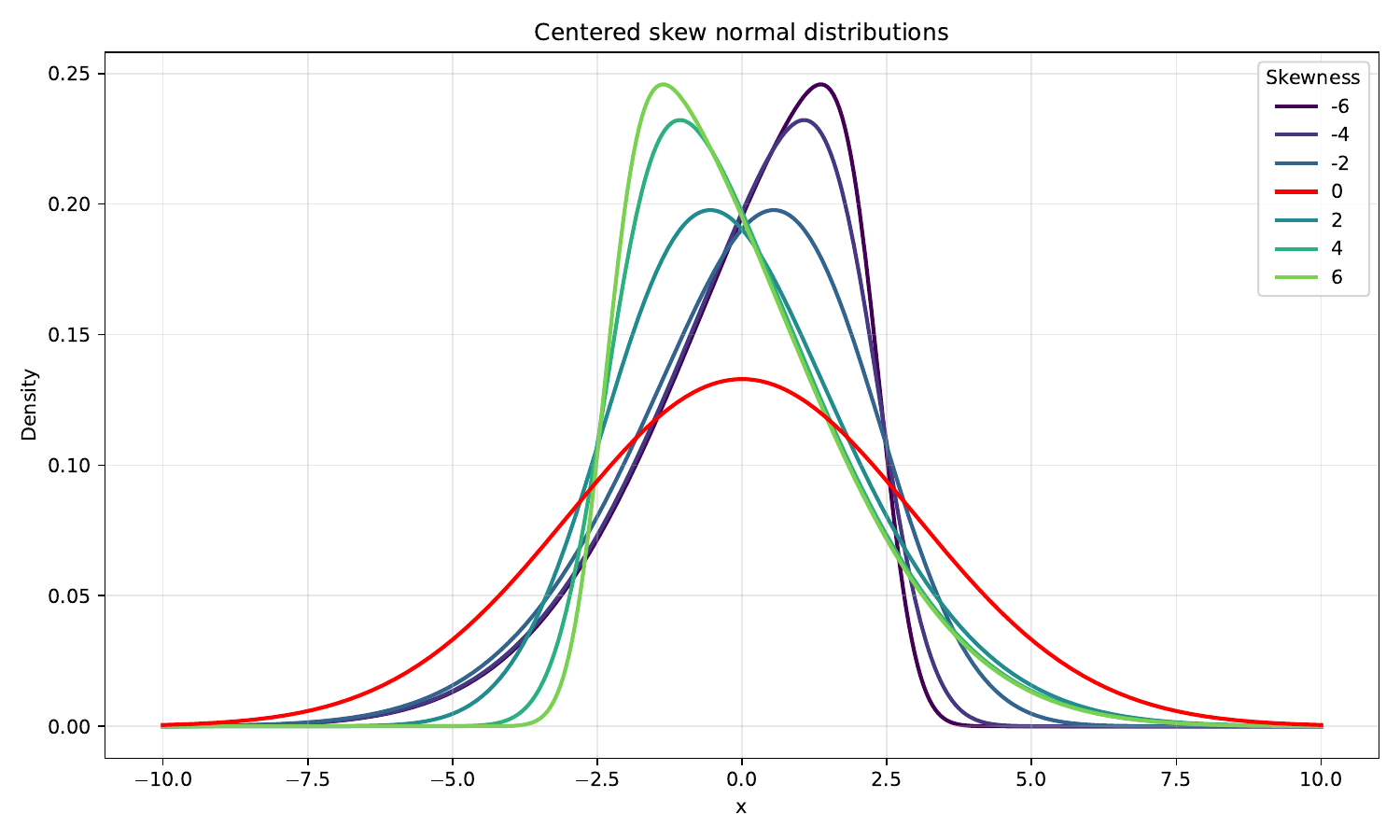}
        % \caption{Subfigure 1}
        % \label{fig:sub1}
    \end{subfigure}
    % \hfill
    % \hspace{-1cm}
    \begin{subfigure}[b]{0.49\textwidth}
        \centering
        \includegraphics[width=\textwidth]{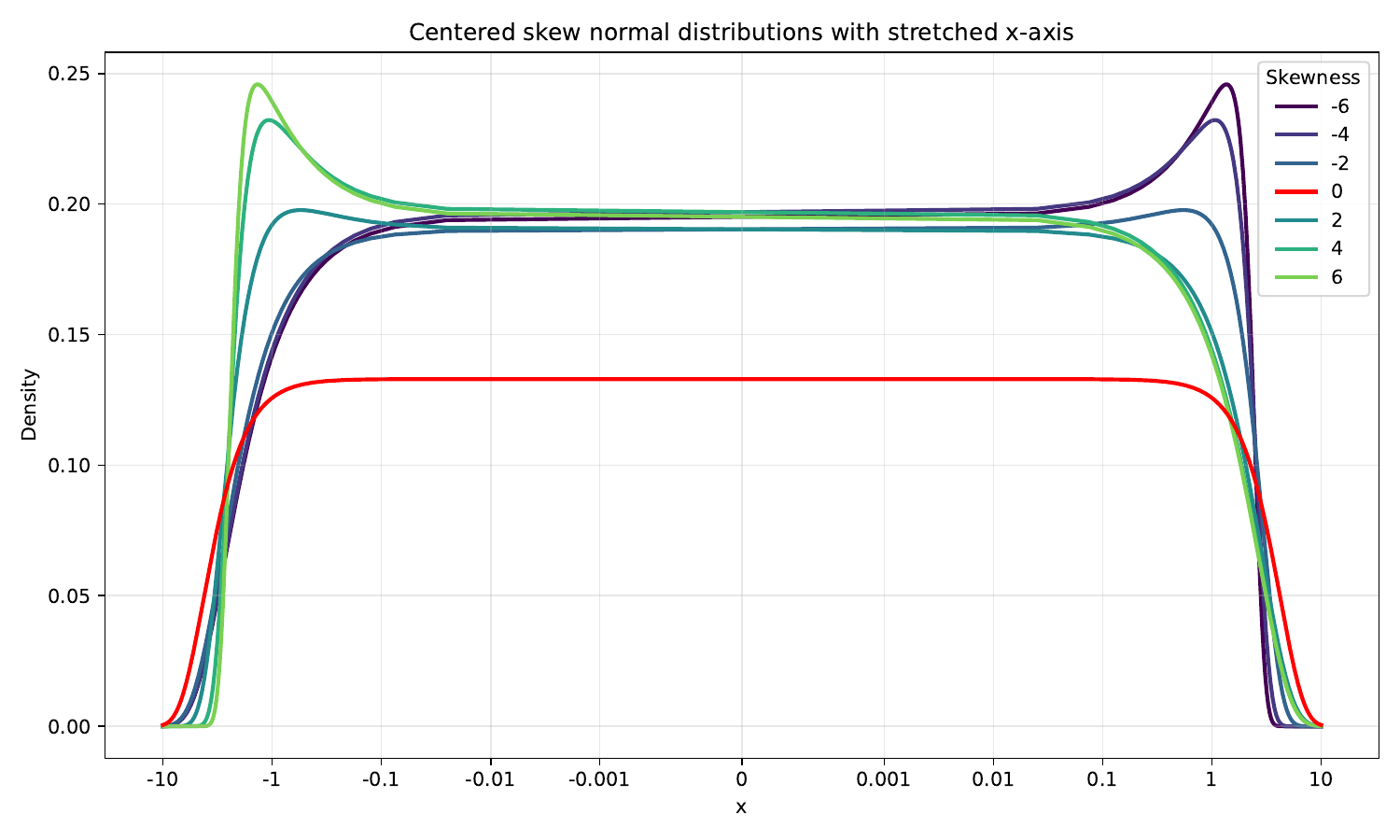}
        % \caption{Subfigure 2}
        % \label{fig:sub2}
    \end{subfigure}
    \caption{
    Plot of $\skewnormalpdf$ where $\alpha \in \{-6, 4, \ldots, 6\}$ and $\omega = 3$ (cf.\ \cref{setting: skew normal}).
    }
    \label{fig:skew normal densities}
\end{figure}

\subsubsection{Simulations for the convergence rates with respect to $n$}

In \cref{fig:Skew Adam training trajectory} we consider \cref{setting: skew normal} and approximately plot
$
  \E\big[ | \mandVcom{0}{M}{n}{i} | \big]
$
and 
$
  \E\big[ | \mandVcom{0}{M}{n}{i} - \vartheta| \big]
$
against 
$n \in \{1, 2, \ldots, \numGDsteps\}$
where 
$M = 1$, $\beta_1= 0.9$, $\beta_2 \in (\cup_{i = 0}^{3} \{ 1-2^{(-4 - 5i/3)} \})$, $\varepsilon=10^{-8}$, $\alpha \in \{-100, 0\}$, and $\omega = 1$ and where $\vartheta$ is a zero of the \Adam\ vector field $f$ (cf.\ \cref{def: Adam vector field numerical}).

\begin{figure}[htbp]
    \centering
    \begin{subfigure}[b]{0.4\textwidth}
        \centering
        \includegraphics[width=\textwidth]{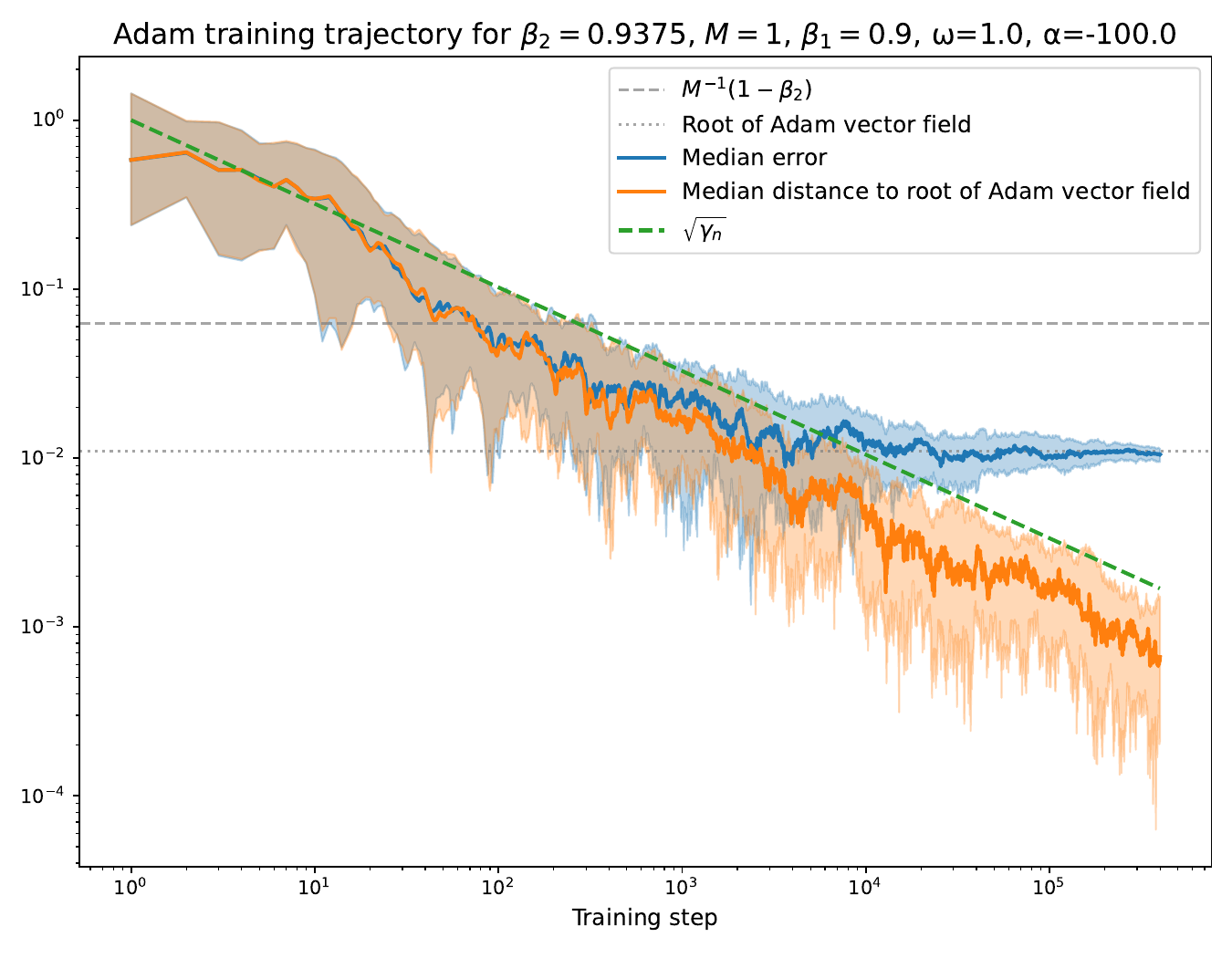}
        % \caption{Subfigure 1}
    \end{subfigure}
    % \hfill
    % \hspace{-1cm}
    \begin{subfigure}[b]{0.4\textwidth}
        \centering
        \includegraphics[width=\textwidth]{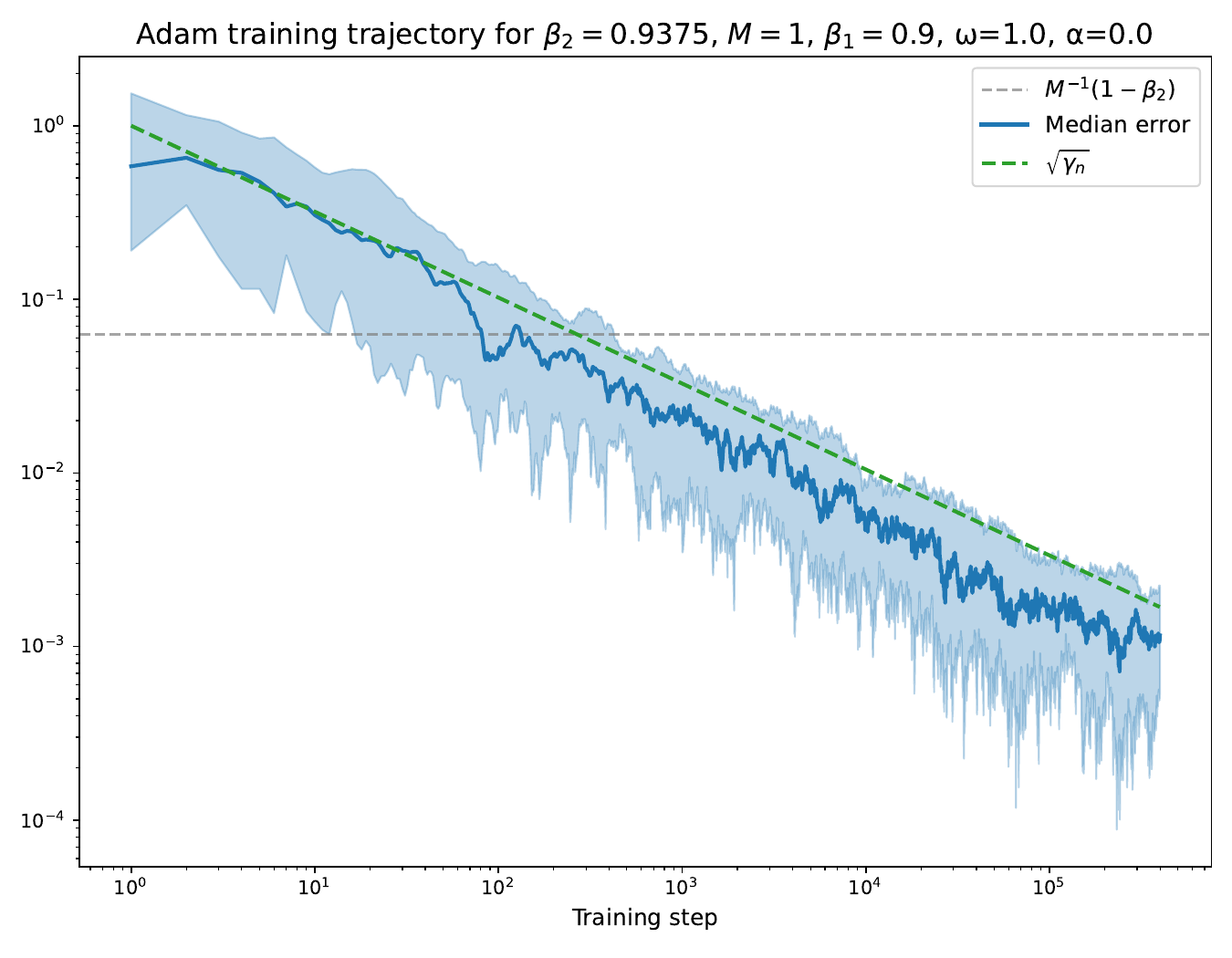}
        % \caption{Subfigure 2}
    \end{subfigure}
    
    % \vspace{-0.5cm}
    
    \begin{subfigure}[b]{0.4\textwidth}
        \centering
        \includegraphics[width=\textwidth]{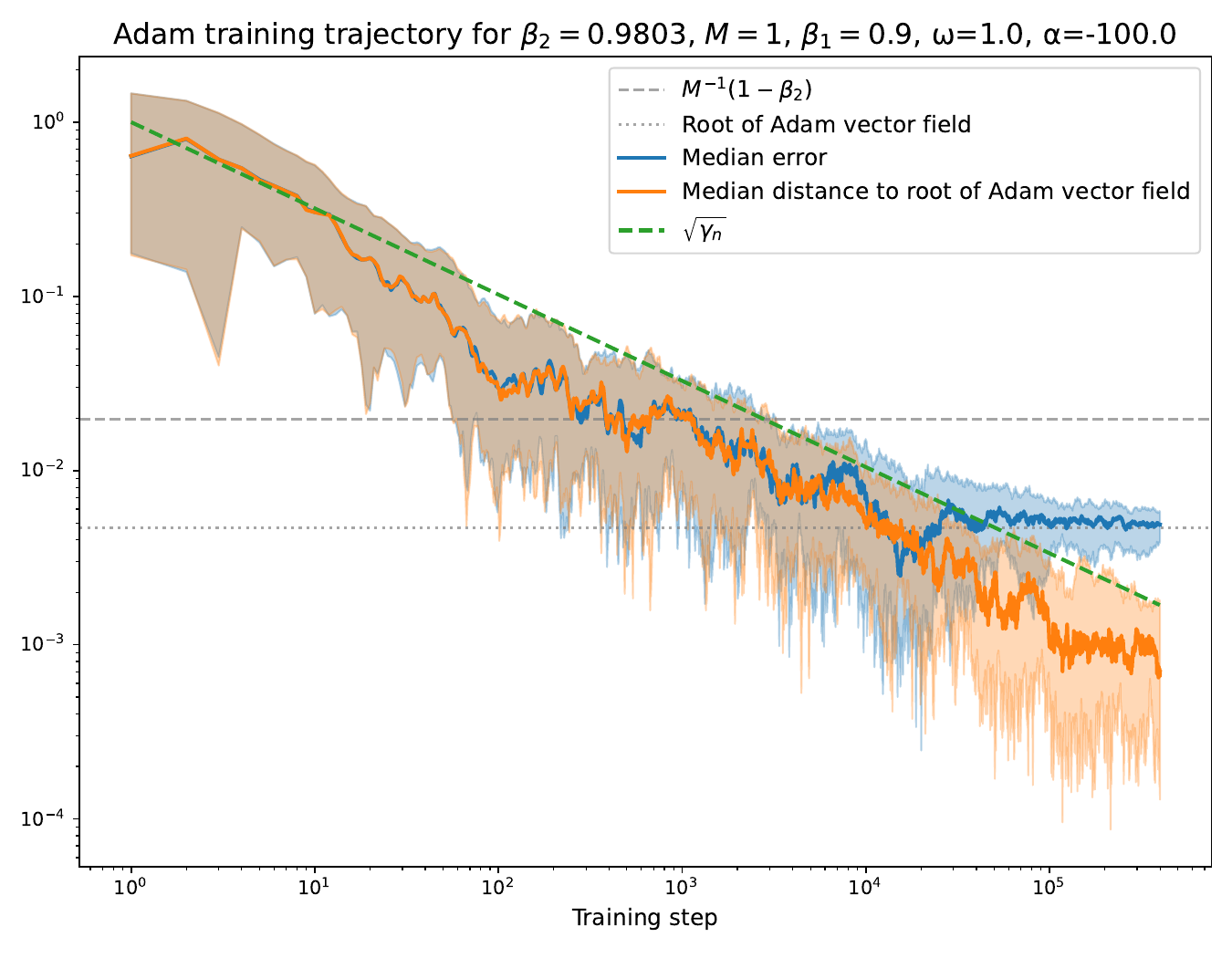}
        % \caption{Subfigure 1}
    \end{subfigure}
    % \hfill
    % \hspace{-1cm}
    \begin{subfigure}[b]{0.4\textwidth}
        \centering
        \includegraphics[width=\textwidth]{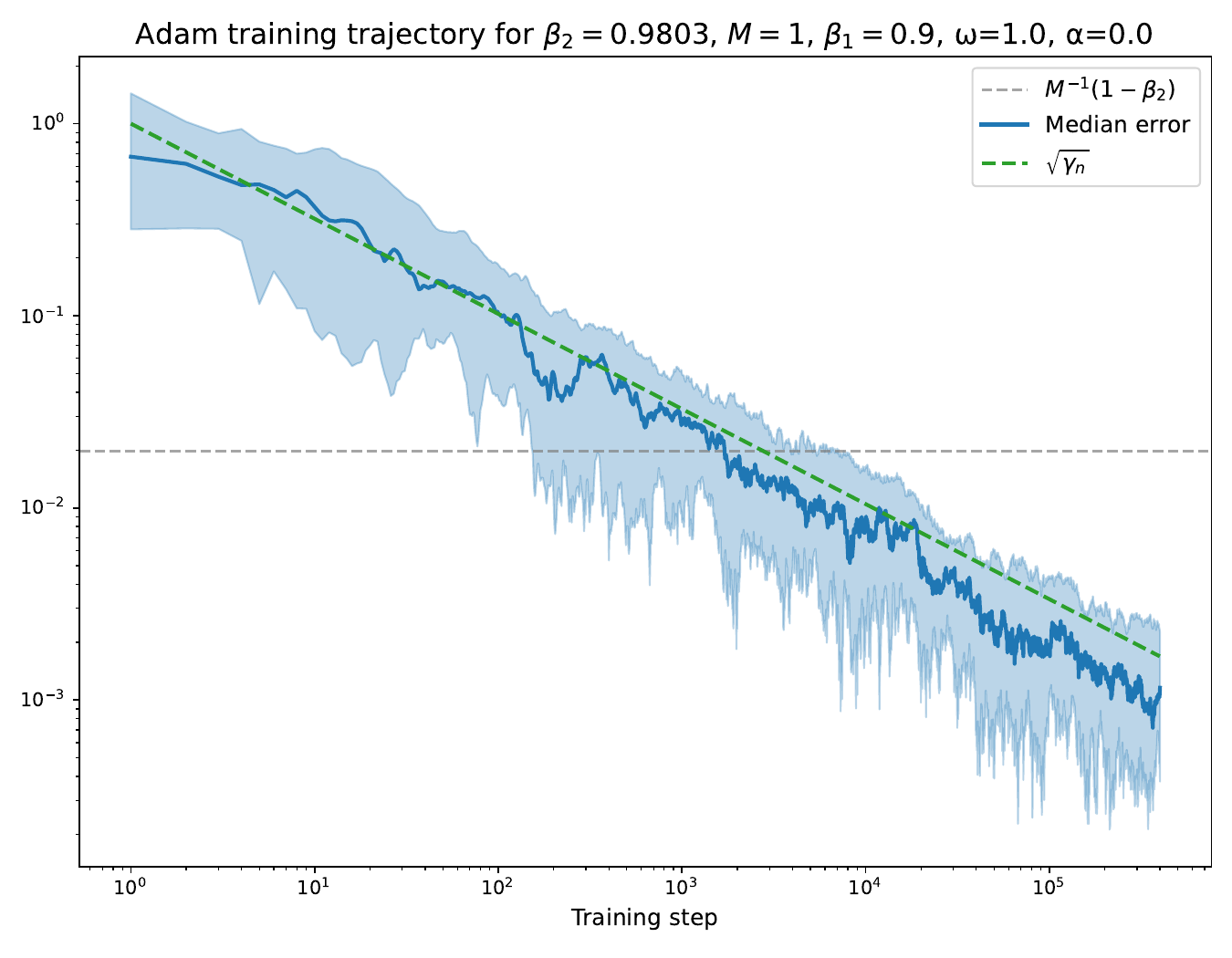}
        % \caption{Subfigure 2}
    \end{subfigure}
    
    % \vspace{-0.5cm}

    \begin{subfigure}[b]{0.4\textwidth}
        \centering
        \includegraphics[width=\textwidth]{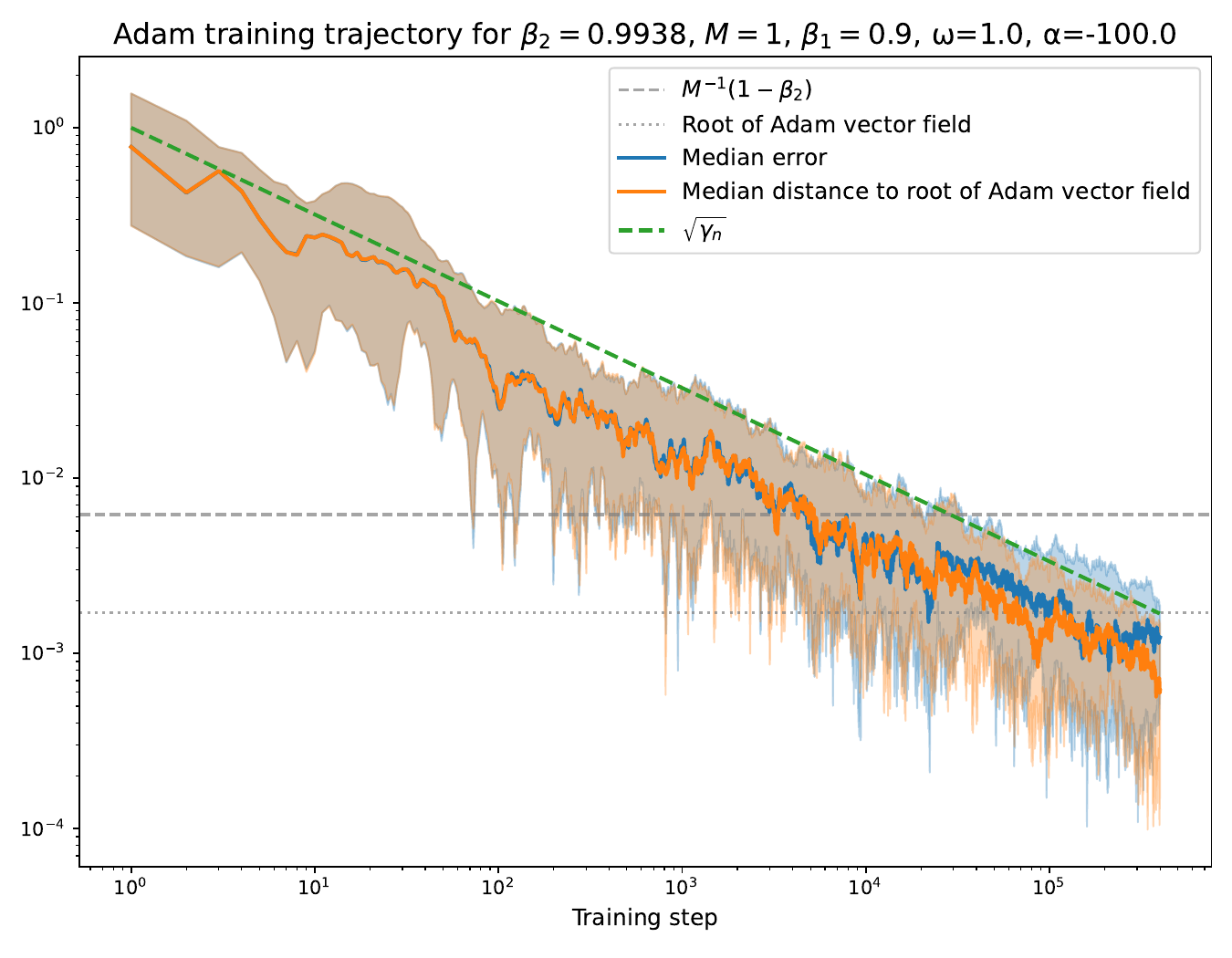}
        % \caption{Subfigure 1}
    \end{subfigure}
    % \hfill
    % \hspace{-1cm}
    \begin{subfigure}[b]{0.4\textwidth}
        \centering
        \includegraphics[width=\textwidth]{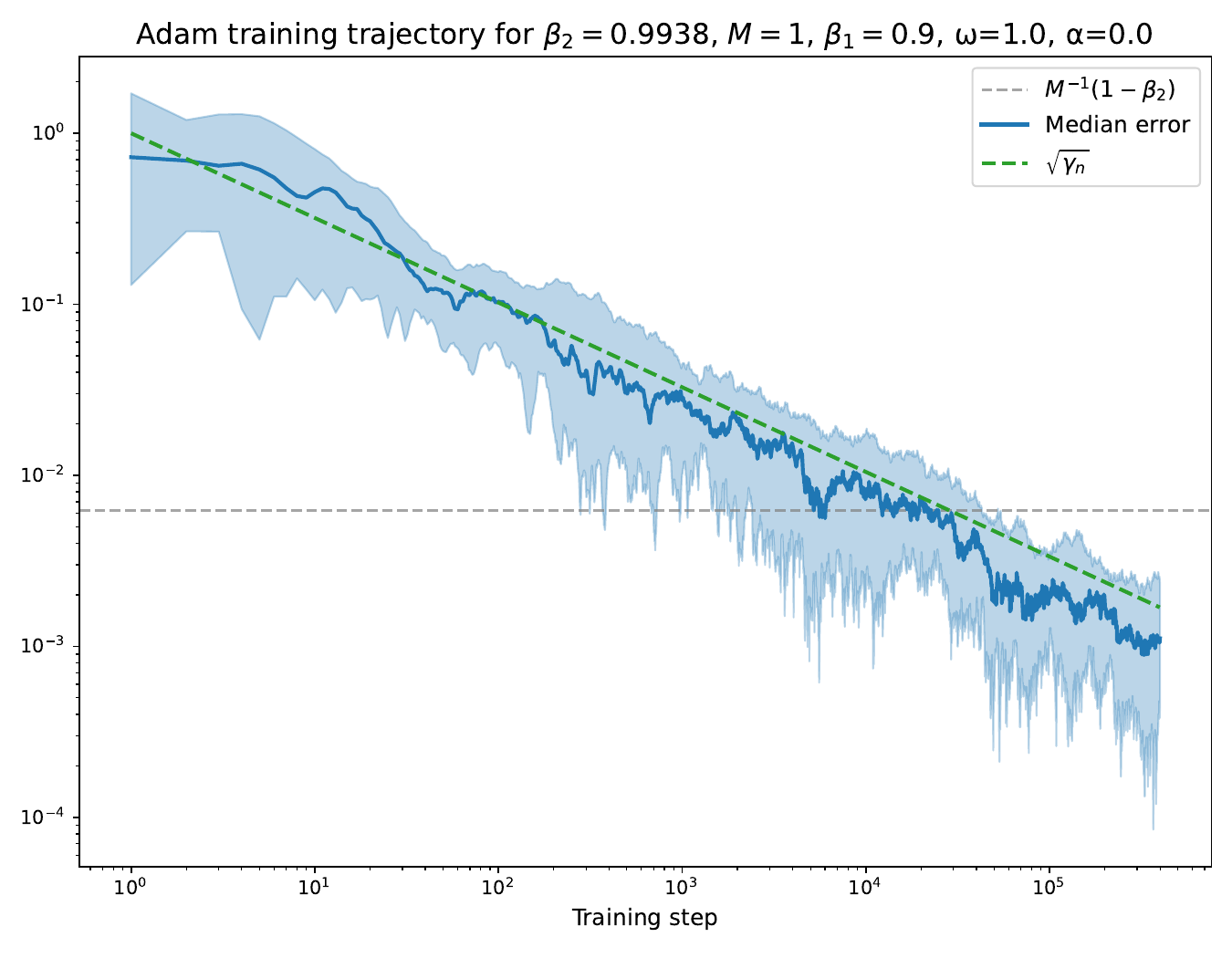}
        % \caption{Subfigure 2}
    \end{subfigure}
    
    % \vspace{-0.5cm}

    \begin{subfigure}[b]{0.4\textwidth}
        \centering
        \includegraphics[width=\textwidth]{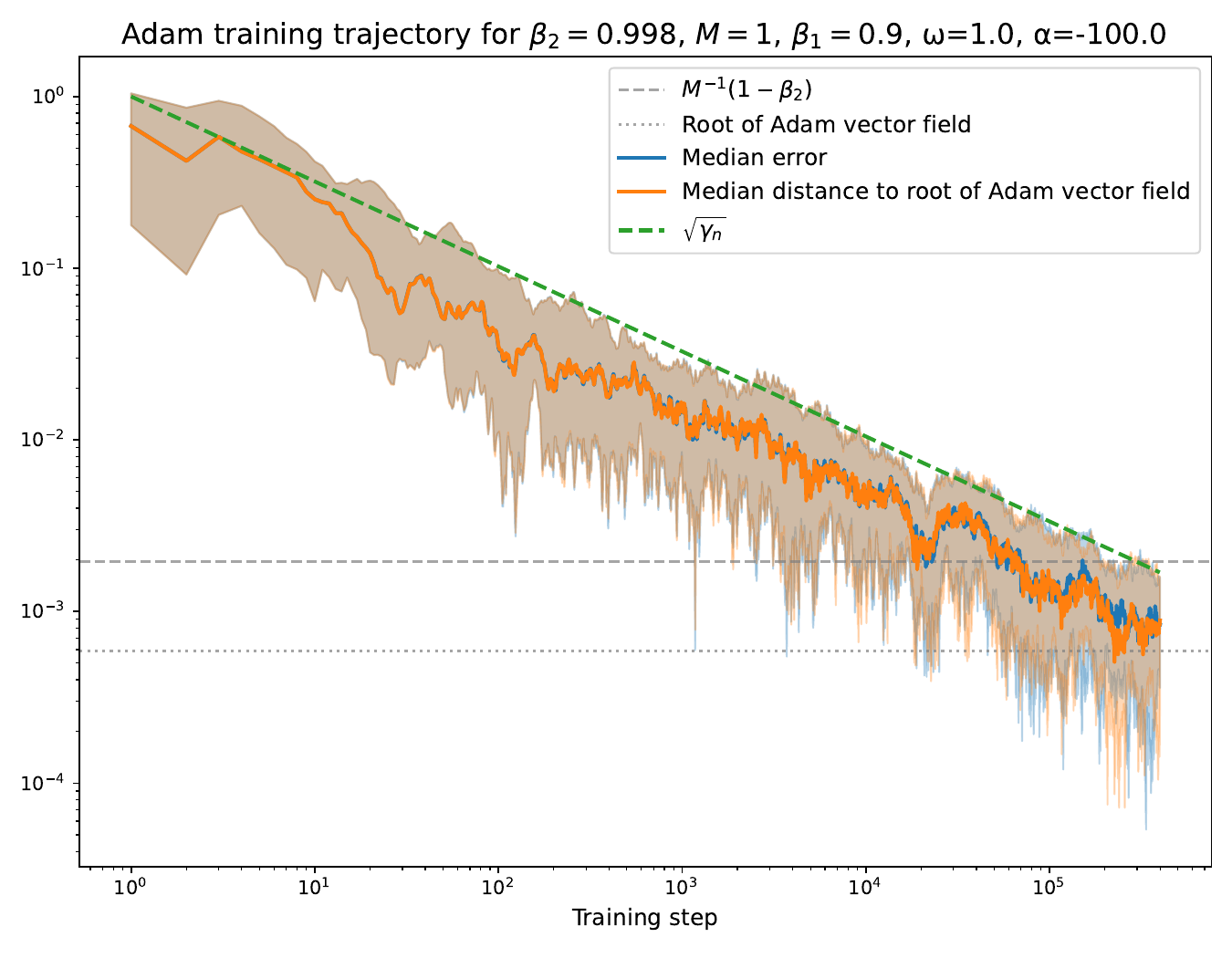}
        % \caption{Subfigure 1}
    \end{subfigure}
    % \hfill
    % \hspace{-1cm}
    \begin{subfigure}[b]{0.4\textwidth}
        \centering
        \includegraphics[width=\textwidth]{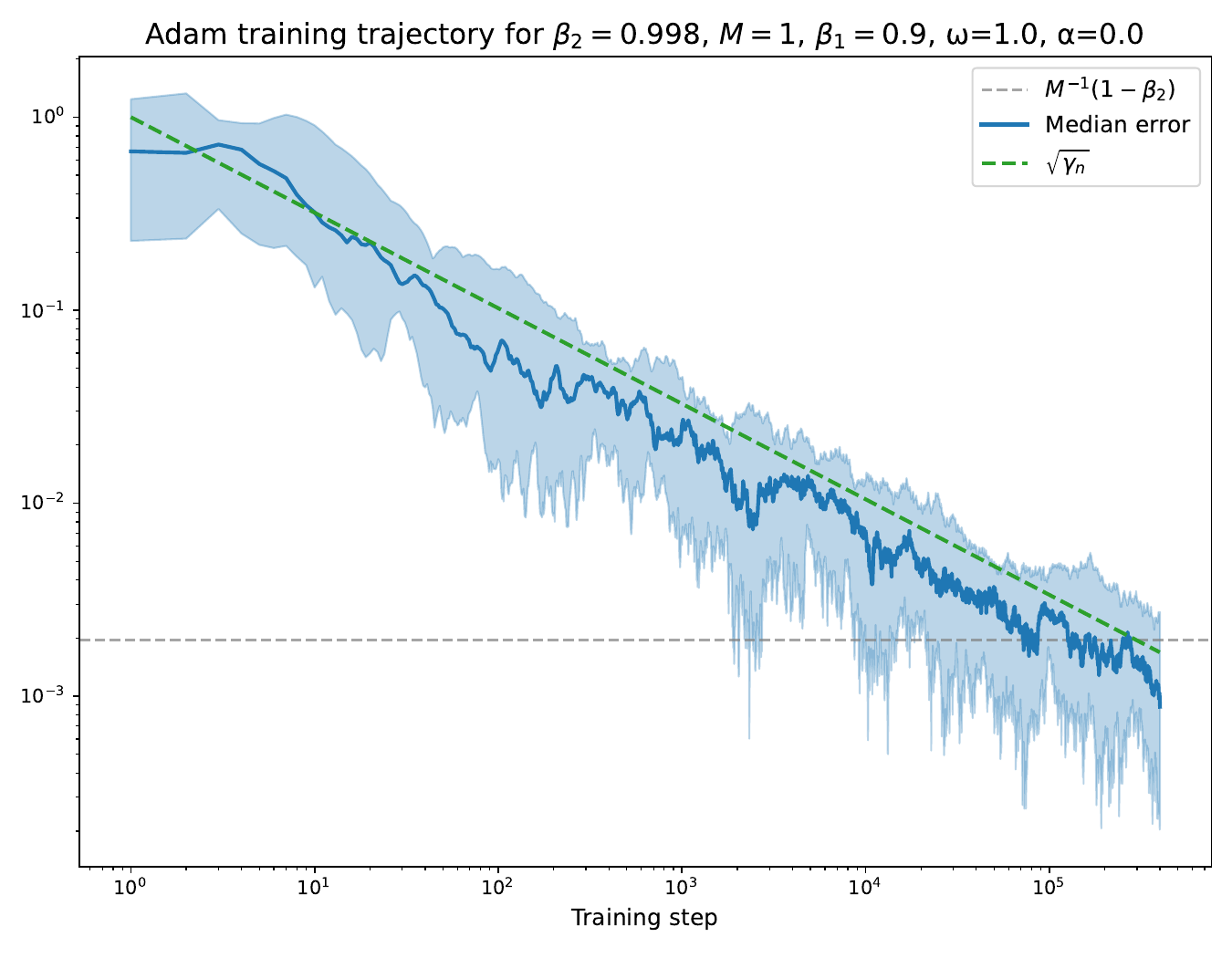}
        % \caption{Subfigure 2}
    \end{subfigure}

    \caption{
    Approximative plot of
    $
      \E\big[ | \mandVcom{0}{M}{n}{i} | \big]
    $
    and 
    $
      \E\big[ | \mandVcom{0}{M}{n}{i} - \vartheta| \big]
    $
    against 
    $n \in \{1, 2, \ldots, \numGDsteps\}$
    where 
    $M = 1$, $\beta_1= 0.9$, $\beta_2 \in (\cup_{i = 0}^{3} \{ 1-2^{(-4 - 5i/3)} \})$, $\varepsilon=10^{-8}$, $\alpha \in \{-100, 0\}$, and $\omega = 1$ and where $\vartheta$ is a zero of the \Adam\ vector field $f$ (cf.\ \cref{def: Adam vector field numerical}).
    \emph{Left:} We consider the case of asymmetrically distributed data ($\alpha = -100$). 
    \emph{Right:} We consider the case of symmetrically distributed data ($\alpha = 0$).
    }
    \label{fig:Skew Adam training trajectory}
\end{figure}

\subsubsection{Simulations for the convergence rates with respect to $ \beta_2 $}

In \cref{fig:skew_cv_rate_beta2} we consider \cref{setting: skew normal} and approximately plot 
$
  \E\big[ | \mandVcom{0}{M}{\numGDsteps}{i} | \big]
$
and
a zero of the Adam vector field $f$ (cf.\ \cref{def: Adam vector field numerical})
against
$\beta_2 \in ( \cup_{i = 0}^{9} \{ 1-2^{(-4 - 5i/9)} \} )$
where $M = 1$, $\beta_1 = 0.9$, $\varepsilon=10^{-8}$, $\alpha \in \{-100, 0\}$, and $\omega = 1$
(cf.\ \cref{setting: skew normal})

\begin{figure}[htbp]
    \centering
    \begin{subfigure}[b]{0.45\textwidth}
        \centering
        \includegraphics[width=\textwidth]{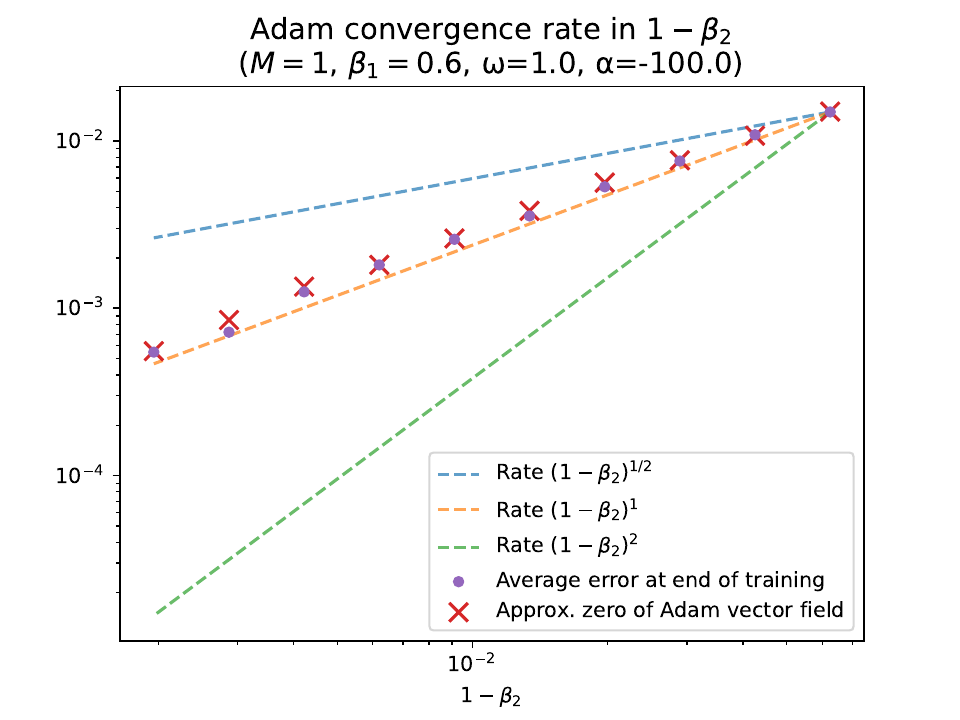}
        % \caption{Subfigure 1}
    \end{subfigure}
    % \hfill
    % \hspace{-1cm}
    \begin{subfigure}[b]{0.45\textwidth}
        \centering
        \includegraphics[width=\textwidth]{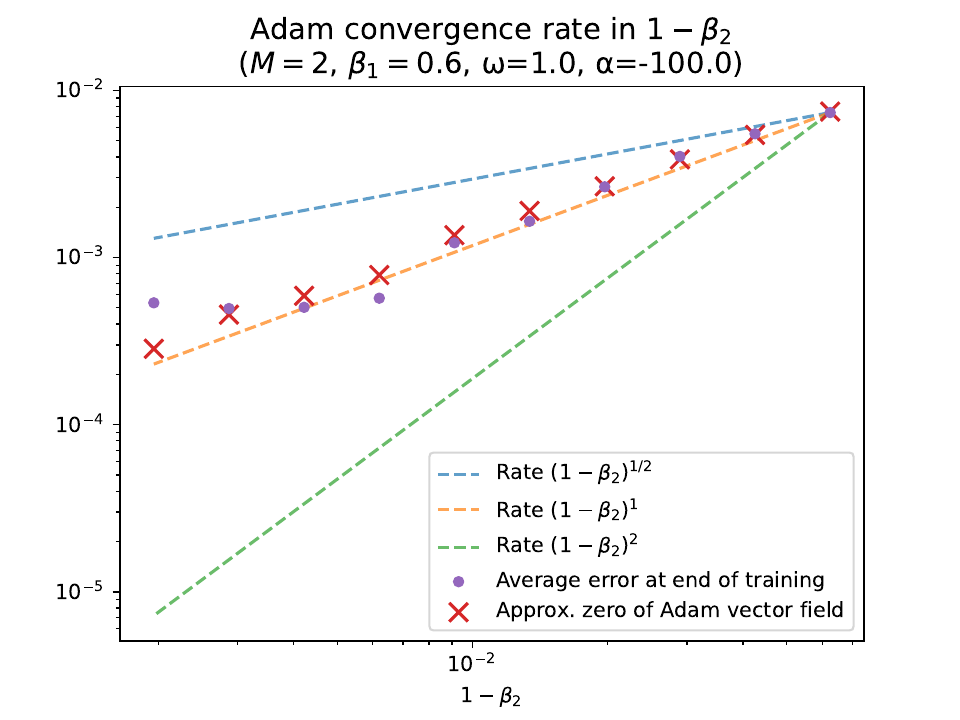}
        % \caption{Subfigure 2}
    \end{subfigure}
    
    % \vspace{-0.2cm}
    
    \begin{subfigure}[b]{0.45\textwidth}
        \centering
        \includegraphics[width=\textwidth]{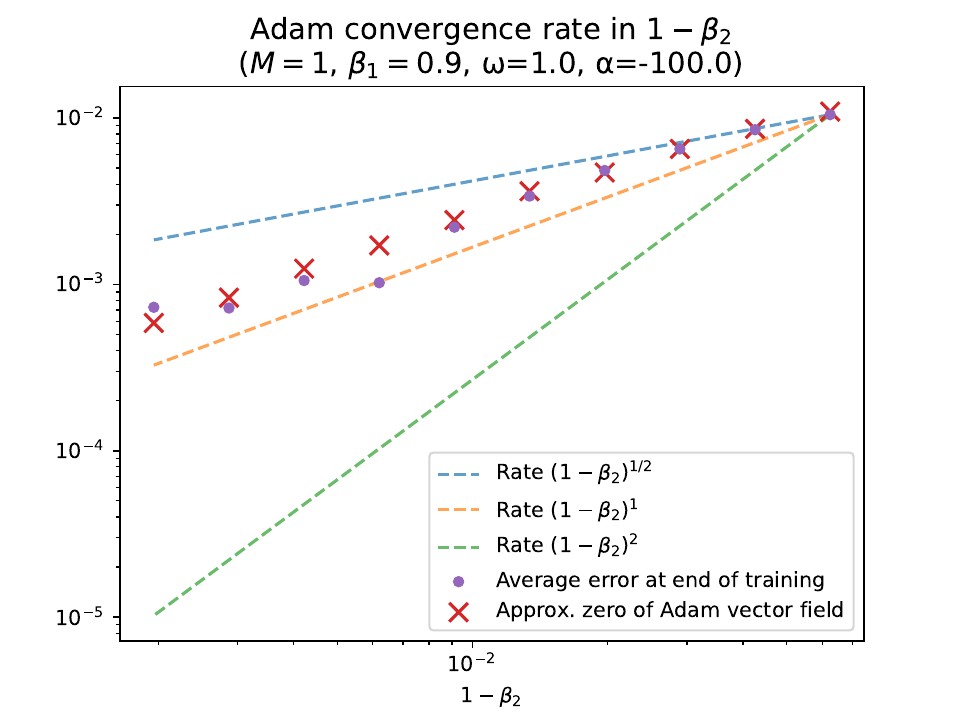}
        % \caption{Subfigure 3}
    \end{subfigure}
    % \hfill
    % \hspace{-1cm}
    \begin{subfigure}[b]{0.45\textwidth}
        \centering
        \includegraphics[width=\textwidth]{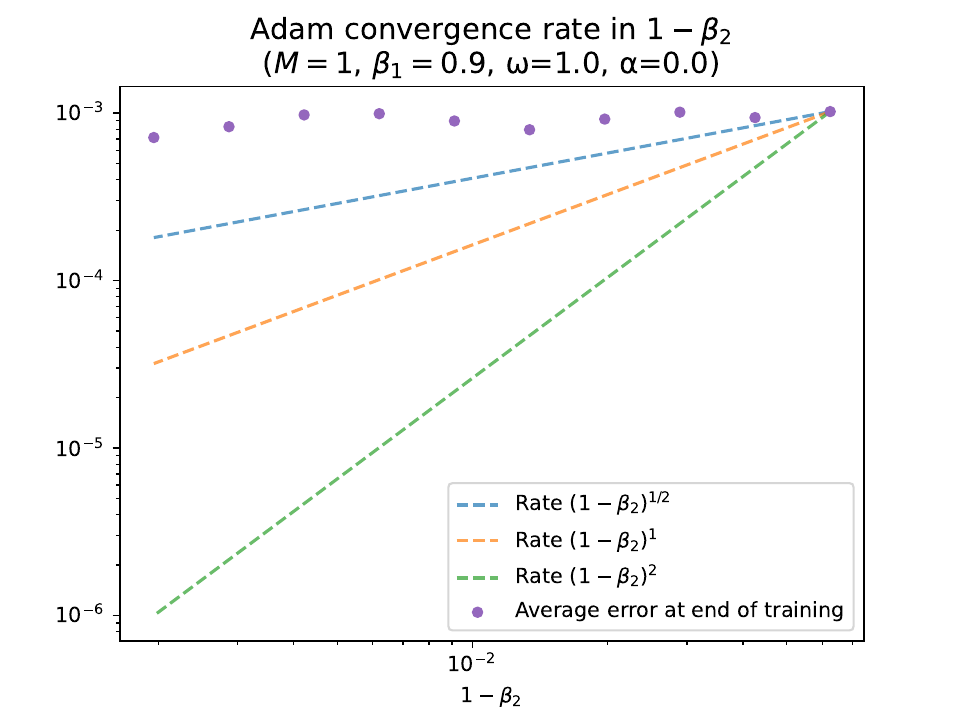}
        % \caption{Subfigure 4}
    \end{subfigure}
    
    \caption{
    Approximative plot 
    of
    $
      \E\big[ | \mandVcom{0}{M}{\numGDsteps}{i} | \big]
    $
    and
    a zero of the Adam vector field $f$ (cf.\ \cref{def: Adam vector field numerical})
    against
    $\beta_2 \in ( \cup_{i = 0}^{9} \{ 1-2^{(-4 - 5i/9)} \} )$
    where $M \in \{1, 2\}$, $\beta_1 = 0.9$, $\varepsilon=10^{-8}$, $\alpha \in \{-100, 0\}$, and $\omega = 1$
    (cf.\ \cref{setting: skew normal}).
    }
    \label{fig:skew_cv_rate_beta2}
\end{figure}

\subsubsection{Simulations for the Adam symmetry theorem}

\newcommand{\numGDstepsSecond}{3000000}
 
In \cref{fig:skew 3d plot 1.1,fig:skew 3d plot 1.2,fig:skew 3d plot 1.3} we consider \cref{setting: skew normal} and approximately plot
$
  \E\big[ \mandVcom{0}{M}{\numGDsteps}{i} \big]
$
and
a zero of the Adam vector field $f$ (cf.\ \cref{def: Adam vector field numerical})
against
$\alpha \in \{-6, 4, \ldots, 6\}$
where $M=1$, $\beta_1 = 0.6$, $\beta_2 = 0.9$, $\varepsilon=10^{-8}$, and $\omega=3$.
Moreover, in \cref{fig:skew 3d plot 2.1,fig:skew 3d plot 2.2,fig:skew 3d plot 2.3} we consider \cref{setting: skew normal} and approximately plot
$
  \E\big[ \mandVcom{0}{M}{\numGDstepsSecond}{i} \big]
$
and
a zero of the Adam vector field $f$ (cf.\ \cref{def: Adam vector field numerical})
against
$\alpha \in \{-6, 4, \ldots, 6\}$
where $M=1$, $\beta_1 = 0.9$, $\beta_2 = 0.95$, $\varepsilon=10^{-8}$, and $\omega=3$.

\begin{figure}[p]
    \centering
    \begin{subfigure}[b]{0.8\textwidth}
        \centering
        \includegraphics[width=\textwidth,page=4]{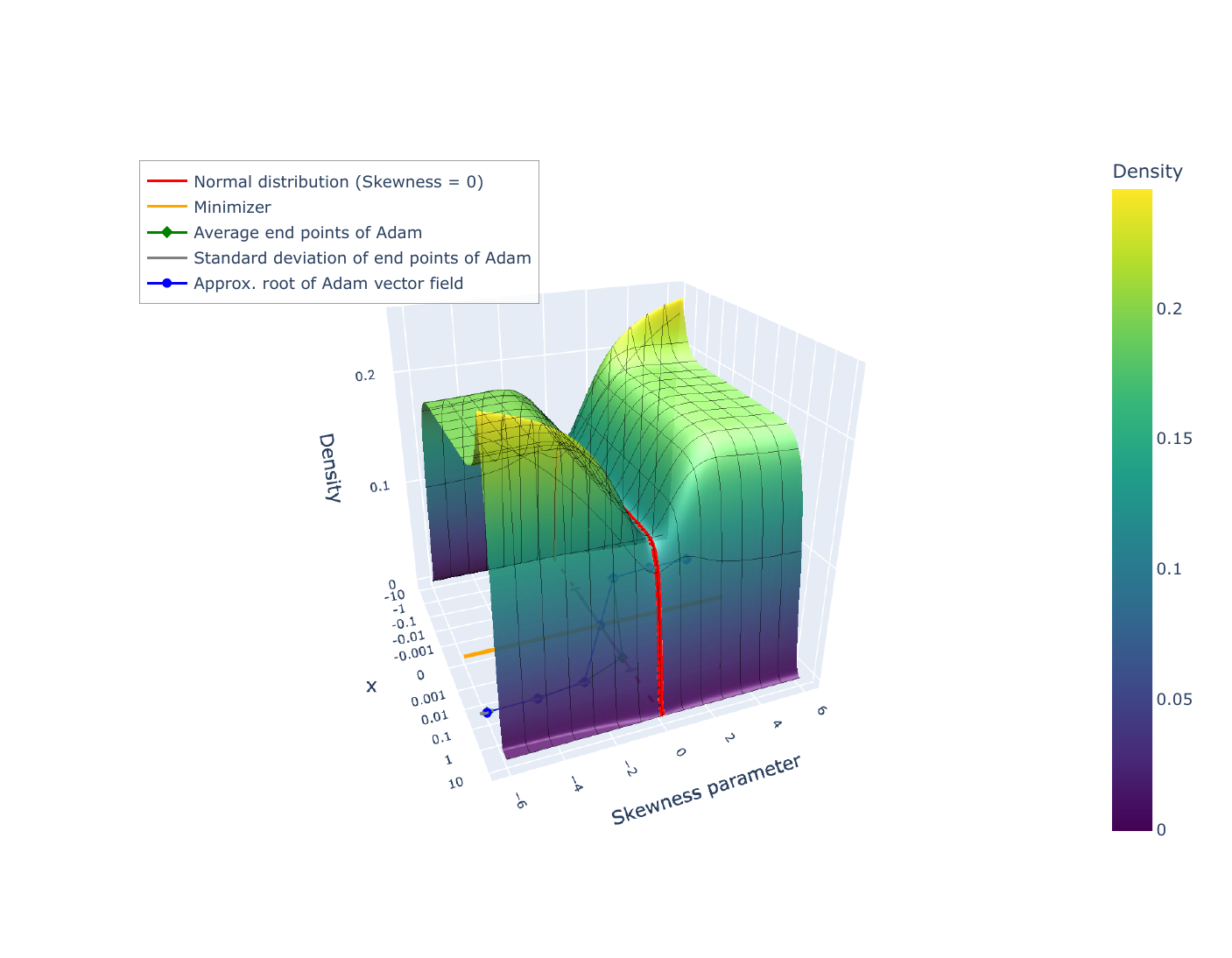}
        % \caption{View 1}
        % \label{fig:sub1}
    \end{subfigure}

    \vspace{-1cm}

      \begin{subfigure}[b]{0.8\textwidth}
        \centering
        \includegraphics[width=\textwidth,page=6]{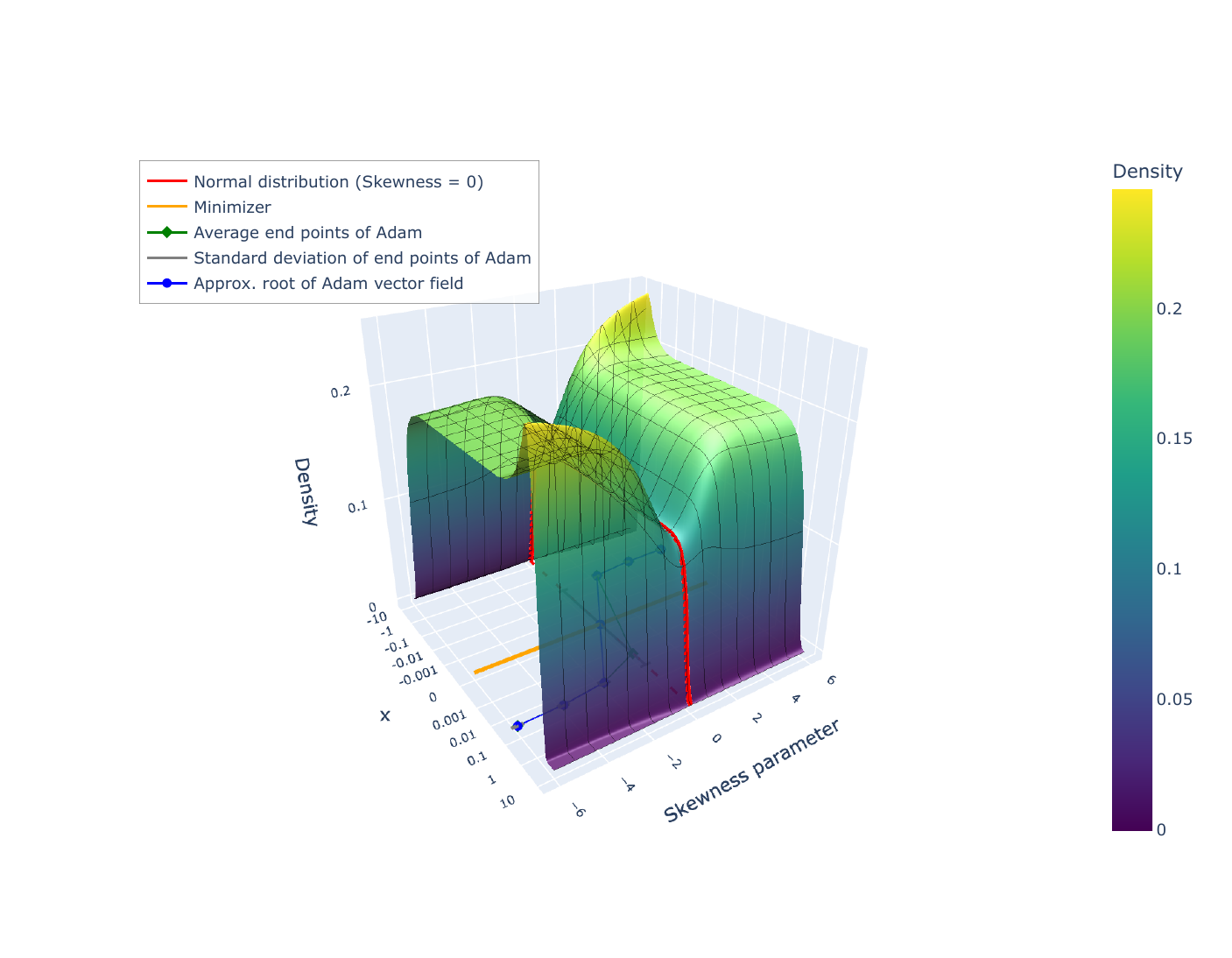}
        % \caption{View 2}
        % \label{fig:sub2}
    \end{subfigure}

    \caption{
    Approximative plot of
    $
      \E\big[ \mandVcom{0}{M}{\numGDsteps}{i} \big]
    $
    and
    a zero of the Adam vector field $f$ (cf.\ \cref{def: Adam vector field numerical})
    against
    $\alpha \in \{-6, 4, \ldots, 6\}$
    where $M=1$, $\beta_1=0.6$, $\beta_2 = 0.9$, $\varepsilon=10^{-8}$, and $\omega=3$.
    The density axis corresponds to $\skewnormalpdf$ (cf.\ \cref{setting: skew normal}).
    All the images in \cref{fig:skew 3d plot 1.1,fig:skew 3d plot 1.2,fig:skew 3d plot 1.3} correspond to the same plot from different angles.
    } 
    \label{fig:skew 3d plot 1.1}
\end{figure}

\begin{figure}[p]
    \centering
    \begin{subfigure}[b]{0.8\textwidth}
        \centering
        \includegraphics[width=\textwidth,page=8]{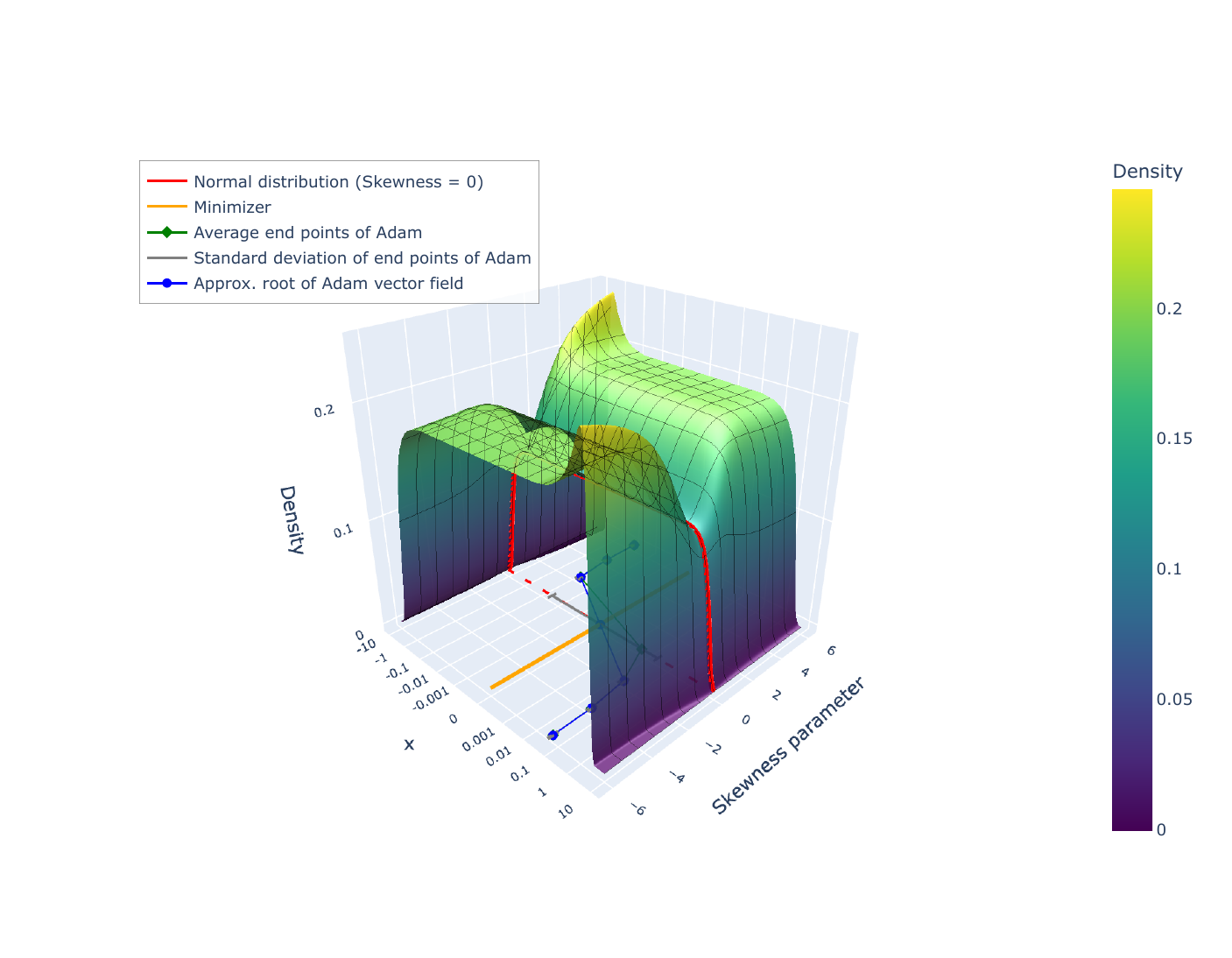}
        % \caption{View 3}
        % \label{fig:sub3}
    \end{subfigure}

    \vspace{-1cm}

    \begin{subfigure}[b]{0.8\textwidth}
        \centering
        \includegraphics[width=\textwidth,page=10]{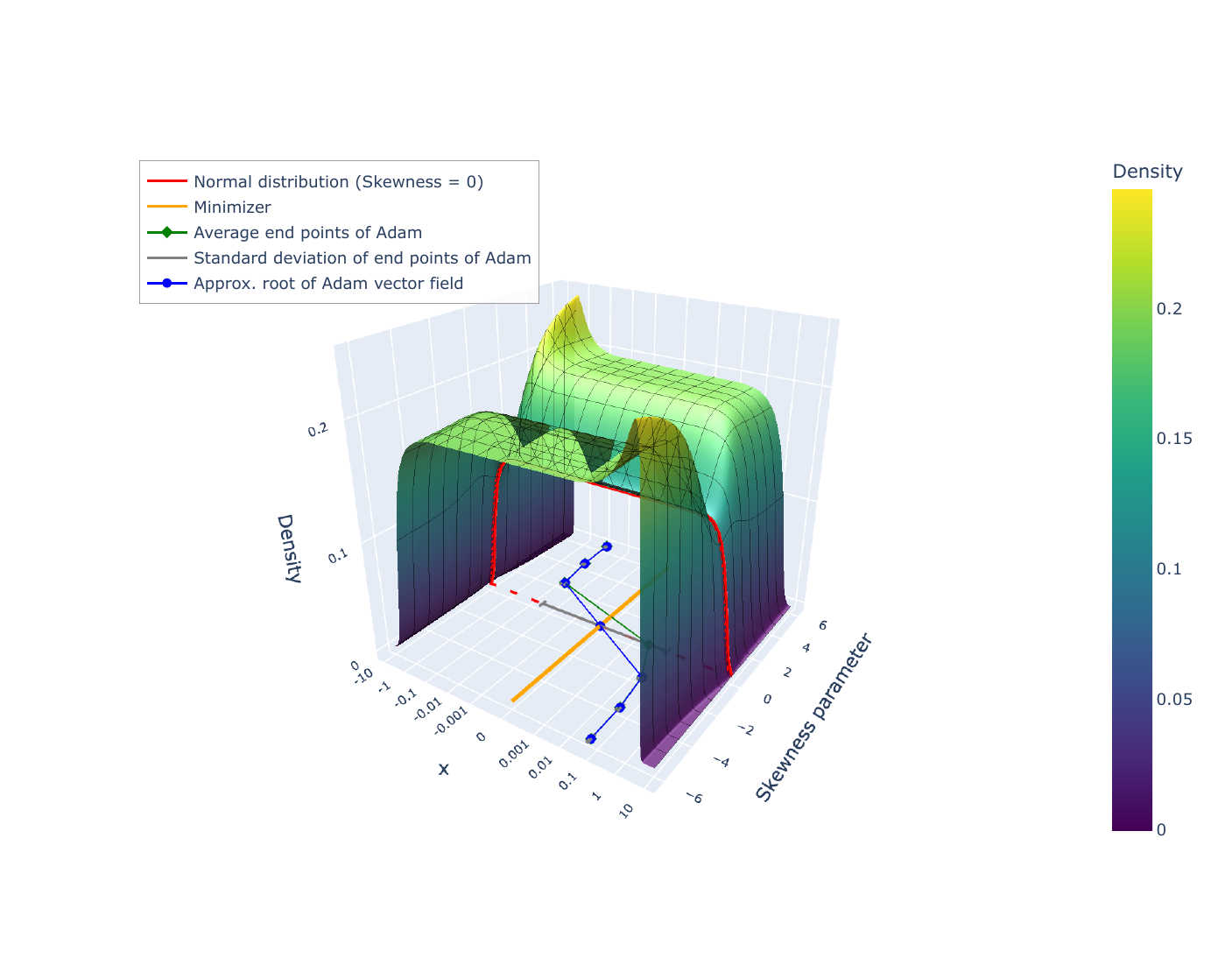}
        % \caption{View 4}
        % \label{fig:sub4}
    \end{subfigure}
    \caption{
    Approximative plot of
    $
      \E\big[ \mandVcom{0}{M}{\numGDsteps}{i} \big]
    $
    and
    a zero of the Adam vector field $f$ (cf.\ \cref{def: Adam vector field numerical})
    against
    $\alpha \in \{-6, 4, \ldots, 6\}$
    where $M=1$, $\beta_1=0.6$, $\beta_2 = 0.9$, $\varepsilon=10^{-8}$, and $\omega=3$.
    The density axis corresponds to $\skewnormalpdf$ (cf.\ \cref{setting: skew normal}).
    All the images in \cref{fig:skew 3d plot 1.1,fig:skew 3d plot 1.2,fig:skew 3d plot 1.3} are the same plot from different angles.
    } 
    \label{fig:skew 3d plot 1.2}
\end{figure}

\begin{figure}[p]
    \centering
    \begin{subfigure}[b]{0.8\textwidth}
        \centering
        \includegraphics[width=\textwidth,page=12]{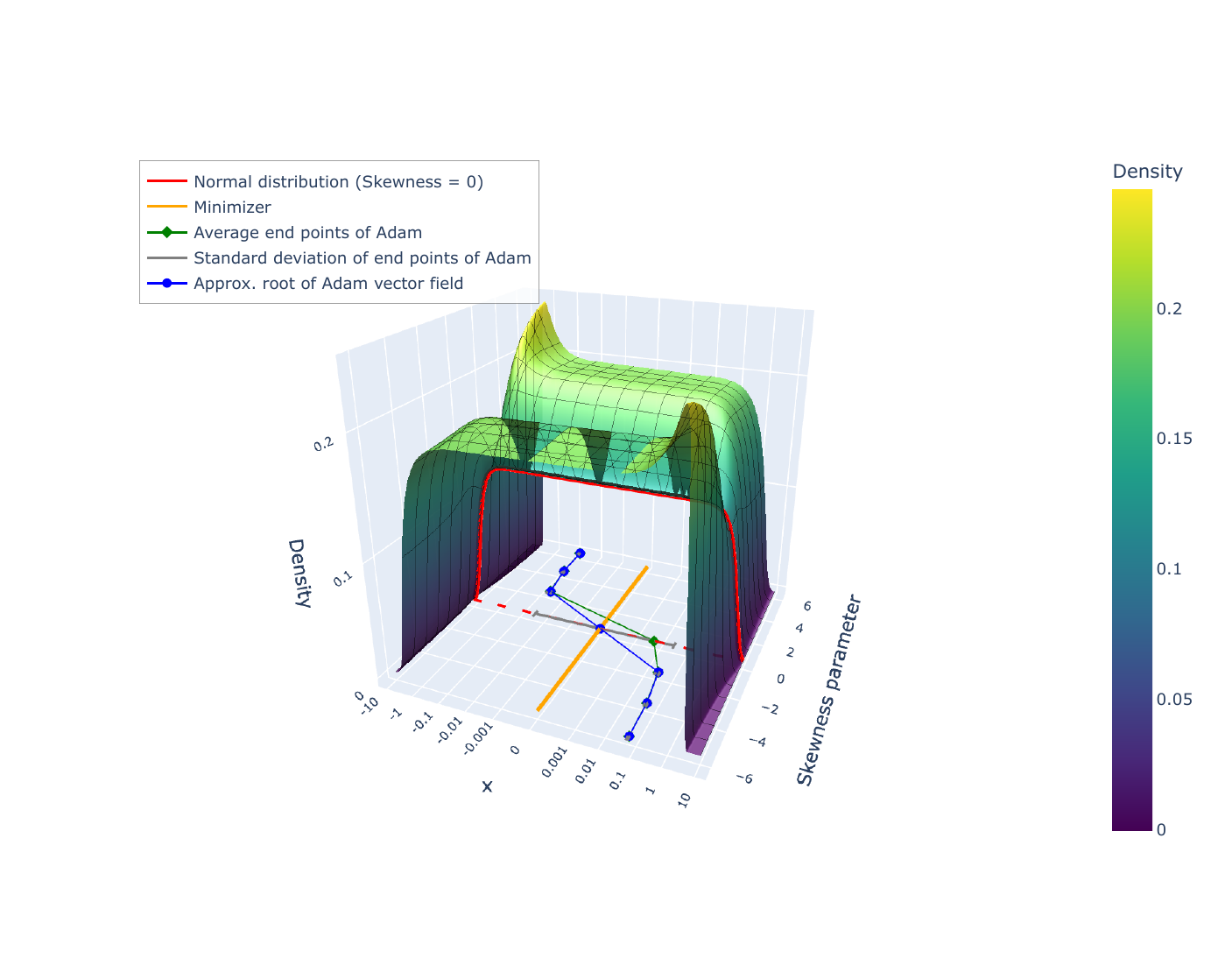}
        % \caption{View 3}
        % \label{fig:sub3}
    \end{subfigure}

    \vspace{-1cm}

    \begin{subfigure}[b]{0.8\textwidth}
        \centering
        \includegraphics[width=\textwidth,page=13]{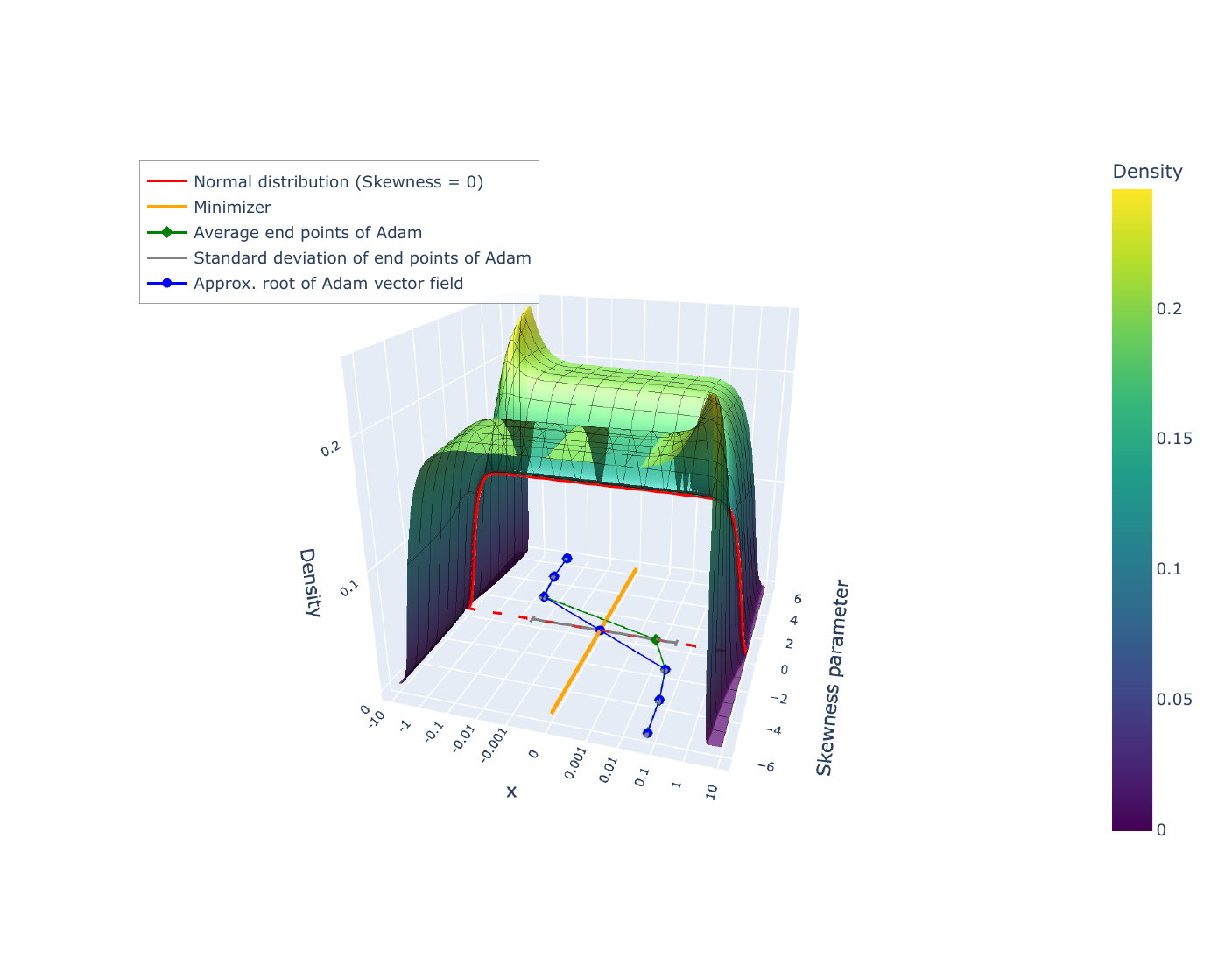}
        % \caption{View 4}
        % \label{fig:sub4}
    \end{subfigure}
    \caption{
    Approximative plot of
    $
      \E\big[ \mandVcom{0}{M}{\numGDsteps}{i} \big]
    $
    and
    a zero of the Adam vector field $f$ (cf.\ \cref{def: Adam vector field numerical})
    against
    $\alpha \in \{-6, 4, \ldots, 6\}$
    where $M=1$, $\beta_1=0.6$, $\beta_2 = 0.9$, $\varepsilon=10^{-8}$, and $\omega=3$.
    The density axis corresponds to $\skewnormalpdf$ (cf.\ \cref{setting: skew normal}).
    All the images in \cref{fig:skew 3d plot 1.1,fig:skew 3d plot 1.2,fig:skew 3d plot 1.3} are the same plot from different angles.
    } 
    \label{fig:skew 3d plot 1.3}
\end{figure}

%%%%%%%%%%%%%%%%%%%%%%%%%%%%%%%%%%%%%%%%

\begingroup
\renewcommand{\numGDsteps}{3000000}

\begin{figure}[p]
    \centering
    \begin{subfigure}[b]{0.8\textwidth}
        \centering
        \includegraphics[width=\textwidth,page=4]{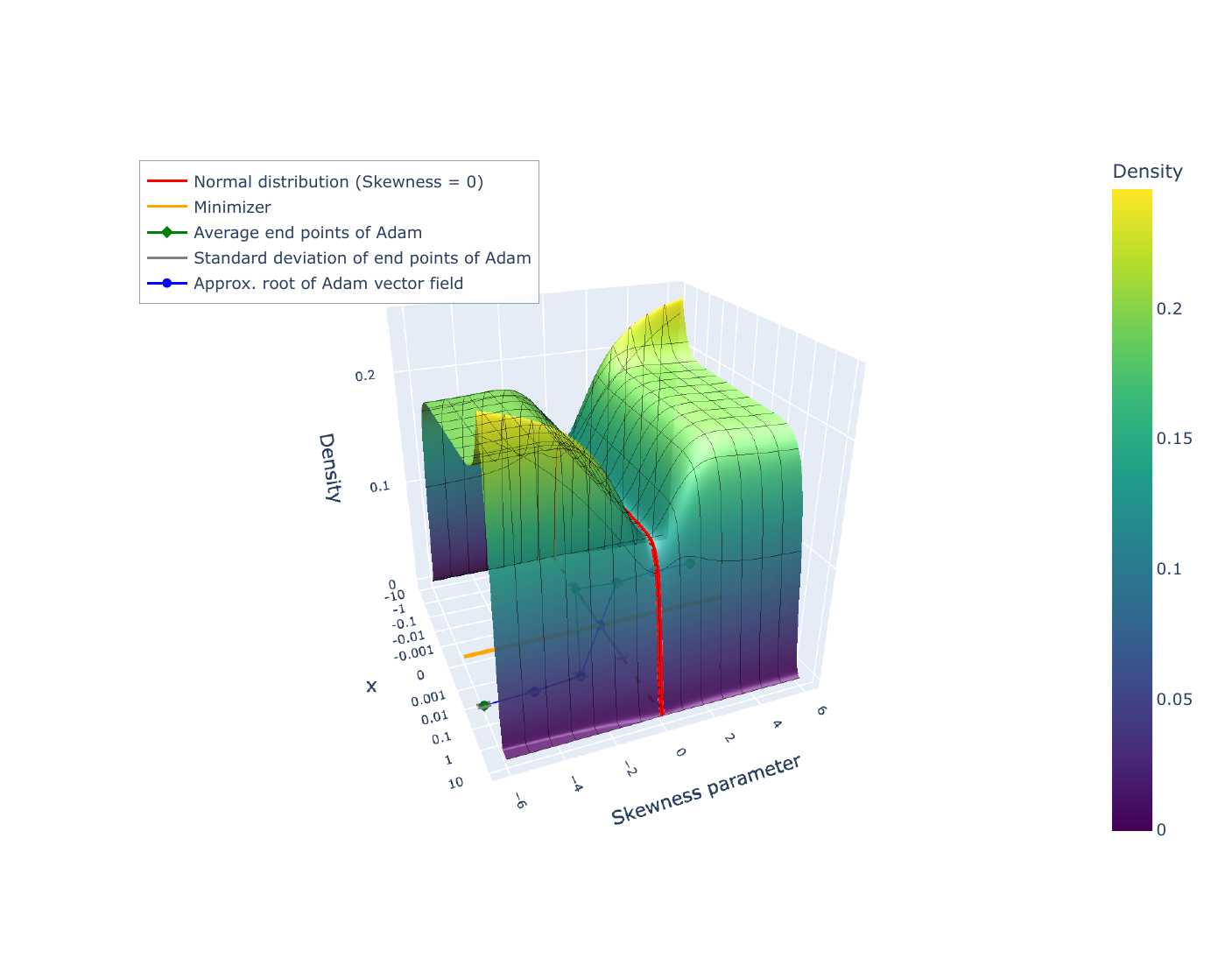}
        % \caption{View 1}
        % \label{fig:sub1}
    \end{subfigure}

    \vspace{-1cm}

      \begin{subfigure}[b]{0.8\textwidth}
        \centering
        \includegraphics[width=\textwidth,page=6]{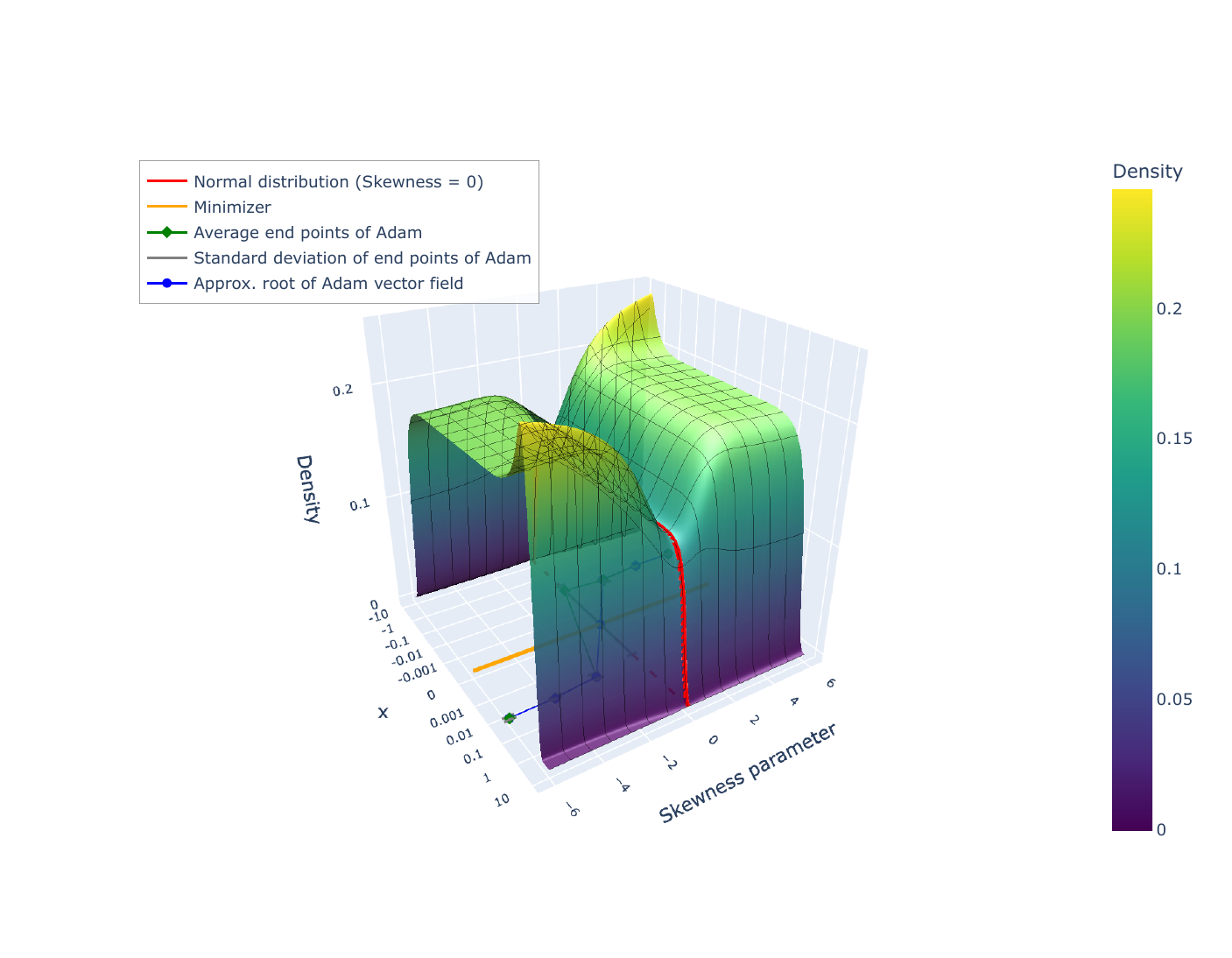}
        % \caption{View 2}
        % \label{fig:sub2}
    \end{subfigure}

    \caption{
    Approximative plot of
    $
      \E\big[ \mandVcom{0}{M}{\numGDsteps}{i} \big]
    $
    and
    a zero of the Adam vector field $f$ (cf.\ \cref{def: Adam vector field numerical})
    against
    $\alpha \in \{-6, 4, \ldots, 6\}$
    where $M=1$, $\beta_1=0.9$, $\beta_2 = 0.95$, $\varepsilon=10^{-8}$, and $\omega=3$.
    The density axis corresponds to $\skewnormalpdf$ (cf.\ \cref{setting: skew normal}).
    All the images in \cref{fig:skew 3d plot 2.1,fig:skew 3d plot 2.2,fig:skew 3d plot 2.3} correspond to the same plot from different angles.
    } 
    \label{fig:skew 3d plot 2.1}
\end{figure}

\begin{figure}[p]
    \centering
    \begin{subfigure}[b]{0.8\textwidth}
        \centering
        \includegraphics[width=\textwidth,page=8]{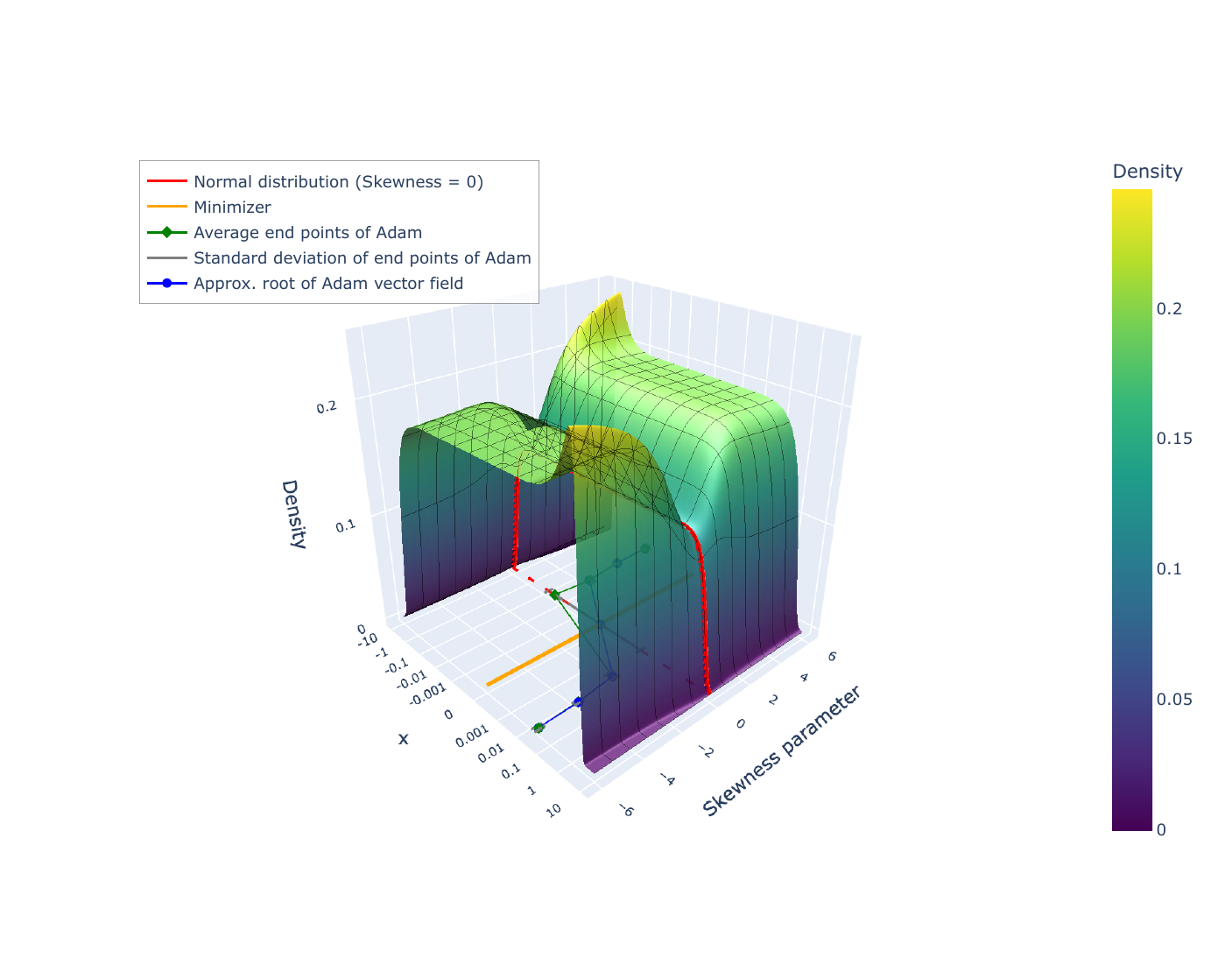}
        % \caption{View 3}
        % \label{fig:sub3}
    \end{subfigure}

    \vspace{-1cm}

    \begin{subfigure}[b]{0.8\textwidth}
        \centering
        \includegraphics[width=\textwidth,page=10]{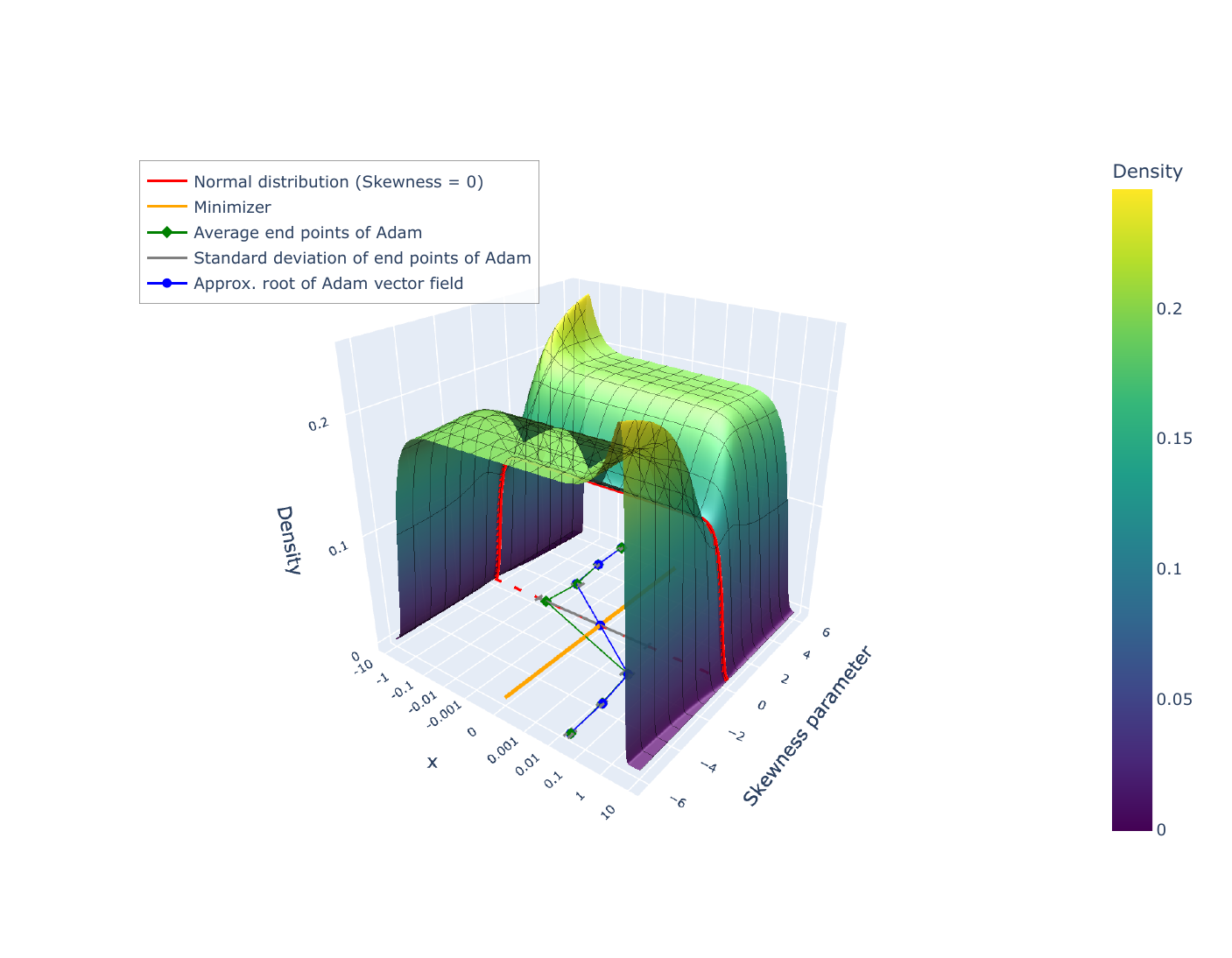}
        % \caption{View 4}
        % \label{fig:sub4}
    \end{subfigure}
    \caption{
    Approximative plot of
    $
      \E\big[ \mandVcom{0}{M}{\numGDsteps}{i} \big]
    $
    and
    a zero of the Adam vector field $f$ (cf.\ \cref{def: Adam vector field numerical})
    against
    $\alpha \in \{-6, 4, \ldots, 6\}$
    where $M=1$, $\beta_1=0.9$, $\beta_2 = 0.95$, $\varepsilon=10^{-8}$, and $\omega=3$.
    The density axis corresponds to $\skewnormalpdf$ (cf.\ \cref{setting: skew normal}).
    All the images in \cref{fig:skew 3d plot 2.1,fig:skew 3d plot 2.2,fig:skew 3d plot 2.3} are the same plot from different angles.
    } 
    \label{fig:skew 3d plot 2.2}
\end{figure}

\begin{figure}[p]
    \centering
    \begin{subfigure}[b]{0.8\textwidth}
        \centering
        \includegraphics[width=\textwidth,page=12]{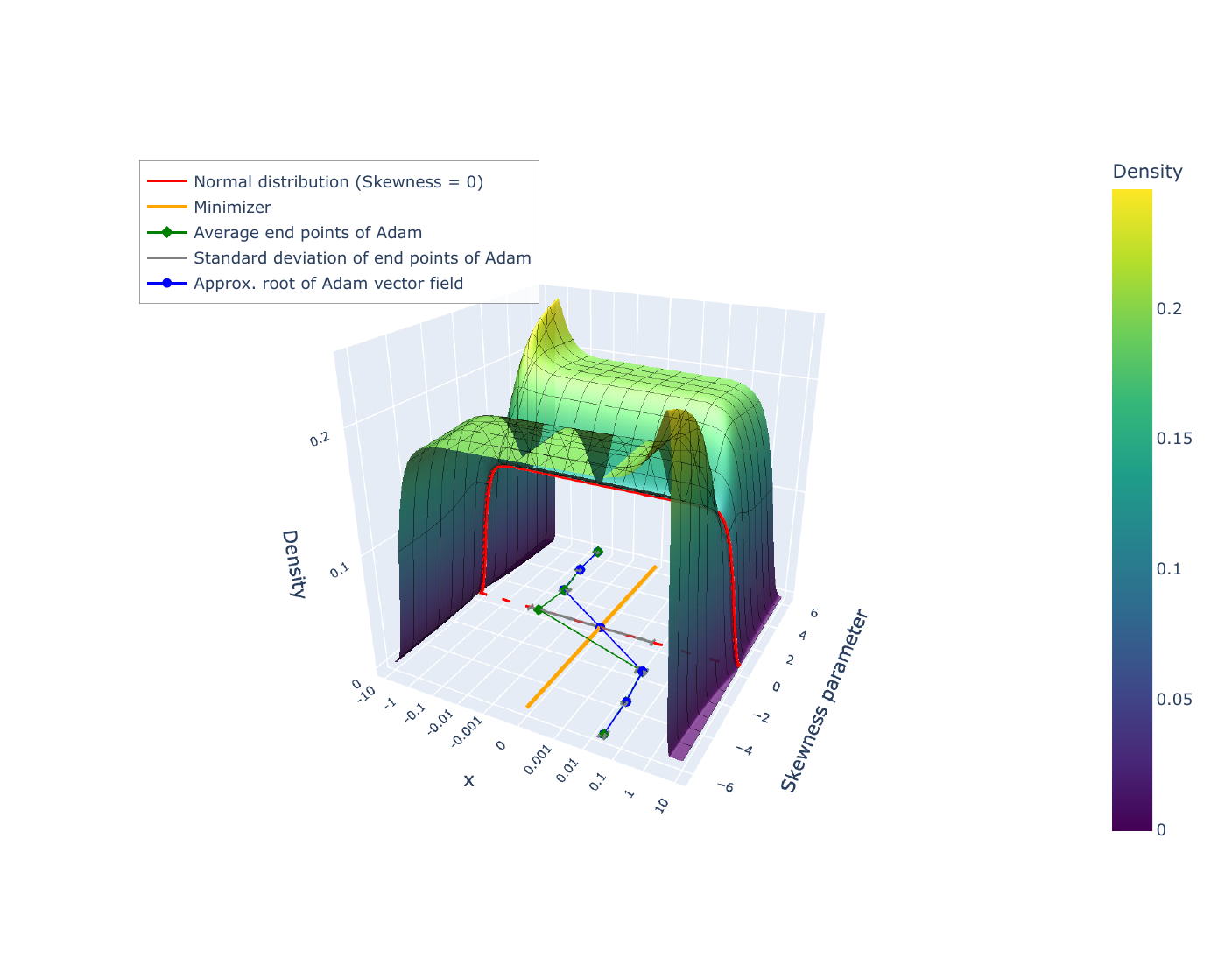}
        % \caption{View 3}
        % \label{fig:sub3}
    \end{subfigure}

    \vspace{-1cm}

    \begin{subfigure}[b]{0.8\textwidth}
        \centering
        \includegraphics[width=\textwidth,page=13]{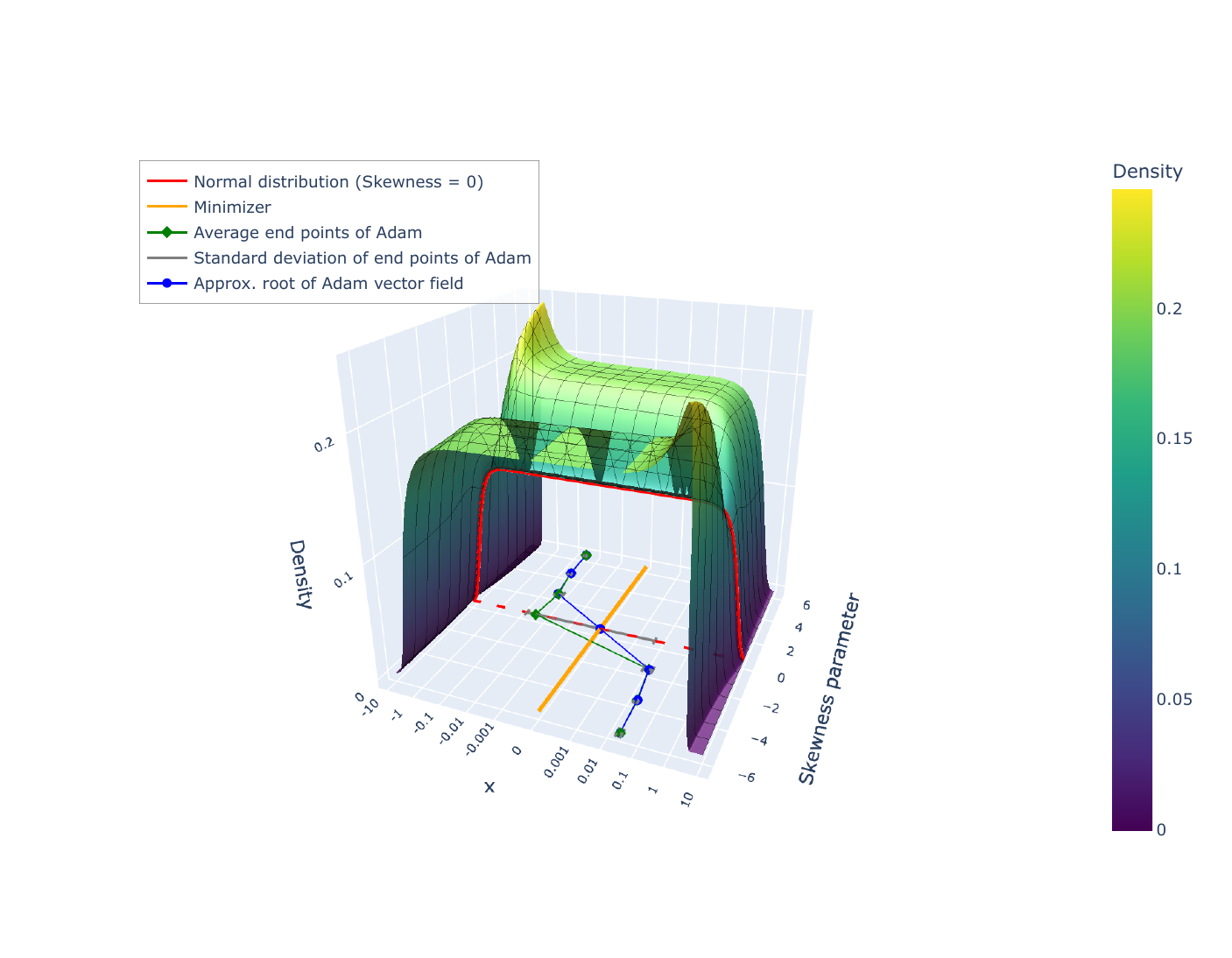}
        % \caption{View 4}
        % \label{fig:sub4}
    \end{subfigure}
    \caption{
    Approximative plot of
    $
      \E\big[ \mandVcom{0}{M}{\numGDsteps}{i} \big]
    $
    and
    a zero of the Adam vector field $f$ (cf.\ \cref{def: Adam vector field numerical})
    against
    $\alpha \in \{-6, 4, \ldots, 6\}$
    where $M=1$, $\beta_1=0.9$, $\beta_2 = 0.95$, $\varepsilon=10^{-8}$, and $\omega=3$.
    The density axis corresponds to $\skewnormalpdf$ (cf.\ \cref{setting: skew normal}).
    All the images in \cref{fig:skew 3d plot 2.1,fig:skew 3d plot 2.2,fig:skew 3d plot 2.3} are the same plot from different angles.
    } 
    \label{fig:skew 3d plot 2.3}
\end{figure}

\endgroup

\endgroup

\subsubsection*{Acknowledgements}
This work has been partially funded by the National Science Foundation of China (NSFC) under grant number W2531010. Moreover, this work has been partially funded by the European Union (ERC, MONTECARLO, 101045811). The views and the opinions expressed in this work are however those of the authors only and do not necessarily reflect those of the European Union or the European Research Council (ERC). Neither the European Union nor the granting authority can be held responsible for them. Moreover, this work has been partially funded by the Deutsche Forschungsgemeinschaft (DFG, German Research Foundation) under Germany’s Excellence Strategy EXC 2044-390685587, Mathematics Münster: Dynamics-Geometry-Structure. Most of the specific formulations in the proofs of this work have been created using \cite{Bennoargumentcommand}.

\pagebreak
\bibliographystyle{acm}
\bibliography{bibfileAdamconvergenceII}

\end{document}